\documentclass[11pt,letterpaper]{article}
\usepackage{float, graphicx}

\usepackage[]{epsfig}
\usepackage{amsmath, amsthm, amssymb,}
\usepackage{epsfig}
\usepackage{verbatim}
\usepackage{multicol}
\usepackage{url}
\usepackage{latexsym}
\usepackage{mathrsfs}
\usepackage[colorlinks, bookmarks=true, citecolor=blue]{hyperref}
\usepackage{amsmath}
\usepackage{enumerate}
\usepackage[normalem]{ulem}
\usepackage{bm}
\usepackage{dsfont}
\usepackage{stmaryrd}
\usepackage{enumitem}
\usepackage[final]{showlabels}
\usepackage{tikz}
\usepackage{cleveref}

\usepackage{cite}
\usepackage{color}
\usepackage{xcolor}

\usepackage{authblk}

\usepackage[margin=0.8in]{geometry}

\makeatletter

\def\@settitle{\begin{center}%
    \bfseries
 \normalfont\LARGE\@title
  \end{center}%
}
\def\@setauthors{\begin{center}%
 \normalsize\@author
  \end{center}%
}

\makeatother

\numberwithin{equation}{section}

\newcommand{\sfb}{{\mathsf b}}

\newcommand{\sfT}{{\mathsf T}}

\renewcommand{\cal}{\mathcal}
\newcommand\cA{{\mathcal A}}
\newcommand\cB{{\mathcal B}}
\newcommand{\cC}{{\cal C}}

\newcommand{\cE}{{\cal E}}
\newcommand{\cF}{{\cal F}}
\newcommand{\cG}{{\cal G}}

\newcommand{\cM}{{\cal M}}

\newcommand{\bme}{{\bm{e}}}

\newcommand{\bmu}{{\bm{u}}}
\newcommand{\bmv}{{\bm{v}}}

\newcommand{\bmx}{{\bm{x}}}
\newcommand{\bmy}{{\bm{y}}}

\newcommand{\bms}{\bm s}

\newcommand{\rd}{{\rm d}}

\newcommand{\ri}{\mathbf{i}}

\newcommand{\bC}{{\mathbb C}}

\newcommand{\bE}{\mathbb{E}}

\newcommand{\bP}{\mathbb{P}}

\newcommand{\bR}{{\mathbb R}}

\newcommand{\bZ}{\mathbb{Z}}

\DeclareMathOperator{\supp}{supp}
\DeclareMathOperator{\dist}{dist}

\DeclareMathOperator{\OO}{O}
\DeclareMathOperator{\oo}{o}

\renewcommand{\Re}{\mathop{\mathrm{Re}}}
\renewcommand{\Im}{\mathop{\mathrm{Im}}}

\renewcommand{\leq}{\leqslant}
\renewcommand{\geq}{\geqslant}

\newcommand{\beq}{\begin{equation}}
\newcommand{\eeq}{\end{equation}}

\theoremstyle{plain} 
\newtheorem{theorem}{Theorem}[section]
\newtheorem*{theorem*}{Theorem}
\newtheorem{lemma}[theorem]{Lemma}
\newtheorem*{lemma*}{Lemma}
\newtheorem{corollary}[theorem]{Corollary}
\newtheorem*{corollary*}{Corollary}
\newtheorem{proposition}[theorem]{Proposition}
\newtheorem*{proposition*}{Proposition}
\newtheorem{assumption}[theorem]{Assumption}
\newtheorem*{assumption*}{Assumption}
\crefname{assumption}{assumption}{assumptions}
\Crefname{assumption}{Assumption}{Assumptions}

\newtheorem{definition}[theorem]{Definition}
\newtheorem*{definition*}{Definition}

\newtheorem*{example*}{Example}
\newtheorem{remark}[theorem]{Remark}

\newtheorem*{remark*}{Remark}
\newtheorem*{remarks*}{Remarks}

\newcommand{\cout}{{\omega_+}}
\newcommand{\cin}{{\omega_-}}

\title{Dynamical edge universality and the Airy$_\beta$ line ensemble}

\author{
Jiaoyang Huang\thanks{Courant Institute of Mathematical Sciences, New York University.}
\quad
Shengjing Xu\thanks{Department of Mathematics, University of Pennsylvania.}
\quad
Yulei Xu\thanks{Department of Mathematical Sciences, Carnegie Mellon University.}
}

\date{}

\providecommand{\E}{\mathbb E}

\providecommand{\D}{\mathscr D}

\begin{document}
\sloppy

\maketitle
\begin{abstract}
We prove dynamical edge universality for three families of interacting particle systems for every $\beta>0$: $\beta$-Dyson Brownian motion with an analytic potential, weighted $\beta$-corners processes, and a class of discrete-time Markov chains. Our results apply to general non-equilibrium initial data and show that the Airy$_\beta$ line ensemble emerges at positive times whenever the limiting density has a regular square-root edge.

More precisely, after centering at the moving deterministic edge and rescaling space and time on the $n^{-2/3}$ and $n^{-1/3}$ scales, respectively, the trajectories of the top $k$ particles converge jointly in distribution to the top $k$ curves of the Airy$_\beta$ line ensemble, for every fixed $k$, in the topology of uniform convergence on compact time intervals.

The proof combines dynamical loop equations with characteristic-flow estimates to establish local laws near the moving edge and upper bounds on the extremal particle. We then identify the stochastic evolution of every subsequential limit of the rescaled Stieltjes transform and establish convergence through the pole-evolution characterization of the Airy$_\beta$ line ensemble.
\end{abstract}
\tableofcontents

\section{Introduction}
\label{sec:introduction}

\subsection{Background and the \texorpdfstring{Airy$_\beta$}{Airy-beta} line ensemble}
\label{subsec:intro-background}

The Airy line ensemble~\cite{prahofer2002scale,corwin2014brownian} is a central object in the study of the Kardar--Parisi--Zhang (KPZ) universality class, which describes the large-scale fluctuations of a broad class of stochastic growth processes; see~\cite{Corwin2012KPZ}. Its top curve is the Airy$_2$ process, whose one-point marginal is the GUE Tracy--Widom distribution~\cite{tracy1994level}. Beyond KPZ growth models, the Airy line ensemble also arises as the edge scaling limit of random surface models, such as random tilings near their arctic boundaries~\cite{ferrari2003step,johansson2005arctic,baik2007discrete,petrov2014asymptotics,duse2018universal,aggarwal2025edge}, and of interacting particle systems, such as Dyson Brownian motion~\cite{corwin2014brownian,landon2017edge,adhikari2020dyson,landon2020edgeequilibrium}.

For general $\beta>0$, the natural analogue of the Airy line ensemble is the Airy$_\beta$ line ensemble, denoted by $\mathrm{ALE}_\beta$, which reduces to the Airy line ensemble when $\beta=2$. It was first introduced in~\cite{landon2020edgeequilibrium} for $\beta\geq 1$, as the edge scaling limit of Dyson Brownian motion at equilibrium, and constructed for all $\beta>0$ in~\cite{GorinXuZhang2024Airy,charofairy}. The Airy$_\beta$ line ensemble is a stationary family of ordered continuous trajectories. At each fixed time, its particles form the Airy$_\beta$ point process, and its top particle follows the Tracy--Widom$_\beta$ distribution~\cite{edelman2007random,MR2813333}. It was characterized through explicit formulas for its multi-time Laplace transform in~\cite{GorinXuZhang2024Airy}, and as the pole evolution of random meromorphic functions satisfying certain stochastic differential equations in~\cite{charofairy}. These works also prove that $\mathrm{ALE}_\beta$ is the edge scaling limit of the G$\beta$E corners process and of Dyson Brownian motion~\cite{GorinXuZhang2024Airy}, and of stationary Dyson Brownian motion with general potentials, stationary Laguerre processes, and stationary Jacobi processes~\cite{charofairy}, for all $\beta>0$.

In this paper, we go beyond the stationary setting. We develop a unified approach to proving convergence to the Airy$_\beta$ line ensemble, based on its characterization in~\cite{charofairy} as the pole evolution of random meromorphic functions. In contrast to~\cite{charofairy}, where all the processes are stationary, the particle systems considered here start from general initial data, so their edges evolve in time. Using this approach, we show that, under regularity assumptions on the limiting family of measures, the edge scaling limits of large families of such non-stationary interacting particle systems are given by the Airy$_\beta$ line ensemble. Our results apply to interacting particle systems whose transition probabilities satisfy the dynamical loop equations introduced in~\cite{Gorin:2022qky}. Examples include $\beta$-Dyson Brownian motion with an analytic potential, weighted $\beta$-corners processes, and discrete-time Markov chains on shifted lattices related to symmetric polynomials.

\subsection{The models}\label{sec:models}

In this article, we focus on the following three families of interacting particle systems, all of which satisfy the dynamical loop equations of~\cite{Gorin:2022qky}. Detailed descriptions of these particle systems are given in \Cref{subsec:models}.

\paragraph{Discrete-time Markov chains.}
The main time-independent examples are \textit{nonintersecting Bernoulli random walks} with $\theta=1$ \cite[Section~3.1]{Gorin:2022qky}, the \textit{ascending Jack process} for any $\theta>0$, obtained as the Jack degeneration of \cite[Claim~3.10]{Gorin:2022qky}, and the \textit{ascending Macdonald process} with principal and dual specializations \cite[Claim~3.10, (61)]{Gorin:2022qky} (see also \cite{borodin2014macdonald}). In the Macdonald case, instead of fixing $q$ and $t$, we take $q_n=e^{-\gamma/n}$ and $t_n=q_n^\theta$ with $\gamma>0$, so that $q_n,t_n\to 1$ as $n\to\infty$. 

More generally, our results apply to analytic, spatially varying weights satisfying the hypotheses of our main theorem.
\Cref{sec:appendix-b-verification} also discusses several time-dependent examples: \textit{uniform lozenge tilings}, \textit{$(q,\kappa)$-distributed lozenge tilings} in the real and imaginary regimes, and lozenge tilings with \textit{trigonometric weights} \cite[Theorem~3.14]{Gorin:2022qky}, \cite{borodin2010q}; \textit{descending Jack} and \textit{descending Macdonald transitions} in particle--hole coordinates \cite[Section~2.2, Claim~3.12, (64)--(65)]{Gorin:2022qky}; and \textit{Koornwinder quasi-branching transitions} \cite[Claim~3.17]{Gorin:2022qky}. 

\paragraph{Weighted $\beta$-corners processes.}
This family consists of the \textit{weighted Dixon--Anderson transitions} of \cite[Section~3.5 and Corollary~3.6]{Gorin:2022qky}. It includes the \textit{unweighted $\beta$-corners process} ($W=0$) and the \textit{Gaussian $\beta$-corners process}, in which the top row follows the Gaussian $\beta$-ensemble \cite{MR3418747,GorinXuZhang2024Airy}. By \Cref{prop:b-corners}, it also includes \textit{quadratic}, \textit{Laguerre}, and \textit{Jacobi potentials}, under the analyticity and derivative conditions stated there and away from the poles of $W'$. Here the potential $W$ changes the conditional law of each row given the row above; the names Laguerre and Jacobi do not refer to the continuous-time Laguerre and Jacobi processes.

\paragraph{Dyson Brownian motion.}
This family includes \textit{additive Dyson Brownian motion} ($V=0$), \textit{Dyson Brownian motion with a quadratic potential} (including the stationary Gaussian, or Hermite, process), and \textit{Dyson Brownian motion with a general analytic potential}, under the hypotheses of \Cref{subsec:prelim-dbm-model,ass:weightfunc}. For $\beta=1$ (resp.\ $\beta=2$), additive Dyson Brownian motion is the eigenvalue process of a real symmetric (resp.\ complex Hermitian) matrix-valued Brownian motion, and the Gaussian case is the eigenvalue process of the corresponding matrix-valued Ornstein--Uhlenbeck process \cite{MR0148397,MR4009708,adhikari2020dyson}. Our results hold for every fixed $\beta>0$; for $\beta<1$, where particles may collide, we use the collision convention in \Cref{subsec:prelim-dbm-model}.

\paragraph{Related continuous-time chains.}
\Cref{sec:appendix-b-verification} also discusses two continuous-time models: \textit{$\theta$-nonintersecting Poisson random walks}, which reduce to ordinary nonintersecting Poisson random walks at $\theta=1$, and walks with \textit{spatially dependent Poisson jump rates} \cite[Corollary~3.2]{Gorin:2022qky}, \cite{huang2017beta}. At $\theta=1$, these walks are related to quantum random walks \cite[Section~5.2]{bufetov2019fourier}: the quantum walk associated with the defining representation is the embedded Poisson chain, while the Bernoulli transition corresponds to a mixture of exterior powers. These continuous-time chains require separate generator and martingale estimates, which we do not pursue here.

\subsection{Main results}
\label{subsec:intro-results}
In this section, we state an informal version of our dynamical edge universality result. All the models discussed in the previous section are interacting particle systems $x_1(t)\geq x_2(t)\geq \cdots \geq x_n(t)$, in discrete or continuous time, living on $\bR$ or on a discrete lattice. For corners processes, the number of particles can change with time.

In all cases, the empirical measure of the particles concentrates around a deterministic family of measures $\{\mu_s\}_{0\leq s\leq T_0}$, where $T_0<1$ for corners processes, with right edge $E_s=\sup\supp\mu_s$ (see \Cref{prop:prelim-deterministic-limit}). Our main assumptions are that the initial configuration concentrates around $\mu_0$ (\Cref{ass:edge-data}\textup{(iii)}), and that $\mu_s$ has \textit{square-root behavior} near its right edge after certain waiting time (see \Cref{ass:edge-data}).

Fix an admissible macroscopic observation time $T\in(s_*,T_0)$, where the square-root behavior applies. There exist scales $\chi\asymp n^{-2/3}$ and $\zeta\asymp n^{-1/3}$, depending on $T$, such that we can rescale the particle system as
\begin{equation}\label{eq:intro-rescaled-particles}
 \widetilde x_i^{(n)}(t)
 =\frac{x_i(T+\zeta t)-E_{T+\zeta t}}{\chi}.
\end{equation}
For discrete-time Markov chains and corners processes, where time is discrete, we extend $\widetilde x_i^{(n)}$ to all $t$ by linear interpolation. Our main result can be summarized as follows.

\begin{theorem*}[Dynamical edge universality]
Fix $\beta>0$, and consider any of the three families of models, together with related continuous-time chains in \Cref{sec:models}, under the hypotheses of \Cref{thm: main thm}. Then for every admissible observation time $T$, every fixed $\sfT>0$, and every $k\geq 1$,
\[
 (\widetilde x_1^{(n)},\ldots,\widetilde x_k^{(n)})
 \Longrightarrow
 (\mathcal A_1^\beta,\ldots,\mathcal A_k^\beta)
 \quad\text{in }C([-\sfT,\sfT],\mathbb R^k),
\]
where $(\mathcal A_i^\beta)_{i\geq 1}$ is the Airy$_\beta$ line ensemble $\mathrm{ALE}_\beta$.
\end{theorem*}

\begin{figure}[t]
    \centering
    \includegraphics[width=0.75\linewidth]{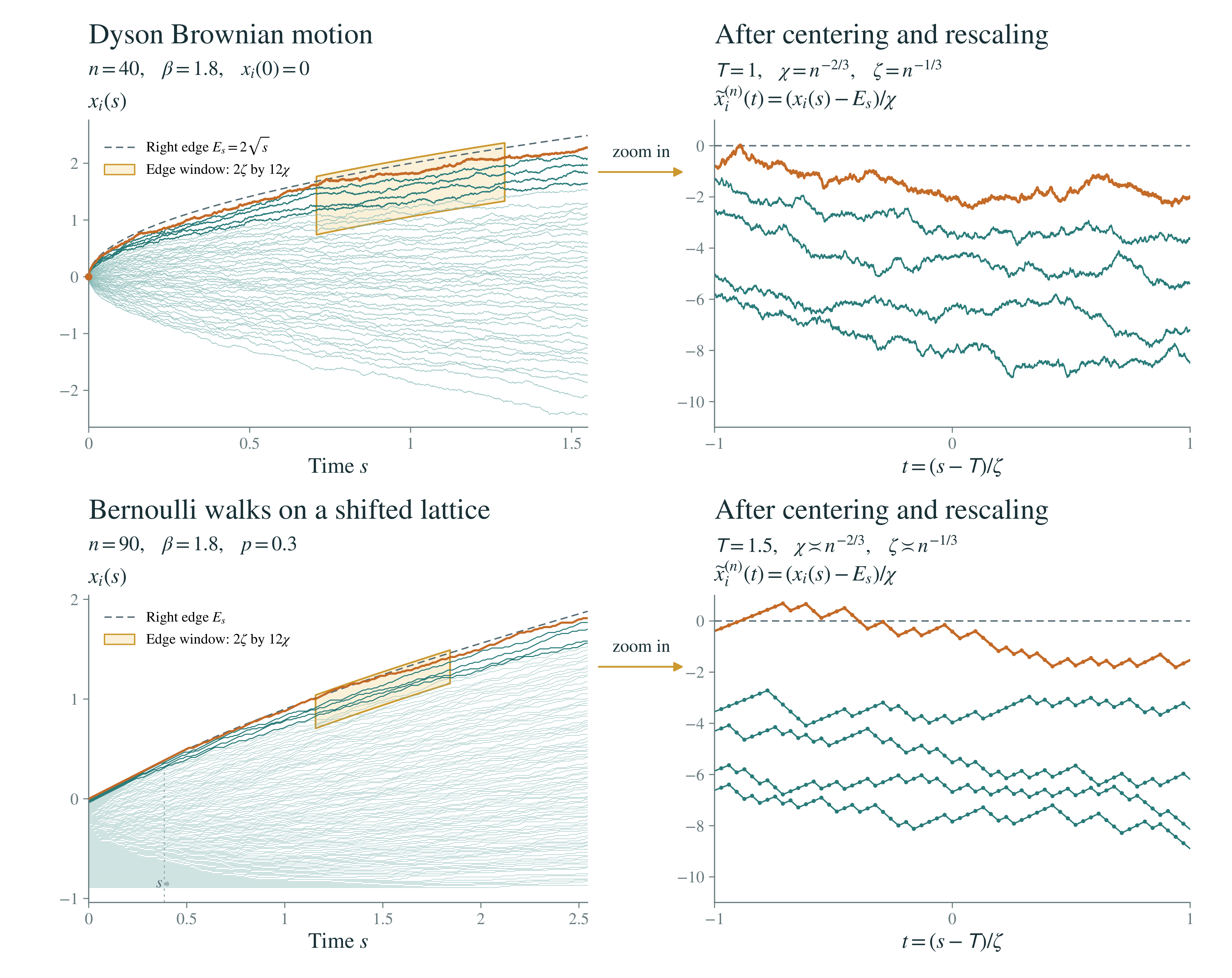}
    \caption{Edge scaling in two particle systems with $\beta=1.8$.
    The top row shows Dyson Brownian motion with $n=40$ particles
    started at zero and observed at $T=1$.
    The bottom row shows $n=90$ Bernoulli walks on a shifted lattice
    with $p=0.3$, started from $x_i(0)=-(i-1)\theta/n$ and observed
    at $T=1.5$.
    In each row, the dashed curve on the left is the deterministic
    right edge $E_s$. The right panel shows the top five trajectories
    in the highlighted window, after centering at $E_s$ and rescaling
    space by $\chi$ and time by $\zeta$.
    The dots in the lower-right panel mark the update times.}
    \label{fig:edge-scaling}
\end{figure}

In other words, the trajectories of any fixed number of top particles
converge jointly in law in the topology of uniform convergence on compact
time intervals. The result applies to a range of initial configurations,
including those out of equilibrium. \Cref{fig:edge-scaling} illustrates
the space--time scaling for Dyson Brownian motion and Bernoulli walks.

\subsection{Proof outline and main novelties}
Our proof proceeds in three steps. First, we analyze a single transition of the particle system using the dynamical loop equations. Second, we combine the resulting estimates with the characteristic flow of the deterministic evolution to control the process near its moving edge and to establish tightness. Finally, we identify the stochastic evolution satisfied by every subsequential limit, and conclude using the pole-evolution characterization of the Airy$_\beta$ line ensemble from~\cite{charofairy}.

\textit{Step 1: One-step analysis through dynamical loop equations.}
Our starting point is the dynamical loop equations of~\cite{Gorin:2022qky}, which are exact identities for conditional expectations of certain observables under a single Markov transition. In~\cite{Gorin:2022qky}, these identities were used to compute the conditional mean, mixed covariances, and higher moments of the increments of the empirical Stieltjes transform at the macroscopic scale. The key new step here is to extend this analysis to the mesoscopic and microscopic scales. For this, we need error bounds that depend explicitly on the distance between the spectral parameter and the support of the particles. These estimates, summarized in \Cref{thm:main}, remain valid for spectral parameters close enough to the support to reach the microscopic edge regime.

\textit{Step 2: Characteristic-flow analysis and tightness.}
Next, we compare the Stieltjes transform $m_s^n$ of the empirical measure with its deterministic limit $m_s$, the Stieltjes transform of $\mu_s$. The function $m_s$ satisfies a transport equation, and evaluating along its characteristic flow
\[
\partial_s z_s=\chi_s(z_s)
\]
cancels the leading transport term in the evolution of $m_s^n(z_s)-m_s(z_s)$. The square-root behavior of $\mu_s$ near its right edge gives quantitative control of these characteristics relative to the moving edge $E_s$. We use this to control the drift errors and martingale fluctuations from Step~1 as they accumulate along the characteristics. A stopping-time argument then yields an upper bound on the rightmost particle, bounds on $m_s^n-m_s$ near the edge, and bounds that hold uniformly over microscopic time windows; see \Cref{large-eigen: DMC}, \Cref{bulktight: DMC}, and \Cref{l: Delta is sml}. After centering at $E_s$ and rescaling by $n^{-2/3}$ in space and $n^{-1/3}$ in time, these bounds imply tightness of the rescaled Stieltjes transforms and of their martingale parts.

\textit{Step 3: Identification of the limiting stochastic evolution.}
Using the estimates from Steps~1 and~2, we show that the drift and martingale decomposition of the evolution of the Stieltjes transform passes to any subsequential limit $Y_t$, which satisfies
\[
dY_t(z)
=
dM_t(z)
+
\left[
\frac{2-\beta}{2\beta}\,\partial_z^2Y_t(z)
+
Y_t(z)\partial_zY_t(z)
-\frac12
\right]dt,
\]
where $\{M_t(z)\}$ is a family of martingales with quadratic covariations
\[
d\langle M(z),M(w)\rangle_t
=
\frac{2}{\beta}\,
\partial_z\partial_w
\frac{Y_t(z)-Y_t(w)}{z-w}\,dt,
\qquad z\ne w.
\]
The characteristic flow estimates from Step~2 also provide the Airy asymptotics of $Y_t$ required by the characterization in~\cite{charofairy}, which identifies the ordered poles of $Y_t$ with the Airy$_\beta$ line ensemble. Thus every subsequential limit has the same law, and together with the tightness from Step~2, this proves the convergence.

Characteristic-flow methods have been used extensively to analyze Dyson Brownian motion and related interacting particle systems. The works~\cite{adhikari2020dyson, MR4009708} combined the characteristics of the McKean--Vlasov equation with stochastic analysis to prove local laws, rigidity estimates, and mesoscopic central limit theorems in the bulk and at the edge. In the discrete setting, tools from stochastic calculus such as It\^o's formula are not available. Instead, Huang~\cite{Jiaoyang:edgestat} combined a mesoscopic analysis of dynamical loop equations with the characteristics of the complex Burgers equation to prove nearly optimal concentration estimates for nonintersecting Bernoulli walks and lozenge tilings. These works form the basis of our approach. However, they all analyze the evolution only at mesoscopic scales. The new challenge here is to control the stochastic evolution itself at the microscopic edge scale, precisely enough to prove tightness at the edge.

The first difficulty is that, to obtain the correct upper bound for the top particle in \Cref{s:eost}, we need to carefully choose a family of characteristics \eqref{e:edge char flow spec} with special end points governed by the control function $f(s)$ in \eqref{e:deff}.
The construction of the control function $f(s)$ is quite nontrivial, as it must yield the correct spatial scaling $n^{-2/3}$ within the $n^{-1/3}$ time window around $T$.
The choice of the stopping time $\tau$ in \eqref{def:stoptime edge DMC} gives a sharp bound on the error $m_s^n-m_s$ near the edge at the microscopic scale, whereas the authors of \cite{edgerigiditydysonbrownian} study only the error at the macroscopic and mesoscopic scale.
The analysis of these characteristics is divided into three parts:
the part before square-root edge behavior;
the part after square-root edge behavior, while the characteristics are still away from the edge; the part close to the edge at the microscopic scale.
The main difficulty is to connect the estimates in these three regimes to obtain \(n^{-2/3}\)-scale bounds on the rightmost particle, uniformly over \(n^{-1/3}\) time windows.

The second difficulty is that, in \Cref{sec:ale-convergence}, it is not obvious that, after suitable scaling, the subleading term $c^n_s$ in \eqref{eq:main-drift} contributes to the second derivative term in the limiting SDE in \Cref{eq:prelim-airy-martingale}.
The subleading terms $c^n_s$ for DMC and weighted $\beta$-corners have completely different expressions.
In \Cref{p:dismark scale sde}, we show that these different subleading terms converge to the same object.
This identifies the common stochastic evolution required by the characterization in~\cite{charofairy}.

Overall, our approach provides a unified framework for proving edge universality in discrete and continuous time particle systems, without requiring explicit correlation kernels or exact formulas.

\subsection{Comparison with previous results}
\label{subsec:intro-comparison}

\paragraph{Gaussian and stationary models.}
Gorin, Xu and Zhang~\cite[Theorems~1.5--1.6]{GorinXuZhang2024Airy} proved convergence to $\mathrm{ALE}_\beta$ for Gaussian $\beta$-corners processes and additive DBM started from zero, for all $\beta>0$. We recover these limits, since at any fixed positive time these models have rescaled Gaussian $\beta$-ensemble marginals; the entrance estimates follow from \cite[Theorem~7.9]{charofairy} and the deterministic regularity from \Cref{thm:app-edge-formation} and \Cref{prop:app-scalar-waiting}. For these models, our contribution is the extension to other initial data and nonconstant weights, and we do not recover their explicit Laplace transform formulas. 

Landon~\cite[Theorems~1.1--1.2]{landon2020edgeequilibrium} proved the edge limit of equilibrium Gaussian DBM for $\beta\geq1$, with quantitative estimates over growing ranks and time ranges; we recover the limit only for fixed ranks and compact time intervals. Huang and Zhang~\cite[Theorem~7.2]{charofairy} proved trajectory convergence for stationary DBM with a general potential and for stationary Laguerre and Jacobi diffusions. Under common assumptions, we extend their DBM result to non-stationary initial data. 

\paragraph{General initial data.}
Adhikari and Huang~\cite{adhikari2020dyson} proved edge rigidity, mesoscopic fluctuations, and short-time edge universality for DBM with a general potential. Aggarwal and Huang~\cite[Theorems~1.3--1.4]{edgerigiditydysonbrownian} treated additive DBM with initial data satisfying a finite-scale lower mass condition, without an initial square-root edge. Under similar assumptions, we upgrade their one-time Tracy--Widom conclusions to trajectory convergence at a fixed positive macroscopic time, but we do not recover their quantitative, $n$-dependent relaxation times. Huang and Xu~\cite{HuangXu2026Support} studied additive DBM with arbitrary initial data, and proved that for all times the support of the empirical measure stays close to the support of the free convolution of the initial empirical measure with the semicircle law.

\paragraph{Discrete ensembles and Jack processes.}
Guionnet and Huang~\cite{GuionnetHuang2019Discrete} proved rigidity and one-time edge universality for discrete $\beta$-ensembles with $\beta\geq1$. We obtain trajectory limits for all $\beta>0$, but only for ensembles that arise as marginals of our transitions. Xu~\cite[Theorems~1.1--1.2 and~4.1]{Xu2026JackAiry} proved multi-time soft-edge moment convergence for a class of Jack processes, including the Jack--Plancherel process, for all $\beta>0$, and for $\beta\geq1$ identified the equal-time limit with the Airy$_\beta$ point process through its Laplace transform. Neither result contains the other. On the one hand, \cite{Xu2026JackAiry} assumes that the power sums of the specialization are nonnegative. This excludes the dual specialization giving our Bernoulli transitions, whose power sums alternate in sign. On the other hand, the Plancherel and geometric transitions in \cite{Xu2026JackAiry} involve continuous-time or unbounded jumps. They are therefore not covered by our discrete-time Markov chains, in which each particle moves by $0$ or $1$ at each step.

\paragraph{Walks and tilings at $\beta=2$.}
Dauvergne, Nica and Vir\'ag~\cite{DauvergneNicaVirag2023Uniform} proved that nonintersecting Bernoulli walks, as well as related integrable models, converge to the Airy line ensemble uniformly on compact sets. For Bernoulli walks, our results recover this limit under our assumptions, with centering at the moving edge, but they do not cover the other jump distributions considered in \cite{DauvergneNicaVirag2023Uniform}. Denisov, FitzGerald and Wachtel~\cite{DenisovFitzGeraldWachtel2026Ordered} proved convergence to the Airy line ensemble for ordered continuous-time random walks with more general increments, when the number of particles grows slowly compared with time; this regime differs from ours. Aggarwal and Huang~\cite[Theorem~2.10]{aggarwal2025edge} proved convergence to the Airy line ensemble near any point of the arctic boundary of a polygonal lozenge tiling that is not a cusp or tangency location. This result does not follow from ours, since the tiling examples in \Cref{sec:appendix-b-verification} do not cover all such polygons.




\paragraph{Organization of the paper.}
\Cref{sec:preliminary} introduces the particle systems, the
deterministic evolution and the characterization of
$\mathrm{ALE}_\beta$. \Cref{sec:dynamical_loop} proves the quantitative
dynamical loop estimates. \Cref{s:eost} develops the
characteristic analysis and the edge bounds.
\Cref{sec:ale-convergence} proves tightness, identifies the
limiting evolution and establishes the convergence theorem.
\Cref{s:dyn-loop-beta-corners} contains the remaining one-step
calculations, \Cref{sec:appendix-b-verification} gives sufficient
conditions for the analytic assumptions, and
\Cref{app:deterministic-proofs} gives sufficient
conditions for the square-root assumptions.

\paragraph{Acknowledgements}
The research of J.H. was supported in part by NSF Grant DMS-2246664, DMS-2331096 and a Sloan Research Fellowship. The research of S.X. is supported in part by NSF Grant DMS-2331096.

\paragraph{Statement on the use of AI.}
The authors proved the math content of this paper one year before
this revision. During the revision, the authors used ChatGPT to
assist with exposition, organization, and literature review, and to
check the arguments and notation for consistency. We also used
ChatGPT to work out and check some of the computations underlying
the examples in Appendices~B and~C. The authors take full
responsibility for the mathematical content and the final text.

\section{Preliminaries and general framework}\label{sec:preliminary}

\subsection{The \texorpdfstring{Airy$_\beta$}{Airy-beta} line ensemble}
\label{subsec:prelim-overview}

Fix $\beta>0$.  The Airy$_\beta$ line ensemble, denoted by
$\mathrm{ALE}_\beta$, is the stationary soft-edge limit of Gaussian Dyson
Brownian motion~\cite[Definition~1.1 and Theorem~1.2]{charofairy}.
At each fixed time its particles have the Airy$_\beta$ point-process law,
and its top curve has the Tracy--Widom--$\beta$ distribution; for
$\beta=1,2$ these are the GOE and GUE edge laws,
respectively~\cite[Theorem~1.1]{MR2813333}.
For $\beta=2$ it is the classical Airy line
ensemble~\cite[Theorem~3.1]{corwin2014brownian}, with the time normalization
fixed below.

Write $\mathbb H$ for the open upper half-plane. A Nevanlinna function is
 a holomorphic map $\mathbb H\to\mathbb H\cup\mathbb R$; it is
\emph{particle-generated} if its Nevanlinna representing measure is a
locally finite sum of unit point masses, with multiplicities
allowed~\cite[Definitions~1.3 and~1.6]{charofairy}. We use the following
characterization.

\begin{theorem}[Huang--Zhang, {\cite[Theorem~1.15]{charofairy}}]
\label{thm:prelim-airychar}
Let $(Y_t)_{t\in\mathbb R}$ be a process of particle-generated Nevanlinna
functions, continuous in the topology of locally uniform convergence on
$\mathbb H$, and extended to the lower half-plane by Schwarz reflection.
Suppose that the following conditions hold.
\begin{enumerate}[label=\textup{(\roman*)},leftmargin=*,itemsep=3pt]
\item There is a deterministic $\mathfrak d>0$, a deterministic sequence
$t_j\to-\infty$, and a tight family of positive random variables
$(C_{*,j})_{j\ge1}$ such that all poles of $Y_{t_j}$ are at most $C_{*,j}$
and, almost surely,
\begin{equation}\label{eq:prelim-airy-normalization}
 \left|Y_{t_j}(z)-\sqrt z\right|
 \le\frac{C_{*,j}\bigl(\operatorname{Im}\sqrt z\bigr)^{1-\mathfrak d}}
          {\operatorname{Im}z},
 \qquad
 \operatorname{Im}z\ge
 C_{*,j}\sqrt{1+\max\{\operatorname{Re}z,0\}}.
\end{equation}
Here the square root is the principal branch.
\item In the integral semimartingale sense,
\begin{equation}\label{eq:prelim-airy-martingale}
 \rd Y_t(z)=\rd M_t(z)+
 \left[
 \frac{2-\beta}{2\beta}\,\partial_z^2Y_t(z)
 +Y_t(z)\partial_zY_t(z)-\frac12
 \right]\rd t,
\end{equation}
where $M$ has joint continuous complex martingale increments on every
finite time interval, relative to a common filtration, with
$M_t(\bar z)=\overline{M_t(z)}$ and quadratic variations
\begin{equation}\label{eq:prelim-airy-bracket}
 \begin{aligned}
 \frac{\rd}{\rd t}\langle M(z),M(w)\rangle_t
 &=\frac2\beta\,\partial_z\partial_w
       \frac{Y_t(z)-Y_t(w)}{z-w}, &&z\ne w,\\
 \frac{\rd}{\rd t}\langle M(z),M(z)\rangle_t
 &=\frac1{3\beta}\,\partial_z^3Y_t(z).
 \end{aligned}
\end{equation}
These identities hold for $z,w\in\mathbb C\setminus\mathbb R$;
jointness is required for every finite set of spectral parameters.
\end{enumerate}
Then the poles of $(Y_t)_{t\in\mathbb R}$, counted with multiplicity and
ordered decreasingly, form a continuous line ensemble with law
$\mathrm{ALE}_\beta$.
\end{theorem}

Writing these poles as $\mathcal A_i^\beta(t)$ and the zeros of
$\operatorname{Ai}$ as $\mathfrak a_1>\mathfrak a_2>\cdots$, the same
hypotheses give the regularized representation
\cite[Corollary~3.2]{charofairy}
\begin{equation}\label{eq:prelim-airy-resolvent}
 Y_t(z)=\sum_{i=1}^{\infty}
 \left(\frac1{\mathcal A_i^\beta(t)-z}-\frac1{\mathfrak a_i}\right)
 -\frac{\operatorname{Ai}'(0)}{\operatorname{Ai}(0)},
 \qquad z\in\mathbb C\setminus\mathbb R.
\end{equation}
The series converges locally uniformly off the real axis. To prove
convergence to $\mathrm{ALE}_\beta$, we will establish tightness and verify
that every subsequential limit satisfies the hypotheses of
\Cref{thm:prelim-airychar}.

\subsection{The three particle systems}\label{subsec:models}

Throughout, $\beta>0$ and $\theta=\beta/2>0$.  Particles are ordered
decreasingly: $\mathbb W_n=\{\bmx\in\mathbb R^n:x_1>\cdots>x_n\}$,
and $\overline{\mathbb W}_n$ is its closure.  Write
$\Delta(\bmx)=\prod_{i<j}(x_i-x_j)$ for the Vandermonde determinant.
The interaction map $\sfb$ is part of the discrete model data; for corners
and Dyson Brownian motion, $\sfb(z)=z$.

\subsubsection{A discrete-time Markov chain on a $\theta$-shifted lattice}
\label{subsec:prelim-dmc-model}

The state space of the Markov chain is the shifted lattice
\[
 \mathbb W_\theta^n=\Bigl\{\bmx\in\mathbb R^n:x_1\in n^{-1}\mathbb Z,
 \quad x_i-x_{i+1}\in n^{-1}(\theta+\mathbb Z_{\ge0}),\quad 1\le i<n\Bigr\}.
\]
From a configuration in $\mathbb W_\theta^n$, the
chain makes the one-step update $\bmy=\bmx+\bme/n$, with $\bme\in\{0,1\}^n$.
Its transition probabilities are
\begin{equation}\label{eq:disc-mc}
 \bP(\bmy=\bmx+\bme/n\mid\bmx)=a_{\bme}(\bmx),\qquad
 a_{\bme}(\bmx)=\frac1{Z(\bmx)}
 \prod_{i<j}\frac{\sfb(x_i+\theta e_i/n)-\sfb(x_j+\theta e_j/n)}
                       {\sfb(x_i)-\sfb(x_j)}
 \prod_{i=1}^n\varphi^{+}(x_i)^{e_i}\varphi^{-}(x_i)^{1-e_i} .
\end{equation}
Here $\varphi^{\pm}$ is positive on the allowed real configurations and $Z(\bmx)$
normalizes the weights.

\subsubsection{The weighted $\beta$-corners transition}
\label{subsec:prelim-corner-model}

The rows are expressed in common macroscopic spatial coordinates, with the
same spatial normalization retained throughout the descending process.
Given an upper row $\bmx\in\mathbb W_{k+1}$, the lower row
$\bmy\in\mathbb W_k$ has density
\begin{equation}\label{eq:beta-corners}
 p_k(\bmy\mid\bmx)=\frac{\mathbf1_{\{x_i>y_i>x_{i+1},\ 1\le i\le k\}}}
                              {Z_k(\bmx)}
 \Delta(\bmy)\prod_{1\le i\le j\le k}
 (y_i-x_{j+1})^{\theta-1}(x_i-y_j)^{\theta-1}
 e^{k\sum_{i=1}^kW(y_i)}.
\end{equation}
Here $W$ is a continuous real-valued potential. This is the weighted
Dixon--Anderson transition in
\cite[Section~3.5 and Corollary~3.6]{Gorin:2022qky}.
For $\theta>0$, its normalization is finite and positive: the unweighted
kernel is integrable, and $W$ is bounded on the compact interlacing region.
The analytic assumptions used below concern its asymptotics.  Iterating \eqref{eq:beta-corners}, at the transitions $k+1\to k$, defines the descending corners process.

\subsubsection{Dyson Brownian motion}\label{subsec:prelim-dbm-model}

Let $V\in C^4(\mathbb R)$ be a real-valued potential such that
$V''(x)\ge-2K$ for all $x\in\mathbb R$, for some constant $K\ge0$.
Let $B_1,\ldots,B_n$ be independent standard Brownian motions.
In the standard macroscopic normalization, the $\beta$-Dyson Brownian
motion with potential $V$ satisfies
\begin{equation}\label{eq:DBM-potential}
 \rd \lambda_i(s)=\frac{\rd B_i(s)}{\sqrt{\theta n}}
 +\frac1n\sum_{j\ne i}\frac{\rd s}{\lambda_i(s)-\lambda_j(s)}
 -\frac12V'(\lambda_i(s))\,\rd s,\qquad 1\le i\le n.
\end{equation}
For initial data in $\overline{\mathbb W}_n$, we use the unique global
ordered strong solution, allowing collisions. The singular interaction
is interpreted through its time integral; collision times have Lebesgue
measure zero~\cite[Theorem~3.1, Proposition~3.2 and Lemma~3.3]{CepaLepingle1997}. The locally Lipschitz potential drift
is handled by localization, and the lower bound on $V''$ prevents explosion;
see \Cref{app:deterministic-proofs}.
The analyticity required below is local to the spectral region.
The case $V=0$ is the unconfined DBM; for $\beta=1,2$, it gives the
eigenvalue dynamics of normalized real symmetric and complex Hermitian
matrix Brownian motion, respectively~\cite{MR0148397}.

\subsection{Empirical Measure and Stieltjes Transform}
\label{subsec:prelim-measures}
\label{subsec:prelim-analytic-setup}

Retain the notations in \Cref{subsec:models}. In this subsection, we recall the empirical Measure and Stieltjes Transform for all three models. Throughout the
rest of this section, $n$ is the particle number for DMC, DBM, and the
top-row size for corners.

\paragraph{Time normalization.}
Fix $T>0$, with $T<1$ for corners. In both discrete models, write
$\bmx_j$ for the particle configuration after $j$ transitions, at grid time
$s_j=j/n$, and $x_{j,i}$ for its $i$--th coordinate. The configuration
contains $n$ particles for DMC and $n-j$ particles for corners.
DBM uses the continuous time $s$ of \eqref{eq:DBM-potential}.

\paragraph{Empirical measures.}
For DMC, define the smoothed empirical measure at grid times by
\[
 \mu_{s_j}^n(\rd r)=\frac1\theta\sum_{i=1}^n
 \mathbf1_{[x_{j,i},x_{j,i}+\theta/n]}(r)\,\rd r.
\]
We also write this measure as $\mu(\rd r;\bmx_j)$.
For corners, use the same top-row normalization at every level:
\begin{equation}\label{eq:corner-level-measures}
 \mu_{s_j}^n(\rd r)=\frac1n\sum_{i=1}^{n-j}
                       \delta_{x_{j,i}}(\rd r).
\end{equation}
The empty terminal row has measure $\mu_1^n=0$.
For both discrete models, define $\mu_s^n$ between grid times by linear
interpolation of $\mu_{s_j}^n$ and $\mu_{s_{j+1}}^n$.
For DBM, set
\[
 \mu_s^n(\rd r)=\frac1n\sum_{i=1}^n\delta_{\lambda_i(s)}(\rd r).
\]

\paragraph{Stieltjes transforms.}
For each time $s\in[0,T]$, with $s=s_j$ for the discrete models,
let $\bm x_s$ be the particle configuration at time $s$, and define
the particle support by
\[
 I^n_s=
 \begin{cases}
 [x_{j,n},x_{j,1}+\theta/n],&\text{DMC},\\
 [x_{j,n-j},x_{j,1}],&\text{corners},\\
 [\lambda_n(s),\lambda_1(s)],&\text{DBM}.
 \end{cases}
\]

For simplicity, we write \(I_s\) whenever there is no ambiguity. 
For each particle support $I_s$, consider a complex neighborhood
$\Lambda_s\subset\bC$ of $I_s$.
For each of the three models, the Stieltjes transform is defined by
\begin{equation}\label{eq:stieltjes-general}
 m_s^n(z)=\int_{I_s}\frac{\sfb'(z)}{\sfb(r)-\sfb(z)}\,\mu_s^n(\rd r),
 \qquad z\in\Lambda_s\setminus I_s.
\end{equation}
For DMC, at the grid times $s=s_j$,
\[
 m_{s_j}^n(z)=\frac1\theta\sum_{i=1}^n
 \int_{x_{j,i}}^{x_{j,i}+\theta/n}
 \frac{\sfb'(z)}{\sfb(r)-\sfb(z)}\,\rd r.
\]
For corners, $\sfb(z)=z$ and
\[
 m_{s_j}^n(z)=\frac1n\sum_{i=1}^{n-j}\frac1{x_{j,i}-z}.
\]
By linearity of \eqref{eq:stieltjes-general}, $m_s^n$ is also the linear
interpolation of its values at successive grid times. For DBM, $\sfb(z)=z$ and
\[
 m_s^n(z)=\frac1n\sum_{i=1}^n\frac1{\lambda_i(s)-z}.
\]
For discrete models, define the discrete increments by
\begin{equation}\label{eq:corners-forward-increment}
D_j^n(z):= n\bigl(m_{s_{j+1}}^n(z)-m_{s_j}^n(z)\bigr).
\end{equation}

A function $f$ satisfies Schwarz symmetry if $f(\bar z)=\overline{f(z)}$.
We choose $\Lambda_s$ invariant under complex conjugation. 

The  three families of models considered in this paper have analytic model data: the potentials, interaction maps and transition weights admit
holomorphic extensions to $\Lambda_s$. We impose the following assumption of common analytic framework, which will be used in the following sections.

\begin{assumption}[Analyticity of the model data]\label{ass:weightfunc}
The potentials $V$ and $W$ are holomorphic on $\Lambda$ and satisfy
Schwarz symmetry. For DMC, $\sfb$ is biholomorphic and $\varphi^{\pm}$ is
holomorphic on $\Lambda=\cup_s \Lambda_s$, and both satisfy Schwarz symmetry, with
$\inf\varphi^{\pm}>c$, for some positive constant on $I=\cup_s I_s$. There is a constant $C\ge1$,
independent of $n$ and time, such that
\[
 |\sfb(z)|\le C,\qquad |\sfb'(z)|\ge C^{-1},\qquad
 |\varphi^{\pm}(z)|\le C,\qquad z\in\Lambda.
\]
\end{assumption}

Although the three models have different microscopic dynamics,
their Stieltjes transforms satisfy similarly evolution rule. We introduce a model-dependent coefficient $\cB_s^n$ and flux $\cF_s^n$ to describe the evolution of $m_{s}^n$. Here $\cF_s^n(p,z)$ is a holomorphic function
in the two complex variables $p$ and $z$
\begin{definition}\label{def: function F}
For each time $s\in[0,T]$, with $s=s_j$ for the discrete models, define the coefficient
$\cB_s^n$ and flux $\cF_s^n$ by
\begin{equation}\label{e:B_function}
 \cB_s^n(z)=\begin{cases}
 \varphi^{-}(z)+\varphi^{+}(z)e^{-\theta m_s^n(z)},&\textnormal{DMC},\\
 (1-s-1/n)W'(z)-\theta m_s^n(z),&\textnormal{corners},\\
 \exp\!\{\frac\theta2(m_s^n(z)+V'(z)/2)^2\},&\textnormal{DBM},
 \end{cases}
\end{equation}
\begin{equation}\label{e:F_CP}
 \cF_s^n(p,z)=\begin{cases}
 \log(\varphi^{-}(z)+\varphi^{+}(z)e^{-\theta p}),&\textnormal{DMC},\\
 -\theta\log((1-s-1/n)W'(z)-\theta p),&\textnormal{corners},\\
 \frac\theta2(p+V'(z)/2)^2,&\textnormal{DBM}.
 \end{cases}
\end{equation}
Define the characteristic velocity to be
\begin{equation}\label{eq:velocity-from-flux}
 \chi_s^n(z)=-\frac1\theta\partial_1\cF_s^n(m_s^n(z),z)
 =\begin{cases}
 \varphi^+(z)e^{-\theta m_s^n(z)}/\cB_s^n(z),&\textnormal{DMC},\\
 -\theta/\cB_s^n(z),&\textnormal{corners},\\
 -m_s^n(z)-V'(z)/2,&\textnormal{DBM}.
 \end{cases}
\end{equation}
\end{definition}

We then give the analytic bounds and logarithm assumption below on $\cB_s^n(z)$. This assumption will play an important role in the one-step analysis in \Cref{sec:dynamical_loop}. In particular, in \cite[Section~4.4]{Gorin:2022qky}, they give the Wiener--Hopf decomposition of $\cB_s^n$, which is a key input in the one-step analysis. The logarithm conditions depend on the model: different from DMC, for the corners process, each transition removes one particle, resulting in winding number $-1$.

\begin{assumption}[Bounds and logarithm branches]\label{ass:prelim-log-branches}
Fix $0<\delta<1$. For each time $s\in[0,T]$, with $s=s_j$ for the
discrete models, let $\mathscr D_s\Subset\Lambda_s\setminus I_s$
be an annulus surrounding $I_s$. We require
\[
 \dist(\mathscr D_s,I_s)\ge n^{-1+\delta}.
\]
There are constants $0<c\le C<\infty$, which may depend on $T$ and
the chosen domains $\Lambda,I$ but are independent of $s$ and $n$, such that
\[
 c\le|\cB_s^n(z)|\le C,
 \qquad z\in\mathscr D_s,
\]
at all such times.

For the discrete models, we additionally impose the following conditions
at every grid time $s=s_j$:
\begin{enumerate}[label=\textup{(\roman*)},leftmargin=*,itemsep=2pt]
\item For DMC, $\cB_s^n$ has a single-valued holomorphic logarithm on
$\Lambda_s\setminus I_s$.
\item For corners, $\operatorname{wind}_{\omega_-}\cB_s^n=-1$, where
$\omega_-\Subset\mathscr D_s$ is a positively oriented contour enclosing $I_s$.
Equivalently, for every $e\in I_s$, the function
$(z-e)\cB_s^n(z)$ has a single-valued holomorphic logarithm on
$\mathscr D_s$.
\end{enumerate}
For DBM, no additional condition on logarithm branches is imposed.
\end{assumption}

As shown in
Appendix~B, the assumptions hold for broad classes of transition weights
and potentials, and in particular for the discrete-time
chains, weighted $\beta$-corners processes and Dyson Brownian
motion considered here. We also verify them for the concrete
examples introduced in the introduction, including the
ascending and descending Jack and Macdonald processes,
lozenge tilings, and weighted corners with quadratic,
Laguerre and Jacobi potentials.

\subsection{Stieltjes Transform evolution Away from the Support}
\label{subsec:prelim-fixed-inputs}

We first record the evolution equations on spectral domains at a fixed
distance from the particle support.
For the discrete models, let $\mathscr F_s=\mathscr F_j$ at $s=s_j$
be the sigma-field generated by the initial configuration and the first
$j$ transitions. Recall $D_j^n$ from
\eqref{eq:corners-forward-increment}, and set
\[
 \Delta\cM_j^n(z)
 =D_j^n(z)-\bE[D_j^n(z)\mid\mathscr F_j].
\]
We use the kernel
\begin{equation}\label{eq:kernel-common}
 K_{\sfb}(w,z)
 =\frac{\sfb'(w)\sfb'(z)}{(\sfb(w)-\sfb(z))^2}.
\end{equation}

Throughout this subsection, $\omega_-\Subset\mathscr D_s$ is positively
oriented, encloses $I_s$, and excludes the spectral parameters.
We work on fixed compact spectral sets in $\Lambda_s$, with $\omega_-$
a fixed positive distance from $I_s$ uniformly in
$s$ and $n$. 

\begin{proposition}[Stieltjes-transform evolution]
\label{thm:loopstudy}
Assume \Cref{ass:weightfunc,ass:prelim-log-branches} and the contour
conditions above. Recall $\chi_s^n$ from
\eqref{eq:velocity-from-flux}.
The following statements hold uniformly in time and in the spectral
parameters.
\begin{enumerate}[label=\textup{(\roman*)},leftmargin=*,itemsep=4pt]
\item \textnormal{Discrete models.}
For DMC and corners, at every grid time $s_j\le T$,
\begin{equation}\label{e:dmg}
 D_j^n(z)
 =-\frac1{2\pi\ri\theta}\oint_{\omega_-}
 \partial_w[\cF_{s_j}^n(m_{s_j}^n(w),w)]
 \frac{\sfb'(z)}{\sfb(w)-\sfb(z)}\,\rd w
 +\Delta\cM_j^n(z)+\OO(n^{-1}),
\end{equation}
and
\begin{equation}\label{e:covT}
 n\bE[\Delta\cM_j^n(z)\Delta\cM_j^n(u)\mid\mathscr F_j]
 =\frac1{2\pi\ri\theta}\oint_{\omega_-}
 \chi_{s_j}^n(w)K_{\sfb}(w,z)K_{\sfb}(w,u)\,\rd w
 +\OO(n^{-1}).
\end{equation}

\item \textnormal{Dyson Brownian motion.}
The Stieltjes transform satisfies the exact identity
\begin{equation}\label{eq:pote_mt}
 \rd m_s^n(z)
 =-\frac1{2\pi\ri\theta}\oint_{\omega_-}
 \partial_w[\cF_s^n(m_s^n(w),w)]
 \frac{\rd w}{w-z}\,\rd s
 +\rd \cM_s^n(z)
 +\frac{1-\theta}{2\theta n}\partial_z^2m_s^n(z)\,\rd s,
\end{equation}
where
\[
 \rd \cM_s^n(z)
 =-\frac1{\sqrt{\theta n^3}}
 \sum_{i=1}^n\frac{\rd B_i(s)}{(\lambda_i(s)-z)^2}.
\]
Its quadratic variation is
\begin{equation}\label{eq:dbm-bracket}
 n^2\frac{\rd}{\rd s}
 \langle\cM^n(z),\cM^n(u)\rangle_s
 =\frac1\theta\int_{I_s}
 \frac{\mu_s^n(\rd r)}{(r-z)^2(r-u)^2}
 =\frac1{2\pi\ri\theta}\oint_{\omega_-}
 \frac{\chi_s^n(w)}{(w-z)^2(w-u)^2}\,\rd w.
\end{equation}
\end{enumerate}
\end{proposition}

\begin{proof}
We defer the proof to \Cref{s:dyn-loop-beta-corners}.
\end{proof}

\subsection{Deterministic limits and characteristic flow}
\label{subsec:prelim-characteristics}

We next identify the deterministic evolution obtained from
\Cref{thm:loopstudy}.
For a family of deterministic measure $(\mu_s)_s$ supported in $I$, write
\[
 m_s(z)=\int_I
 \frac{\sfb'(z)}{\sfb(r)-\sfb(z)}\,\mu_s(\rd r).
\]
Define $\cB_s$, $\cF_s$ and $\chi_s$ by
\Cref{def: function F}, replacing $m_s^n$ by $m_s$ and setting
$n^{-1}=0$. In particular, for corners,
$\cB_s(z)=(1-s)W'(z)-\theta m_s(z)$.

\begin{proposition}[Deterministic evolution]
\label{prop:prelim-deterministic-limit}
Assume $\operatorname{supp}\mu_0^n\subset I$ and
$\mu_0^n\Rightarrow\mu_0$ weakly in probability, where $I$ is compact
and $\mu_0$ is a deterministic probability measure.
For the discrete models, take $I$ to contain the empirical supports
throughout $[0,T]$, and use common fixed working domains and contours
satisfying the conditions of \Cref{thm:loopstudy}.

Under \Cref{ass:weightfunc,ass:prelim-log-branches}, there is a
deterministic continuous path of finite positive measures $(\mu_s)$
such that, for every bounded continuous function $f$ on $\mathbb R$,
\[
 \sup_{0\le s\le T}
 \left|\int f\,\rd \mu_s^n-\int f\,\rd \mu_s\right|
 \xrightarrow[n\to\infty]{\mathbb P}0.
\]
The limiting measure has mass $1$ for DMC and DBM, and mass $1-s$
for corners.

On any time interval on which $\operatorname{supp}\mu_s\subset I$,
its Stieltjes transform satisfies
\begin{equation}\label{e:burger}
 \begin{aligned}
 -\partial_s m_s(z)
 &=\frac1{2\pi\ri\theta}\oint_{\omega_-}
 \partial_w[\cF_s(m_s(w),w)]
 \frac{\sfb'(z)}{\sfb(w)-\sfb(z)}\,\rd w,\\
 m_0(z)
 &=\int_I\frac{\sfb'(z)}{\sfb(r)-\sfb(z)}\,\mu_0(\rd r),
 \end{aligned}
\end{equation}
for $z\in\Lambda$ outside $\omega_-$.
This equation uniquely determines the limiting path among continuous
positive-measure paths supported in $I$ and satisfying the limiting
conditions of \Cref{ass:prelim-log-branches}.
\end{proposition}
\begin{proof}
We defer the proof to Appendix~C.
\end{proof}

For DBM with potential, define
\begin{equation}\label{eq:def_g_func}
 g(z,r)
 =\frac{V'(r)-V'(z)-(r-z)V''(z)}{2(r-z)^2},
 \qquad g(z,z)=\frac14V'''(z).
\end{equation}
Then \eqref{e:burger} becomes
\begin{equation}\label{equ:MVequ}
 \partial_s m_s(z)
 =m_s'(z)\bigl(m_s(z)+V'(z)/2\bigr)
 +\frac12m_s(z)V''(z)
 +\int_I g(z,r)\,\mu_s(\rd r).
\end{equation}
In particular, for DBM with $V=0$,
$\partial_s m_s=m_sm_s'$ and $\chi_s=-m_s$.
For corners with $W=0$,
$\partial_s m_s=-m_s'/m_s$ and $\chi_s=1/m_s$.

\paragraph{Characteristic flow.}
We follow the velocity $\chi_s$ to cancel the spatial derivative in
\eqref{e:burger}, as in
\cite[Section~2, (2.25)--(2.26)]{MR4009708}.
Fix positively oriented contours $\omega_-$ and $\omega_+$ in a common
working annulus, with $\omega_+$ enclosing $\omega_-$, and assume
$\cB_s$ is nonzero between them.

For $z$ between $\omega_-$ and $\omega_+$, applying Cauchy's theorem, together
with the chain rule, yield
\begin{equation}\label{eq:flux-chain-rule}
 \partial_z[\cF_s(m_s(z),z)]
 =-\theta\chi_s(z)\partial_zm_s(z)
 +\partial_2\cF_s(m_s(z),z).
\end{equation}
Consequently,
\begin{equation}\label{e:burger2}
 \partial_s m_s(z)+\chi_s(z)\partial_zm_s(z)=F_s(z),
\end{equation}
where
\begin{equation}\label{def of F_s}
 F_s(z)
 =\frac1\theta\partial_2\cF_s(m_s(z),z)
 -\frac1{2\pi\ri\theta}\oint_{\omega_+}
 \partial_w[\cF_s(m_s(w),w)]
 \frac{\sfb'(z)}{\sfb(w)-\sfb(z)}\,\rd w.
\end{equation}

To solve the ODE \eqref{e:burger2}, we use the
method of characteristics flow. 
In \Cref{s:eost}, we estimate $m_s^n-m_s$ along the characteristics flow to give a finer estimate which help to establish tightness and Airy condition in \Cref{thm:prelim-airychar}. We introduce the characteristic
flow as follows.

\begin{definition}[Characteristic flow]
\label{def:prelim-characteristic-flow}
For $u$ between $\omega_-$ and $\omega_+$, the characteristic flow is
the unique solution of
\begin{equation}\label{e:stij characteristic flow}
 \partial_s z_s(u)=\chi_s(z_s(u)),\qquad z_0(u)=u,
\end{equation}
up to its first exit from the region between the contours.
Along this flow,
\[
 \frac{\rd}{\rd s}m_s(z_s(u))=F_s(z_s(u)).
\]
\end{definition}

\section{Dynamical loop equations}
\label{sec:dynamical_loop}
The goal of this section is to establish the enhanced one-step estimates for DMC and
weighted $\beta$-corners, uniformly over all grid times in $[0,T]$. The estimate holds in the mesoscopic and microscopic scales.
We extend \Cref{thm:loopstudy} from fixed-distance spectral domains to domains
approaching the particle support $I_s$ at any time $s$, with quantitative bounds
on the drift error, the covariance error and higher moments of the
martingale increment. These estimates provide the drift and fluctuation
bounds needed for the studying of edge-scaling limit and build tightness in \Cref{sec:ale-convergence}.

\subsection{Main one-step estimates}\label{sec:main}

For any grid time $0\le s=s_j=j/n\le T$, with $T<1$ for corners, recall
\begin{equation}\label{eq:increment-common}
 \begin{aligned}
 D_s^n(z)&=n\bigl(m_{s+1/n}^n(z)-m_s^n(z)\bigr),\\
 \Delta\mathcal M_s^n(z)&=D_s^n(z)-\E_sD_s^n(z),
 \qquad \E_s=\E[\,\cdot\mid\mathscr F_s].
 \end{aligned}
\end{equation}
Recall that $\mathscr D_s\Subset\Lambda_s\setminus I_s$ is an annulus
surrounding $I_s$ and satisfies
\(
 \dist(\mathscr D_s,I_s)\ge n^{-1+\delta}
\)
for a fixed constant $0<\delta<1$.

Recall the kernel $K_{\mathsf b}(w,z)$ from \eqref{eq:kernel-common}
and the characteristic velocity $\chi_s^n$ from
\eqref{eq:velocity-from-flux}. Denote the leading drift in
\eqref{e:dmg} by
\begin{equation}\label{eq:leading-common}
 a_s^n(z)=-\frac1{2\pi\ri\theta}\oint_{\omega_{-}}
 \partial_w[\cF_s^n(m_s^n(w),w)]
 \frac{\mathsf b'(z)}{\mathsf b(w)-\mathsf b(z)}\,\rd w.
\end{equation}
Denote
\begin{equation}\label{eq:main-P}
 P_s^n(z)=1+\frac{|\Im\chi_s^n(z)|}{|\Im z|},
 \qquad z\in\D_s\setminus\mathbb R.
\end{equation}
At real points, use the continuous extension
\[
 P_s^n(x)=1+|(\chi_s^n)'(x)|,
 \qquad x\in\D_s\cap\mathbb R.
\]

\begin{theorem}[One-step estimates]\label{thm:main}
Assume \Cref{ass:weightfunc,ass:prelim-log-branches} and the domain
and contour conditions above. The following estimates hold uniformly
over all grid times $s=s_j\le T$.
For every $z\in\D_s\setminus\mathbb R$, the increment admits the expansion
\begin{equation}\label{eq:main-drift}
 D_s^n(z)=a_s^n(z)+\frac1n c_s^n(z)
          +\Delta\mathcal M_s^n(z)+\mathcal R_s^n(z).
\end{equation}
For every fixed integer $p\ge0$, the remainder satisfies
\begin{equation}\label{eq:main-second}
 |\partial_z^p\mathcal R_s^n(z)|\le
 \frac{C P_s^n(z)}{n^2\dist(z,I_s)^{p+2}}.
\end{equation}
For $z,u\in\D_s$ with
$\dist(z,I_s),\dist(u,I_s)\asymp d$, the conditional covariance satisfies
\begin{equation}\label{eq:main-cov}
 \begin{aligned}
 \E_s[\Delta\mathcal M_s^n(z)\Delta\mathcal M_s^n(u)]
 ={}&\frac1{2\pi\ri\theta n}\oint_{\omega_{-}}
 \chi_s^n(w)K_{\mathsf b}(w,z)K_{\mathsf b}(w,u)\,\rd w+\OO\!\left(\frac{P_s^n(z)+P_s^n(u)}{n^2d^3}\right),
 \end{aligned}
\end{equation}
where $\omega_{-}$ doesn't enclose $z,u$.
The same formula holds with $u$ replaced by $\bar z$. For every fixed integer $p\ge1$,
\begin{equation}\label{eq:main-moment}
 \E_s|\Delta\mathcal M_s^n(z)|^{2p}\le C \left[
 \left(\frac{P_s^n(z)}{n\dist(z,I_s)^2}\right)^p
 +\frac{P_s^n(z)}{n^{2p-1}\dist(z,I_s)^{4p-2}}\right].
\end{equation}
All the constants may depend on $\delta$, $T$, but are independent of $n$, $s$.
\end{theorem}

\begin{remark}
The second-order coefficient in \eqref{eq:main-drift} is
\[
 c_s^n(z)=
 \begin{cases}
 -\displaystyle\frac1{2\pi\ri\theta}\oint_{\omega_-}
 \cE_2(w,z)K_{\sfb}(w,z)\,\rd w,
 &\textnormal{DMC},\\[8pt]
 \displaystyle\frac1{2\pi\ri\theta}\oint_{\omega_-}
 \frac{\chi_s^n(w)a_s^n(w)-\frac12\partial_w\chi_s^n(w)}
      {(w-z)^2}\,\rd w,
 &\textnormal{corners}.
 \end{cases}
\]
Here $\cE_2$ is given by \eqref{err:secondder}, evaluated at time
$s$ with the time index suppressed. The contour $\omega_-$ encloses
$I_s$ but not $z$. For corners, the contour in
\eqref{eq:leading-common} defining $a_s^n(w)$ lies strictly inside
$\omega_-$.
These formulas are established in \eqref{t:loopstudy2} and
\eqref{eq:corners-second-order-mean}, respectively.
\end{remark}

\subsection{Discrete Markov chain}\label{s:DMC bulkedgeest}

In the part, we establish \Cref{thm:main} for the discrete Markov
chain, including the drift approximation, conditional covariance,
and moment bounds. In the rest of the sections, we always fix the grid time $s=s_j=j/n$. We retain the notation introduced in
Sections~2 and~3.1 and, for simplicity, write $I=I_s$,
$\mathscr D=\mathscr D_s$, $m=m_s^n$, $\cB=\mathcal B_s^n$,
and $\chi=\chi_s^n$. Without loss of generality, we normalize $\varphi^-=1$ by replacing
$\varphi^+$ with $\varphi:=\varphi^+/\varphi^-$ and absorbing
$\prod_{i=1}^n\varphi^-(x_i)$ into the normalizing constant.

\subsubsection{Deformation and preliminary estimates}
\label{s:Bulkdiscrete}
We define
\begin{equation}\label{e:A_def}
\cA(z):=\bE\left[
\prod_{i=1}^n
\frac{\sfb(z)-\sfb(x_i+\theta e_i/n)}
     {\sfb(z)-\sfb(x_i)}
\right],
\quad
\cG(z):=e^{-\theta m(z)}.
\end{equation}
Recall $\cB$ from \eqref{e:B_function}, the functions $\cA$, $\cG$,
and $\cB$ are holomorphic on $\Lambda\setminus I$.
For $z\in\Lambda\setminus I$,
we denote
\begin{align}\label{e:C-factor}
\cC(z)
&=\bE\Biggl[
\prod_{i=1}^n
\frac{\sfb(z)-\sfb(x_i+\theta e_i/n)}
     {\sfb(z)-\sfb(x_i)}
\Biggl(
\varphi(z)\prod_{i=1}^n
\frac{\sfb(z+\theta/n)-\sfb(x_i+\theta e_i/n)}
     {\sfb(z)-\sfb(x_i+\theta e_i/n)}
+1
\Biggr)
\Biggr].
\end{align}
By \eqref{e:A_def}, the first summand in \eqref{e:C-factor} becomes
\begin{align}
\begin{split}
&\phantom{{}={}}\bE\left[
\prod_{i=1}^n
\frac{\sfb(z)-\sfb(x_i+\theta e_i/n)}
     {\sfb(z)-\sfb(x_i)}
\cdot
\varphi(z)\prod_{i=1}^n
\frac{\sfb(z+\theta/n)-\sfb(x_i+\theta e_i/n)}
     {\sfb(z)-\sfb(x_i+\theta e_i/n)}
\right]\\
&=
\varphi(z)\,\cA(z+\theta/n)\prod_{i=1}^n
\frac{\sfb(z+\theta/n)-\sfb(x_i)}
     {\sfb(z)-\sfb(x_i)}.
\end{split}
\end{align}
Therefore,
\begin{equation}\label{e:C-AG}
\cC(z)=\cA(z)
+\varphi(z)\,\cA(z+\theta/n)\prod_{i=1}^n
\frac{\sfb(z+\theta/n)-\sfb(x_i)}
     {\sfb(z)-\sfb(x_i)}.
\end{equation}
Define
\begin{equation}\label{e:hz}
h_-(z):=\frac{1}{2\pi\ri}\oint_{\cin}
\frac{\log\cB(w)\sfb'(w)}{\sfb(w)-z}\,\rd w,
\quad
h_+(z):=\frac{1}{2\pi\ri}\oint_{\cout}
\frac{\log\cB(w)\sfb'(w)}{\sfb(w)-z}\,\rd w,
\end{equation}
where the contours $\cout$ enclosed $\sfb^{-1}(z)$, while $\cin$ not.
Since $h_+$ is holomorphic inside $\sfb(\cout)$, it satisfies
$|\partial_z^kh_+(z)|\le C$ for every fixed $k\ge0$.
$h_-$ is holomorphic outside $\sfb(\cin)$ and satisfies
$h_-(u)=\OO(u^{-1})$, hence $e^{h_-(u)}\to1$, as $u\to\infty$.
\medskip

\paragraph{Deformation of the transition weights}
Fix a positive integer $R\ge1$.
We introduce the deformed transition weights $a^{(\bmv,\bms)}_{\bme}$
and the deformed weight function $\varphi^{(\bmv,\bms)}$ by
\begin{align}\label{e:phi-tilde-def}
&\varphi^{(\bmv,\bms)}(z)
:=\varphi(z)\prod_{j=1}^r
\left(
\frac{\sfb(v_j)-\sfb(z)}
     {\sfb(v_j)-\sfb(z+\theta/n)}
\right)^{s_j},\\
&a^{(\bmv,\bms)}_{\bme}(\bmx)
=\frac{1}{Z^{(\bmv,\bms)}(\bmx)}
\prod_{1\le i<j\le n}
\frac{\sfb(x_i+\theta e_i/n)-\sfb(x_j+\theta e_j/n)}
     {\sfb(x_i)-\sfb(x_j)}
\prod_{i=1}^n\varphi^{(\bmv,\bms)}(x_i)^{e_i}.
\end{align}
for either
\begin{align}\label{e:term1}
0\le r\le2R-1,\quad
\bmv=(v_1,v_2,\ldots,v_r),\quad
\bms=(s_1,s_2,\ldots,s_r),\quad
v_j\in\mathscr D,\quad s_j=\pm1,
\end{align}
or
\begin{align}\label{e:term2a}
r=2R,\quad
\bmv=(v_1,v_2,\ldots,v_r),\quad
\bms=(s_1,s_2,\ldots,s_r),\quad
v_{2j-1}=\overline v_{2j}\in\mathscr D,\quad s_j=1.
\end{align}
Let $\mathbb E^{(\bmv,\bms)}$ denote the conditional expectation
under the deformed transition weights $a^{(\bmv,\bms)}_{\bme}(\bmx)$.
By \Cref{def: function F}
\begin{align}\label{e: def tilde A}
\cA^{(\bmv,\bms)}(z)
&:=\mathbb E^{(\bmv,\bms)}\left[
\prod_{i=1}^n
\frac{\sfb(z)-\sfb(x_i+\theta e_i/n)}
     {\sfb(z)-\sfb(x_i)}
\right],\notag\\
\cB^{(\bmv,\bms)}(z)
&:=1+\varphi^{(\bmv,\bms)}(z)\cG(z),\notag\\
\chi^{(\bmv,\bms)}(z)
&:=\frac{\varphi^{(\bmv,\bms)}(z)\cG(z)}
        {\varphi^{(\bmv,\bms)}(z)\cG(z)+1}.
\end{align}
We recall the following fixed-distance estimates from
\cite{Gorin:2022qky}.

\begin{proposition}[\cite{Gorin:2022qky}, Propositions~4.11 and~4.13]
Fix $\mathfrak r>0$ independently of $n$.
Let $z\in\mathscr D$ and $\bmv\in\mathscr D^r$ satisfy
\eqref{e:term1} or \eqref{e:term2a}, and assume
\[
\dist(\{z,v_1,\ldots,v_r\},I)\ge2\mathfrak r.
\]
Then
\begin{equation}\label{e:initial bootstrap}
\bigl|\cA^{(\bmv,\bms)}(z)e^{-h_-(\sfb(z))}\bigr|
=e^{\OO(1/n)}.
\end{equation}
The implicit constant may depend on $\mathfrak r$ and $R$,
but is independent of $n$ and the spectral parameters.
Moreover,
\begin{equation}\label{e: bulk error}
\log\cA^{(\bmv,\bms)}(z)
=h_-(\sfb(z))
+\frac{1}{2\pi\ri}\oint_{\cin}
\frac{
\log\!\left(1+\frac1n\cE_1^{(\bmv,\bms)}(w)\right)\sfb'(w)}
{\sfb(w)-\sfb(z)}\,\rd w,
\end{equation}
where $\cin\subset\mathscr D$ encloses $I$ but none of
$z,v_1,\ldots,v_r$, and stays a fixed positive distance
from $I$ and these parameters. Here
\begin{equation}\label{e:error-term-1-clean: deform}
\frac1n\cE_1^{(\bmv,\bms)}(z)
=
\frac{\cA^{(\bmv,\bms)}(z+\theta/n)}
     {\cA^{(\bmv,\bms)}(z)}
\cdot
\frac{
\displaystyle\varphi^{(\bmv,\bms)}(z)
\prod_{i=1}^n
\frac{\sfb(z+\theta/n)-\sfb(x_i)}
     {\sfb(z)-\sfb(x_i)}
}{\cB(z)}
-\chi(z).
\end{equation}
\end{proposition}

Let $\cC^{(\bmv,\bms)}$ be defined by \eqref{e:C-AG}
with $\cA$ and $\varphi$ replaced by their deformed versions.
The same pole cancellation as in
\cite[Theorem~1.1]{Gorin:2022qky} shows that
$\cC^{(\bmv,\bms)}$ is holomorphic inside $\cin$,
provided that $\cin$ excludes all $v_j$ and $v_j-\theta/n$.
By \eqref{e:error-term-1-clean: deform},
\begin{equation}\label{e:C-factorization-prop}
\cC^{(\bmv,\bms)}(z)
=\cA^{(\bmv,\bms)}(z)\cB(z)
\left(1+\frac1n\cE_1^{(\bmv,\bms)}(z)\right).
\end{equation}

By \Cref{l: derivative bdd}, the winding number of
$\cA^{(\bmv,\bms)}$ remains zero under continuation
from the fixed-distance regime.
For $nd$ sufficiently large, later on \eqref{e:cEbound} gives
$|\cE_1^{(\bmv,\bms)}/n|<1/2$ on $\cin$.
Since $\cB$ admits a holomorphic logarithm,
\eqref{e:C-factorization-prop} yields
\[
\operatorname{wind}_{\cin}\cC^{(\bmv,\bms)}
=
\operatorname{wind}_{\cin}\cA^{(\bmv,\bms)}
+\operatorname{wind}_{\cin}\cB
=0.
\]
Thus $\cC^{(\bmv,\bms)}$ has no zeros inside $\cin$
by the argument principle.
Taking logarithms with the continued branches and applying
Cauchy's formula gives \eqref{e: bulk error}.

We next record estimates for the Stieltjes transform.

\begin{lemma}\label{p:DMC prelimestimate}
Assume \Cref{ass:weightfunc,ass:prelim-log-branches}.
Recall $m$ and $\chi$ from \eqref{eq:stieltjes-general}
and \eqref{eq:velocity-from-flux}.
For every $z\in\mathscr D$,
\begin{equation}\label{e:prelimestimate}
\Re m(z)\ge-C,\qquad
|\Im m(z)|,\;|\chi(z)|\le C.
\end{equation}
\end{lemma}
\begin{proof}
Taylor expansion gives uniformly
\[
 \frac{\sfb'(z)}{\sfb(x)-\sfb(z)}
 =\frac1{x-z}+\OO(1).
\]
Thus, writing $z=u+\ri v$ with $v\ne0$,
\[
 |\Im m(z)|
 \lesssim 1+\int_I\frac{|v|}{(x-u)^2+v^2}\,\mu(\rd x;\bmx)
 \lesssim 1+\int_{\mathbb R}
              \frac{|v|}{(x-u)^2+v^2}\,\rd x
 \lesssim 1.
\]
For real $z\notin I$, Schwarz symmetry gives $\Im m(z)=0$.

By \Cref{ass:prelim-log-branches}, we have $|\cB(z)|\asymp 1$ uniformly for $z\in\mathscr D$. In particular,
\(
\Bigl|\varphi(z)e^{-\theta m(z)}\Bigr|\lesssim 1,
\)
and hence $\Re m(z) \geq  -C$ on $\mathscr D$. Moreover,
\[
|\chi(z)|
=\left|\frac{\varphi(z)\cG(z)}{\cB(z)}\right|
=\left|1-\frac{1}{\cB(z)}\right|
\lesssim 1,
\]
again uniformly for $z\in\mathscr D$.
\end{proof}

\begin{lemma}\label{lem:stieltjes-bounds}
Under \Cref{ass:weightfunc} and \Cref{ass:prelim-log-branches}, let
$z\in\mathscr D\setminus\mathbb R$, with
$d=\dist(z,I)$ .
For every fixed integer $p\ge1$,
\begin{equation}\label{eq:stieltjes-higher-derivative-bounds}
 |\partial_z^p m(z)|
 \lesssim_p
 \min\left\{
 \frac1{d^p},\
 \frac1{d^{p-1}}
 \left(1+\frac{|\Im m(z)|}{|\Im z|}\right)
 \right\}.
\end{equation}
The first bound also holds for real $z\notin I$.
\end{lemma}

\begin{proof}
Write
\[
 m_0(z)=\int_I\frac{\mu(\rd s;\bmx)}{s-z},
 \qquad q(z)=m(z)-m_0(z).
\]
Since $q$ extends
holomorphically across $I$ and is uniformly bounded on a fixed
neighborhood. One has
\[
 |q^{(p)}(z)|\lesssim 1,
 \qquad |\Im q(z)|\lesssim|\Im z|.
\]
Therefore, it yields
\begin{equation}\label{eq:mp-basic}
 \begin{aligned}
 |\partial_z^p m(z)|
 &\lesssim_p
 1+\int_I\frac{\mu(\rd s;\bmx)}{|s-z|^{p+1}}\le
 1+\frac1{d^{p-1}}
       \int_I\frac{\mu(\rd s;\bmx)}{|s-z|^2}\\
 &=1+\frac1{d^{p-1}}
       \frac{|\Im m_0(z)|}{|\Im z|} \lesssim_p
 \frac1{d^{p-1}}
 \left(1+\frac{|\Im m(z)|}{|\Im z|}\right).
 \end{aligned}
\end{equation}
On the other hand, Let \(u\) be the point of \(I\) nearest to \(z\)., 
then $|s-z|^2\ge d^2+(s-u)^2$ for $s\in I$.
We obtain
\[
 |\partial_z^p m(z)|
 \lesssim_p
 1+\int_I\frac{\rd s}{|s-z|^{p+1}}
 \le
 1+\int_{\mathbb R}
       \frac{\rd t}{(d^2+t^2)^{(p+1)/2}}
 \lesssim_p d^{-p}.
\]
This argument also applies to real $z\notin I$.
\end{proof}

\begin{lemma}\label{lem:dmc-velocity-derivatives}
Under \Cref{ass:weightfunc} and \Cref{ass:prelim-log-branches}, recall $P$ from
\eqref{eq:main-P}. For $z\in\mathscr D\setminus\mathbb R$, denote
$d=\dist(z,I)$. Then, for every fixed integer $p\ge1$,
\begin{equation}\label{eq:dmc-velocity-derivatives}
 P(z)\lesssim\frac1d,\qquad
 |\partial_z^p\chi(z)|+|\partial_z^p\log\cB(z)|
 \lesssim_p\frac{P(z)}{d^{p-1}}.
\end{equation}
Moreover, $P(w)\asymp P(z)$ whenever
$w,z\in\mathscr D\setminus\mathbb R$ and $|w-z|\le d/2$.
\end{lemma}

\begin{proof}
Write $m=m_0+q$ as in the preceding proof and set
$t_i=x_i+\theta/n$. One has
\[
 e^{-\theta m_0(z)}
 =\prod_{i=1}^n\frac{z-x_i}{z-t_i}
 =1+\sum_{i=1}^n\frac{r_i}{z-t_i},
 \qquad
 r_i=\frac{\theta}{n}
      \prod_{j\ne i}\frac{t_i-x_j}{t_i-t_j}\ge0.
\]
The nonnegativity follows from interlacing, and comparison at infinity
gives $\sum_i r_i=\theta$.

Set $f=\varphi e^{-\theta q}$. Then $f$ is Schwarz symmetric,
$f(t_i)\ge0$, and its derivatives are uniformly bounded. Thus
\begin{equation}\label{eq:dmc-velocity-positive-representation}
 \cB(z)=1+a(z)+\sum_i\frac{\nu_i}{z-t_i},
 \qquad
 \nu_i=r_i f(t_i)\ge0,
\end{equation}
where
\[
 a(z)=f(z)+\sum_i r_i\frac{f(z)-f(t_i)}{z-t_i}.
\]
Since $a(z)$ is holomorphic, $\sum_i r_i=\theta$, Cauchy's estimates and Schwarz symmetry give
$|a^{(k)}(z)|\lesssim_k1$ and $|\Im a(z)|\lesssim|\Im z|$.

Put $S(z)=\sum_i\nu_i/|z-t_i|^2$. Since $\chi=1-\cB^{-1}$,
taking imaginary parts in
\eqref{eq:dmc-velocity-positive-representation} yields
\[
 P(z)
 =1+\frac{|\Im\cB(z)|}{|\cB(z)|^2|\Im z|}
 \asymp1+S(z).
\]
Here we used $|\cB|\asymp1$ and
$\Im\cB(z)/\Im z=\Im a(z)/\Im z-S(z)$.

Let \(u\) be the point of \(I\) nearest to \(z\). For every $t_i\in I$,
\[
 \Re\frac{z-u}{z-t_i}
 =\frac{d^2+(\Re z-u)(u-t_i)}{|z-t_i|^2}
 \ge\frac{d^2}{|z-t_i|^2}.
\]
Consequently,
\[
 d^2S(z)
 \le\Re\bigl[(z-u)(\cB(z)-1-a(z))\bigr]
 \le d\,|\cB(z)-1-a(z)|
 \lesssim d.
\]
This proves $P(z)\lesssim d^{-1}$.

Differentiating
\eqref{eq:dmc-velocity-positive-representation} gives, for $j\ge1$,
\[
 |\partial_z^j\cB(z)|
 \lesssim_j
 1+\sum_i\frac{\nu_i}{|z-t_i|^{j+1}}
 \le 1+\frac{S(z)}{d^{j-1}}
 \lesssim_j\frac{P(z)}{d^{j-1}}.
\]
The chain rule, $|\cB|\asymp1$, and $dP(z)\lesssim1$ imply
\[
 |\partial_z^p\chi(z)|+|\partial_z^p\log\cB(z)|
 \lesssim_p
 \sum_{k=1}^p\frac{P(z)^k}{d^{p-k}}
 \lesssim_p\frac{P(z)}{d^{p-1}}.
\]
Finally, $|w-z|\le d/2$ gives
$|w-t_i|\asymp|z-t_i|$ for every $i$, and hence
\[
 P(w)\asymp1+S(w)\asymp1+S(z)\asymp P(z).
\]
\end{proof}

The next lemma quantifies how the deformation affects $\varphi$
and $\chi$.

\begin{lemma}\label{lem:phi-tilde-comparison}
Assume \Cref{ass:weightfunc,ass:prelim-log-branches}.
Let $r$ be fixed, $\bms\in\{\pm1\}^r$, and let
$w\in\mathscr D$ and $\bmv\in\mathscr D^r$ satisfy
$\dist(\bmv,w)\asymp\dist(\{\bmv,w\},I)\asymp d$.
For every fixed integer $p\ge0$,
\begin{align}
\begin{split}\label{e:phideform}
\partial_w^p\varphi^{(\bmv,\bms)}(w)
&=\partial_w^p\varphi(w)
+\OO\left(\frac1{nd^{p+1}}\right),\\
\partial_w^p\chi^{(\bmv,\bms)}(w)
&=\partial_w^p\chi(w)
+\OO\left(\frac1{nd^{p+1}}\right),\\
\chi^{(\bmv,\bms)}(w+\theta/n)-\chi^{(\bmv,\bms)}(w)
&=\chi(w+\theta/n)-\chi(w)
+\OO\left(\frac1{n^2d^2}\right).
\end{split}
\end{align}
Here the derivatives are taken with $\bmv,\bms$ fixed.
\end{lemma}

\begin{proof}
Choose a sufficiently small fixed $c>0$ so that $B(w,cd)$
and its translate by $\theta/n$ lie in the fixed working
neighborhood, and the disk stays a distance $\gtrsim d$
from $I$ and every $v_j,v_j-\theta/n$.
By the positive-residue representation in the proof of
\Cref{lem:dmc-velocity-derivatives}, for $|\zeta-w|\le cd$,
\[
 |\cB(\zeta)-\cB(w)|
 \lesssim |\zeta-w|
 \left(1+\sum_i\frac{\nu_i}{|w-t_i|^2}\right)
 \lesssim\frac{|\zeta-w|}{d}.
\]
Shrinking $c$ if necessary, we obtain
$|\cB(\zeta)|\asymp1$ throughout $B(w,cd)$.

For each $j$, uniformly on this disk,
\[
 \frac{\sfb(v_j)-\sfb(\zeta)}
      {\sfb(v_j)-\sfb(\zeta+\theta/n)}
 =
 1+\frac{\sfb(\zeta+\theta/n)-\sfb(\zeta)}
         {\sfb(v_j)-\sfb(\zeta+\theta/n)}
 =
 1+\OO\left(\frac1{nd}\right),
\]
since the numerator of the second term is $\OO(n^{-1})$
and its denominator has modulus $\gtrsim d$.
By \eqref{e:phi-tilde-def},
\begin{equation}\label{e:phidiff1comp-clean}
 \varphi^{(\bmv,\bms)}(\zeta)-\varphi(\zeta)
 =\varphi(\zeta)\,\OO\left(\frac1{nd}\right).
\end{equation}
Consequently,
\[
 \cB^{(\bmv,\bms)}(\zeta)-\cB(\zeta)
 =(\cB(\zeta)-1)\,\OO\left(\frac1{nd}\right),
\]
so $|\cB^{(\bmv,\bms)}(\zeta)|\asymp1$. Hence
\begin{equation}\label{e:phidiff2comp-clean}
 \chi^{(\bmv,\bms)}(\zeta)-\chi(\zeta)
 =
 \frac{(\varphi^{(\bmv,\bms)}(\zeta)-\varphi(\zeta))\cG(\zeta)}
      {\cB^{(\bmv,\bms)}(\zeta)\cB(\zeta)}
 =\OO\left(\frac1{nd}\right).
\end{equation}

Cauchy's formula applied to \eqref{e:phidiff1comp-clean} and
\eqref{e:phidiff2comp-clean} on $B(w,cd)$ gives the first two
bounds in \eqref{e:phideform}.
Integrating the $p=1$ bound for
$\chi^{(\bmv,\bms)}-\chi$ over $[w,w+\theta/n]$
gives the last bound.
\end{proof}

\subsubsection{Weak Estimates}

In this section, we prove a weak estimate for $\log  {\cA^{(\bmv,\bms)}}(z)$, when $\operatorname{dist}\bigl(\{\bmv,z\}, I\bigr)$ is small. 

\begin{proposition}\label{p: weak estimate}
Assume \Cref{ass:weightfunc} and \Cref{ass:prelim-log-branches}.
Let $(\bmv,\bms)$ belong to either
\eqref{e:term1}--\eqref{e:term2a}, with
$\bmv\in\mathscr D^r$ and
\(
 \dist(\{v_1,\ldots,v_r\},I)\ge d,
\)
with $dn$ is big enough.
Then, uniformly for $w\in\mathscr D$ satisfying
\[
 \dist\bigl(w,I\cup\{v_1,\ldots,v_r\}\bigr)\ge d/2,
\]
one has
\(
 \left|\cE_1^{(\bmv,\bms)}(w)\right|\lesssim d^{-1}.
\)
Moreover, for every $z\in\mathscr D$ with
$\dist(z,I)\ge d$,
\begin{equation}\label{e: bulk error 2}
 \left|
 \log\cA^{(\bmv,\bms)}(z)-h_-(\sfb(z))
 \right|
 \lesssim\frac1{nd}.
\end{equation}
\end{proposition}

The proof of \Cref{p: weak estimate} follows the argument of \cite[Theorem~6.7]{Jiaoyang:edgestat}. 
The main input is the lemma stated below, which shows that a weak bound on \(\cA^{(\bmv,\bms)}(z)\) in the region
\(
\operatorname{dist}\bigl(\{z,\bmv\}, I\bigr)\ge d
\)
can be bootstrapped to an improved bound in the larger region
\(
\operatorname{dist}\bigl(\{z,\bmv\}, I\bigr)\ge \Bigl(1-\frac{1}{K}\Bigr)d,
\)
for some constant \(K\ge 1\).

\begin{lemma}\label{p:one step estimate}
Retain the notation and assumptions of \Cref{p: weak estimate}.
Fix $\mathfrak C>1$. There exists $K=K(\mathfrak C,R)$ such
that the following holds for
\(
 Mn^{-1}\le d\le 2\mathfrak r,
\)
where $M$ is sufficiently large depending on $\mathfrak C$ and $R$.
Suppose that, for $(\bmv,\bms)$ in 
\eqref{e:term1}--\eqref{e:term2a},
\begin{equation}\label{e:C bound}
 \mathfrak C^{-1}
 \le
 \left|\cA^{(\bmv,\bms)}(z)e^{-h_-(\sfb(z))}\right|
 \le \mathfrak C
\end{equation}
for $z,v_1,\ldots,v_r\in\mathscr D$ satisfying
\(
 \dist(\{z,v_1,\ldots,v_r\},I)\ge d.
\)
Then the following statements hold.
\begin{itemize}
\item
For $w,v_1,\ldots,v_r\in\mathscr D$ satisfying
\[
 \dist(\{w,v_1,\ldots,v_r\},I)
 \ge \left(1-\frac1K\right)d,
 \qquad
 |w-v_j|\ge\frac{d}{2K},
 \quad 1\le j\le r,
\]
one has
\begin{equation}\label{e:cEbound}
 \left|\cE_1^{(\bmv,\bms)}(w)\right|
 \lesssim d^{-1}.
\end{equation}

\item
For $z,v_1,\ldots,v_r\in\mathscr D$ satisfying
\(
 \dist(\{z,v_1,\ldots,v_r\},I)
 \ge \left(1-\frac1{2K}\right)d,
\)
one has
\begin{align}\label{e:improved estimate bulk}
 \partial_z\log\cA^{(\bmv,\bms)}(z)
 &=
 \partial_z h_-(\sfb(z))
 +\OO\left(\frac1{nd^2}\right) =
 \frac1{2\pi\ri}\oint_{\cin}
 \frac{\log\cB(w)\,\sfb'(w)\sfb'(z)}
      {(\sfb(w)-\sfb(z))^2}\,\rd w
 +\OO\left(\frac1{nd^2}\right),
\end{align}
where $\cin$ encloses $I$ but none of $z,v_1,\ldots,v_r$,
and is chosen in the region
\[
 \left\{
 w\in\mathscr D:
 \dist(w,I)\ge\left(1-\frac1K\right)d,\quad
 \dist(w,\{z,v_1,\ldots,v_r\})\ge\frac{d}{2K}
 \right\}.
\]
Moreover,
\begin{equation}\label{e: C bound improved}
 \left|\cA^{(\bmv,\bms)}(z)e^{-h_-(\sfb(z))}\right|
 =
 1+\OO\left(\frac1{nd}\right).
\end{equation}
\end{itemize}
\end{lemma}

The following lemma is from \cite[Proposition 6.10]{Jiaoyang:edgestat} which gives estimates on the derivatives of $\log \cA^{(\bmv,\bms)}(z)$ in terms of $z,\bm v$.

\begin{lemma}\label{l: derivative bdd}
Under the condition of \Cref{p:one step estimate}.
There exists a constant $K=K(\mathfrak C,R)$ such that the
following holds for $z\in\mathscr D$ and
$\bmv=(v_1,\ldots,v_r)\in\mathscr D^r$ satisfying
\[
 \dist(\{z,v_1,\ldots,v_r\},I)\ge (1-1/K)d.
\]
For the case in \eqref{e:term1}, one has
\[
 \left|\partial_z\log\cA^{(\bmv,\bms)}(z)\right|,
 \quad
 \left|\partial_{v_j}\log\cA^{(\bmv,\bms)}(z)\right|
 \lesssim d^{-1},
 \qquad 1\le j\le r,
\]
together with
\[
 (2\mathfrak C)^{-1}
 \le
 \left|\cA^{(\bmv,\bms)}(z)e^{-h_-(\sfb(z))}\right|
 \le 2\mathfrak C.
\]
For the case in \eqref{e:term2a},
\(
 \left|\partial_z\log\cA^{(\bmv,\bms)}(z)\right|
 \lesssim d^{-1}
\)
also holds.
\end{lemma}

\begin{proof}
In the coordinate $u=\sfb(z)$, the shifted particle intervals
remain disjoint real intervals, while \eqref{e:hz} and
\Cref{lem:dmc-velocity-derivatives} give
$|\partial_z h_-(\sfb(z))|\lesssim d^{-1}$.
Arguments of
\cite[Proposition~6.10]{Jiaoyang:edgestat} establish nonvanishing
and the asserted logarithmic derivative bounds for
\eqref{e:term1}.
For \eqref{e:term2a}, positivity of the deformed weights and
interlacing give the spectral derivative bound directly.
Integrating the derivative bounds along the continuation
paths of length $\OO(d/K)$ and taking $K$ sufficiently large
yields the two-sided estimate.
\end{proof}

\begin{lemma}\label{l:factor b equal G}
Assume \Cref{ass:weightfunc}, let $w\in\mathscr D\setminus\mathbb R$ satisfy
$\dist(w,I)\gtrsim d$, with $n d$ is
sufficiently large. Then
\begin{align}\label{e:interstep}
\begin{split}
\prod_{i=1}^n
\frac{\sfb(w+\theta/n)-\sfb(x_i)}
     {\sfb(w)-\sfb(x_i)}
&=\cG(w+\theta/n)L_G(w)\\
&=\cG(w)\exp\left[
\OO\left(
\min\left\{
\frac1{nd},
\frac{|\Im m(w)|+|\Im w|}{n|\Im w|}
\right\}
\right)\right].
\end{split}
\end{align}
Here $L_G$ is
holomorphic in $\Lambda$,
with a well defined logarithm satisfying
\[
|\partial_w^p\log L_G(w)|\lesssim_p n^{-1},
\qquad p\ge0.
\]
\end{lemma}

\begin{proof}
Recall $q=m-m_0$ from the proof of
\Cref{lem:stieltjes-bounds}. Then
\[
L_G(w)
=
e^{\theta q(w+\theta/n)}
\prod_{i=1}^n
\frac{(\sfb(w+\theta/n)-\sfb(x_i))(w-x_i)}
     {(\sfb(w)-\sfb(x_i))(w+\theta/n-x_i)}.
\]
The apparent singularities are removable, so $L_G$ is
holomorphic and nonvanishing on $\Lambda$.
Taylor expansion gives
\[
\begin{aligned}
\log L_G(w)
&=\theta q(w+\theta/n)
+\frac{\theta}{n}\sum_{i=1}^n
\left[
\frac{\sfb'(w)}{\sfb(w)-\sfb(x_i)}
-\frac1{w-x_i}
\right]+\OO(n^{-1})\\
&=\theta\bigl(q(w+\theta/n)-q(w)\bigr)+\OO(n^{-1})
=\OO(n^{-1}).
\end{aligned}
\]
Similar estimates yield,
for every fixed integer $p\ge0$,
\[
 |\partial_w^p\log L_G(w)|\lesssim_p n^{-1}.
\]

It remains to compare $\cG(w+\theta/n)$ with $\cG(w)$.
For $0\le t\le\theta/n$ and $s\in I$, we have
$|w+t-s|\asymp|w-s|$. By \Cref{lem:stieltjes-bounds},
\[
\theta|m(w+\theta/n)-m(w)|
\lesssim
\min\left\{
\frac1{nd},
\frac{|\Im m(w)|+|\Im w|}{n|\Im w|}
\right\}.
\]
The bound
$\log L_G=\OO(n^{-1})$ can be absorbed into either estimate:
\[
\left|
\log L_G(w)-\theta\bigl(m(w+\theta/n)-m(w)\bigr)
\right|
\lesssim
\min\left\{
\frac1{nd},
\frac{|\Im m(w)|+|\Im w|}{n|\Im w|}
\right\}.
\]
The desired expansion now follows from the exact identity
\[
\prod_{i=1}^n
\frac{\sfb(w+\theta/n)-\sfb(x_i)}
     {\sfb(w)-\sfb(x_i)}
=
\cG(w)\exp\left[
\log L_G(w)-\theta\bigl(m(w+\theta/n)-m(w)\bigr)
\right].
\]
For real $w$, the same argument uses only
$|m'(w+t)|\lesssim d^{-1}$.
\end{proof}

\begin{proof}[Proof of \Cref{p:one step estimate}]
Let $K_0$ be the constant in \Cref{l: derivative bdd}.
Choose $K\ge2K_0$ and then take $M$ sufficiently large.
By \Cref{l: derivative bdd},
\[
 \frac{\cA^{(\bmv,\bms)}(w+\theta/n)}
      {\cA^{(\bmv,\bms)}(w)}
 =
 \exp\left\{
 \int_w^{w+\theta/n}
 \partial_u\log\cA^{(\bmv,\bms)}(u)\,\rd u
 \right\}
 =
 1+\OO\left(\frac1{nd}\right).
\]
The deformation formula gives
\[
 \prod_{j=1}^r
 \left(
 \frac{\sfb(v_j)-\sfb(w)}
      {\sfb(v_j)-\sfb(w+\theta/n)}
 \right)^{s_j}
 =
 1+\OO\left(\frac1{nd}\right).
\]
Combining this with \Cref{l:factor b equal G}, we obtain
\[
 \frac{\varphi^{(\bmv,\bms)}(w)}{\cB(w)}
 \prod_{i=1}^n
 \frac{\sfb(w+\theta/n)-\sfb(x_i)}
      {\sfb(w)-\sfb(x_i)}
 =
 \chi(w)\left(1+\OO\left(\frac1{nd}\right)\right).
\]
Substitution into \eqref{e:error-term-1-clean: deform},
together with $|\chi(w)|=|1-\cB(w)^{-1}|\lesssim1$,
yields
\[
 \frac1n\cE_1^{(\bmv,\bms)}(w)
 =
 \OO\left(\frac1{nd}\right).
\]
This proves \eqref{e:cEbound}. The same calculation applies
with the separation $d/(2K)$ replaced by $d/(8K)$.

For the remaining assertions, let $c_1$ be the oriented
boundary of the union of disks of radius $d/(8K)$ centered
at $z,v_1,\ldots,v_r$, and let $c_2$ be the fixed outer
contour in $\Lambda$. Thus
\[
 \operatorname{length}(c_1)\lesssim_R d,
 \qquad
 \dist(c_1,\{z,v_1,\ldots,v_r\})\ge\frac{d}{8K},
\]
whereas $c_2$ has bounded length and stays a fixed positive
distance from $\{z,\bm v \}$.

Since 
$\log(1+\cE_1^{(\bmv,\bms)}(w)/n)$ is holomorphic between
$\cin$ and $c_2$ outside these disks and is
$\OO((nd)^{-1})$. On $c_2$, the preceding deformation bounds and the fixed-distance
argument give $\cE_1^{(\bmv,\bms)}(w)=\OO(1)$.
Hence, for $k=1,2$, Cauchy's theorem gives
\begin{align*}
 &\left|
 \oint_{\cin}
 \frac{\log(1+\cE_1^{(\bmv,\bms)}(w)/n)\,\sfb'(w)}
      {(\sfb(w)-\sfb(z))^k}\,\rd w
 \right|=
 \left|
 \left(\oint_{c_2}-\oint_{c_1}\right)
 \frac{\log(1+\cE_1^{(\bmv,\bms)}(w)/n)\,\sfb'(w)}
      {(\sfb(w)-\sfb(z))^k}\,\rd w
 \right|
 \lesssim
 \frac1{nd}\left(1+\frac{d}{d^k}\right)
 \lesssim\frac1{nd^k}.
\end{align*}
Applying the case $k=1$ to \eqref{e: bulk error} yields
\[
 \log\cA^{(\bmv,\bms)}(z)-h_-(\sfb(z))
 =
 \OO\left(\frac1{nd}\right),
\]
and exponentiation proves \eqref{e: C bound improved}.
Applying the case $k=2$ to the $z$ derivative of
\eqref{e: bulk error}, taken with the original inner
contour held fixed locally, proves
\eqref{e:improved estimate bulk}.
\end{proof}

\begin{proof}[Proof of \Cref{p: weak estimate}]
We first consider $d\le2\mathfrak r$.
By \eqref{e:initial bootstrap}, there exists
$\mathfrak C>1$ such that
\begin{equation}\label{e:prior bdd}
 \mathfrak C^{-1}
 \le
 \left|\cA^{(\bmv,\bms)}(z)e^{-h_-(\sfb(z))}\right|
 \le\mathfrak C
\end{equation}
whenever
$\dist(\{z,v_1,\ldots,v_r\},I)\ge2\mathfrak r$,
uniformly over \eqref{e:term1}--\eqref{e:term2a}.
Let $K=K(\mathfrak C,R)$ be the constant in
\Cref{p:one step estimate}, and set
\[
 \eta_0=2\mathfrak r,
 \qquad
 \eta_{i+1}=\left(1-\frac1{2K}\right)\eta_i.
\]
Take $M$ sufficiently large that the one-step lemma
applies whenever $\eta_i\ge M/(2n)$ and that
\eqref{e: C bound improved} restores
\eqref{e:prior bdd} with the same constant $\mathfrak C$.
Induction therefore propagates \eqref{e:prior bdd}
to all these scales.

Choose $i$ such that $\eta_{i+1}\le d/2<\eta_i$.
Applying \Cref{p:one step estimate} at scale $\eta_i\asymp d$
gives the proof.
\end{proof}

\subsubsection{Improved Estimate}\label{s:DMCedge estimate}
In this subsection, we refine \Cref{p: weak estimate} in order to establish
the second-order mean and covariance estimates in \Cref{thm:main}.
We retain the standing assumptions and spectral-domain conventions,
and suppress the fixed time index as above.
Recall $P(z)=P_s^n(z)$ defined in \eqref{eq:main-P}.

\begin{proposition}\label{p:edge impr}
Assume \Cref{ass:weightfunc} and \Cref{ass:prelim-log-branches}.
In \eqref{e:phi-tilde-def}, take $r=1$, $\bmv=(v)$ and $\bms=(1)$.
For $v\in\mathscr D\setminus\mathbb R$, set $d=\dist(v,I)$
with $d\ge n^{-1+\delta}$.
Then, for every fixed integer $p\ge1$,
\begin{equation}\label{e:partialzlnA}
 \left|
 \left.
 \partial_w^p
 \left(
 \log\cA^{(\bmv,\bms)}(w)-h_-(\sfb(w))
 \right)
 \right|_{w=v}
 \right|
 \lesssim_p \frac{P(v)}{nd^p}.
\end{equation}
\end{proposition}

For the single deformation in \Cref{p:edge impr}, define
\begin{align*}
 L_A(w)
 &:=
 \frac{
 \cA^{(\bmv,\bms)}(w+\theta/n)
 e^{-h_-(\sfb(w+\theta/n))}
 }{
 \cA^{(\bmv,\bms)}(w)e^{-h_-(\sfb(w))}
 },\\
 L_G(w)
 &:=
 \frac{1}{\cG(w+\theta/n)}
 \prod_{i=1}^n
 \frac{\sfb(w+\theta/n)-\sfb(x_i)}
      {\sfb(w)-\sfb(x_i)},\\
 L(w)
 &:=
 L_A(w)L_G(w)
 \frac{\varphi(w)\cG(w+\theta/n)}
      {\cB(w+\theta/n)}
 \exp\left(
 h_+(\sfb(w+\theta/n))-h_+(\sfb(w))
 \right).
\end{align*}
Since $\log\cB=h_+\circ\sfb-h_-\circ\sfb$,
\eqref{e:error-term-1-clean: deform} becomes
\begin{equation}\label{e:error-term-2-clean: deform}
\begin{aligned}
 \frac1n\cE_0^{(\bmv,\bms)}(w)
 &:=
 L(w)-\chi(w),\\
 \frac1n\cE_1^{(\bmv,\bms)}(w)
 &=
 \frac1n\cE_0^{(\bmv,\bms)}(w)
 +L(w)\left(
 \frac{\sfb(v)-\sfb(w)}
      {\sfb(v)-\sfb(w+\theta/n)}-1
 \right).
\end{aligned}
\end{equation}
By \Cref{l:factor b equal G}, for every fixed integer $p\ge0$,
\[
 |\partial_w^p\log L_G(w)|
 +
 |\partial_w^p(L_G(w)-1)|
 \lesssim_p n^{-1}.
\]
The proof of \Cref{p:edge impr} uses the following two lemmas.

\begin{lemma}\label{l:bdd for err and L(z)}
Retain the notation and assumptions of \Cref{p:edge impr}.
For a sufficiently small fixed $c>0$, set
\begin{equation}\label{eq:dmc-Mv-definition}
 M_v:=
 \sup_{|\zeta-v|\le2cd}
 \left|\partial_\zeta
 \bigl(\log\cA^{(\bmv,\bms)}(\zeta)-h_-(\sfb(\zeta))\bigr)
 \right|.
\end{equation}
Then, uniformly for $|w-v|\le cd$,
\begin{equation}\label{e: higher derivative of err_0}
\begin{aligned}
 \left|\partial_w^p
       \left(\frac1n\cE_0^{(\bmv,\bms)}(w)\right)\right|
 &\lesssim_p\frac{P(v)+M_v}{nd^p},&&p\ge0,\\
 |L(w)|&\lesssim1,\\
 |\partial_w^pL(w)|
 &\lesssim_p\frac{P(v)}{d^{p-1}}+\frac{M_v}{nd^p},&&p\ge1.
\end{aligned}
\end{equation}
\end{lemma}

\begin{proof}
By \Cref{p: weak estimate},
\[
 \left|\log\cA^{(\bmv,\bms)}(\zeta)-h_-(\sfb(\zeta))\right|
 \lesssim\frac1{nd}
\]
uniformly for $|\zeta-v|\le3cd$, since
$\dist(\zeta,I)\asymp d$ there. Cauchy's formula therefore gives
$M_v\lesssim(nd^2)^{-1}$.

We first prove the $p=0$ estimates on the slightly larger
neighborhood $|w-v|\le3cd/2$.
Since $\theta/n\ll d$, the definition of $L_A$ gives
\[
\begin{aligned}
 L_A(w)
 &=\exp\left\{\int_w^{w+\theta/n}
 \partial_u\bigl(\log\cA^{(\bmv,\bms)}(u)-h_-(\sfb(u))\bigr)
 \,\rd u\right\} =\exp\{\OO(M_v/n)\}=1+\OO(M_v/n),
\end{aligned}
\]
where $M_v/n\lesssim(nd)^{-2}\ll1$.
The uniform lower bound on $\varphi$ and its analytic regularity give
\[
 \frac{\varphi(w)}{\varphi(w+\theta/n)}-1
 =\frac{\varphi(w)-\varphi(w+\theta/n)}{\varphi(w+\theta/n)}
 =\OO(n^{-1}).
\]
Likewise, regularity of $h_+\circ\sfb$ gives
$h_+(\sfb(w+\theta/n))-h_+(\sfb(w))=\OO(n^{-1})$.
Using $\varphi\cG/\cB=\chi$ and \Cref{l:factor b equal G},
we obtain
\[
\begin{aligned}
 \frac{L(w)}{L_A(w)}
 &=\chi(w+\theta/n)L_G(w)
 \frac{\varphi(w)}{\varphi(w+\theta/n)}
 e^{h_+(\sfb(w+\theta/n))-h_+(\sfb(w))}\\
 &=\chi(w+\theta/n)\bigl(1+\OO(n^{-1})\bigr).
\end{aligned}
\]
In particular, $|L/L_A|\lesssim1$ and $|L|\lesssim1$.

\Cref{lem:dmc-velocity-derivatives} and local comparability of $P$ give
\[
 |\chi(w+\theta/n)-\chi(w)|
 \le\int_0^{\theta/n}|\chi'(w+t)|\,\rd t
 \lesssim\frac{P(v)}n.
\]
Together with $|\chi|\lesssim1$ and $P(v)\ge1$, this yields
$|L/L_A-\chi|\lesssim P(v)/n$. Consequently,
\[
\begin{aligned}
 \left|\frac1n\cE_0^{(\bmv,\bms)}(w)\right|
 &=|L(w)-\chi(w)|\\
 &=\left|L_A(w)\left(\frac{L(w)}{L_A(w)}-\chi(w)\right)
       +(L_A(w)-1)\chi(w)\right|
 \lesssim\frac{P(v)+M_v}{n}.
\end{aligned}
\]
For $|w-v|\le cd$, Cauchy's formula on
$|\zeta-w|=cd/2$ now gives
\[
 \left|\partial_w^p\left(\frac1n\cE_0^{(\bmv,\bms)}(w)\right)\right|
 \lesssim_p d^{-p}
 \sup_{|\zeta-w|=cd/2}\left|\frac1n\cE_0^{(\bmv,\bms)}(\zeta)\right|
 \lesssim_p\frac{P(v)+M_v}{nd^p}.
\]
This proves the estimates for $\cE_0/n$.
Finally, differentiating $L=\chi+\cE_0^{(\bmv,\bms)}/n$ and using
\Cref{lem:dmc-velocity-derivatives}, we obtain
\[
 |\partial_w^pL(w)|
 \lesssim_p\frac{P(v)}{d^{p-1}}+
              \frac{P(v)+M_v}{nd^p}
 \lesssim_p\frac{P(v)}{d^{p-1}}+\frac{M_v}{nd^p},
 \qquad p\ge1,
\]
where the last step uses $nd\gg1$.
\end{proof}

\begin{lemma}\label{lem:dmc-Mv-improvement}
Under the assumptions of \Cref{p:edge impr}, one has
\begin{equation}\label{eq:dmc-Mv-improved}
 M_v\lesssim\frac{P(v)}{nd}.
\end{equation}
\end{lemma}

\begin{proof}
By \Cref{p: weak estimate},
$\log\cA^{(\bmv,\bms)}(\zeta)-h_-(\sfb(\zeta))
=\OO((nd)^{-1})$ uniformly for $|\zeta-v|\le4cd$.
Cauchy's formula therefore gives
\[
 \sup_{|\zeta-v|\le3cd}
 \left|
 \partial_\zeta
 \bigl(\log\cA^{(\bmv,\bms)}(\zeta)
       -h_-(\sfb(\zeta))\bigr)
 \right|
 \lesssim\frac1{nd^2}.
\]

To improve this bound, recall for any $|\zeta-v|\le3cd$,
\[
 \partial_\zeta
 \bigl(\log\cA^{(\bmv,\bms)}(\zeta)
       -h_-(\sfb(\zeta))\bigr)
 =
 \frac1{2\pi\ri}\oint_{\cin}
 \frac{
 \log(1+\cE_1^{(\bmv,\bms)}(w)/n)
 \,\sfb'(w)\sfb'(\zeta)}
 {(\sfb(w)-\sfb(\zeta))^2}\,\rd w.
\]
choose $\cin$ enclosing $I$ and leaving
the closed disk $|\zeta-v|\le4cd$ in its exterior.

Fix $2c\le r<r'\le3c$, and let $\omega_v$ be the positively
oriented circle $|w-v|=(r+r')d/2$.
For $|\zeta-v|\le rd$, both $\zeta$ and the shifted pole
$v-\theta/n$ lie inside $\omega_v$. Cauchy's theorem deforms $\omega_-$ to
$\omega_+-\omega_v$, where $\omega_+$ is
the fixed outer contour; see
\Cref{fig:contour-deformation}.
The outer-contour contribution is $\OO(n^{-1})$ by
\Cref{p: weak estimate}.
It remains to estimate the local integral over $\omega_v$.

\begin{figure}[t]
    \centering
    \includegraphics[width=0.8\linewidth]
        {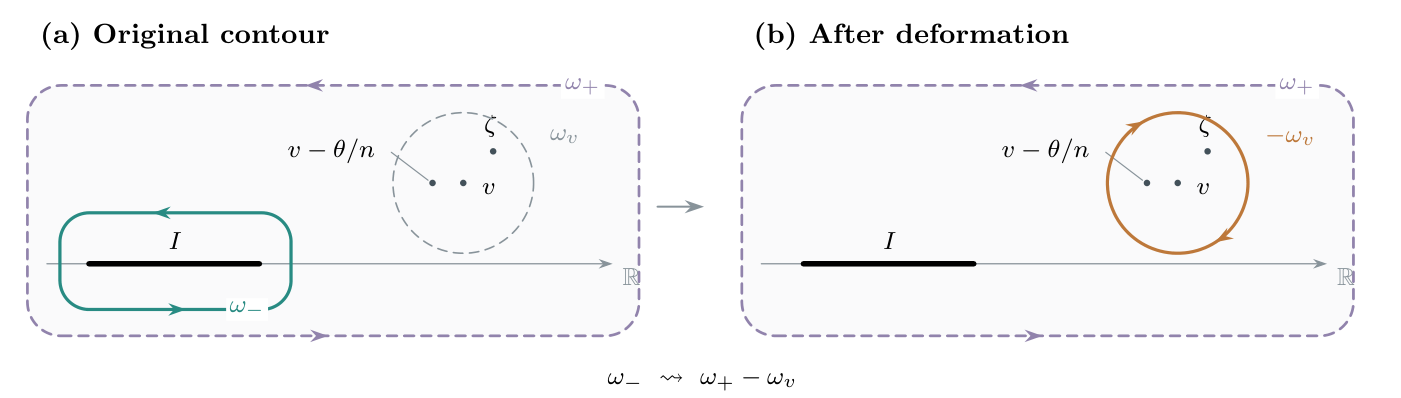}
    \caption{Contour deformation in the improved DMC estimate.
    The contour $\omega_-$ encloses $I$ and leaves the spectral
    points outside. Deforming it to the fixed outer contour $\omega_+$
    produces the local circle $\omega_v$ around $\zeta$ and
    $v-\theta/n$. The orange clockwise contour in panel (b)
    represents the subtraction of the positively oriented $\omega_v$.}
    \label{fig:contour-deformation}
\end{figure}

For this purpose, using
$L=\chi+\cE_0^{(\bmv,\bms)}/n$ in
\eqref{e:error-term-2-clean: deform}, we decompose
$\cE_1^{(\bmv,\bms)}/n$ into three terms:
\[
\begin{aligned}
 \frac1n\cE_1^{(\bmv,\bms)}(w)
 ={}&
 \frac1n\cE_0^{(\bmv,\bms)}(w)
 \frac{\sfb(v)-\sfb(w)}
      {\sfb(v)-\sfb(w+\theta/n)}
 \\
 &+
 \bigl(\chi(w)-\chi(v)\bigr)
 \left(
 \frac{\sfb(v)-\sfb(w)}
      {\sfb(v)-\sfb(w+\theta/n)}-1
 \right)
 \\
 &+
 \chi(v)
 \left(
 \frac{\sfb(v)-\sfb(w)}
      {\sfb(v)-\sfb(w+\theta/n)}-1
 \right).
\end{aligned}
\]
We estimate these three terms on $\omega_v$.

For the first term, the zeroth-order estimate in
\Cref{l:bdd for err and L(z)}, with the derivative
controlled on $|\zeta-v|\le r'd$, together with
\[
 \frac{\sfb(v)-\sfb(w)}
      {\sfb(v)-\sfb(w+\theta/n)}
 =1+\OO\left(\frac1{nd}\right),
\]
gives
\[
 \left|
 \frac1n\cE_0^{(\bmv,\bms)}(w)
 \frac{\sfb(v)-\sfb(w)}
      {\sfb(v)-\sfb(w+\theta/n)}
 \right|
 \lesssim
 \frac1n\left(
 P(v)+
 \sup_{|\zeta-v|\le r'd}
 \left|
 \partial_\zeta
 \bigl(\log\cA^{(\bmv,\bms)}(\zeta)
       -h_-(\sfb(\zeta))\bigr)
 \right|
 \right).
\]
For the second term, \Cref{lem:dmc-velocity-derivatives}
gives $|\chi(w)-\chi(v)|\lesssim dP(v)$, and hence
\[
 \left|
 \bigl(\chi(w)-\chi(v)\bigr)
 \left(
 \frac{\sfb(v)-\sfb(w)}
      {\sfb(v)-\sfb(w+\theta/n)}-1
 \right)
 \right|
 \lesssim dP(v)\frac1{nd}
 =\frac{P(v)}n.
\]
For the last term, since $\theta/n\ll d$, Taylor expansion gives
\begin{equation}\label{eq:dmc-shifted-principal-part}
 \frac{\sfb(v)-\sfb(w)}
      {\sfb(v)-\sfb(w+\theta/n)}-1
 =
 -\frac{\theta\sfb'(v)}
 {n\bigl(\sfb(w)-\sfb(v-\theta/n)\bigr)}
 +\OO(n^{-1}).
\end{equation}
Combining these estimates and using
$|\chi(v)|\lesssim1$ and $P(v)\ge1$, we obtain
\begin{equation}\label{eq:dmc-pole-subtraction-bound}
\begin{aligned}
 &\left|
 \frac1n\cE_1^{(\bmv,\bms)}(w)
 +\frac{\theta\chi(v)\sfb'(v)}
 {n\bigl(\sfb(w)-\sfb(v-\theta/n)\bigr)}
 \right|
 \lesssim
 \frac1n\left(
 P(v)+
 \sup_{|\zeta-v|\le r'd}
 \left|
 \partial_\zeta
 \bigl(\log\cA^{(\bmv,\bms)}(\zeta)
       -h_-(\sfb(\zeta))\bigr)
 \right|
 \right).
\end{aligned}
\end{equation}
Both terms in \eqref{eq:dmc-pole-subtraction-bound} are
$\OO((nd)^{-1})$.
Since the logarithm is uniformly Lipschitz near $1$, we have
\[
\begin{aligned}
 &\left|
 \log\left(1+\frac1n\cE_1^{(\bmv,\bms)}(w)\right)
 -
 \log\left(
 1-\frac{\theta\chi(v)\sfb'(v)}
 {n\bigl(\sfb(w)-\sfb(v-\theta/n)\bigr)}
 \right)
 \right| \lesssim
 \left|
 \frac1n\cE_1^{(\bmv,\bms)}(w)
 +\frac{\theta\chi(v)\sfb'(v)}
 {n\bigl(\sfb(w)-\sfb(v-\theta/n)\bigr)}
 \right|,
\end{aligned}
\]
with both logarithms chosen near $1$.

For $|z-v|\le rd$,
\begin{equation}\label{eq:dmc-pole-cancellation}
 \oint_{\omega_v}
 \log\left(
 1-\frac{\theta\chi(v)\sfb'(v)}
 {n\bigl(\sfb(w)-\sfb(v-\theta/n)\bigr)}
 \right)
 \frac{\sfb'(w)\sfb'(z)}
      {(\sfb(w)-\sfb(z))^2}\,\rd w=0.
\end{equation}
Indeed, expanding the logarithm in its uniformly convergent
series on $\omega_v$ and setting $u=\sfb(w)$, each term is
a constant multiple of
\[
 \oint_{\sfb(\omega_v)}
 \frac{\rd u}
 {\bigl(u-\sfb(v-\theta/n)\bigr)^j
  \bigl(u-\sfb(z)\bigr)^2},
 \qquad j\ge1.
\]
This integral vanishes because all finite poles are enclosed
and the integrand is $\OO(u^{-j-2})$ at infinity.

Subtracting \eqref{eq:dmc-pole-cancellation} from the local
integral obtained above, we can therefore apply
\eqref{eq:dmc-pole-subtraction-bound}.
Since $\omega_v$ has length $\OO(d)$ and its derivative
kernel is $\OO_{r,r'}(d^{-2})$, combining the local estimate
with the outer-contour contribution gives
\begin{equation}\label{eq:dmc-local-contraction}
\begin{aligned}
 &\sup_{|z-v|\le rd}
 \left|
 \partial_z
 \bigl(\log\cA^{(\bmv,\bms)}(z)-h_-(\sfb(z))\bigr)
 \right|
 \lesssim_{r,r'}
 \frac{P(v)}{nd}
 +\frac1{nd}
 \sup_{|\zeta-v|\le r'd}
 \left|
 \partial_\zeta
 \bigl(\log\cA^{(\bmv,\bms)}(\zeta)
       -h_-(\sfb(\zeta))\bigr)
 \right|,
\end{aligned}
\end{equation}
where we used $P(v)\ge1$ to absorb the outer contribution.

Choose a fixed integer $J$ such that $\delta(J+1)>2$.
Iterating \eqref{eq:dmc-local-contraction} over $J$ equally
spaced steps from $2c$ to $3c$, and applying the initial
bound on the outermost disk, yields
\[
 M_v
 \lesssim_J
 \frac{P(v)}{nd}
 \sum_{j=0}^{J-1}\frac1{(nd)^j}
 +\frac1{(nd)^J}\frac1{nd^2}
 \lesssim_J
 \frac{P(v)}{nd}
 +\frac1{(nd)^J}\frac1{nd^2}.
\]
Since $nd\ge n^\delta$ and $P(v)\ge1$,
\[
 dM_v
 \lesssim_J
 \frac{P(v)}n+(nd)^{-J-1}
 \le
 C_J\left(\frac{P(v)}n+n^{-\delta(J+1)}\right)
 \lesssim\frac{P(v)}n.
\]
This proves \eqref{eq:dmc-Mv-improved}.
\end{proof}

Combining the two lemmas, we obtain $|L(w)|\lesssim1$ and,
uniformly for $|w-v|\le cd$,
\begin{equation}\label{eq:dmc-local-coefficients-improved}
\begin{aligned}
 \left|\partial_w^p\left(\frac1n\cE_0^{(\bmv,\bms)}(w)\right)\right|
 &\lesssim_p\frac{P(v)}{nd^p},&&p\ge0,\\
 |\partial_w^pL(w)|
 &\lesssim_p\frac{P(v)}{d^{p-1}},&&p\ge1,
\end{aligned}
\end{equation}

\begin{proof}[Proof of \Cref{p:edge impr}]
By \Cref{lem:dmc-Mv-improvement}, $M_v\lesssim P(v)/(nd)$.
Cauchy's formula applied to the spectral derivative in the
neighborhood defining $M_v$ gives, for every fixed $p\ge1$,
\[
 \left|\left.\partial_z^p
 \bigl(\log\cA^{(\bmv,\bms)}(z)-h_-(\sfb(z))\bigr)
 \right|_{z=v}\right|
 \lesssim_p d^{1-p}M_v
 \lesssim_p\frac{P(v)}{nd^p}.
\]
This is \eqref{e:partialzlnA}.
\end{proof}

\subsubsection{Proof of \Cref{thm:main} for DMC}
\paragraph{Second-order mean}
In this part, we finish the proof of \Cref{thm:main} for the discrete Markov chain models.

For simplicity, we suppress the time index, define
\begin{equation}\label{err:secondder}
\begin{aligned}
\cE_2(w,z)
={}&\chi(w)\Biggl(
\frac{\theta\sfb'(w)}{\sfb(z)-\sfb(w)}
+\theta\partial_w h_-(\sfb(w))\\
&\qquad+\frac{\theta^2}{2}\int_I
\frac{\sfb''(w)(\sfb(w)-\sfb(x))-(\sfb'(w))^2
      -\sfb'(w)\sfb'(x)}
     {(\sfb(w)-\sfb(x))^2}\,\mu(\rd x;\bmx)
\Biggr)-\frac{\cB'(w)}{2\cB(w)}.
\end{aligned}
\end{equation}

\begin{corollary}\label{cor:one step result}
Assume \Cref{ass:weightfunc,ass:prelim-log-branches} and the domain and contour conditions above.
For $s=s_j\le T$ and $z\in\mathscr D_s\setminus\mathbb R$, with $d=\dist(z,I_s)\ge n^{-1+\delta}$, we have
\begin{equation}\label{t:loopstudy2}
 \E_sD_s^n(z)=a_s^n(z)+\frac1n c_s^n(z)+\mathcal R_s^n(z),
 \qquad
 c_s^n(z)=-\frac1{2\pi\ri\theta}\oint_{\cin}\cE_2(w,z)K_{\sfb}(w,z)\,\rd w.
\end{equation}
Here $\cE_2$ is defined in \eqref{err:secondder}, with the time index suppressed, and $\cin$ encloses $I_s$ but not $z$.
For every fixed integer $p\ge0$,
\[
 |\partial_z^p\mathcal R_s^n(z)|\lesssim_p\frac{P_s^n(z)}{n^2d^{p+2}}.
\]
\end{corollary}

\begin{proof}
Fix $s=s_j$ and suppress the time index. Denote
\[
 Q_{\bme}(v):=\prod_{i=1}^n
 \frac{\sfb(v)-\sfb(x_i+\theta e_i/n)}{\sfb(v)-\sfb(x_i)}.
\]
By \eqref{e:phi-tilde-def}, $a_{\bme}^{(z,1)}=a_{\bme}Q_{\bme}(z)^{-1}/\E Q_{\bme}(z)^{-1}$.
This gives, with the deformation parameter held fixed when differentiating in $\zeta$,
\[
\begin{aligned}
 \left.\partial_\zeta\log\cA^{(z,1)}(\zeta)\right|_{\zeta=z}
 &=\E\partial_z\log Q_{\bme}(z)\\
 &=\E\sum_i e_i\left[
 \frac{\sfb'(z)}{\sfb(z)-\sfb(x_i+\theta/n)}
 -\frac{\sfb'(z)}{\sfb(z)-\sfb(x_i)}
 \right].
\end{aligned}
\]
To compute $\left.\partial_\zeta\log\cA^{(z,1)}(\zeta)\right|_{\zeta=z}$, we differentiate \eqref{e: bulk error} and deform the contour integral to $\cout-\omega_z$, with $\omega_z=\{|w-z|=cd\}$ as in the proof of \Cref{lem:dmc-Mv-improvement}.

We first compute the second-order expansion of the local integral
\[
 \oint_{\omega_z}\log\left(1+\frac{\cE_1^{(z,1)}(w)}n\right)K_{\sfb}(w,z)\,\rd w.
\]
Recall $L,L_A,L_G$ defined above \eqref{e:error-term-2-clean: deform}, with $v=z$. By that identity,
\[
\begin{aligned}
 \frac1n\cE_1^{(z,1)}(w)
 ={}&\frac1n\cE_0^{(z,1)}(w)
 +\frac1n\cE_0^{(z,1)}(w)\left(
 \frac{\sfb(z)-\sfb(w)}{\sfb(z)-\sfb(w+\theta/n)}-1
 \right)+\chi(w)\left(
 \frac{\sfb(z)-\sfb(w)}{\sfb(z)-\sfb(w+\theta/n)}-1
 \right).
\end{aligned}
\]
By \eqref{eq:dmc-Mv-improved}, $L_A-1=\OO(P(z)/(n^2d))$ locally.
Taylor expansion of the regular kernel, as in the proof of \Cref{lem:stieltjes-bounds}, gives
\[
 \log L_G(w)
 =-\frac{\theta^2}{2n}\int_I(\partial_w+\partial_x)
 \frac{\sfb'(w)}{\sfb(w)-\sfb(x)}\,\mu(\rd x;\bmx)+\OO(n^{-2}).
\]
Indeed, the contributions of $(w-x)^{-1}$ cancel exactly.
Substituting this into $L$, using $\log\cB=h_+\circ\sfb-h_-\circ\sfb$ and \Cref{lem:dmc-velocity-derivatives}, Taylor expansion yields
\[
 \frac1n\cE_0^{(z,1)}(w)=L(w)-\chi(w)
 =\frac1n\left(
 \cE_2(w,z)+\frac{\cB'(w)}{2\cB(w)}
 -\frac{\theta\chi(w)\sfb'(w)}{\sfb(z)-\sfb(w)}
 \right)+\OO\left(\frac{P(z)}{n^2d}\right).
\]
Here the integral terms combine into
$\frac{\theta^2\chi(w)}2\int_I(\partial_w-\partial_x)
[\sfb'(w)/(\sfb(w)-\sfb(x))]\,\mu(\rd x;\bmx)$,
as in \eqref{err:secondder}.
By \eqref{eq:dmc-local-coefficients-improved} and the $\OO((nd)^{-1})$ bound on the deformation factor minus $1$,
\[
 \frac1n\cE_0^{(z,1)}(w)\left(
 \frac{\sfb(z)-\sfb(w)}{\sfb(z)-\sfb(w+\theta/n)}-1
 \right)=\OO\left(\frac{P(z)}{n^2d}\right).
\]
Combining these expansions gives
\[
\begin{aligned}
 &\frac{\cE_1^{(z,1)}(w)}n
 -\frac1n\left(\cE_2(w,z)+\frac{\cB'(w)}{2\cB(w)}\right)\\
 &\qquad=\chi(w)\left[
 \frac{\sfb(z)-\sfb(w)}{\sfb(z)-\sfb(w+\theta/n)}
 -1-\frac{\theta\sfb'(w)}{n(\sfb(z)-\sfb(w))}
 \right]+\OO\left(\frac{P(z)}{n^2d}\right).
\end{aligned}
\]
Combining these estimates, and using \eqref{eq:dmc-pole-subtraction-bound}, \eqref{eq:dmc-Mv-improved}, and \eqref{eq:dmc-pole-cancellation} for the nonlinear logarithmic remainder, we obtain
\[
 \left|\oint_{\omega_z}\left[
 \log\left(1+\frac{\cE_1^{(z,1)}(w)}n\right)
 -\frac1n\left(\cE_2(w,z)+\frac{\cB'(w)}{2\cB(w)}\right)
 \right]K_{\sfb}(w,z)\,\rd w\right|
 \lesssim\frac{P(z)}{n^2d^2}.
\]
Here we use the following estimate, with $[\cdots]$ denoting the bracketed deformation remainder above:
\[
\begin{aligned}
 \left|\oint_{\omega_z}\chi(w)[\cdots]K_{\sfb}(w,z)\,\rd w\right|
 \le{}&
 \left|\oint_{\omega_z}(\chi(w)-\chi(z))[\cdots]K_{\sfb}(w,z)\,\rd w\right|+\left|\oint_{\omega_z}\chi(z)[\cdots]K_{\sfb}(w,z)\,\rd w\right|
 \lesssim\frac{P(z)}{n^2d^2}.
\end{aligned}
\]
In particular, \Cref{lem:dmc-velocity-derivatives} gives $|\chi(w)-\chi(z)|\lesssim dP(z)$.
Since the bracketed term is $\OO(n^{-2}d^{-2})$ on $\omega_z$, $\operatorname{length}(\omega_z)\lesssim d$, and $|K_{\sfb}(w,z)|\lesssim d^{-2}$, we obtain
\[
 \left|\oint_{\omega_z}(\chi(w)-\chi(z))[\cdots]K_{\sfb}(w,z)\,\rd w\right|
 \lesssim d\cdot dP(z)\cdot\frac1{n^2d^2}\cdot d^{-2}
 =\frac{P(z)}{n^2d^2}.
\]
For the constant term, Cauchy's theorem allows us to enlarge $\omega_z$ to a fixed-radius contour $\Gamma_z$ in the analytic coordinate neighborhood, enclosing both $z$ and $z-\theta/n$.
On $\Gamma_z$, the denominators are uniformly bounded away from zero; Taylor expansion therefore bounds the bracketed term by $\OO(n^{-2})$, while $|K_{\sfb}(w,z)|\lesssim1$.
Since $\Gamma_z$ has bounded length and $|\chi(z)|\lesssim1$, we obtain
\[
 \left|\oint_{\omega_z}\chi(z)[\cdots]K_{\sfb}(w,z)\,\rd w\right|
 \lesssim\frac1{n^2}\lesssim\frac{P(z)}{n^2d^2}.
\]

We next estimate the outer-contour contribution.
We show that the same expansion holds on $\cout$ with remainder $\OO(n^{-2})$.
By the insertion identities, Cauchy--Schwarz, and \Cref{p: weak estimate} with $R\ge2$,
\[
 \sum_{\bme}|a_{\bme}^{(z,1)}|
 \le\left|
 \frac{\cA^{(z,1)}(z)}{\cA^{((z,\bar z),(1,1))}(\bar z)}
 \right|^{1/2}\lesssim1.
\]
At fixed distance from $I$ and $z$, the ratios $Q_{\bme}(w)$ and their derivatives are uniformly bounded.
Together with $|\cA^{(z,1)}(w)|\asymp1$ from the weak estimate, this gives $|\partial_w\log\cA^{(z,1)}(w)|\lesssim1$.
Substituting into \eqref{e:error-term-1-clean: deform} and using \Cref{l:factor b equal G}, we obtain $\cE_1^{(z,1)}(w)=\OO(1)$ on fixed outer neighborhoods separated from $z$.
Similarly we have along $\cout$,
\[
 \log\left(1+\frac{\cE_1^{(z,1)}(w)}n\right)
 =\frac1n\left(\cE_2(w,z)+\frac{\cB'(w)}{2\cB(w)}\right)+\OO(n^{-2}).
\]
Combining these estimate, we obtain
\begin{equation}\label{eq:dmc-second-order-insertion}
\begin{aligned}
 \left.\partial_\zeta\log\cA^{(z,1)}(\zeta)\right|_{\zeta=z}
 ={}&\partial_z h_-(\sfb(z)) +\frac1{2\pi\ri n}\oint_{\cin}
 \left(\cE_2(w,z)+\frac{\cB'(w)}{2\cB(w)}\right)
 K_{\sfb}(w,z)\,\rd w
 +\OO\left(\frac{P(z)}{n^2d^2}\right).
\end{aligned}
\end{equation}

We now convert this logarithmic derivative into the mean $\E D(z)$.
By definition,
\[
 -\theta\E D(z)
 =n\int_0^{1/n}\E\sum_i e_i\left[
 \frac{\sfb'(z)}{\sfb(z)-\sfb(x_i+\theta/n+t)}
 -\frac{\sfb'(z)}{\sfb(z)-\sfb(x_i+t)}
 \right]\rd t.
\]
Following the argument of \cite[Proposition~4.16, (147)--(150)]{Gorin:2022qky}, with the remainders controlled by \eqref{e:partialzlnA} and \Cref{lem:dmc-velocity-derivatives}, we obtain
\begin{equation}\label{eq:dmc-physical-mean-conversion}
 -\E D(z)
 =\frac1\theta\left.\partial_\zeta\log\cA^{(z,1)}(\zeta)\right|_{\zeta=z}
 +\frac12\int_{\mathbb R}K_{\sfb}(x,z)\,
 \E[\mu(\rd x;\bmy)-\mu(\rd x;\bmx)]
 +\OO\left(\frac{P(z)}{n^2d^2}\right).
\end{equation}
The same argument gives, to first order,
\[
\begin{aligned}
 \int_{\mathbb R}K_{\sfb}(x,z)\,
 \E[\mu(\rd x;\bmy)-\mu(\rd x;\bmx)]
 &=-\frac1{2\pi\ri n}\oint_{\cin}K_{\sfb}(w,z)\E D(w)\,\rd w\\
 &=-\frac1{2\pi\ri\theta n}\oint_{\cin}
 \frac{\cB'(w)}{\cB(w)}K_{\sfb}(w,z)\,\rd w
 +\OO\left(\frac{P(z)}{n^2d^2}\right),
\end{aligned}
\]
where $\cin$ encloses both the current and next-step supports but not $z$.

Substituting this identity and \eqref{eq:dmc-second-order-insertion} into \eqref{eq:dmc-physical-mean-conversion}, the $\cB'/(2\cB)$ terms cancel.
Since $a(z)=-\theta^{-1}\partial_z h_-(\sfb(z))$, we conclude that
\[
 \E D(z)=a(z)-\frac1{2\pi\ri\theta n}\oint_{\cin}
 \cE_2(w,z)K_{\sfb}(w,z)\,\rd w
 +\OO\left(\frac{P(z)}{n^2d^2}\right).
\]
This proves \eqref{t:loopstudy2} for $p=0$.
The remainder is holomorphic on a neighborhood of $z$ of radius comparable to $d$, where the same estimate holds.
Local comparability of $P$ and Cauchy's formula therefore give the asserted derivative bounds for every fixed $p\ge0$.
\end{proof}

\paragraph{Conditional covariance}

\begin{proposition}\label{p:enhance martingale}
Assume \Cref{ass:weightfunc,ass:prelim-log-branches} and the domain
and contour conditions above. For $s=s_j\le T$ and
$z,u\in\mathscr D_s\setminus\mathbb R$ with
$\dist(z,I_s)\asymp\dist(u,I_s)\asymp d\ge n^{-1+\delta}$,
\begin{equation}\label{e:quadramartin dmc}
 \E_s[\Delta\mathcal M_s^n(z)\Delta\mathcal M_s^n(u)]
 =\frac1{2\pi\ri\theta n}\oint_{\cin}
 \chi_s^n(w)K_{\sfb}(w,z)K_{\sfb}(w,u)\,\rd w
 +\OO\left(\frac{P_s^n(z)+P_s^n(u)}{n^2d^3}\right).
\end{equation}
Here $\cin$ encloses $I_s$ but neither spectral parameter.
For fixed integers $a,b\ge0$, the $\partial_z^a\partial_u^b$ derivative
of the remainder is bounded by
$C_{a,b}(P_s^n(z)+P_s^n(u))/(n^2d^{a+b+3})$.
Coincident and conjugate parameters are allowed; in particular,
$\E_s|\Delta\mathcal M_s^n(z)|^2\lesssim P_s^n(z)/(nd^2)$.
\end{proposition}

\begin{proof}
We suppress the time index and follow
\cite[Proposition~4.18, (160)--(166)]{Gorin:2022qky}.
Recall $Q_{\bme}$ from the proof of \Cref{cor:one step result}.
For $(\bmv,\bms)=((z,u,v),(1,1,-1))$, differentiating in
$\zeta$ with the insertion parameters fixed gives
\[
 \left.\partial_\zeta\log\cA^{(\bmv,\bms)}(\zeta)\right|_{\zeta=z}
 =\frac{\E[\partial_z\log Q_{\bme}(z)\,Q_{\bme}(v)/Q_{\bme}(u)]}
        {\E[Q_{\bme}(v)/Q_{\bme}(u)]}.
\]
By \Cref{p: weak estimate}, the denominator is nonzero
for $|v-u|\le cd$. Differentiating at $v=u$ gives the covariance.

We first extend \eqref{eq:dmc-Mv-improved} to these insertions.
Write $[\xi_0,\ldots,\xi_q]_{\sfb}f$ for the divided difference
of $f\circ\sfb^{-1}$ at $\sfb(\xi_0),\ldots,\sfb(\xi_q)$,
with coincident nodes defined by continuity.
Uniformly for $\xi_0,\ldots,\xi_q$ in sufficiently small
$cd$-neighborhoods of $z$ or $u$, we have
\begin{equation}\label{eq:dmc-mixed-divdiff}
 |[\xi_0,\ldots,\xi_q]_{\sfb}\chi^j|
 \le C^{j+q}\frac{P(z)+P(u)}{d^{q-1}},
 \qquad j,q\ge1.
\end{equation}
For $q=1$, this follows from Cauchy--Schwarz in
\eqref{eq:dmc-velocity-positive-representation} and
$\chi=1-\cB^{-1}$.
Higher orders follow from the divided-difference recursion
for separated nodes and Cauchy's formula for nearby nodes.
The estimate also holds after multiplication by analytic
factors with uniformly bounded derivatives.

Define $L_A,L_G,L$ as above
\eqref{e:error-term-2-clean: deform}, using the present
$\cA^{(\bmv,\bms)}$.
Then \eqref{e:error-term-1-clean: deform} gives
\[
 1+\frac1n\cE_1^{(\bmv,\bms)}(w)
 =
 1-\chi(w)+L(w)\prod_{i=1}^3
 \left(
 \frac{\sfb(v_i)-\sfb(w)}
      {\sfb(v_i)-\sfb(w+\theta/n)}
 \right)^{s_i}.
\]
Fix an integer $J$ such that $\delta(J+1)>2$.
For $j=0,\ldots,J$, set
\[
 M_j=
 \sup_{\min\{|w-z|,|w-u|\}\le (2+j/J)cd}
 \left|
 \partial_w\bigl(
 \log\cA^{(\bmv,\bms)}(w)-h_-(\sfb(w))
 \bigr)
 \right|.
\]

In the logarithm of the previous expression, first replace
$L$ by $\chi$.
The pole decomposition in the proof of
\Cref{cor:one step result} and \eqref{eq:dmc-mixed-divdiff}
bound the contribution of order $n^{-k}$ to the derivative
Cauchy integral by
$C^k(P(z)+P(u))/(nd)^k$.
The difference of the two logarithms is bounded by
$C(P(z)+P(u)+M_{j+1})/n$.
Since the local contours have length $\OO(d)$ and the
derivative Cauchy kernel is $\OO(d^{-2})$,
the argument of \Cref{lem:dmc-Mv-improvement} gives
\[
 M_j\le\frac{C(P(z)+P(u))}{nd}
       +\frac{C}{nd}M_{j+1},
 \qquad
 M_j\lesssim\frac1{nd^2}.
\]
Here the second bound follows from \Cref{p: weak estimate}
and Cauchy's formula.
On $\cout$, the argument in the proof of
\Cref{cor:one step result}, using \Cref{p: weak estimate}
with $R\ge4$, gives a contribution $\OO(n^{-1})$
and $L_A-1=\OO(n^{-2})$.
Iterating this inequality $J$ times,
the last term $(C/(nd))^J/(nd^2)$ is absorbed by
$(P(z)+P(u))/(nd)$.
Thus, on the smaller neighborhoods,
\begin{equation}\label{eq:dmc-multi-insertion-improved}
 \left|
 \partial_w\bigl(\log\cA^{(\bmv,\bms)}(w)-h_-(\sfb(w))\bigr)
 \right|
 \lesssim\frac{P(z)+P(u)}{nd}.
\end{equation}

Repeating the expansion leading to
\eqref{eq:dmc-second-order-insertion} and subtracting the
formula for $\cA^{(z,1)}$, the terms independent of $u,v$
cancel. Uniformly for $|v-u|\le cd$, we obtain
\begin{equation}\label{eq:dmc-paired-response}
\begin{aligned}
 &\left.\partial_\zeta
 \left(\log\cA^{(\bmv,\bms)}(\zeta)
       -\log\cA^{(z,1)}(\zeta)\right)\right|_{\zeta=z}\\
 &\quad=\frac{\theta}{2\pi\ri n}\oint_{\cin}
 \chi(w)\sfb'(w)
 \left(
 \frac1{\sfb(u)-\sfb(w)}
 -\frac1{\sfb(v)-\sfb(w)}
 \right)
 K_{\sfb}(w,z)\,\rd w
 +\OO\left(\frac{P(z)+P(u)}{n^2d^2}\right).
\end{aligned}
\end{equation}
Indeed, \eqref{eq:dmc-multi-insertion-improved} gives
$L_A-1=\OO((P(z)+P(u))/(n^2d))$ locally.
Using $d(P(z)+P(u))\lesssim1$, the higher-order terms
containing $L-\chi$ are also
$\OO((P(z)+P(u))/(n^2d))$ before integration, as in the
proof of \Cref{cor:one step result}.
The terms from the insertion factors and the shifts by
$\theta/n$ are bounded after integration by
\eqref{eq:dmc-mixed-divdiff}.
On $\cout$, the error is $\OO(n^{-2})$.

The error term in \eqref{eq:dmc-paired-response} is
holomorphic in $v$ for $|v-u|<cd$.
Differentiating at $v=u$, with the contours fixed,
and applying Cauchy's formula give
\begin{equation}\label{eq:dmc-logratio-cov}
\begin{aligned}
 &\operatorname{Cov}\bigl(
 \partial_z\log Q_{\bme}(z),
 \partial_u\log Q_{\bme}(u)\bigr) =
 \frac{\theta}{2\pi\ri n}\oint_{\cin}
 \chi(w)K_{\sfb}(w,z)K_{\sfb}(w,u)\,\rd w
 +\OO\left(\frac{P(z)+P(u)}{n^2d^3}\right).
\end{aligned}
\end{equation}
Taking $u=\bar z$, \eqref{eq:dmc-mixed-divdiff} bounds
the integral by $CP(z)/d^2$.
Cauchy's formula then gives, for $q=0,1$,
\[
 \bigl\|
 \partial_z^{q+1}\log Q_{\bme}(z)
 -\E\partial_z^{q+1}\log Q_{\bme}(z)
 \bigr\|_2
 \lesssim\frac{\sqrt{P(z)}}{\sqrt n\,d^{q+1}}.
\]

Finally, the computation in the proof of
\Cref{cor:one step result} and the residue theorem give
\begin{equation}\label{eq:dmc-exact-step-transfer}
\begin{aligned}
 -\theta D(z)=n\int_0^{1/n}\Bigg[
 &\left.\partial_\zeta\log Q_{\bme}(\zeta)
   \right|_{\zeta=z-t} +\frac1{2\pi\ri}\oint_{\omega_+}
 \partial_w\log Q_{\bme}(w)\,\sfb'(z)
 \left(
 \frac1{\sfb(z)-\sfb(w+t)}
 -\frac1{\sfb(z)-\sfb(w)}
 \right)\rd w
 \Bigg]\rd t.
\end{aligned}
\end{equation}
The kernel in the contour integral has residues $-1$ and $1$
at $w=z-t$ and $w=z$, respectively, and is $\OO(t)$
on $\omega_+$.
Subtracting expectations and using the preceding $L^2$
bounds and Minkowski's inequality, we obtain
\[
 \left\|
 \Delta\mathcal M(z)+\frac1\theta
 \bigl(\partial_z\log Q_{\bme}(z)
       -\E\partial_z\log Q_{\bme}(z)\bigr)
 \right\|_2
 \lesssim\frac{\sqrt{P(z)}}{n^{3/2}d^2}.
\]
Cauchy--Schwarz and \eqref{eq:dmc-logratio-cov} prove
\eqref{e:quadramartin dmc} and the variance bound.
The mixed derivative bounds follow from Cauchy's formula
on slightly larger neighborhoods of radius comparable to $d$
and local comparability of $P$.
\end{proof}

\paragraph{Higher conditional moments}

\begin{proposition}\label{prop:dmc-higher-moments}
Assume \Cref{ass:weightfunc,ass:prelim-log-branches} and the domain
and contour conditions above. For $s=s_j\le T$,
$z\in\mathscr D_s\setminus\mathbb R$, and fixed integers $p\ge1$, $q\ge0$,
\begin{equation}\label{eq:dmc-higher-conditional-moments}
 \E_s|\partial_z^q\Delta\mathcal M_s^n(z)|^{2p}
 \lesssim_{p,q}
 \left(\frac{P_s^n(z)}{nd^{2q+2}}\right)^p
 +\frac{P_s^n(z)}{n^{2p-1}d^{2pq+4p-2}},
 \qquad d=\dist(z,I_s)\ge n^{-1+\delta}.
\end{equation}
\end{proposition}

\begin{proof}
We suppress the time indices and repeat the insertion argument in
the proofs of \Cref{cor:one step result} and \Cref{p:enhance martingale}.
Taking $R\ge2p$ in \Cref{p: weak estimate} allows all the required
insertions and ensures nonvanishing on the local parameter
neighborhoods.

For $2\le r\le2p$ and $u_1,\ldots,u_r\in\{z,\bar z\}$, we obtain
\[
 \left|
 \left.\partial_{v_1}\cdots\partial_{v_r}
 \log\E\prod_{j=1}^r
 \frac{Q_{\bme}(v_j)}{Q_{\bme}(u_j)}
 \right|_{\bm v=\bm u}
 \right|
 \lesssim_r\frac{P(z)}{n^{r-1}d^{2r-2}}.
\]
Indeed, differentiate the exact deformed loop equation
\eqref{e:error-term-2-clean: deform} before setting $\bm v=\bm u$.
Each derivative of an explicit deformation factor is
$\OO((nd^2)^{-1})$, and the cancellation in
\eqref{eq:dmc-pole-cancellation} retains $P(z)$ after projection.
The highest derivative entering through $\log L_A$ is integrated
over $[w,w+\theta/n]$, so its projection carries the factor
$\OO((nd)^{-1})$; all remaining terms involve lower-order
derivatives. Induction, starting from the weak estimate and using
the same iteration as in the covariance proof, gives the bound above.

These mixed logarithmic derivatives are the joint cumulants of
$\partial_{u_j}\log Q_{\bme}(u_j)$. The moment--cumulant formula
for the centered variables therefore gives
\[
\begin{aligned}
 \E\left|\partial_z\log Q_{\bme}(z)
          -\E\partial_z\log Q_{\bme}(z)\right|^{2p}
 &\lesssim_p
 \sum_{a=1}^p\frac{P(z)^a}{n^{2p-a}d^{4p-2a}} \lesssim_p
 \left(\frac{P(z)}{nd^2}\right)^p
 +\frac{P(z)}{n^{2p-1}d^{4p-2}}.
\end{aligned}
\]
By \eqref{eq:dmc-exact-step-transfer} and Minkowski's inequality,
the same estimate holds for $\Delta\mathcal M(z)$.
Indeed, the fixed outer-contour contribution is $\OO(n^{-1})$
in $L^{2p}$, since its kernel is $\OO(t)$ and
$\partial_w\log Q_{\bme}(w)=\OO(1)$ there; its $2p$-th moment
is absorbed by the first term.
Finally, Cauchy's formula on $|w-z|=cd$ and local comparability
of $P$ give \eqref{eq:dmc-higher-conditional-moments}
for every fixed $q\ge0$.
\end{proof}

\begin{proof}[Proof of \Cref{thm:main}]
The first-order estimate in the mean proof, \Cref{cor:one step result},
and \Cref{p:enhance martingale,prop:dmc-higher-moments} prove
the \Cref{thm:main}.
\end{proof}

\subsection{Weighted \texorpdfstring{$\beta$}{beta}--Corner process}\label{sec:corners}
We prove the corners case of \Cref{thm:main}, following the
arguments for DMC in Section~3.2 and recording the changes for corners.
We retain \Cref{ass:weightfunc,ass:prelim-log-branches} and the
domain and contour conditions of \Cref{sec:main}.

Fix $s=j/n\le T<1$ and suppress the time indices, writing
$I=I_s$, $\Lambda=\Lambda_s$ and $\mathscr D=\mathscr D_s$.
We use the normalization \eqref{eq:corner-level-measures}:
the current row $\bmx$ has $k+1=n-j$ particles and the next row
$\bmy$ has $k=n-j-1$ particles.
All expectations below are conditional on $\bmx$.
Since $\sfb(z)=z$, \eqref{e:B_function},
\eqref{eq:velocity-from-flux} and \eqref{eq:kernel-common} give
\[
 \cB(z)=\frac{k}{n}W'(z)-\theta m(z),\qquad
 \chi(z)=-\frac{\theta}{\cB(z)},\qquad
 K_{\sfb}(w,z)=\frac1{(w-z)^2}.
\]
All estimates below are uniform in $s,n$.

\paragraph{Preliminary estimates and deformation.}
We use the characteristic-polynomial ratio corresponding to
$Q_{\bme}$ in the proof of \Cref{cor:one step result}:
\begin{equation}\label{eq:corners-ratio}
 R_{\bmy}(z)=\frac{\prod_{a=1}^{k+1}(z-x_a)}
                  {\prod_{i=1}^{k}(z-y_i)},
 \qquad D(z)=\partial_z\log R_{\bmy}(z).
\end{equation}
The second identity follows directly from \eqref{eq:corner-level-measures} and \eqref{eq:corners-forward-increment}.
We abbreviate $R_{\bmy}$ to $R$.

For the contour estimates below, we require a positively oriented
outer contour $\Gamma$ enclosing $I$ and the spectral neighborhoods
of \Cref{sec:main}, with
$\overline{\operatorname{Int}\Gamma}\subset\Lambda$.
The contour stays a fixed positive distance from $I$ and these
neighborhoods, and satisfies $|\cB|\asymp1$ on a fixed neighborhood
of $\Gamma$ and $\operatorname{wind}_{\Gamma}\cB=-1$.
All these bounds and separation constants are uniform in $s,n$
and the admissible current particle configuration.

The coefficient $\cB$ has $k+1$ simple poles at $x_a$, each with
residue $\theta/n$, and opposite signs at the endpoints of each gap.
The winding condition therefore gives exactly one simple zero
$\lambda_i\in(x_{i+1},x_i)$ and no other zeros inside $\Gamma$.
Define
\begin{equation}\label{eq:corners-reference-factor}
 H(z)=\frac{\prod_{a=1}^{k+1}(z-x_a)}
           {\prod_{i=1}^{k}(z-\lambda_i)}.
\end{equation}
By \eqref{e:F_CP} and \eqref{eq:leading-common} and the residue theorem,
\[
 a(z)=\frac1{2\pi\ri}\oint_{\cin}
       \frac{\cB'(w)}{\cB(w)(w-z)}\,\rd w
     =\frac{H'(z)}{H(z)}.
\]
This agrees with \eqref{eq:corners-factor-derivative};
in the notation there, $H=(z-e)e^{h_-}$.

Interlacing gives $|H|\asymp1$ near $\Gamma$.
Since $H$ cancels the zeros and poles of $\cB$ inside $\Gamma$,
the maximum principle applied to $\cB H$ and its reciprocal gives
$|\cB H|\asymp1$ there.
Together with \Cref{ass:prelim-log-branches}, this yields
$|H|\asymp1$ on the spectral neighborhoods.
Moreover, $a+\cB'/\cB=\partial_z\log(\cB H)$ has bounded derivatives
of every fixed order on a fixed smaller domain.

The residues $-\theta/\cB'(\lambda_i)$ of $\chi$ are positive.
The positive-residue argument in the proof of
\Cref{lem:dmc-velocity-derivatives} therefore gives
\begin{equation}\label{eq:corners-velocity-jets}
 P(z)\lesssim d^{-1},\qquad
 |\partial_z^q\chi(z)|\lesssim_q\frac{P(z)}{d^{q-1}}
 \quad(q\ge1),\qquad d=\dist(z,I),
\end{equation}
and local comparability of $P$ on $|w-z|\le cd$.
At real points outside $I$, use the continuous extension of $P$.

For a fixed number of insertions $\bmv=(v_1,\ldots,v_r)$ and
$\bms\in\{\pm1\}^r$, define
\begin{equation}\label{e:phi-tilde-def-Corner}
 \cA^{(\bmv,\bms)}(z)
 =\frac{\E[R(z)\prod_iR(v_i)^{-s_i}]}
        {\E\prod_iR(v_i)^{-s_i}},\qquad
 W^{(\bmv,\bms)}(w)
 =W(w)+\frac1k\sum_i s_i\log(v_i-w).
\end{equation}
For all sufficiently large $n$, the normalization factor is nonzero
by \Cref{prop:corners-source-expansion}.
By \eqref{e:B_function}, the deformed coefficient is
$\cB^{(\bmv,\bms)}(w)=\cB(w)+n^{-1}\sum_i s_i/(w-v_i)$.
Applying \Cref{lem:dyn-loop-one-step} to the deformed potential,
as in the proof of \Cref{prop:one-step-mean}, gives
\begin{equation}\label{eq:corners-deformed-loop}
 \cC^{(\bmv,\bms)}(w)
 =\left(\cB(w)+\frac1n\sum_i\frac{s_i}{w-v_i}\right)
    \cA^{(\bmv,\bms)}(w)
 -\frac1n\partial_w\cA^{(\bmv,\bms)}(w).
\end{equation}
Here $\cC^{(\bmv,\bms)}$ is holomorphic inside $\cin$, which
encloses $I$ and excludes the insertions.
\paragraph{Weak and improved estimates.}
The following proposition gives the corners analogues of
\Cref{p: weak estimate,p:edge impr} and the higher-order expansion
used in the moment estimates.

\begin{proposition}\label{prop:corners-source-expansion}
Assume \Cref{ass:weightfunc,ass:prelim-log-branches} and the contour
conditions above. Fix integers $r\ge0$, $p\ge1$, and
$\bms\in\{\pm1\}^r$.
Suppose that $z,v_1,\ldots,v_r$ lie in the admissible domain and
have distances to $I$ comparable to $d\ge n^{-1+\delta}$.
Put $P_*=\max\{P(z),P(v_1),\ldots,P(v_r)\}$.
For all sufficiently large $n$, uniformly as these parameters
vary independently within distance $cd$ of their given values,
\[
 \left|\E\prod_iR(v_i)^{-s_i}\right|\asymp1,\qquad
 \log\frac{\cA^{(\bmv,\bms)}(z)}{H(z)}
 =\sum_{j=1}^p\frac{\ell_j^{(\bmv,\bms)}(z)}{n^j}
  +\OO\left(\frac{P_*}{n^{p+1}d^p}\right).
\]
The logarithm is chosen near zero, and
$|\ell_j^{(\bmv,\bms)}(z)|\lesssim_j P_*d^{1-j}$.
The first coefficient is
\begin{equation}\label{eq:corners-first-coefficient}
 \ell_1^{(\bmv,\bms)}(z)
 =\frac1{2\pi\ri\theta}\oint_{\cin}
 \frac{\chi(w)}{w-z}
 \left(a(w)-\sum_i\frac{s_i}{w-v_i}\right)\,\rd w.
\end{equation}
Here $\ell_j$ is a sum of products of iterated contour integrals
containing at most $j$ factors of $\sum_i s_i/(w-v_i)$.

Each derivative in $z$ or $v_j$
costs a factor $d^{-1}$ in the remainder bounds. In particular,
\begin{equation}\label{eq:corners-improved-insertion}
 \left|\partial_z
 \log\frac{\cA^{(\bmv,\bms)}(z)}{H(z)}\right|
 \lesssim\frac{P_*}{nd}.
\end{equation}
\end{proposition}

\begin{proof}
Put $M=\E\prod_iR(v_i)^{-s_i}$ and
$F(z)=H(z)^{-1}\E[R(z)\prod_iR(v_i)^{-s_i}]$.
Apply the loop identity of \Cref{lem:dyn-loop-one-step} with
the factor $\prod_iR(v_i)^{-s_i}$ inserted.
The Cauchy projection used in \eqref{eq:corners-factor-exact}
then gives
\begin{equation}\label{eq:corners-integral-equation}
 F(z)=M+\frac1{2\pi\ri\theta n}\oint_{\cin}
 \frac{\chi(w)}{w-z}
 \left[F'(w)+\left(a(w)-\sum_i\frac{s_i}{w-v_i}\right)F(w)\right]
 \,\rd w.
\end{equation}

We successively replace $F$ in the integral term of
\eqref{eq:corners-integral-equation} by the right-hand side
of that equation, obtaining a finite expansion in powers of $n^{-1}$.
The integral evaluated at $F=1$ is $\OO(P_*)$.
Indeed, the identity
$\chi a=\chi(a+\cB'/\cB)+\chi'$ and
\eqref{eq:corners-velocity-jets} control the term containing $a$.
For the insertion terms, apply Cauchy--Schwarz to the positive
residues of $\chi$, as in the proof of
\eqref{eq:dmc-mixed-divdiff}. For each subleading term, integration by parts and the
contour deformation used in the proof of \Cref{p:one step estimate}
give a bound of $C/d$ times the supremum of its input on slightly
larger spectral neighborhoods.
Thus the $j$th iterated integral applied to $1$ is
$\OO_j(P_*d^{1-j})$, and its contribution to the expansion is
$\OO_j(P_*/(n^jd^{j-1}))$.

Interlacing and $|H|\asymp1$ give $|F|\lesssim_r d^{-r-1}$.
Each substitution reduces the terminal bound by a factor
$C/(nd)$. Since $nd\ge n^\delta$, a fixed number of substitutions,
depending only on $p,r,\delta$, yields
\[
 F(z)=M\left(1+\OO\left(\frac{P_*}{n}\right)\right)
       +\OO\left(\frac{P_*}{n^{p+1}d^p}\right).
\]
If $s_i=1$, then
$H(v_i)F(v_i)=\E\prod_{j\ne i}R(v_j)^{-s_j}$.
If $s_i=-1$, apply this expansion to the shorter insertion list
at $z=v_i$.
Since $|H(v_i)|\asymp1$ and
$P_*/n\lesssim(nd)^{-1}=\oo(1)$, induction from the empty list
gives $|M|\asymp1$.

The finite expansion of $F/M=\cA^{(\bmv,\bms)}/H$ therefore has
remainder $\OO(P_*/(n^{p+1}d^p))$ when truncated at order $p$.
Expanding its logarithm gives the asserted expansion.
The first integral gives \eqref{eq:corners-first-coefficient},
and $dP_*\lesssim1$ preserves the coefficient bounds under
products arising from the logarithm.
Each iteration introduces at most one factor of the insertion
sum, so the coefficient of $n^{-j}$ has total insertion degree
at most $j$.

The estimates hold on slightly larger parameter neighborhoods.
Local comparability of $P$ and Cauchy's formula give the derivative
bounds, as in the proof of \Cref{p:edge impr}, and imply
\eqref{eq:corners-improved-insertion}.
\end{proof}

\paragraph{One-step estimates.}
We use the insertion identities from
\Cref{lem:corners-inverse-tilt} and the proofs of
\Cref{p:enhance martingale,prop:dmc-higher-moments},
with the remainders supplied by
\Cref{prop:corners-source-expansion}.

\begin{proposition}\label{prop:corners-one-step}
Under the preceding assumptions, for $s=j/n\le T<1$ and
$d=\dist(z,I_s)\ge n^{-1+\delta}$,
\begin{equation}\label{eq:corners-second-order-mean}
 \E_sD_s^n(z)=a_s^n(z)+\frac1n c_s^n(z)+\mathcal R_s^n(z),
 \qquad
 c_s^n(z)=\frac1{2\pi\ri\theta}\oint_{\cin}
 \frac{\chi_s^n(w)a_s^n(w)-\tfrac12\partial_w\chi_s^n(w)}
      {(w-z)^2}\,\rd w.
\end{equation}
The contour defining $a_s^n(w)$ lies strictly inside $\cin$.
For every fixed integer $q\ge0$,
\[
 |\partial_z^q\mathcal R_s^n(z)|
 \lesssim_q\frac{P_s^n(z)}{n^2d^{q+2}},\qquad
 |\partial_z^qc_s^n(z)|
 \lesssim_q\frac{P_s^n(z)}{d^{q+1}}.
\]
For $\dist(z,I_s)\asymp\dist(u,I_s)\asymp d$,
\begin{equation}\label{eq:corners-sharp-covariance}
 \E_s[\Delta\mathcal M_s^n(z)\Delta\mathcal M_s^n(u)]
 =\frac1{2\pi\ri\theta n}\oint_{\cin}
 \frac{\chi_s^n(w)}{(w-z)^2(w-u)^2}\,\rd w
 +\OO\left(\frac{P_s^n(z)+P_s^n(u)}{n^2d^3}\right).
\end{equation}
For fixed integers $a,b\ge0$, the mixed derivative
$\partial_z^a\partial_u^b$ of the remainder is bounded by
$C_{a,b}(P_s^n(z)+P_s^n(u))/(n^2d^{a+b+3})$.
Coincident and conjugate parameters are allowed.
For fixed integers $p\ge1$ and $q\ge0$,
\begin{equation}\label{eq:corners-higher-moment-bound}
 \E_s|\partial_z^q\Delta\mathcal M_s^n(z)|^{2p}
 \lesssim_{p,q}
 \left(\frac{P_s^n(z)}{nd^{2q+2}}\right)^p
 +\frac{P_s^n(z)}{n^{2p-1}d^{2pq+4p-2}}.
\end{equation}
\end{proposition}

\begin{proof}
We suppress the time indices.
For the mean, \eqref{eq:corners-tilt-diagonal} gives
$\E D(z)=
\left.\partial_\zeta\log\cA^{(z,1)}(\zeta)\right|_{\zeta=z}$,
with the insertion fixed during differentiation.
Apply \Cref{prop:corners-source-expansion} with $p=1$.
Differentiating \eqref{eq:corners-first-coefficient} yields
\[
 \E D(z)=a(z)+\frac1{2\pi\ri\theta n}\oint_{\cin}
 \left[\frac{\chi(w)a(w)}{(w-z)^2}
       -\frac{\chi(w)}{(w-z)^3}\right]\rd w
 +\OO\left(\frac{P(z)}{n^2d^2}\right).
\]
Integration by parts gives
\eqref{eq:corners-second-order-mean}, in agreement with
\eqref{eq:corners-explicit-mean}.
For the coefficient bound, use
$\chi a-\chi'/2=\chi(a+\cB'/\cB)+\chi'/2$
and \eqref{eq:corners-velocity-jets}.
Cauchy's formula gives the remainder derivative bounds.

For covariance, use the insertions
$((z,u,v),(1,1,-1))$ as in the proof of
\Cref{p:enhance martingale}.
By definition, it gives
\[
 \left.
 \partial_v\partial_\zeta
 \log\cA^{((z,u,v),(1,1,-1))}(\zeta)
 \right|_{\zeta=z,\,v=u}
 =\operatorname{Cov}(D(z),D(u)).
\]
The insertion-dependent term in
\eqref{eq:corners-first-coefficient} gives the integral in
\eqref{eq:corners-sharp-covariance}.
The two derivatives of the remainder cost $d^{-2}$,
giving the stated error.

For higher moments, first differentiate
$\log\E\prod_{i=1}^rR(v_i)/R(u_i)$ in $v_1$ and set $v_1=u_1$.
Cancellation as in \eqref{eq:corners-tilt-diagonal} expresses
the result as a logarithmic spectral derivative of $\cA$
with $2r-1$ insertions.
Apply \Cref{prop:corners-source-expansion} with $p=r-1$.
The remaining $r-1$ parameter derivatives must act on $r-1$
factors of the insertion sum, so all coefficients with
$j<r-1$ vanish.
The coefficient and remainder bounds therefore give, for $r\ge2$,
\[
 \left|
 \left.\partial_{v_1}\cdots\partial_{v_r}
 \log\E\prod_{i=1}^r\frac{R(v_i)}{R(u_i)}
 \right|_{\bmv=\bmu}
 \right|
 \lesssim_r\frac{P_*}{n^{r-1}d^{2r-2}}.
\]
The remainder is smaller by a factor $(nd)^{-1}$.
Taking $u_i\in\{z,\bar z\}$ and expanding the centered moment
as in the proof of \Cref{prop:dmc-higher-moments} gives
\eqref{eq:corners-higher-moment-bound} for $q=0$.
Since $D=R'/R$ exactly, no further conversion is required.
Cauchy's formula on the local $d$-scale neighborhoods and local
comparability of $P$ give the remaining derivative estimates.
\end{proof}

Together with Section~3.2, this proves \Cref{thm:main}
under the stated contour conditions.

\section{Evolution of Stieltjes Transform}\label{s:eost}
In this section, for three particle systems: DMC, weighted $\beta$-corners,
and Dyson Brownian motion, we establish an upper bound on the rightmost
particle $x_1(s)$ relative to the edge $E_s$, together with bulk
estimates for $m_s^n-m_s$.
The proof combines the one-step estimates of Section~3 with
the characteristic flow introduced in
\Cref{def:prelim-characteristic-flow}.
These estimates will be used in Section~5 to prove tightness
of the rescaled Nevanlinna field.

The main idea is that we estimate $m_s^n-m_s$ separately: away from the support
and near the edge.
Away from the support, we have the
rough bounds in \Cref{s4:outer-log}, which give the
initial input for the local characteristics.
Near the edge, we use the square-root behavior of $\mu_s$
to estimate the characteristic flow and obtain more precise
bounds on the error terms and martingales.
Combining these estimates with Gr\"onwall's inequality gives
\Cref{lem:s4-edge-improvement,lem:s4-bulk-characteristic}.

Fix a any time $T_0=\OO(1)>0$, with $T_0<1$ for the corners process.
Throughout this section, the estimates of Sections~2 and~3 are
applied with $T=T_0$.
In the below statements, time variables are
restricted to $\mathbb Z_n=n^{-1}\mathbb Z_{\ge0}$ for DMC
and corners, and range continuously for DBM.

\subsection{Main results and assumptions}\label{subsec:edge-main-results}
We retain the notation $\mu_s^n,m_s^n$ from
\Cref{subsec:prelim-measures}, and $\mu_s,m_s,\cF_s$ and $z_s(u)$
from \Cref{subsec:prelim-characteristics}.
Write $x_1(s)$ for the rightmost particle, with
$x_1(s)=\lambda_1(s)$ for DBM, and let
$E_s=\sup\operatorname{supp}\mu_s$ denote the deterministic
right edge.

\begin{assumption}[Regular edge and initial data]
\label[assumption]{ass:edge-data}
There is $s_*\in[0,T_0)$ such that the following hold.
\begin{enumerate}[label=\textup{(\roman*)},leftmargin=*,itemsep=4pt]
\item For every $a\in(s_*,T_0)$, there are $r_*,c_*,C_*>0$
such that, uniformly for $s\in[a,T_0]$, $\mu_s$ has a density on
$(E_s-r_*,E_s)$ satisfying
\begin{equation}\label{def: classical measure}
 \rho_s(E_s-x)=\frac{\sqrt{\cA(s)}}{\pi}\sqrt{x}\,(1+\OO(x)),
 \qquad 0<x<r_*,\qquad c_*\le\cA(s)\le C_*.
\end{equation}
The constants may depend on $a$, but not on $n$.
\item For every such $a$, the flux $\cF_s$ is holomorphic and
uniformly bounded on a fixed complex neighborhood of
$\{(m_s(E_s),E_s):a\le s\le T_0\}$, and
\begin{equation}\label{eq:edge-flux-nondegeneracy}
 \partial_{11}\cF_s(m_s(E_s),E_s)\ge c_*,\qquad a\le s\le T_0.
\end{equation}
\item The initial configurations may be random. There is a fixed
compact interval $[\ell_0,r_0]$ such that
$\mathbb P(I_0\subset[\ell_0,r_0])\to1$, and
\begin{equation}\label{eq:edge-initial-preparation}
 (\sup I_0-E_0)_+\xrightarrow{\mathbb P}0,
 \qquad
 \sup_{z\in\mathcal K}|m_0^n(z)-m_0(z)|
 =\OO_{\mathbb P}(n^{-1})
\end{equation}
for every fixed compact set
$\mathcal K\Subset\Lambda\setminus[\ell_0,E_0]$.
\end{enumerate}
\end{assumption}

\begin{remark}
For DBM, we also assume the corresponding conditions at
$L_s=\inf\operatorname{supp}\mu_s$ on each $[a,T_0]$, $a>s_*$,
and
\(
 (L_0-\inf I_0)_+\xrightarrow{\mathbb P}0.
\)
In \Cref{ass:edge-data}\textup{(iii)}, replace
$[\ell_0,E_0]$ by $[L_0,E_0]$. The characteristic-flow and support estimates below also
apply at the left edge by the same argument. For simplicity, later on, we will not mention the analysis near left edge.

For DMC, $\sup I_0=x_1(0)+\theta/n$, so the first condition in
\eqref{eq:edge-initial-preparation} is equivalent to
$(x_1(0)-E_0)_+\to0$ in probability.
\end{remark}

In \Cref{app:deterministic-proofs}, we give sufficient conditions on
the model data and the initial measure $\mu_0$ under which
parts~\textup{(i)} and~\textup{(ii)} of \Cref{ass:edge-data} hold
at the right edge after a waiting time.
We apply \Cref{thm:app-edge-formation} to concrete models from the
introduction and give examples of admissible initial measures in
\Cref{rem:app-initial-tests}.

\begin{theorem}[Upper bound at the edge]\label{large-eigen: DMC}
Under \Cref{ass:weightfunc,ass:prelim-log-branches,ass:edge-data},
fix $s_*<T<T_0$ and a sufficiently large constant $M'>0$.
There are constants $C,c,\delta>0$ such that, for every sufficiently
large fixed $M\ge M'$ and all $n\ge n_0(M)$,
\begin{equation}\label{eq:edge-window-bound}
 \mathbb P\left(
 \sup_{|s-T|\le M'n^{-1/3}}(x_1(s)-E_s)
 \le CM^2n^{-2/3}\right)\ge1-Ce^{-cM^\delta}.
\end{equation}
The constants may depend on $T,M'$ and the fixed model and
regularity bounds, but not on $M,n$.
The threshold $n_0(M)$ may also depend on the initial preparation.
\end{theorem}

For the bulk estimate, set
\begin{equation}\label{def: bulk region}
 R(K)=\left\{w\in\mathbb H:|w|\le n^{1/10},\quad
 \Im w\ge K^2\sqrt{1+\max\{\Re w,0\}}\right\},\qquad K\ge1.
\end{equation}

\begin{theorem}[Bulk estimate]\label{bulktight: DMC}
Under \Cref{ass:weightfunc,ass:prelim-log-branches,ass:edge-data},
fix a compact interval $J\Subset(s_*,T_0)$.
There are $C,c,K_0>0$ and $0<\alpha\le\alpha'<1$ such that,
for every fixed $K\ge K_0$, all $n\ge n_0(K)$, and every
deterministic time $s\in J$, with probability at least
$1-Ce^{-cK^\alpha}$,
\begin{equation}\label{e:bulk estimate}
 \left|m_s^n(E_s+n^{-2/3}w)-m_s(E_s+n^{-2/3}w)\right|
 \le\frac{C(\Im\sqrt w)^{\alpha'}}{n^{1/3}\Im w},
 \qquad w\in R(K).
\end{equation}
The constants may depend on $J$ and the fixed model and regularity
bounds, but not on $s,K,n$. The event may depend on $s$.
The threshold $n_0(K)$ may also depend on the initial preparation.
\end{theorem}

For fixed $s_*<a<T$, let $\cC(s)$, $a\le s\le T_0$,
be the exterior contours constructed in
\Cref{subsec:s4-characteristics}.

\begin{proposition}[Uniform bounds on a time window]
\label{l: Delta is sml}
Under \Cref{ass:weightfunc,ass:prelim-log-branches,ass:edge-data},
fix $s_*<T<T_0$.
For every $M'>0$, compact $\mathcal K\Subset\mathbb H$ and
$\varepsilon>0$, there are $C<\infty$ and $n_0$ such that,
for all $n\ge n_0$, with probability at least $1-\varepsilon$,
\begin{equation}\label{eq:edge-compact-window}
 \sup_{|s-T|\le M'n^{-1/3}}\sup_{w\in\mathcal K}
 \left|m_s^n(E_s+n^{-2/3}w)-m_s(E_s+n^{-2/3}w)\right|
 \le\frac C{n^{1/3}},
\end{equation}
and
\begin{equation}\label{eq:edge-exterior-window}
 \sup_{|s-T|\le M'n^{-1/3}}\sup_{z\in\cC(s)}
 |m_s^n(z)-m_s(z)|\le\frac Cn.
\end{equation}
The constants may depend on $\varepsilon,M',\mathcal K,T,a$,
the fixed working domains, the regularity bounds and the tightness
bounds in \eqref{eq:edge-initial-preparation}, but not on $n$.
\end{proposition}

\subsection{Deterministic limit and regular edge}\label{subsec:s4-macroscopic}
In this part, we study how deterministic Stieltjes transform behaves near the edge when we have the square root behavior assumption.

By \Cref{prop:prelim-deterministic-limit}, $\mu_s^n$ converge to $\mu_s$ as $n \to \infty$.
Recall from \Cref{subsec:prelim-characteristics} that
\begin{equation}\label{s4:det-eq}
 \partial_s m_s(z)+\chi_s(z)\partial_zm_s(z)=F_s(z),\qquad
 \chi_s(z)=-\frac1\theta\partial_1\cF_s(m_s(z),z).
\end{equation}
The function $F_s$ is defined in \eqref{def of F_s}.  Recall the kernel decomposition from the proof of \Cref{lem:stieltjes-bounds}:
\begin{equation}\label{s4:kernel-split}
 \frac{\mathsf b'(z)}{\mathsf b(r)-\mathsf b(z)}
 =\frac1{r-z}+g_{\mathsf b}(r,z).
\end{equation}
The function
$g_{\mathsf b}$ is Schwarz symmetric and has uniformly bounded
fixed-order derivatives on smaller working domains. In particular,
$|\Im g_{\mathsf b}(r,z)|\le C\Im z$ for real $r$.

For convenience, we write
\begin{equation}\label{s4:ordinary-transforms}
 \widehat m_s^n(z)=\int\frac{\mu_s^n(\rd r)}{r-z},\qquad
 \widehat m_s(z)=\int\frac{\mu_s(\rd r)}{r-z}.
\end{equation}
These are ordinary-kernel transforms of the same measures. They agree
with $m_s^n,m_s$ for corners; for DMC their differences are the analytic
terms in \eqref{s4:kernel-split}.

For the following estimates, fix $a\in(s_*,T_0)$ and put
\begin{equation}\label{s4:working-time}
 a_n=\lceil na\rceil/n\quad\text{for DMC and corners},
 \qquad a_n=a\quad\text{for DBM}.
\end{equation}
The edge estimates below are used on $[a,T_0]$, with constants
allowed to depend on $a$. 

\begin{lemma}\label[lemma]{lem:imastij}
Under \Cref{ass:weightfunc,ass:prelim-log-branches,ass:edge-data},
there exists $r>0$ such that, for every $s\in[a,T_0]$ and
$\zeta=\kappa+\ri \eta$ with $0<\eta\le\kappa$ and $|\zeta|\le r$,
\begin{equation}\label{s4:edge-expansion}
 m_s(E_s+\zeta)
 =m_s(E_s)+\sqrt{\cA(s)}\sqrt\zeta+\OO(|\zeta|).
\end{equation}
The error has derivative $\OO(1)$ with respect to $\zeta$
and imaginary part $\OO(\eta)$.
For every $\zeta\in\mathbb H$ with $|\zeta|\le r$,
\begin{equation}\label{s4:im-edge}
 \Im m_s(E_s+\zeta)\asymp\Im\sqrt\zeta.
\end{equation}
The constants in these estimates are independent of $s$.

For $s\in[a,T_0]$, we have $E_s'=\chi_s(E_s)$ and
\begin{equation}\label{s4:edge-time}
 \frac{\rd}{\rd s}m_s(E_s)
 =F_s(E_s)
 +\frac{\partial_{11}\cF_s(m_s(E_s),E_s)}{2\theta}\cA(s).
\end{equation}
The functions $\sqrt{\cA(s)}$, $m_s(E_s)$, and $E_s'$ are
Lipschitz in $s$ on this interval.
\end{lemma}

\begin{proof}

By \eqref{s4:kernel-split} and the identity
\[
 \zeta\int_0^\infty\frac{\rd x}{\sqrt{x}(x+\zeta)}
 =\pi\sqrt\zeta,
\]
we have
\[
 m_s(E_s+\zeta)-m_s(E_s)
 =\sqrt{\cA(s)}\sqrt\zeta
  +\mathrm I_1+\mathrm I_2+\mathrm I_3+\mathrm I_4,
\]
where
\begin{align*}
 \mathrm I_1
 &=\zeta\int_0^{r_*}
   \frac{\rho_s(E_s-x)-\pi^{-1}\sqrt{\cA(s)}\sqrt{x}}
        {x(x+\zeta)}\,\rd x,
\quad
 \mathrm I_2
 =-\frac{\sqrt{\cA(s)}}{\pi}\,
   \zeta\int_{r_*}^{\infty}\frac{\rd x}{\sqrt{x}(x+\zeta)},\\
 \mathrm I_3
 &=\int_{(-\infty,E_s-r_*]}
   \left(\frac1{r-E_s-\zeta}-\frac1{r-E_s}\right)\mu_s(\rd r),
 \quad
 \mathrm I_4
 =\int_{\mathbb R}
   \bigl[g_{\mathsf b}(r,E_s+\zeta)
         -g_{\mathsf b}(r,E_s)\bigr]\mu_s(\rd r).
\end{align*}
Here $\mathrm I_1$ is the contribution of the density
remainder in \eqref{def: classical measure}, and
$\mathrm I_2$ subtracts the tail introduced by extending
the leading integral to infinity.
The terms $\mathrm I_3$ and $\mathrm I_4$ come from the
integral over $(-\infty,E_s-r_*]$ and the kernel correction,
respectively.

By \eqref{def: classical measure}, for $\kappa\ge0$,
\[
 |\mathrm I_1|
 \le C|\zeta|\int_0^{r_*}x^{-1/2}\,\rd x,
 \qquad
 |\partial_\zeta\mathrm I_1|
 \le C\int_0^{r_*}x^{-1/2}\,\rd x.
\]
The terms $\mathrm I_2,\mathrm I_3,\mathrm I_4$ are analytic
for sufficiently small $|\zeta|$, vanish at $\zeta=0$,
and have bounded first derivatives. Hence
$|\mathrm I_j|\le C|\zeta|$ and
$|\partial_\zeta\mathrm I_j|\le C$ for $j=1,\ldots,4$.
Each term is real for real $\zeta\ge0$.
Integrating its derivative along the vertical segment
from $\kappa$ to $\kappa+\ri \eta$ therefore gives
$|\Im\mathrm I_j|\le C\eta$.
This proves \eqref{s4:edge-expansion} and the remainder bounds.

For general $\zeta=\kappa+\ri \eta\in\mathbb H$ near zero,
the same decomposition gives
\[
 \Im m_s(E_s+\zeta)
 =\eta\int_0^{r_*}
   \frac{\rho_s(E_s-x)}{(x+\kappa)^2+\eta^2}\,\rd x+\OO(\eta),
\]
where the error comes from the imaginary part of $\mathrm I_3+\mathrm I_4$.
By \eqref{def: classical measure}, the density is comparable
to $\sqrt{x}$, while
\[
 \eta\int_0^{r_*}\frac{\sqrt{x}}{(x+\kappa)^2+\eta^2}\,\rd x
 =\pi\Im\sqrt\zeta+\OO(\eta).
\]
Thus \eqref{s4:im-edge} follows by shrinking the neighborhood,
since $\eta/\Im\sqrt\zeta=2\Re\sqrt\zeta$ tends to zero.

We next establish continuity at the edge.
By \eqref{def: classical measure},
\(
 \mu_s((E_s-h,E_s))\asymp h^{3/2}
\)
for sufficiently small $h>0$, with constants independent of $s$.
If $s_j\to s$, weak convergence from
\Cref{prop:prelim-deterministic-limit} rules out both
$E_{s_j}\ge E_s+2h$ and $E_{s_j}\le E_s-2h$ along a subsequence:
apply the Portmanteau theorem to $[E_s+h,\infty)$ and
$(E_s-h,\infty)$, respectively.
Hence $E_s$ is continuous.
For fixed $\kappa>0$, $s\mapsto m_s(E_s+\kappa)$ is therefore
continuous. The uniform limits
\[
 m_s(E_s+\kappa)\longrightarrow m_s(E_s),
 \qquad
 \frac{m_s(E_s+\kappa)-m_s(E_s)}{\sqrt\kappa}
 \longrightarrow\sqrt{\cA(s)}
\]
from \eqref{s4:edge-expansion} then give continuity of
$m_s(E_s)$ and $\cA(s)$.

Differentiating the Stieltjes integral and using
\eqref{def: classical measure} gives
\[
 \begin{aligned}
 m_s'(E_s+\kappa)
 &=\frac{\sqrt{\cA(s)}}{2\sqrt\kappa}+\OO(1),\\
 m_s''(E_s+\kappa)
 &=-\frac{\sqrt{\cA(s)}}{4\kappa^{3/2}}+\OO(\kappa^{-1/2}).
 \end{aligned}
\]
Thus, for small $\epsilon>0$, there is a unique nearby
$x_s^\epsilon>E_s$ with $m_s'(x_s^\epsilon)=\epsilon^{-1}$,
and $x_s^\epsilon-E_s\asymp\epsilon^2$.
The implicit function theorem and \eqref{s4:det-eq} yield
\[
 \frac{\rd}{\rd s}x_s^\epsilon
 =\chi_s(x_s^\epsilon)
 +\frac{\chi_s'(x_s^\epsilon)m_s'(x_s^\epsilon)
              -F_s'(x_s^\epsilon)}
        {m_s''(x_s^\epsilon)}
 =\chi_s(E_s)+\OO(\epsilon).
\]
Here $\chi_s',F_s'=\OO(\epsilon^{-1})$ at $x_s^\epsilon$,
by the definition of $\chi_s$ and \eqref{def of F_s},
whereas $-m_s''(x_s^\epsilon)\asymp\epsilon^{-3}$.
Integrating and letting $\epsilon\downarrow0$ proves
$E_s'=\chi_s(E_s)$.

Finally, \eqref{s4:det-eq} and the edge expansion give
\[
 \begin{aligned}
 \frac{\rd}{\rd s}m_s(E_s+\epsilon)
 &=F_s(E_s+\epsilon)
   +[\chi_s(E_s)-\chi_s(E_s+\epsilon)]m_s'(E_s+\epsilon)\\
 &=F_s(E_s)
   +\frac{\partial_{11}\cF_s(m_s(E_s),E_s)}{2\theta}\cA(s)
   +\OO(\sqrt\epsilon).
 \end{aligned}
\]
Integrating and letting $\epsilon\to 0$ proves
\eqref{s4:edge-time}. By \eqref{s4:edge-time} and the preceding estimate,
\[
 \left|\frac{\rd}{\rd s}
 \bigl(m_s(E_s+\epsilon)-m_s(E_s)\bigr)\right|
 \le C\sqrt\epsilon.
\]
Integrating between $s$ and $t$, dividing by $\sqrt\epsilon$,
and letting $\epsilon\to 0$, we obtain from
\eqref{s4:edge-expansion}
\[
 |\sqrt{\cA(t)}-\sqrt{\cA(s)}|\le C|t-s|.
\]
The bounded right-hand sides give Lipschitz continuity of
$E_s$ and $m_s(E_s)$.
The formula
$E_s'=-\theta^{-1}\partial_1\cF_s(m_s(E_s),E_s)$
and the smooth time dependence of the coefficients in
\Cref{def: function F} then show that $E_s'$ is Lipschitz.
\end{proof}

For the bulk characteristics, we also need the square-root
expansion on the upper half-plane to the left side of the right edge $E_s$.

\begin{lemma}\label{lem:edge-interior-expansion}
There exist constants $r,C>0$ such that, for every
$s\in[a,T_0]$ and $\zeta=\kappa+\ri \eta\in\mathbb H$ with
$|\zeta|\le r$, writing $\sqrt\zeta=u+\ri v$, when $u\le4v$, we have
\begin{equation}\label{s4:interior-expansion}
 \left|m_s(E_s+\zeta)-m_s(E_s)
       -\sqrt{\cA(s)}\sqrt\zeta\right|
 \le C|\zeta|
       \left[1+v\log\left(2+\frac{v^2}{\eta}\right)\right]
\end{equation}
\end{lemma}

\begin{proof}
In the decomposition from the proof of \Cref{lem:imastij},
$\mathrm I_2,\mathrm I_3,\mathrm I_4$ remain $\OO(|\zeta|)$,
so only $\mathrm I_1$ requires an estimate.
Take $r<r_*/2$ and split its integral at $x=2v^2$.
By \eqref{def: classical measure},
\[
 |\mathrm I_1|
 \le C|\zeta|\left[
 v\int_0^{2v^2}\frac{\rd x}{\sqrt{(x+\kappa)^2+\eta^2}}
 +\int_{2v^2}^{r_*}\frac{\rd x}{\sqrt{x}}
 \right],
\]
where we used $\kappa=u^2-v^2$ and
$|x+\zeta|\ge x/2$ for $x\ge2v^2$.
The two integrals are bounded by
$C\log(2+v^2/\eta)$ and $C$, respectively.
Combining the four terms proves
\eqref{s4:interior-expansion}.
\end{proof}

\subsection{Characteristic flow analysis}\label{subsec:s4-characteristics}
\begin{figure}[htbp]
 \centering
 \includegraphics[width=0.7\linewidth]{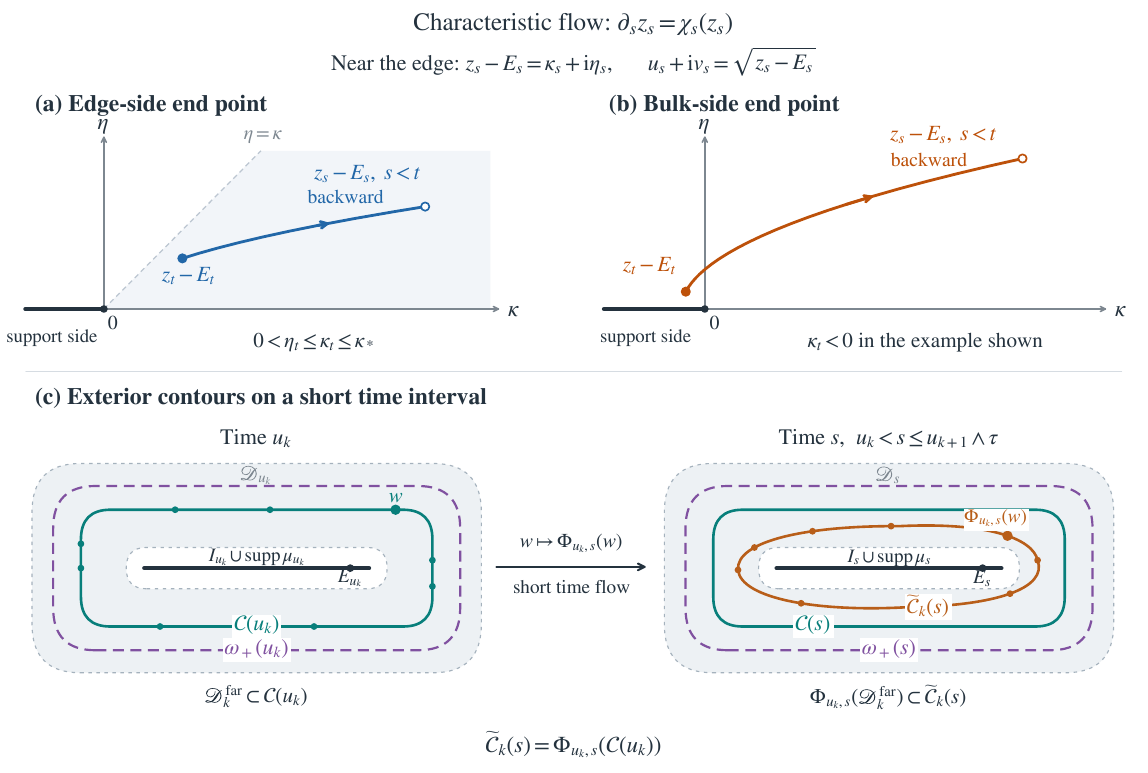}
 \caption{Characteristic paths near and away from the right edge.
Panels (a) and (b) show paths in the coordinates
$z_s-E_s=\kappa_s+\ri\eta_s$, with arrows indicating backward time.
In panel (c), $\mathcal C(s)$ and $\omega_+(s)$ are the exterior
contours. On each short time interval $[u_k,u_{k+1}]$,
the characteristic flow maps $\mathcal C(u_k)$ to
$\widetilde{\mathcal C}_k(s)$, which lies inside $\mathcal C(s)$.}
 \label{fig:s4-characteristics}
\end{figure}

Recall the characteristic flow $z_s(u)$ from
\Cref{def:prelim-characteristic-flow}. In this subsection,
we first estimate the real and imaginary parts of
$z_s(u)-E_s$ near the right edge.
We then choose exterior contours and estimate $m_s^n-m_s$
along characteristic flows on short time intervals.
The characteristic paths and contours are illustrated in
\Cref{fig:s4-characteristics}.

\paragraph{Velocity bounds.}
Recall the characteristic velocity
$\chi_s^n(z)=-\theta^{-1}\partial_1\cF_s^n(m_s^n(z),z)$
from \eqref{eq:velocity-from-flux}, and its deterministic
counterpart $\chi_s$ from \Cref{subsec:prelim-characteristics}.
The following estimates will be used to control the imaginary
part of the characteristic flow. 

\begin{lemma}\label{lem:velocity-bounds}
Under \Cref{ass:weightfunc,ass:prelim-log-branches}, let
$s\in[0,T_0]$ and $z\in\mathscr D_s\cap\mathbb H$, and
denote $\eta=\Im z$. For each of the three models,
the characteristic velocities $\chi$ satisfy
\begin{equation}\label{s4:one-sided-velocity}
 \begin{aligned}
  \Im\chi_s^n(z)&\le C\eta,
  \qquad \Im\chi_s(z)\le C\eta,\\
  |\partial_z\chi_s^n(z)|
  &\le-\frac{\Im\chi_s^n(z)}{\eta}+C.
 \end{aligned}
\end{equation}
The constant $C$ depends only on the fixed model, spectral
domains and assumption constants, and is independent of
$s,n,z$.
\end{lemma}

\begin{proof}
Fix $s$ and suppress the indices $s,n$.
For DMC, for simplicity, use the normalization $\varphi^-=1$ from
\Cref{s:DMC bulkedgeest}.
By \eqref{eq:dmc-velocity-positive-representation} and
the bounds on its analytic term in the proof of
\Cref{lem:dmc-velocity-derivatives},
\[
 \Im\cB(z)\le C\eta,\qquad
 |\cB'(z)|\le-\frac{\Im\cB(z)}{\eta}+C.
\]
For corners, the same inequalities follow by taking
imaginary parts and differentiating
\[
 \cB(z)=(1-s-1/n)W'(z)
       +\frac{\theta}{n}\sum_i\frac1{z-x_i},
\]
using Schwarz symmetry and the boundedness of $W''$.
Since $\chi=1-\cB^{-1}$ for DMC and
$\chi=-\theta\cB^{-1}$ for corners, the lower bound
$|\cB|\ge c$ in \Cref{ass:prelim-log-branches} gives
\eqref{s4:one-sided-velocity}.

For DBM, \eqref{eq:velocity-from-flux} gives
\[
 \chi(z)=-\frac12V'(z)
         +\frac1n\sum_{i=1}^n\frac1{z-\lambda_i}.
\]
Both estimates follow directly by taking imaginary parts
and differentiating, using Schwarz symmetry and the
boundedness of $V''$.
Finally, the deterministic estimate follows by passing
to the limit using \Cref{prop:prelim-deterministic-limit}.
\end{proof}

\paragraph{Characteristics flow near the edge.}
Write a characteristic as $z_s=E_s+\kappa_s+\ri \eta_s$, and
denote $u_s+\ri v_s=\sqrt{z_s-E_s}$.
We study these paths backwards from their end points.
We first consider end points to the right of the edge,
following the characteristic-flow argument in
\cite[Proposition~3.5]{adhikari2020dyson} and
\cite[Section~3]{edgerigiditydysonbrownian}.
We then extend the estimates to end points in the bulk
near the edge.

\begin{lemma}\label{c:kappa:DMC}
Under \Cref{ass:weightfunc,ass:prelim-log-branches,ass:edge-data},
there exist constants $\kappa_*,c>0$ such that the
following holds.
Fix $t\in[a,T_0]$ and a characteristic with endpoint satisfying
$0<\eta_t\le\kappa_t\le\kappa_*$.
Trace it backwards until time $a$ or the first time
$\kappa_s=2\kappa_*$, and denote this interval by $[b,t]$.
Then $0<\eta_s\le\kappa_s\le2\kappa_*$ throughout $[b,t]$, and
\begin{equation}\label{e:sharp real part DMC}
 -c\le\frac{\rd}{\rd s}\sqrt{\kappa_s}\le-\frac{1}{c}<0,
 \qquad
 -\Im\chi_s(z_s)\asymp\frac{\eta_s}{\sqrt{\kappa_s}}.
\end{equation}
Moreover, $\eta_s/\kappa_s$ is strictly increasing in $s$.
Writing $q=-\Im\chi_t(z_t)>0$, we have
\begin{equation}\label{e:Im part flow DMC}
 -\Im\chi_s(z_s)\asymp q,\qquad
 \eta_s\asymp\eta_t+q(t-s),\qquad
 \sqrt{\kappa_s}\asymp\sqrt{\kappa_t}+t-s.
\end{equation}
All constants depend only on $T_0$ and the fixed model and
assumption constants.
\end{lemma}

\begin{proof}
By \Cref{lem:imastij}, $E_s'=\chi_s(E_s)$.
Taylor expansion of the velocity using
\eqref{s4:edge-expansion} gives
\begin{equation}\label{s4:velocity-edge}
 \chi_s(E_s+\zeta)-\chi_s(E_s)
 =-\frac{\partial_{11}\cF_s(m_s(E_s),E_s)}{\theta}
       \sqrt{\cA(s)}\sqrt\zeta+\OO(|\zeta|).
\end{equation}
The remainder has bounded spatial derivative and imaginary
part $\OO(\Im\zeta)$, by the remainder bounds in
\Cref{lem:imastij} and Schwarz symmetry.
Taking real and imaginary parts, and using
\eqref{eq:edge-flux-nondegeneracy}, proves
\eqref{e:sharp real part DMC} while
$0<\eta_s\le\kappa_s\le2\kappa_*$, for sufficiently small
$\kappa_*$.

Since $\eta_s=2u_sv_s$ and $\kappa_s=u_s^2-v_s^2$,
the same expansion gives
\[
 \frac{\rd}{\rd s}\log\frac{\eta_s}{\kappa_s}
 =\frac{\partial_{11}\cF_s(m_s(E_s),E_s)\sqrt{\cA(s)}}
        {\theta}
   \frac{u_s^2+v_s^2}{2u_s\kappa_s}+\OO(1)
 \ge\frac{c}{\sqrt{\kappa_s}}>0.
\]
Thus $\eta_s\le\kappa_s$ is preserved backwards, and the
path continues until time $a$ or $\kappa_s=2\kappa_*$.

To compare the imaginary velocities $\Im\chi_s(z_s)$ along the path,
differentiate using \Cref{def:prelim-characteristic-flow}:
\begin{equation}\label{s4:time-velocity}
 \frac{\rd}{\rd s}\chi_s(z_s)
 =-\frac1\theta\bigl(
 \partial_s\partial_1\cF_s
 +\partial_{11}\cF_s F_s(z_s)
 +\partial_{12}\cF_s\chi_s(z_s)\bigr),
\end{equation}
where the derivatives of $\cF_s$ are evaluated at
$(m_s(z_s),z_s)$.
By \eqref{def of F_s} and the expressions for $\cF_s$
in \Cref{def: function F}, and Schwarz symmetry,
\[
 \left|\frac{\rd}{\rd s}\Im\chi_s(z_s)\right|
 \le C\bigl(\Im m_s(z_s)+\eta_s\bigr)
 \le C\bigl(-\Im\chi_s(z_s)\bigr).
\]
Here the last inequality follows from
\eqref{s4:im-edge} and \eqref{e:sharp real part DMC}.
Gr\"onwall's inequality gives
$-\Im\chi_s(z_s)\asymp q$.
Integrating $\eta_s'=\Im\chi_s(z_s)$ and
\eqref{e:sharp real part DMC} proves
\eqref{e:Im part flow DMC}.
\end{proof}

For bulk endpoints, we use \eqref{s4:interior-expansion}
to estimate the flow in the coordinates $u_s,v_s$.

\begin{lemma}\label{lem:bulk-characteristics}
Under \Cref{ass:weightfunc,ass:prelim-log-branches,ass:edge-data},
let $\kappa_*>0$ be as in \Cref{c:kappa:DMC},
decreased if necessary.
Fix $K\ge1$ and $t\in[a,T_0]$.
Consider a characteristic $z_s$.
Suppose that
\[
 |z_t-E_t|\le n^{-2/3+1/10},
 \qquad \eta_t\ge K^2n^{-2/3}.
\]
Trace the characteristic backwards until time $a$ or the
first time $\kappa_s=2\kappa_*$, and denote the resulting
interval by $[b,t]$.
\begin{equation}\label{c:bulkkappa DMC}
 -C\le u_s'\le-c<0,\qquad
 |v_s'|\le Cv_s,\qquad
 -\Im\chi_s(z_s)\asymp v_s.
\end{equation}
Consequently,
\begin{equation}\label{s4:uv-flow}
 u_s\asymp u_t+t-s,\qquad
 v_s\asymp v_t,\qquad
 \eta_s=2u_sv_s.
\end{equation}
Moreover, $u_s/v_s$ is strictly decreasing in $s$, so the
path crosses $u_s=4v_s$ at most once.
All constants depend only on $T_0$ and the fixed model and
assumption constants.
\end{lemma}

\begin{proof}
Since $E_s'=\chi_s(E_s)$, the characteristic equation gives
\[
 (u_s+\ri v_s)'
 =\frac{\chi_s(z_s)-\chi_s(E_s)}{2(u_s+\ri v_s)}.
\]
Continue the path backwards under the bounds
$v_s\le2e^{CT_0}v_t$ and $\eta_s\ge\eta_t/2$,
with $C$ a fixed sufficiently large constant.
The endpoint assumptions imply
$v_t\le n^{-1/3+1/20}$ and $v_t^2/\eta_t\le n^{1/10}$.
Hence
\[
 v_s\log\left(2+\frac{v_s^2}{\eta_s}\right)
 \le Cv_t\log n=\oo(1).
\]

When $u_s\le4v_s$, Taylor expansion of the velocity using
\eqref{s4:interior-expansion} therefore gives
\[
 (u_s+\ri v_s)'
 =-\frac{\partial_{11}\cF_s(m_s(E_s),E_s)\sqrt{\cA(s)}}
         {2\theta}+\OO(v_s).
\]
When $u_s\ge4v_s$, use \eqref{s4:velocity-edge}.
Its remainder is $\OO(u_s^2)$ and has imaginary part
$\OO(u_sv_s)$; after division by $u_s+\ri v_s$, its absolute
value is $\OO(u_s)$ and its imaginary part is $\OO(v_s)$.
Thus, on both parts of the path,
\[
 -C\le u_s'\le-c<0,\qquad |v_s'|\le Cv_s,
\]
after choosing $\kappa_*$ sufficiently small.
Since $u_s^2=\kappa_s+v_s^2\le2\kappa_*+\oo(1)$,
the identity $\eta_s'=2(v_su_s'+u_sv_s')$ also gives
$cv_s\le-\eta_s'\le Cv_s$.

Integration improves the assumed bounds to
$v_s\le e^{C(t-s)}v_t$ and $\eta_s\ge\eta_t$.
These estimates therefore continue until time $a$ or
$\kappa_s=2\kappa_*$, proving
\eqref{c:bulkkappa DMC}.
Integrating the bounds for $u_s'$ and $v_s'$ gives
\eqref{s4:uv-flow}.
Finally,
\[
 \frac{\rd}{\rd s}\frac{u_s}{v_s}
 =\frac{u_s'}{v_s}-\frac{u_sv_s'}{v_s^2}
 \le\frac{-c+Cu_s}{v_s}<0,
\]
after decreasing $\kappa_*$ if necessary.
\end{proof}

For clarity, we refer to the characteristic flows satisfying the
endpoint conditions in \Cref{c:kappa:DMC,lem:bulk-characteristics}
as \textit{edge and bulk characteristics}, respectively. We call their
restrictions to the time interval $[b,t]$ \textit{local edge and bulk characteristics}.

\paragraph{Characteristic flows away from the edge.}
We choose exterior contours for all times and estimate the
transform along characteristic flows on short time intervals.

For a sufficiently small fixed $\varepsilon>0$, put
\[
 \ell(s)=
 \begin{cases}
  \ell_0-\varepsilon,&\text{DMC and corners},\\
  L_s-\kappa_*/4,&\text{DBM}.
 \end{cases}
\]
Fix outer contours $\omega_+(s)$ satisfying the contour
conditions of Sections~2 and~3, enclosing the $ \cC(s)$. Let $C$ be a uniform constant in
\eqref{s4:one-sided-velocity}.
Choose $K_0>C_1$ and a sufficiently small
$0<h_*<\kappa_*/4$, and set
\[
 h_s=h_*e^{K_0(s-T_0)},\qquad
 \cC(s)=
 \partial\left\{
 z:\ell(s)<\Re z<E_s+\kappa_*/4,\quad |\Im z|<h_s
 \right\},
\]
with positive orientation.

We show that characteristics starting on $\cC(s)$ move
into the rectangle as time increases.
For DMC, positivity of the weights gives
$\chi_s(\ell(s))>0$.
For corners, \eqref{eq:corners-reference-factor} gives
$\cB_s(\ell(s))<0$, hence $\chi_s(\ell(s))>0$.
For DBM, the same argument as at the right edge gives
$L_s'=\chi_s(L_s)$ and
\(
 \chi_s(L_s-\kappa_*/4)-L_s'\ge c\sqrt{\kappa_*}.
\)
Together with \eqref{s4:velocity-edge}, analyticity and
\eqref{s4:one-sided-velocity},
for sufficiently small $h_*$, one has
\begin{equation}\label{s4:contour-inward}
 \begin{aligned}
 \Re\chi_s(\ell(s)+\ri y)-\ell'(s)&\ge c,\\
 \Re\chi_s(E_s+\kappa_*/4+\ri y)-E_s'&\le-c,
       && |y|\le h_s,\\
 \Im\chi_s(x+\ri h_s)-h_s'
       &\le-(K_0-C_1)h_s<0,
       && \ell(s)\le x\le E_s+\kappa_*/4.
 \end{aligned}
\end{equation}
The similar bound works for the bottom side of $\cC(s)$ via Schwarz symmetry.
By definition, $\cC(s)$ encloses $\operatorname{supp}\mu_s$, with
\[
 \inf_{0\le s\le T_0}
 \operatorname{dist}
 \bigl(\cC(s),\operatorname{supp}\mu_s\cup\omega_+(s)\bigr)>0.
\]

For $u\in[a,T_0]$ and $w$ in a fixed neighborhood of $\cC(a)$,
define $\Phi_{u,s}$ as the characteristic flow map starting at time $a$ by
\[
 \partial_s\Phi_{u,s}(w)=\chi_s(\Phi_{u,s}(w)),
 \qquad \Phi_{u,u}(w)=w.
\]
It is clear that these maps
are defined for $u\le s\le\min\{u+h,T_0\}$, with uniformly
bounded derivatives. 

Choose a deterministic partition of time intervel $[0,T_0]: a_n=u_0<\cdots<u_J$ with a fixed number
$J$ and $u_{k+1}-u_k\le h$.
For the discrete models, take $u_k\in\mathbb Z_n$ and
$u_J=\lfloor nT_0\rfloor/n$; for DBM, take $u_J=T_0$.
Set
\[
 \widetilde{\cC}_k(s)
 =\Phi_{u_k,s}\bigl(\cC(u_k)\bigr),
 \qquad u_k\le s\le u_{k+1}.
\]
By \eqref{s4:contour-inward}, one has
\[
 \widetilde{\cC}_k(s)\subset\operatorname{Int}\cC(s),
 \qquad u_k<s\le u_{k+1}.
\]
After suitable decreasing $h$, $\widetilde{\cC}_k(s)$ still encloses
$\operatorname{supp}\mu_s$ and satisfies $\operatorname{dist}
 \bigl(\widetilde{\cC}_k(s),\operatorname{supp}\mu_s\cup\omega_+(s)\bigr)>0.$

The following lemma compares the transform difference
on the exterior contours.

\begin{lemma}\label{lem:exterior-cauchy}
Under \Cref{ass:weightfunc}, for the contours chosen
as above,
fix a time $s$ and suppose that $\cC(s)$ encloses
$\operatorname{supp}\mu_s^n\cup\operatorname{supp}\mu_s$.
Then, for every $z$ outside $\cC(s)$ in the spectral domain,
\begin{equation}\label{e:cauchy int for stij}
 m_s^n(z)-m_s(z)
 =
 -\frac{\mathsf b'(z)}{2\pi\ri}
 \oint_{\cC(s)}
 \frac{m_s^n(w)-m_s(w)}
      {\mathsf b(w)-\mathsf b(z)}\,\rd w.
\end{equation}
Consequently,
\begin{equation}\label{s4:outer-cauchy}
 \sup_{z\in\omega_+(s)}|m_s^n(z)-m_s(z)|
 \le C\sup_{z\in\cC(s)}|m_s^n(z)-m_s(z)|.
\end{equation}
The same bound also holds for any fixed spatial derivative.

If $u_k\le s\le u_{k+1}$ and $\widetilde{\cC}_k(s)$
also encloses both supports, then
\begin{equation}\label{s4:contour-transfer}
 \sup_{z\in\cC(s)}
 \left|\frac{m_s^n(z)-m_s(z)}{\mathsf b'(z)}\right|
 \le
 \sup_{z\in\widetilde{\cC}_k(s)}
 \left|\frac{m_s^n(z)-m_s(z)}{\mathsf b'(z)}\right|.
\end{equation}
The constants depend only on the fixed working domains,
contour separations and model bounds, and are independent
of $s$ and $n$.
\end{lemma}

\begin{proof}
By \eqref{eq:stieltjes-general}, the function
\[
 q\longmapsto
 \int\frac{(\mu_s^n-\mu_s)(\rd r)}{\mathsf b(r)-q}
\]
is holomorphic outside
$\mathsf b(\operatorname{supp}\mu_s^n
\cup\operatorname{supp}\mu_s)$ and vanishes at infinity.
Its value at $q=\mathsf b(z)$ is
$(m_s^n(z)-m_s(z))/\mathsf b'(z)$.
The exterior Cauchy formula on $\mathsf b(\cC(s))$
therefore gives \eqref{e:cauchy int for stij}.
The fixed contour separation gives
\eqref{s4:outer-cauchy}, and differentiation gives
the stated derivative bounds.

By the contour construction and biholomorphicity of
$\mathsf b$, the curve $\mathsf b(\cC(s))$ lies outside
or on $\mathsf b(\widetilde{\cC}_k(s))$.
Applying the maximum principle to the same function
in the exterior of $\mathsf b(\widetilde{\cC}_k(s))$
gives \eqref{s4:contour-transfer}.
\end{proof}

\subsection{Stochastic evolution of the Stieltjes transform}
\label{subsec:s4-stochastic}
We now combine the estimates of \Cref{thm:main} with the
characteristic-flow analysis in \Cref{subsec:s4-characteristics}
to prove the main results of this section.
We first derive the stochastic evolution of $m_s^n-m_s$ along
the deterministic characteristics and estimate its drift,
martingale fluctuations, and accumulated errors.
A stopping-time argument then gives the upper bound for the
rightmost particle in \Cref{large-eigen: DMC} and the bulk
estimate in \Cref{bulktight: DMC}.
Finally, we obtain the uniform time-window bounds near the edge
and on the exterior contours stated in \Cref{l: Delta is sml}.

\subsubsection{Evolution along characteristics flow}
\label{subsec:s4-increments}

We first derive the evolution of $m_s^n-m_s$ along the
deterministic characteristics and record the estimates used
in the edge and bulk arguments below.
Recall $D_s^n$, $\Delta\mathcal M_s^n$ and $\mathbb E_s$
from \eqref{eq:increment-common}, the drift $a_s^n$ from
\eqref{eq:leading-common}, and $P_s^n$ from \eqref{eq:main-P}.
We write
\[
 \delta_s(z)=m_s^n(z)-m_s(z),
 \qquad
 \xi_s(z)=\frac1n\Delta\mathcal M_s^n(z).
\]
For DMC and weighted $\beta$-corners, time takes values in
$[0,T_0]\cap\mathbb Z_n$, where
$\mathbb Z_n=n^{-1}\mathbb Z_{\ge0}$;
for DBM, time ranges continuously over $[0,T_0]$.

For DMC and weighted $\beta$-corners,
\Cref{thm:main} gives, for every $z$ in the spectral domain with
$d=\operatorname{dist}(z,I_s)\ge n^{-1+\delta}$,
\begin{align}
 |\mathbb E_sD_s^n(z)-a_s^n(z)|
 &\le \frac{CP_s^n(z)}{nd},
 \label{s4:mean-input}\\
 \mathbb E_s|\Delta\mathcal M_s^n(z)|^{2p}
 &\le C_p\left[
 \left(\frac{P_s^n(z)}{nd^2}\right)^p
 +\frac{P_s^n(z)}{n^{2p-1}d^{4p-2}}
 \right],
 \label{s4:moment-input}
\end{align}
for each fixed integer $p\ge1$. For the first, the same contour deformation and pole subtraction
as in the proof of \Cref{p:edge impr} give
$|c_s^n(z)|\le CP_s^n(z)/d$ for DMC; for corners, this bound
is included in \Cref{prop:corners-one-step}.
Combining this with \eqref{eq:main-drift} and
\eqref{eq:main-second} proves \eqref{s4:mean-input}.

We first derive the one step time evolution of $\delta_{s}(z_s)$ along the char flow
for DMC and weighted $\beta$-corners.
Fix $t\in[0,T_0]\cap\mathbb Z_n$ and a deterministic
characteristic flow $(z_s)_{0\le s\le t}$.
For $s\in[0,t)\cap\mathbb Z_n$, write $s'=s+1/n$.
By \Cref{def:prelim-characteristic-flow},
\begin{equation}\label{sde for DMC1}
 \begin{aligned}
 \delta_{s'}(z_{s'})-\delta_s(z_s)
 ={}&
 \bigl(m_{s'}^n(z_{s'})-m_s^n(z_{s'})\bigr)
 +\bigl(m_s^n(z_{s'})-m_s^n(z_s)\bigr)-\int_s^{s'}F_r(z_r)\,\rd r.
 \end{aligned}
\end{equation}
For the first term, \eqref{eq:main-drift} in
\Cref{thm:main} gives
\[
 m_{s'}^n(z_{s'})-m_s^n(z_{s'})
 =
 \frac{a_s^n(z_{s'})}{n}
 +\frac{c_s^n(z_{s'})}{n^2}
 +\xi_s(z_{s'})
 +\frac{\mathcal R_s^n(z_{s'})}{n}.
\]
For the second term, the characteristic equation gives
\[
 m_s^n(z_{s'})-m_s^n(z_s)
 =
 \int_s^{s'}\chi_r(z_r)\partial_zm_s^n(z_r)\,\rd r.
\]

Define $F_s^n$ by replacing $m_s,\cF_s$ with $m_s^n,\cF_s^n$
in \eqref{def of F_s}, using the same outer contour.
By \eqref{eq:leading-common} and \eqref{e:burger2},
we obtain
\(
 a_s^n(z)
 =-\chi_s^n(z)\partial_zm_s^n(z)+F_s^n(z).
\)
We first collect the leading term in \eqref{sde for DMC1}. Evaluating these terms
at the initial point gives
\begin{equation}\label{e:first derivative DMC}
 \begin{aligned}
 \mathcal E_1(s)
 &:=
 \frac{a_s^n(z_s)}n
 +\frac{\chi_s(z_s)}n\partial_zm_s^n(z_s)
 -\frac{F_s(z_s)}n =
 \frac{\chi_s(z_s)-\chi_s^n(z_s)}n
       \partial_zm_s^n(z_s)
 +\frac{F_s^n(z_s)-F_s(z_s)}n.
 \end{aligned}
\end{equation}
Collecting the remaining terms, we set
\begin{equation}\label{e:second derivative DMC}
 \begin{aligned}
 \mathcal E_2(s)
 &:=
 \xi_s(z_{s'})
 +\frac{c_s^n(z_{s'})}{n^2}
 +\frac{\mathcal R_s^n(z_{s'})}{n}+
 \frac{a_s^n(z_{s'})-a_s^n(z_s)}n\\
 &\quad+
 \int_s^{s'}
 \left[
 \chi_r(z_r)\partial_zm_s^n(z_r)
 -\chi_s(z_s)\partial_zm_s^n(z_s)
 \right]\rd r-
 \int_s^{s'}\left[F_r(z_r)-F_s(z_s)\right]\rd r.
 \end{aligned}
\end{equation}
Therefore we has
\begin{equation}\label{sde for DMC 2}
 \delta_{s'}(z_{s'})-\delta_s(z_s)
 =\mathcal E_1(s)+\mathcal E_2(s).
\end{equation}
For DBM, \eqref{eq:pote_mt} and
\Cref{def:prelim-characteristic-flow} give the corresponding
continuous-time decomposition
\begin{equation}\label{s4:dbm-characteristic-evolution}
 \rd \delta_s(z_s)
 =\rd \mathcal E_1(s)+\rd \mathcal E_2(s),
\end{equation}
where
\begin{equation}\label{s4:dbm-error-decomposition}
 \begin{aligned}
 \rd \mathcal E_1(s)
 &=
 \left[
 (\chi_s-\chi_s^n)\partial_zm_s^n+F_s^n-F_s
 \right](z_s)\,\rd s,\\
 \rd \mathcal E_2(s)
 &=
 \frac{1-\theta}{2\theta n}
 \partial_z^2m_s^n(z_s)\,\rd s
 +\rd \mathcal M_s^n(z_s).
 \end{aligned}
\end{equation}
Here $\rd \mathcal M_s^n(z_s)$ is the martingale
differential in \eqref{eq:pote_mt}, evaluated at $z_s$. We use the same definition of $P_s^n$ as in
\eqref{eq:main-P} for DBM.

The following lemma estimates $\mathcal E_1(s)$ and
$\mathcal E_2(s)$ along the local edge and bulk
characteristics of
\Cref{c:kappa:DMC,lem:bulk-characteristics},
and the exterior characteristics.

\begin{lemma}\label{lem:s4-characteristic-evolution}
For DMC and weighted $\beta$-corners, assume
\Cref{ass:weightfunc,ass:prelim-log-branches,ass:edge-data}. Assume that the char flow
$(z_r)_{s\le r\le s'}$ lies in the spectral
domain with
\[
 \inf_{s\le r\le s'}
 \operatorname{dist}(z_r,I_s)\ge n^{-1+\delta}.
\]
For the edge and bulk characteristics in
\Cref{c:kappa:DMC,lem:bulk-characteristics}, 
recall that $\eta_s=\Im z_s>0$. For DMC, assume additionally that
\(
 \sup_{s\le r\le s'}
 |m_s^n(z_r)-m_r(z_r)|=\oo(1).
\)
Then we have
\begin{equation}\label{s4:sharp-linear}
 \begin{aligned}
 |\mathcal E_1(s)|
 \le{}&
 \frac{|\delta_s(z_s)|}{n}
 \left(
 -\frac{\Im\chi_s(z_s)}{\eta_s}
 +C+C\frac{|\delta_s(z_s)|}{\eta_s}
 \right)+\frac Cn\sup_{w\in\omega_+(s)}|\delta_s(w)|
 +\frac{CP_s^n(z_s)}{n^2},
 \end{aligned}
\end{equation}
and
\begin{equation}\label{s4:step-remainder}
 |\mathcal E_2(s)-\xi_s(z_{s'})|
 \le
 \frac{CP_s^n(z_s)}
      {n^2\operatorname{dist}(z_s,I_s)}.
\end{equation}
For a characteristic segment $(z_r)_{s\le r\le s'}$ satisfying
\[
 \inf_{s\le r\le s'}
 \operatorname{dist}
 \bigl(z_r,I_s\cup\operatorname{supp}\mu_r\bigr)\ge c>0,
\]
where $c$ is independent of $n$, the estimates become
\[
 \begin{aligned}
 |\mathcal E_1(s)|
 &\le \frac Cn
 \left(
 |\delta_s(z_s)|
 +\sup_{w\in\omega_+(s)}|\delta_s(w)|
 \right)+\frac C{n^2},\\
 |\mathcal E_2(s)-\xi_s(z_{s'})|
 &\le\frac C{n^2}.
 \end{aligned}
\]
The final $n^{-2}$ term in the estimates for
$\mathcal E_1(s)$ is absent for DMC.
The constants depend only on the fixed model,
spectral domains and assumption constants, and additionally
on $c$.
\end{lemma}

\begin{proof}
Recall
\Cref{c:kappa:DMC,lem:bulk-characteristics}.
Differentiating \eqref{eq:velocity-from-flux}, and using
\eqref{eq:dmc-velocity-derivatives} for DMC and
\eqref{eq:corners-velocity-jets} for corners, gives,
uniformly for $s\le r\le s'$,
\[
 \begin{aligned}
 |\partial_zm_s^n(z_r)|
 &\le CP_s^n(z_s),\\
 |\partial_z^2m_s^n(z_r)|
 +|\partial_za_s^n(z_r)|
 &\le
 \frac{CP_s^n(z_s)}{\operatorname{dist}(z_s,I_s)}.
 \end{aligned}
\]
Here the required lower bound on
$|\partial_{11}\cF_s^n(m_s^n(z_r),z_r)|$
follows for corners from
\[
 |\partial_{11}\cF_s^n(m_s^n(z_r),z_r)|
 =\frac{\theta^3}{|\cB_s^n(z_r)|^2}\ge c>0
\]
by \Cref{ass:prelim-log-branches}.
For DMC, the additional smallness condition and the local
edge expansion place $(m_s^n(z_r),z_r)$ in the fixed
neighborhood of the macroscopic edge where
\Cref{ass:edge-data}(ii) gives the same lower bound,
provided that $\kappa_*$ is sufficiently small and $n$ is
sufficiently large.
These lower bounds allow us to recover the derivatives
of $m_s^n$ from those of $\chi_s^n$.
The bound
$P_s^n(z_s)\le C/\operatorname{dist}(z_s,I_s)$
absorbs the quadratic derivative terms.
For the drift derivative, we also use
$a_s^n=-\chi_s^n\partial_zm_s^n+F_s^n$
and the contour separation in \eqref{def of F_s}.
Local comparability $P_s^n(z_r)\asymp P_s^n(z_s)$
holds uniformly since
$|z_r-z_s|=\OO(n^{-1})
=\oo(\operatorname{dist}(z_s,I_s))$.

We first estimate $\mathcal E_1(s)$.
In \eqref{e:first derivative DMC}, Taylor-expand
$\partial_1\cF_s^n$ and $\partial_2\cF_s^n$
in their first variable around $m_s^n(z_s)$ and integrate
the contour difference from \eqref{def of F_s}
by parts.
For corners, this expansion follows from the exact
difference formula for reciprocals and the denominator
bounds at the two end points.
The linear local terms combine into
$-\partial_z\chi_s^n(z_s)\delta_s(z_s)$, giving
\[
 \begin{aligned}
 &\left|
 n\mathcal E_1(s)
 +\partial_z\chi_s^n(z_s)\delta_s(z_s)
 \right| \le
 C\sup_{w\in\omega_+(s)}|\delta_s(w)|
 +C|\delta_s(z_s)|^2
       \bigl(1+|\partial_zm_s^n(z_s)|\bigr)
 +\frac{CP_s^n(z_s)}n.
 \end{aligned}
\]
The contour estimate follows by transferring the
derivative to the kernel, which is uniformly bounded
together with its derivatives by contour separation.
The last term comes from the finite-$n$ flux correction
in \eqref{e:F_CP} and is absent for DMC.

The Stieltjes integral estimate in \eqref{eq:mp-basic}
and boundedness of the transforms give
$|\partial_zm_s^n(z_s)|\le C(1+\eta_s^{-1})$.
Moreover, \eqref{eq:velocity-from-flux} gives
\[
 |\Im(\chi_s^n(z_s)-\chi_s(z_s))|
 \le C|\delta_s(z_s)|
      +\frac{C\eta_sP_s^n(z_s)}n.
\]
For corners, \eqref{e:B_function} and
$\chi_s^n=-\theta/\cB_s^n$ imply
\[
 |\Im m_s^n(z_s)|
 \le C\bigl(\eta_s+|\Im\chi_s^n(z_s)|\bigr)
 \le C\eta_sP_s^n(z_s).
\]
Thus the change from $1-s$ to $1-s-1/n$ contributes
at most $Cn^{-1}(\eta_s+|\Im m_s^n(z_s)|)$,
by Schwarz symmetry.
For DMC, this finite-$n$ contribution is absent.
Applying \eqref{s4:one-sided-velocity} and using
$P_s^n(z_s)/n\le C$ therefore yields
\[
 |\partial_z\chi_s^n(z_s)|
 \le
 -\frac{\Im\chi_s(z_s)}{\eta_s}
 +C+C\frac{|\delta_s(z_s)|}{\eta_s}.
\]
Combining these estimates with boundedness of
$\delta_s(z_s)$ proves \eqref{s4:sharp-linear}.

We next estimate $\mathcal E_2(s)-\xi_s(z_{s'})$.
By \eqref{eq:main-drift},
\eqref{s4:mean-input}, and local comparability,
\[
 \left|
 \frac{c_s^n(z_{s'})}{n^2}
 +\frac{\mathcal R_s^n(z_{s'})}{n}
 \right|
 \le
 \frac{CP_s^n(z_s)}
      {n^2\operatorname{dist}(z_s,I_s)}.
\]
It remains to control the drift difference and the two integrals in \eqref{e:second derivative DMC}.
By \eqref{s4:time-velocity},
\eqref{def of F_s}, and
\Cref{def:prelim-characteristic-flow},
\[
 \left|\frac{\rd}{\rd r}\chi_r(z_r)\right|
 +\left|\frac{\rd}{\rd r}F_r(z_r)\right|\le C.
\]
For the local terms, we use
$\frac{\rd}{\rd r}m_r(z_r)=F_r(z_r)$.
For the contour terms, we fix the contour locally
in time by holomorphic deformation and apply
\eqref{e:burger2} at its fixed positive distance
from the deterministic support.

Since the empirical transform is held at time $s$,
\[
 \frac{\rd}{\rd r}
 \bigl[\chi_r(z_r)\partial_zm_s^n(z_r)\bigr]
 =
 \left(\frac{\rd}{\rd r}\chi_r(z_r)\right)
       \partial_zm_s^n(z_r)
 +\chi_r(z_r)^2\partial_z^2m_s^n(z_r).
\]
Applying these derivative bounds to the drift difference
and the two integrals in \eqref{e:second derivative DMC},
and using $s'-s=1/n$, bounds each term by
\[
 \frac C{n^2}
 \left(
 \frac{P_s^n(z_s)}{\operatorname{dist}(z_s,I_s)}
 +1
 \right).
\]
Since $P_s^n\ge1$ and the spectral domain is bounded,
this proves \eqref{s4:step-remainder}.

Finally, on the characteristics far from the support,
the fixed distance from the empirical and deterministic
supports bounds the transforms and their spatial
derivatives uniformly.
The local flux differences are Lipschitz,
and the contour difference has the same bound obtained
above by integration by parts.
Applying these bounds to \eqref{e:first derivative DMC}
and \eqref{e:second derivative DMC} proves the two exterior
estimates.
\end{proof}

Similar estimates apply to DBM, where for simplicity, we skip the proof.
\begin{lemma}\label{lem:s4-dbm-characteristic-evolution}
For DBM, retain the same notations as above. For the edge and bulk characteristics, we have
\begin{equation}\label{s4:dbm-sharp-linear}
 \begin{aligned}
 |\rd \mathcal E_1(s)|
 \le{}&
 \left[
 |\delta_s(z_s)|
 \left(
 -\frac{\Im\chi_s(z_s)}{\eta_s}
 +C+C\frac{|\delta_s(z_s)|}{\eta_s}
 \right) 
 +C\sup_{w\in\omega_+(s)}|\delta_s(w)|
 \right]\rd s,
 \end{aligned}
\end{equation}
and
\begin{equation}\label{s4:dbm-second-order}
 \left|
 \rd \mathcal E_2(s)-\rd \mathcal M_s^n(z_s)
 \right|
 \le
 \frac{CP_s^n(z_s)}
      {n\operatorname{dist}(z_s,I_s)}\,\rd s.
\end{equation}
To control the accumulated martingale term, we use
\eqref{eq:dbm-bracket}, which gives
\begin{equation}\label{s4:dbm-martingale-bound}
 \frac1{\theta n^2}
 \int\frac{\mu_s^n(\rd x)}{|x-z_s|^4}
 \le
 \frac{CP_s^n(z_s)}
      {n^2\operatorname{dist}(z_s,I_s)^2}.
\end{equation}
Similarly, for the exterior characteristic segments
\((z_r)_{s\le r\le s'},\)
with $\inf_{s\le r\le s'}\operatorname{dist}(z_r,I_s)\ge c$,
the estimates become
\[
 \begin{aligned}
 |\rd\mathcal E_1(s)|
 &\le C
 \left(
 |\delta_s(z_s)|
 +\sup_{w\in\omega_+(s)}|\delta_s(w)|
 \right)\rd s,\\
 \left|
 \rd\mathcal E_2(s)-\rd\mathcal M_s^n(z_s)
 \right|
 &\le\frac Cn\,\rd s.
 \end{aligned}
\]
The quadratic-variation rate in
\eqref{s4:dbm-martingale-bound} is bounded by $C/n^2$.
The constants depend only on the fixed model,
spectral domains and assumption constants, and additionally
on $c$.
\end{lemma}

\begin{remark}
    All the exterior estimates hold for all $0\le s<s'=s+1/n\le T_0$ and do not require
\Cref{ass:edge-data}.
\end{remark}

\subsubsection{Preparation in time interval $[0,a]$}
\label{subsec:s4-preparation}
In this part, we prove \Cref{prop:s4-positive-time-preparation}, which
shows that \Cref{ass:edge-data}(iii) still holds at any fixed time
$a>s_*$, without a square-root assumption before $a$. \Cref{prop:s4-positive-time-preparation} provides the initial bounds for the local edge
and bulk analysis on $[a_n,T_0]$, without requiring square-root behavior throughout $[0,a]$. 

For each sufficiently small fixed $\varepsilon>0$, we construct
a deterministic function $r(s)>E_s$
on $[0,a_n]$, with $r(a_n)=E_{a_n}+\varepsilon$.
Using the characteristics estimates, we show that
$\sup I_s<r(s)$ throughout $[0,a_n]$ with probability tending to one,
and obtain $\OO_{\mathbb P}(n^{-1})$ bounds for
$m_{a_n}^n-m_{a_n}$.

Let
\(
 \mathcal G_n^0=\{I_0\subset[\ell_0,r_0]\}
 \in\mathscr F_0.
\)
By \Cref{ass:edge-data}\textup{(iii)},
$\mathbb P((\mathcal G_n^0)^c)=\oo(1)$.
On this event, \Cref{lem:b-bounds} gives
\[
 I_s\subset[\ell_0,r_0+s]\quad\text{for DMC},\qquad
 I_s\subset[\ell_0,r_0]\quad\text{for corners}.
\]
For DBM, choose $R$ sufficiently large, and set
$\sigma_R=\inf\{s\le T_0:\max_i|\lambda_i(s)|\ge R\}$.
The support comparison in \cite[Proposition~2.5]{MR4009708},
\Cref{lem:dbm-beta-comparison}, and the exponential martingale
inequality give
\begin{equation}\label{s4:dbm-confinement}
 \mathbb P\bigl(\mathcal G_n^0\cap\{\sigma_R\le T_0\}\bigr)
 \le Ce^{-cn}.
\end{equation}
Set $\sigma_R=\infty$ for the discrete models. Recall the common spectral domains and outer contour mentioned above.

We first construct a deterministic function $r(s)>E_s$
on $[0,t]$, with $r(t)=E_t+\varepsilon$.
The bounds in the following lemma will allow us to control
$\sup I_s$ throughout $[0,t]$.
\begin{lemma}\label{lem:s4-exterior-comparison}
Fix $t\in(s_*,T_0)$ and a sufficiently small $\varepsilon>0$.
There are a deterministic Lipschitz function $r:[0,t]\to\mathbb R$
and constants $d,\sigma>0$ such that
\begin{equation}\label{s4:exterior-comparison}
 r(t)=E_t+\varepsilon,\qquad r(s)-E_s\ge d,\qquad
 r'(s)-\chi_s(r(s))\ge\sigma.
\end{equation}
For DBM, same requirement applies to the the left curve, that is
$l(t)=L_t-\varepsilon$, $L_s-l(s)\ge d$ and
$l'(s)-\chi_s(l(s))\le-\sigma$.
The constants are uniform for $t$ in a compact subinterval of
$(s_*,T_0)$, with $\varepsilon$ fixed.
\end{lemma}

\begin{proof}
We first construct a path staying to the right of $E_s$, then
perturb its velocity to obtain the strict inequality in
\eqref{s4:exterior-comparison}.

\smallskip
\noindent\emph{Step 1: Continue the characteristic away from the support.}
Starting from $x_t=E_t+\varepsilon$, solve
$x_s'=\chi_s(x_s)$ backwards and put $v_s=\chi_s(x_s)$.
By \eqref{s4:time-velocity} and \eqref{eq:velocity-from-flux},
while $x_s$ remains in a fixed compact interval inside the
working domain,
\[
 |v_s|\le C\quad\text{for DBM},\qquad
 c\le v_s\le1-c\quad\text{for DMC},\qquad
 c\le|v_s|\le C\quad\text{for corners}.
\]
The equations for $(x_s,v_s)$ consequently have uniformly
bounded derivatives on a fixed complex neighborhood.

For $u$ near $x_t$, let $\Phi_{t,s}(u)$ be the position at time
$s$ of the characteristic ending at $u$ at time $t$.
By $\widehat m_s'(x_s)\ge0$, \eqref{s4:kernel-split} and the
velocity formulas, $\chi_s'(x_s)\le C$. Hence
\begin{equation}\label{s4:exterior-jacobian}
 \left.\partial_u\Phi_{t,s}(u)\right|_{u=x_t}
 =
 \exp\left\{-\int_s^t\chi_q'(x_q)\,\rd q\right\}
 \ge e^{-CT_0}.
\end{equation}
Analytic dependence on the terminal data gives inverse position
maps on uniform neighborhoods.
For the discrete models, nonvanishing of $\cB_s$ in the upper
working domain follows from
\Cref{ass:prelim-log-branches},
\eqref{eq:corners-reference-factor}, and Hurwitz's theorem.
The characteristic equation and \eqref{eq:velocity-from-flux}
therefore extend $m_s$ analytically across $x_s$, with real
values on the real axis. By \eqref{s4:kernel-split} and
Stieltjes inversion,
\(
 \mu_s((x_s-d_0,x_s+d_0))=0
\)
for some fixed $d_0>0$. The uniform inverse neighborhoods allow
this continuation across any proposed first contact with the
support.

\smallskip
\noindent\emph{Step 2: Exclude mass to the right of $x_s$.}
Choose a positively oriented contour $\gamma_s$ enclosing the
part of the support to the right of $x_s$, with its left real
intersection in the preceding zero-mass interval and its right
intersection beyond all the supports.
The residue formula and \eqref{s4:det-eq} give
\begin{equation}\label{s4:exterior-tail-mass}
 \frac{\rd}{\rd s}\mu_s((x_s,\infty))
 =
 \begin{cases}
 0,&\text{DBM},\\
 -\theta^{-1}\operatorname{wind}_{\gamma_s}\cB_s,
     &\text{DMC},\\
 \operatorname{wind}_{\gamma_s}\cB_s,
     &\text{corners}.
 \end{cases}
\end{equation}
The coefficient $\cB_s$ is nonzero on $\gamma_s$, by Step~1,
the velocity formulas at $x_s$, and the fixed outer-contour
bounds at the other real intersection.
Its winding number is therefore constant and is zero near $t$,
where $x_s>E_s$. Thus $\mu_s((x_s,\infty))=0$.
Together with Step~1, this gives
\(
 x_s-E_s\ge d_0.
\)
The same argument applies when the terminal point lies beyond
all the supports.

\smallskip
\noindent\emph{Step 3: Continue the path to time zero.}
Fix $R_+>\sup_{0\le s\le T_0}E_s$, with
$E_t+\varepsilon<R_+$ and fixed separation from the boundary
of the working domain.
Continue the backward path with upper reflection at $R_+$.
On each interval without reflection, Steps~1--2 apply; at
$R_+$ the distance from the support is fixed.
The estimates are uniform over these intervals, so the reflected
path is defined on $[0,t]$ and satisfies
\[
 x_t=E_t+\varepsilon,\qquad x_s-E_s\ge d_0.
\]

\smallskip
\noindent\emph{Step 4: Obtain the velocity inequality.}
Repeat the reflected construction with velocity $\chi_s+\sigma$
and the same terminal value, and denote the resulting path by
$r(s)$. The Lipschitz property of the reflection map and
Gronwall's inequality give
\(
 \sup_{0\le s\le t}|r(s)-x_s|\le C\sigma.
\)
The reflection term is nonnegative in forward time, so
$r'(s)\ge\chi_s(r(s))+\sigma$ almost everywhere.
Choosing $C\sigma<d_0/2$ gives
\[
 r(t)=E_t+\varepsilon,\qquad
 r(s)-E_s\ge d_0/2,\qquad
 r'(s)-\chi_s(r(s))\ge\sigma
\]
almost everywhere, proving \eqref{s4:exterior-comparison}.
For DBM, use lower reflection at a fixed point to the left of
all supports and replace $\chi_s$ by $\chi_s-\sigma$.
All constants are uniform for $t$ in a fixed compact subinterval
of $(s_*,T_0)$, with $\varepsilon$ fixed.
\end{proof}

In the following, using the characteristic estimates, we show that
$\sup I_s<r(s)$ throughout $[0,a_n]$ with probability tending
to one, and obtain $\OO_{\mathbb P}(n^{-1})$ bounds for
$m_{a_n}^n-m_{a_n}$ on the fixed compact sets $\mathcal K$.
\begin{proposition}\label{prop:s4-positive-time-preparation}
Under
\Cref{ass:edge-data}, fix $a\in(s_*,T_0)$ and use
\eqref{s4:working-time}. Then
\begin{equation}\label{s4:propagated-preparation}
 (\sup I_{a_n}-E_{a_n})_+\xrightarrow{\mathbb P}0,
 \qquad
 \sup_{z\in\mathcal K}|\delta_{a_n}(z)|
 =\OO_{\mathbb P}(n^{-1})
\end{equation}
for every fixed compact
$\mathcal K\Subset\Lambda\setminus[\ell_0,E_a]$.
For DBM, $(L_{a_n}-\inf I_{a_n})_+\to0$ in probability, and use
$\mathcal K\Subset\Lambda\setminus[L_a,E_a]$.
The transform estimate holds for every fixed spatial derivative
on a smaller compact set.
\end{proposition}

\begin{proof}
We give the proof for DMC.

\emph{Step 1. Comparison curves and stopping.}
Fix a small $\varepsilon>0$ and let $r$ be given by
\Cref{lem:s4-exterior-comparison} with $t=a_n$.
Choose a small fixed $h>0$ and put $\eta_n=cn^{-1/2}$.
Take grid times $\mathfrak t_j$ covering $[0,a_n]$, with spacing
comparable to $c_1\eta_n$, except for the last interval.
Trace the characteristic from
\[
 w_{\mathfrak t_j}^j=r(\mathfrak t_j)+\ri\eta_n
\]
back to $\alpha_j=(\mathfrak t_j-h)\vee0$, rounded up to the grid.
By \eqref{s4:exterior-comparison}, for sufficiently small $h$,
\[
 \Im w_s^j\asymp\eta_n,\qquad
 c_2(\mathfrak t_j-s)\le\Re w_s^j-r(s)
 \le C(\mathfrak t_j-s).
\]
Choose $0<\Delta<c_2h/4$ and a contour $\Gamma(s)$ with vertical
sides $\ell_0-\Delta$ and $r(s)+\Delta$.
Using \eqref{s4:exterior-comparison} and
\eqref{s4:one-sided-velocity}, construct its heights and
short-time transports as in \Cref{subsec:s4-characteristics}.
This gives a partition $0=u_0<\cdots<u_{N_{\mathrm c}}=a_n$,
with $N_{\mathrm c}=\OO(1)$, satisfying the contour comparisons
of \Cref{lem:exterior-cauchy} before a support crossing.
If $\alpha_j>0$, then $w_{\alpha_j}^j$ lies outside
$\Gamma(\alpha_j)$ at a fixed distance.

Let $\mathcal K_0$ contain fixed neighborhoods of $\Gamma(0)$
and all $w_0^j$ with $\alpha_j=0$.
By \eqref{eq:edge-initial-preparation}, on an
$\mathscr F_0$-measurable event contained in $\mathcal G_n^0$
and of probability $1-\oo(1)$,
\[
 \sup I_0<r(0),\qquad
 \sup_{z\in\mathcal K_0}|\delta_0(z)|\le\frac{\log n}{n}.
\]
On this event, stop at $a_n$, the first support crossing, or the
first failure of
\[
 \sup_{z\in\Gamma(s)}|\delta_s(z)|<\frac{(\log n)^{10}}n,
 \qquad
 |\delta_s(w_s^j)|<n^{-5/8},
 \quad \alpha_j\le s\le\mathfrak t_j.
\]
Denote this time by $\tau_{\mathrm c}$.
Set it and all stopped errors and martingales to zero on the
complementary initial event.
Before $\tau_{\mathrm c}$,
\begin{equation}\label{s4:coarse-distance}
 \operatorname{dist}(w_s^j,I_s)
 \ge c\sqrt{(\mathfrak t_j-s)^2+\eta_n^2}.
\end{equation}
The fixed contour separations persist through the stopping
update, since each particle moves by at most $1/n$.

\emph{Step 2. Exterior estimates.}
Let $X_k$ be the supremum of $|\delta_s|$ on $\Gamma(s)$ for
$0\le s\le u_k\wedge\tau_{\mathrm c}$, and let $N_k^*$ be the
stopped contour-martingale supremum on $[u_k,u_{k+1}]$.
The fixed-distance evolution, contour comparison and
Sobolev--Doob estimate give
\[
 \mathbb E(N_k^*)^2\le Cn^{-2},\qquad
 X_{k+1}\le C(X_k+N_k^*+n^{-1}).
\]
Here we use \eqref{s4:moment-input} and
\Cref{lem:exterior-cauchy,lem:s4-characteristic-evolution},
followed by Gronwall's inequality.
Since $X_0=\OO_{\mathbb P}(n^{-1})$ and
$N_{\mathrm c}=\OO(1)$, we obtain
$X_{N_{\mathrm c}}=\OO_{\mathbb P}(n^{-1})$.
In particular, with probability $1-\oo(1)$,
\begin{equation}\label{s4:coarse-outer-log}
 \sup_{s\le\tau_{\mathrm c}}
 \sup_{z\in\Gamma(s)\cup\omega_+(s)}|\delta_s(z)|
 \le\frac{C(\log n)^2}{n}.
\end{equation}

\emph{Step 3. Estimates along $w_{\mathfrak t_j}^j$.}
Put $z=w_s^j$ and $\eta=\Im z$.
Before the stop, the deterministic support distance and
\eqref{s4:kernel-split} give
\[
 \int\frac{\mu_s^n(\rd x)}{|x-z|^2}
 =\frac{\Im\widehat m_s^n(z)}{\eta}
 \le C+\frac{|\delta_s(z)|}{\eta}\le C.
\]
Hence $|\partial_zm_s^n(z)|+P_s^n(z)\le C$ and
$|\partial_z^2m_s^n(z)|\le C/\operatorname{dist}(z,I_s)$.
Kernel comparison gives the same bounds throughout one step.
Substituting into \eqref{e:first derivative DMC} and
\eqref{e:second derivative DMC}, using \eqref{s4:mean-input},
bounds the leading drift by
\(
 \frac Cn\left(
 |\delta_s(w_s^j)|+\sup_{\omega_+(s)}|\delta_s|
 \right),
\)
with accumulated remainder at most $C(\log n)/n$ by
\eqref{s4:coarse-distance}.

Let $N^j$ be the martingale along $w_s^j$, started at $\alpha_j$
and stopped at $\tau_{\mathrm c}$, including the stopping update.
By \eqref{s4:moment-input} and \eqref{s4:coarse-distance},
\begin{equation}\label{s4:coarse-bracket}
 \langle\Re N^j\rangle+\langle\Im N^j\rangle
 \le Cn^{-3/2}.
\end{equation}
The eighth-moment estimate, the Burkholder--Rosenthal inequality
and a union bound over $\OO(n^{1/2})$ paths give
\[
 \mathbb P\left(\sup_{j,s}|N_s^j|>n^{-2/3}\right)
 \le Cn^{-1/6}.
\]

If $\alpha_j=0$, use the initial transform bound.
If $0<\alpha_j\le\tau_{\mathrm c}$, use
\eqref{e:cauchy int for stij} on $\Gamma(\alpha_j)$ and
\eqref{s4:coarse-outer-log}.
Thus every entrance error is at most $C(\log n)^2/n$.
Gronwall's inequality gives, with probability $1-\oo(1)$,
\begin{equation}\label{s4:coarse-improvement}
 \sup_{j,\,\alpha_j\le s\le\mathfrak t_j\wedge\tau_{\mathrm c}}
 |\delta_s(w_s^j)|
 \le C\bigl(n^{-2/3}+(\log n)^2/n\bigr)
 =\oo(n^{-5/8}).
\end{equation}
This includes the stopping update and paths entering at the stop.

\emph{Step 4. Excluding a support crossing.}
Work on the common event of the preceding estimates.
Suppose that the support first reaches or crosses $r(s)$.
Choose the first $\mathfrak t_j\ge s$.
Since the mesh is $\oo(h)$, this path is defined at $s$, and,
with $z=w_s^j$ and $\eta=\Im z$,
\[
 |\Re z-r(s)|\le C\eta_n,\qquad \eta\asymp\eta_n.
\]
The rightmost interval has mass $1/n$ and length
$\theta/n=\oo(\eta_n)$.
Its distance from $\Re z$ is $\OO(\eta_n)$, including the
$\OO(n^{-1})$ overshoot at a discrete crossing. Therefore
\[
 n\eta\,\Im\widehat m_s^n(z)
 =n\int\frac{\eta^2}{|x-z|^2}\,\mu_s^n(\rd x)
 \ge c_0>0.
\]
On the other hand, the deterministic support distance,
\eqref{s4:kernel-split} and \eqref{s4:coarse-improvement} give
\[
 n\eta\,\Im\widehat m_s^n(z)
 \le Cn\eta^2+n\eta|\delta_s(z)|
 \le Cc^2+Cc n^{-1/8}.
\]
Taking $c$ sufficiently small and then $n$ sufficiently large
gives a contradiction.

Thus no support crossing occurs up to $\tau_{\mathrm c}$.
The bounds \eqref{s4:coarse-outer-log} and
\eqref{s4:coarse-improvement} also exclude either transform
threshold. Hence, with probability $1-\oo(1)$,
\[
 \tau_{\mathrm c}=a_n,\qquad
 \sup I_{a_n}<r(a_n)=E_{a_n}+\varepsilon.
\]
Since $\varepsilon>0$ was arbitrary, the support assertion follows.

For a prescribed $\mathcal K$, choose the terminal margins and
contour height so that $\mathcal K$ lies outside $\Gamma(a_n)$
at a fixed distance.
Step 2 and $\mathbb P(\tau_{\mathrm c}=a_n)\to1$ yield
\(
 \sup_{z\in\Gamma(a_n)}|\delta_{a_n}(z)|
 =\OO_{\mathbb P}(n^{-1}).
\)
Apply \eqref{e:cauchy int for stij} to obtain the estimate on
$\mathcal K$, and Cauchy's formula on a slightly larger compact
set for the derivative estimates.

For corners, use interlacing and the rightmost atom in Step 4.
For DBM, include the left comparison curve and $\sigma_R$ in the
stopping argument, and use \eqref{s4:dbm-error-decomposition},
\eqref{eq:dbm-bracket} and the BDG inequality.
The same estimates follow, and \eqref{s4:dbm-confinement}
controls the additional stopping event.
\end{proof}

\subsubsection{Upper bound at the edge}\label{subsec:s4-upper}

In this subsection, we prove \Cref{large-eigen: DMC} for all three models.
We give the proofs for DMC and weighted $\beta$-corners and
briefly give the continuous-time modification for DBM, using
\Cref{lem:s4-dbm-characteristic-evolution}.
For the estimates of discrete model, fix $0<\delta<1/3$ in
\Cref{thm:main}.

\paragraph{Control function and characteristics near edge.}
We first define the control function and choose the char flows that is useful for \Cref{c:kappa:DMC}.
Fix $T\in(s_*,T_0)$ in \Cref{large-eigen: DMC} and choose
$s_*<a<T$. Recall $a_n$ from \eqref{s4:working-time}.
All local estimates below are used on $[a_n,T_0]$.
Constants may depend on $a$ and $T-a$, but not on $M$ or $n$.

For a sufficiently large fixed $M$, define
\begin{equation}\label{e:deff}
 \sqrt{f(s)}=Mn^{-1/3}+v|s-T|,
 \qquad a_n \le s\le T_0,
\end{equation}
where $v>0$ is a sufficiently small fixed constant,
depending only on $T_0$ and the constants in
\Cref{c:kappa:DMC}.

To choose the time mesh and the end points of the
characteristics, fix a sufficiently small constant
$c>0$ and set
\begin{equation}\label{s4:edge-probe}
 \eta(t)=c n^{-1/2}f(t)^{1/4},
 \qquad a_n \le t\le T_0.
\end{equation}

Starting from $t_0=a_n$, define
\begin{equation}\label{s4:edge-time-mesh}
 t_{j+1}=
 \begin{cases}
 \displaystyle
 \min\left\{
 \frac{\lfloor nT_0\rfloor}{n},
 t_j+\frac1n
 \left\lfloor cn\frac{\eta(t_j)}{\sqrt{f(t_j)}}\right\rfloor
 \right\},
 &\text{DMC and corners},\\[6pt]
 \displaystyle
 \min\left\{
 T_0,t_j+c\frac{\eta(t_j)}{\sqrt{f(t_j)}}
 \right\},
 &\text{DBM}.
 \end{cases}
\end{equation}
Stop the recursion when the right end point is reached,
and denote the last mesh time by $t_N$.
For DMC and corners, this gives a finite mesh for all sufficiently large $n$
\[
 a_n=t_0<t_1<\cdots<t_N=\frac{\lfloor nT_0\rfloor}{n},
 \qquad t_j\in\mathbb Z_n.
\]
Thus $t_j\in\mathbb Z_n$ and
$t_N=\lfloor nT_0\rfloor/n$ for DMC and corners,
whereas $t_N=T_0$ for DBM.
Except possibly for the final truncated interval,
$t_{j+1}-t_j\asymp\eta(t_j)/\sqrt{f(t_j)}$.
The final interval satisfies the same upper bound and is
treated separately whenever a lower bound is needed.

\begin{figure}
    \centering
    \includegraphics[width=0.5\linewidth]{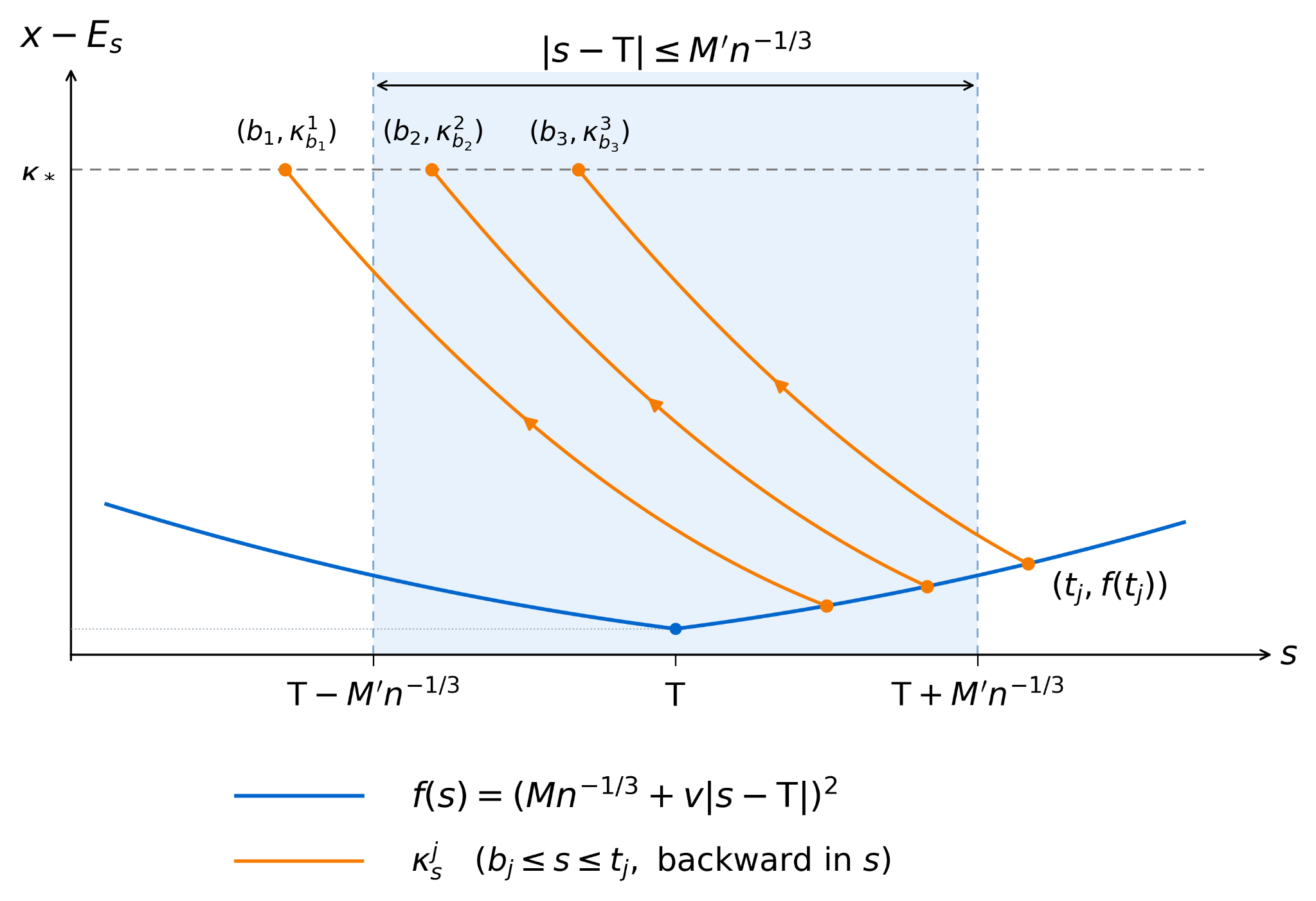}
    \caption{{Control function and local characteristics.}
    The blue curve is $f(s)$ from \eqref{e:deff}.
    The orange curves show
    $\kappa_s^j=\operatorname{Re}z_s^j-E_s$
    for $b_j\le s\le t_j$, with end points
    $(b_j,\kappa_{b_j}^j)$ and $(t_j,f(t_j))$.
    The arrows indicate decreasing $s$.
    The displayed paths have $b_j>a_n$, with
    $\kappa_{b_j}^j=\kappa_*+\OO(n^{-1})$.
    By \Cref{lem:s4-edge-path-geometry},
    $\kappa_s^j\ge f(s)$ on $[b_j,t_j]$.
    The shaded region is
    $|s-T|\le M'n^{-1/3}$.}
    \label{fig:edge-control}
\end{figure}

For each mesh time $t_j$, consider the deterministic
characteristic with end point
\begin{equation}\label{e:edge char flow spec}
z_{t_j}^j=E_{t_j}+f(t_j)+\ri\eta(t_j),
 \qquad z_s^j=E_s+\kappa_s^j+\ri\eta_s^j,
\end{equation}
Trace this characteristic backwards until time $a_n$
or its first contact with $\kappa_s^j=\kappa_*$,
where $\kappa_*$ is from \Cref{c:kappa:DMC}.
For DBM, let $b_j$ be this time.
For the discrete models, let $b_j$ be the smallest
element of $\mathbb Z_n$ greater than or equal to this time.
We use $(z_s^j)_{b_j\le s\le t_j}$ only on this interval.
The control function and the selected backward characteristics
are illustrated in \Cref{fig:edge-control}.

The following lemma estimates $\kappa_s^j$ and $\eta_s^j$
on $[b_j,t_j]$.
When $b_j>a_n$, it also places $z_{b_j}^j$ outside
$\cC(b_j)$ at a fixed positive distance.

\begin{lemma}\label{lem:s4-edge-path-geometry}
For sufficiently large fixed $M$ and $n\ge n_0(M)$,
the control function $f$ in \eqref{e:deff} satisfies
\[
 |\sqrt{f(s)}-\sqrt{f(u)}|\le v|s-u|,
 \qquad f(s)\le\kappa_*/16,
 \qquad f(a_n)\ge\frac{v^2}{4}(T-a)^2>0,
 \qquad s,u\in[a_n,T_0].
\]
For each $j$, on $[b_j,t_j]$ we have
$0<\eta_s^j\le\kappa_s^j\le\kappa_*$ and
\begin{equation}\label{e: realimage}
 \begin{aligned}
 \sqrt{\kappa_s^j}
 &\asymp \sqrt{f(t_j)}+t_j-s,\\
 \frac{\eta_s^j}{\sqrt{\kappa_s^j}}
 &\asymp-\Im\chi_s(z_s^j)
 \asymp\frac{\eta(t_j)}{\sqrt{f(t_j)}}.
 \end{aligned}
\end{equation}
Moreover,
\begin{equation}\label{s4:edge-barrier-distance}
 \kappa_u^j-f(u)
 \ge\frac{s-u}{C}
       \bigl(\sqrt{\kappa_s^j}+s-u\bigr),
 \qquad b_j\le u\le s\le t_j.
\end{equation}
If $b_j>a_n$, then
\[
 0\le \kappa_*-\kappa_{b_j}^j\le \frac{C}{n},
 \qquad
 \eta_{b_j}^j\le CM^{-1/2}n^{-1/3},
\]
with $\kappa_{b_j}^j=\kappa_*$ for DBM. Consequently,
$z_{b_j}^j$ lies outside $\mathcal C(b_j)$ and
\(\operatorname{dist}\bigl(z_{b_j}^j,\mathcal C(b_j)\bigr)
 \ge \frac{\kappa_*}{2}.
\)
The constants depend only on the fixed parameters and working
domains, and are independent of $j$, $M$ and $n$.
\end{lemma}

\begin{proof}
The properties of $f$ follow from \eqref{e:deff},
taking $v$ sufficiently small and then $n\ge n_0(M)$.
At each end point,
\[
 \kappa_{t_j}^j=f(t_j),
 \qquad
 \frac{\eta_{t_j}^j}{\kappa_{t_j}^j}
 =cn^{-1/2}f(t_j)^{-3/4}\le cM^{-3/2}.
\]
Thus \Cref{c:kappa:DMC} applies and gives the sector
condition and \eqref{e: realimage} on $[b_j,t_j]$.

We next compare the characteristic with the control
function. By \eqref{e:sharp real part DMC},
\(
 -C\le\frac{\rd}{\rd s}\sqrt{\kappa_s^j}
 \le-C^{-1}.
\)
Together with the Lipschitz bound on $\sqrt f$,
this yields
\[
 \sqrt{\kappa_u^j}-\sqrt{f(u)}
 \ge
 \sqrt{\kappa_s^j}-\sqrt{f(s)}
 +(C^{-1}-v)(s-u),
 \qquad b_j\le u\le s\le t_j.
\]
Choose $v<(2C)^{-1}$.
Taking $s=t_j$ gives $\kappa_u^j\ge f(u)$
on $[b_j,t_j]$.
Multiplying the preceding inequality by
$\sqrt{\kappa_u^j}+\sqrt{f(u)}$ and using
\(
 \sqrt{\kappa_u^j}
 \ge\sqrt{\kappa_s^j}+C^{-1}(s-u)
\)
proves \eqref{s4:edge-barrier-distance}.

Suppose now that $b_j>a_n$.
For DBM, $\kappa_{b_j}^j=\kappa_*$.
For the discrete models, rounding the contact time
up to the next grid time gives
$0\le\kappa_*-\kappa_{b_j}^j\le C/n$.
Also, \eqref{e: realimage} implies
\[
 \eta_{b_j}^j
 \le Cn^{-1/2}f(t_j)^{-1/4}
 \le CM^{-1/2}n^{-1/3}=\oo(1).
\]
By the definition of $\cC(s)$, its right side at time
$b_j$ has real part $E_{b_j}+\kappa_*/4$.
Hence $z_{b_j}^j$ lies outside $\cC(b_j)$ and
\(
 \operatorname{dist}(z_{b_j}^j,\cC(b_j))
 \ge \frac{\kappa_*}{2}
\)
for sufficiently large $n$.
\end{proof}

\paragraph{Exterior characteristics.}
Use the partition $(u_k)$ from $u_0=a_n$ to $u_J=t_N$
and the characteristic maps $\Phi_{u_k,s}$ from
\Cref{subsec:s4-characteristics}.
For each $k$, choose a finite $1/n$-mesh
$\mathscr D_k^{\mathrm{far}}\subset\cC(u_k)$ satisfying
\begin{equation}\label{def of Dfar}
 \sup_{z\in\cC(u_k)}
 \min_{w\in\mathscr D_k^{\mathrm{far}}}|z-w|
 \le\frac1n,
 \qquad
 \#\mathscr D_k^{\mathrm{far}}=\OO(n).
\end{equation}
For $w\in\mathscr D_k^{\mathrm{far}}$, consider the
characteristic $\Phi_{u_k,s}(w)$ on $[u_k,u_{k+1}]$.
The spatial derivative bound for $\Phi_{u_k,s}$ shows that
\[
 \sup_{z\in\widetilde{\cC}_k(s)}
 \min_{w\in\mathscr D_k^{\mathrm{far}}}
 |z-\Phi_{u_k,s}(w)|
 \le\frac Cn,
 \qquad u_k\le s\le u_{k+1}.
\]

\paragraph{Stopping time.}
We first record the initial bounds at $a_n$.
By \Cref{lem:s4-edge-path-geometry} and \eqref{e:deff},
for sufficiently large $n$,
\[
 \kappa_{a_n}^j\ge f(a_n)
 \ge\frac{v^2}{4}(T-a)^2>0,
 \qquad b_j=a_n.
\]
These starting points remain in the fixed spectral domains.
Together with the contour separation, this allows us to choose
a fixed compact set
$\mathcal K_a\Subset\Lambda\setminus[\ell_0,E_a]$
containing fixed neighborhoods of $\cC(a_n)$,
$\omega_+(a_n)$ and all $z_{a_n}^j$ with $b_j=a_n$.
For DBM, replace $[\ell_0,E_a]$ by $[L_a,E_a]$ and include
the corresponding left-edge points.

Choose $\varepsilon>0$ smaller than $\kappa_*/16$,
$v^2(T-a)^2/8$ and the fixed contour separations at $a_n$.
By \Cref{prop:s4-positive-time-preparation} and
\eqref{s4:dbm-confinement}, choose events
$A_n\in\mathscr F_{a_n}$ with
\[
 A_n\subset\mathcal G_n^0\cap\{\sigma_R>a_n\},
 \qquad q_n:=\mathbb P(A_n^c)=\oo(1),
\]
such that, on $A_n$,
\begin{equation}\label{s4:initial-good-event}
 \sup I_{a_n}<E_{a_n}+\varepsilon,
 \qquad
 \sup_{z\in\mathcal K_a}|\delta_{a_n}(z)|
 \le\frac{\log n}{n}.
\end{equation}
For DBM, also require $\inf I_{a_n}>L_{a_n}-\varepsilon$.
The set $\mathcal K_a$, the choice of $\varepsilon$ and the events
$A_n$ can be chosen independently of $M$.
We use the same event for all local edge characteristics.

On $A_n$, define the stopping time below on $[a_n,t_N]$.
On $A_n^c$, set $\tau=a_n$.

\begin{definition}\label{def:s4-edge-stopping}
On $A_n$, with $\inf\varnothing=\infty$, define
\begin{equation}\label{def:stoptime edge DMC}
 \begin{aligned}
 \tau=t_N\wedge\sigma_R\wedge\inf\Bigl\{s\in[a_n,t_N]:\;&
 \sup_{z\in\cC(s)}|\delta_s(z)| \ge\frac{(\log n)^{10}}n,\quad\text{or}\\
 &|\delta_s(z_s^j)|\ge
 \frac{(n^{1/3}\sqrt{\kappa_s^j})^{1/8}}
      {n\sqrt{\kappa_s^j\eta_s^j}}
 \quad\text{for some }j\text{ with }b_j\le s\le t_j
 \Bigr\}.
 \end{aligned}
\end{equation}
For DMC and corners, the infimum is over grid times;
For DBM, we also include the corresponding left-edge
conditions in the definition of the stopping time. 
\end{definition}

\begin{remark}
For DBM, the same argument applies at the left edge, so we omit the details.
\end{remark}

For each local characteristic, write its accumulated stopped martingale as
\begin{equation}\label{s4:edge-stopped-martingale}
 N_s^j=
 \begin{cases}
 \displaystyle\mathbf1_{A_n}\sum_{u\in[b_j,s)\cap\mathbb Z_n}
 \mathbf1_{\{u<\tau\}}\xi_u(z_{u+1/n}^j),
 &\text{DMC and corners},\\[6pt]
 \displaystyle \mathbf1_{A_n}\int_{b_j}^{s}
 \mathbf1_{\{u\le\tau\}}\,\rd \mathcal M_u^n(z_u^j),
 &\text{DBM},
 \end{cases}
 \qquad b_j\le s\le t_j.
\end{equation}
These are martingales by \eqref{eq:increment-common} and
\eqref{eq:pote_mt}. In the discrete case, the update reaching
$\tau$ is included. We use the same stopping convention on the exterior characteristics.

\begin{proposition}\label{prop:s4-edge-full-time}
For each of the three models, for every sufficiently large fixed
$M$ and all $n\ge n_0(M)$, with probability at least
\[
 1-q_n-Ce^{-cM^{1/8}}-Cn^{-20}-Ce^{-c(\log n)^2},
\]
we have $\tau=t_N$, and
\begin{equation}\label{s4:full-support}
 \sup I_s<E_s+f(s),\qquad a_n\le s\le t_N.
\end{equation}
For DBM, $\inf I_s>L_s-f(s)$ on the same event and interval.
The constants are independent of $M$ and $n$.
\end{proposition}

We first deduce the support bound before $\tau$, then estimate
the exterior char flows error, the local drift, and the martingales.

\begin{proposition}
    \label{l:beforestop DMC}
For each of the three models, on $A_n$, for $a_n\le s<\tau$,
\[
 \sup I_s<E_s+f(s).
\]
\end{proposition}

\begin{proof}
Work on $A_n$. At $s=a_n$, the assertion follows from
\eqref{s4:initial-good-event} and the choice of $\varepsilon$.
Fix $a_n<s<\tau$ and choose the first mesh time $t_j\ge s$.
By \eqref{e:deff}, \eqref{s4:edge-probe}, and
\eqref{s4:edge-time-mesh},
\[
 t_j-t_{j-1}
 \le Cn^{-1/2}f(t_{j-1})^{-1/4}
 \le CM^{-1/2}n^{-1/3}=\oo(1).
\]
If $b_j>a_n$, then $\kappa_{b_j}^j=\kappa_*+O(n^{-1})$.
By \eqref{e:sharp real part DMC} and
$f(t_j)\le\kappa_*/16$,
\[
 t_j-b_j
 \ge C^{-1}
 \bigl(\sqrt{\kappa_{b_j}^j}-\sqrt{f(t_j)}\bigr)
 \ge c>0.
\]
Since $t_j-t_{j-1}=o(1)$, we have
$b_j\le t_{j-1}\le s$ for sufficiently large $n$.
This also holds when $b_j=a_n$.

We can therefore apply \eqref{e: realimage} on $[s,t_j]$.
Together with \eqref{s4:edge-time-mesh} and the Lipschitz
bound on $\sqrt f$, this gives
\[
 t_j-s\le\frac{Cc\eta(t_j)}{\sqrt{f(t_j)}},
 \qquad
 \kappa_s^j\asymp f(t_j)\asymp f(s),
 \qquad
 \eta_s^j\asymp\eta(t_j).
\]
Integrating
$|\partial_r\kappa_r^j|\le C\sqrt{\kappa_r^j}$
from \eqref{e:sharp real part DMC}, we obtain
\[
 0\le\kappa_s^j-f(s)
 \le C\sqrt{f(t_j)}(t_j-s)
 \le Cc\eta_s^j.
\]

Since $s<\tau$, the local transform bound and
\eqref{s4:im-edge} imply
\[
 \begin{aligned}
 n\eta_s^j\Im m_s^n(z_s^j)
 &\le n\eta_s^j\Im m_s(z_s^j)
       +n\eta_s^j|\delta_s(z_s^j)|\\
 &\le Cn\frac{(\eta_s^j)^2}{\sqrt{\kappa_s^j}}
 +(n^{1/3}\sqrt{\kappa_s^j})^{1/8}
       \sqrt{\frac{\eta_s^j}{\kappa_s^j}}\\
 &\le Cc^2+C\sqrt c\,M^{-5/8}.
 \end{aligned}
\]
Suppose now that
\[
 |\sup I_s-E_s-f(s)|\le\eta_s^j.
\]
Then $|\sup I_s-\Re z_s^j|\le C\eta_s^j$.
For corners and DBM, the rightmost atom has mass $1/n$,
so its contribution to the positive Stieltjes kernel gives
\[
 n\eta_s^j\Im m_s^n(z_s^j)\ge c_0,
 \qquad c_0>0.
\]
For DMC, the rightmost interval has mass $1/n$ and
length $\theta/n=\oo(\eta_s^j)$.
Integrating the ordinary kernel over this interval and
using \eqref{s4:kernel-split} therefore gives
\begin{equation}\label{s4:one-particle-signal}
 n\eta_s^j\Im m_s^n(z_s^j)
 \ge c_0-Cn(\eta_s^j)^2.
\end{equation}
Since $n(\eta_s^j)^2\le Cc^2$, both cases contradict
the preceding upper bound when $c$ is sufficiently
small and then $M$ is sufficiently large.
Consequently, for every $s<\tau$,
\[
 \sup I_s\notin
 [E_s+f(s)-\eta_s^j,\ E_s+f(s)+\eta_s^j].
\]

It remains to use the dynamics to exclude a crossing.
By \eqref{s4:initial-good-event} and the choice of $\varepsilon$,
$\sup I_{a_n}<E_{a_n}+\varepsilon<E_{a_n}+f(a_n)$.
Suppose that the support first reaches or crosses the
control curve at a time $a_n<b<\tau$.
For DBM, continuity gives
$\sup I_b=E_b+f(b)$, contradicting the exclusion
just proved.
For the discrete models, the bounded upward steps of
DMC, the interlacing of corners, and the Lipschitz
continuity of $E_s+f(s)$ give
\[
 0\le\sup I_b-E_b-f(b)\le C/n.
\]
Choosing the first mesh time $t_j\ge b$, we have
$n^{-1}=\oo(\eta_b^j)$, so this also places $\sup I_b$
in the excluded interval.
Thus no crossing occurs before $\tau$, proving the lemma.
\end{proof}

\paragraph{Exterior estimates.}
We estimate $\delta_s$ along the exterior characteristics before the stopping time $\tau$
and transfer the bounds to $\cC(s)$ and $\omega_+(s)$
using \Cref{lem:exterior-cauchy}.

On $A_n$, the same argument as in \Cref{l:beforestop DMC}
gives, for DBM,
\[
 \inf I_s>L_s-f(s),\qquad a_n\le s<\tau.
\]
Together with \Cref{lem:b-bounds,l:beforestop DMC} and
the bound $f(s)\le\kappa_*/16$ from
\Cref{lem:s4-edge-path-geometry}, this yields
\[
 I_s\subset
 \begin{cases}
  [\ell_0,E_s+\kappa_*/8],&\text{DMC and corners},\\
  [L_s-\kappa_*/8,E_s+\kappa_*/8],&\text{DBM},
 \end{cases}
 \qquad a_n\le s<\tau.
\]
These inclusions also hold at $\tau$ for sufficiently large $n$.
If $\tau=a_n$, use \eqref{s4:initial-good-event}.
If $\tau>a_n$, this follows for DBM by continuity.
For the discrete models, we have
\[
 \sup I_\tau\le E_\tau+f(\tau)+C/n,
 \qquad
 \inf I_\tau\ge\ell_0,
\]

On $A_n$, for $u_k\le\tau$ and $s\in[u_k,u_{k+1}\wedge\tau]$,
the contours $\cC(s)$ and $\widetilde{\cC}_k(s)$ enclose
$\operatorname{supp}\mu_s^n\cup\operatorname{supp}\mu_s$
and remain a fixed positive distance from this union.
Thus \Cref{lem:exterior-cauchy} applies.

\begin{lemma}\label{s4:outer-log}
For each of the three models, outside an event of probability
$Ce^{-c(\log n)^2}$, the following estimates hold on $A_n$:
\begin{equation}\label{s4:outer-log-bound}
 \sup_{a_n\le s\le\tau}
 \left(
 \sup_{z\in\cC(s)}|\delta_s(z)|
 +\sup_{z\in\omega_+(s)}|\delta_s(z)|
 \right)
 \le\frac{C (\log n)^2}{n}.
\end{equation}
Moreover, for all $j$ with $b_j\le\tau$,
\begin{equation}\label{s4:entrance-bound}
 |\delta_{b_j}(z_{b_j}^j)|
 \le\frac{C (\log n)^2}{n}.
\end{equation}
The constants are independent of $M$ and $n$.
The same estimates hold with $\tau$ replaced by any
stopping time taking values in $[a_n,\tau]$.
\end{lemma}

\begin{proof}
On $A_n$, the support bounds and contour separation established
above imply that both transforms and their spatial derivatives are uniformly bounded
on fixed neighborhoods of the contours up to $\tau$.
Since $|\mathsf b'|\asymp1$, \eqref{s4:contour-transfer} and
\eqref{s4:outer-cauchy} give
\[
 \sup_{z\in\cC(s)}|\delta_s(z)|
 +\sup_{z\in\omega_+(s)}|\delta_s(z)|
 \le C\sup_{z\in\widetilde{\cC}_k(s)}|\delta_s(z)|,
 \qquad u_k\le s\le u_{k+1}\wedge\tau.
\]

We first consider the discrete models.
For $w\in\mathscr D_k^{\mathrm{far}}$, the fixed-distance case
of \eqref{s4:moment-input} gives, on $\{u<\tau\}$,
\[
 \mathbb E_u
 \left|\xi_u\bigl(\Phi_{u_k,u+1/n}(w)\bigr)\right|^2
 \le \frac C{n^3},
 \qquad
 \left|\xi_u\bigl(\Phi_{u_k,u+1/n}(w)\bigr)\right|
 \le \frac Cn.
\]
The second bound follows from the bounded particle steps for DMC.
For corners, it follows by integration by parts in the interlacing
distribution functions, including the boundary term from the loss
of mass $1/n$.
Since $I_{s+1/n}\subset I_s+[0,1/n]$ for DMC and
$I_{s+1/n}\subset I_s$ for corners, the same contour
bounds hold at time $s+1/n$.

Using the stopping convention in \eqref{s4:edge-stopped-martingale},
Freedman's maximal inequality therefore gives
\[
 \mathbb P\left(
 \sup_{u_k\le t\le u_{k+1}}
 \left|
 \sum_{u\in[u_k,t)\cap\mathbb Z_n}
 \mathbf1_{\{u<\tau\}}
 \xi_u\bigl(\Phi_{u_k,u+1/n}(w)\bigr)
 \right|
 >\frac{C(\log n)^2}{n}
 \right)
 \le Ce^{-c(\log n)^2}.
\]
Since $J=\OO(1)$ and $\#\mathscr D_k^{\mathrm{far}}=\OO(n)$,
a union bound gives the same estimate simultaneously for all
$k,w$, outside an event of the stated probability.

Work on $A_n$ and this event, and fix $k$ with $u_k\le\tau$.
Sum \eqref{sde for DMC 2} along $\Phi_{u_k,s}(w)$ and apply
the exterior estimates in \Cref{lem:s4-characteristic-evolution}.
By \eqref{def of Dfar} and the spatial derivative bound for the
flow, the images of the mesh form a $C/n$-mesh of
$\widetilde{\cC}_k(t)$.
The contour comparison above and the uniform spatial derivative
bound thus yield
\begin{equation}\label{s4:outer-gronwall}
 \begin{aligned}
 \sup_{z\in\widetilde{\cC}_k(t)}|\delta_t(z)|
 \le{}&
 \sup_{z\in\cC(u_k)}|\delta_{u_k}(z)|
 +\frac{C(\log n)^2}{n} +\frac Cn\sum_{u\in[u_k,t)\cap\mathbb Z_n}
       \sup_{z\in\widetilde{\cC}_k(u)}|\delta_u(z)|,
 \quad u_k\le t\le u_{k+1}\wedge\tau.
 \end{aligned}
\end{equation}
Here the accumulated $\OO(n^{-2})$ drift errors and the mesh
approximation each contribute $\OO(n^{-1})$.

Discrete Gronwall applied to \eqref{s4:outer-gronwall} gives
\begin{equation}\label{s4:weighted-closure}
 \sup_{u_k\le t\le u_{k+1}\wedge\tau}
 \sup_{z\in\widetilde{\cC}_k(t)}|\delta_t(z)|
 \le C\left(
 \sup_{z\in\cC(u_k)}|\delta_{u_k}(z)|
 +\frac{(\log n)^2}{n}
 \right).
\end{equation}
On $A_n$, \eqref{s4:initial-good-event} gives
$\sup_{z\in\cC(a_n)}|\delta_{a_n}(z)|\le(\log n)/n$.
At each subsequent time $u_{k+1}\le\tau$,
\eqref{s4:contour-transfer} bounds the error on $\cC(u_{k+1})$
by the preceding estimate on $\widetilde{\cC}_k(u_{k+1})$.
Iterating over the fixed number $J$ of intervals and using
\eqref{s4:outer-cauchy} proves \eqref{s4:outer-log-bound}.

For DBM, use the stopped martingale integral along
$\Phi_{u_k,s}(w)$ in place of the discrete sum.
By \Cref{lem:s4-dbm-characteristic-evolution}, its
quadratic-variation rate is bounded by $C/n^2$, and the
additional drift is bounded by $C/n$.
The continuous exponential martingale inequality gives the
same simultaneous martingale bound.
Replacing the sums in \eqref{s4:outer-gronwall} by integrals,
the preceding argument proves \eqref{s4:outer-log-bound}
for DBM as well.

Finally, if $b_j>a_n$, \Cref{lem:s4-edge-path-geometry} places
$z_{b_j}^j$ outside $\cC(b_j)$ at a fixed positive distance.
Thus \eqref{e:cauchy int for stij} and
\eqref{s4:outer-log-bound} give
\[
 |\delta_{b_j}(z_{b_j}^j)|
 \le C\sup_{z\in\cC(b_j)}|\delta_{b_j}(z)|
 \le\frac{C(\log n)^2}{n},
 \qquad a_n<b_j\le\tau.
\]
If $b_j=a_n$, then $z_{a_n}^j\in\mathcal K_a$, and
\eqref{s4:initial-good-event} gives
$|\delta_{a_n}(z_{a_n}^j)|\le(\log n)/n$.
This proves \eqref{s4:entrance-bound}.
The stopping convention includes the update reaching $\tau$,
and the same event controls every earlier stopping time.
\end{proof}

\paragraph{Accumulated error terms.}
We estimate the accumulated nonmartingale terms in
$\mathcal E_2$.
We also bound the predictable quadratic variation of
the stopped martingales in \eqref{s4:edge-stopped-martingale},
which will be used to estimate the accumulated martingale
in the next paragraph.
The support bound in \Cref{l:beforestop DMC} supplies
the spectral distances needed for these estimates.

\begin{lemma}\label{lem:s4-edge-errors}
For each of the three models and every $j$, when
$b_j\le u<t_j\wedge\tau$, we have
\begin{equation}\label{s4:edge-P}
 |\delta_u(z_u^j)|\le CM^{-5/8}[-\Im\chi_u(z_u^j)],
 \qquad P_u^n(z_u^j)\le\frac C{\sqrt{\kappa_u^j}}.
\end{equation}
For $b_j\le u<s\le t_j$ with $u<\tau$,
\begin{equation}\label{s4:edge-distance}
 \operatorname{dist}(z_u^j,I_u)
 \ge\frac1C(\sqrt{\kappa_s^j}+s-u)
 \sqrt{(s-u)^2+\frac{(\eta_s^j)^2}{\kappa_s^j}}.
\end{equation}

The stopped martingale in \eqref{s4:edge-stopped-martingale} satisfies
\begin{equation}\label{s4:edge-bracket}
 \langle\Re N^j\rangle_s+\langle\Im N^j\rangle_s
 \le\frac C{n^2\kappa_s^j\eta_s^j},
 \qquad b_j\le s\le t_j.
\end{equation}
For $b_j\le s\le t_j\wedge\tau$, the discrete models satisfy
\begin{equation}\label{s4:edge-drift-error}
 \sum_{u\in[b_j,s)\cap\mathbb Z_n}
 \left(
 |\mathcal E_2(u)-\xi_u(z_{u+1/n}^j)|
 +\frac{P_u^n(z_u^j)}{n^2}
 \right)
 \le\frac{C\log(2+n^{1/3}\sqrt{\kappa_s^j})}{n\kappa_s^j}.
\end{equation}
For DBM, the same upper bound holds for
\(
 \frac{|1-\theta|}{2\theta n}
 \int_{b_j}^s|\partial_z^2m_u^n(z_u^j)|\,\rd u.
\)
The constants are independent of $j,M,n$.
\end{lemma}

\begin{proof}
Fix $j$ and omit it from the path notation.
For $u<\tau$, \eqref{def:stoptime edge DMC} and
\eqref{e: realimage} give
\[
 \frac{|\delta_u(z_u)|}{-\Im\chi_u(z_u)}
 \le C\bigl(n^{1/3}\sqrt{f(t_j)}\bigr)^{3/4}
       \bigl(n^{1/3}\sqrt{\kappa_u}\bigr)^{-11/8}
 \le CM^{-5/8}.
\]
They also give
$|m_u^n(z_u)-m_u(z_u)|\le Cn^{-1/3}M^{-9/8}=\oo(1)$.
On $\omega_+(u)$, \Cref{lem:exterior-cauchy} and
the definition of $\tau$ imply
\[
 \sup_{w\in\omega_+(u)}|m_u^n(w)-m_u(w)|
 \le \frac{C(\log n)^{10}}n=\oo(1).
\]
We can therefore Taylor-expand the flux in its first
variable, as in the proof of
\Cref{lem:s4-characteristic-evolution}.
The resulting comparison of $\chi_u^n$ and $\chi_u$,
together with \eqref{eq:main-P}, gives
\[
 P_u^n(z_u)
 \le C+C\frac{-\Im\chi_u(z_u)+|\delta_u(z_u)|}{\eta_u}
       +\frac CnP_u^n(z_u).
\]
For DBM, this follows directly from its velocity formula
in \eqref{eq:velocity-from-flux}.
The last term is present only for corners and can be
absorbed for sufficiently large $n$.
Combining the first estimate above with
\eqref{e: realimage} proves \eqref{s4:edge-P}.

We next estimate the distance from $z_u$ to $I_u$.
By \Cref{l:beforestop DMC}, $\sup I_u<E_u+f(u)$
for $u<\tau$. Hence
\[
 \operatorname{dist}(z_u,I_u)
 \ge \sqrt{(\kappa_u-f(u))^2+\eta_u^2}.
\]
For $b_j\le u<s\le t_j$, the geometric estimates
\eqref{s4:edge-barrier-distance} and \eqref{e: realimage}
give
\[
 \kappa_u-f(u)
 \ge \frac{s-u}{C}\bigl(\sqrt{\kappa_s}+s-u\bigr),
 \qquad
 \eta_u\asymp
 \frac{\eta_s}{\sqrt{\kappa_s}}
       \bigl(\sqrt{\kappa_s}+s-u\bigr).
\]
Substituting these bounds proves \eqref{s4:edge-distance}.

For the discrete models, we also check the estimates
along each one-step segment.
By \eqref{s4:edge-probe} and \eqref{e: realimage},
\(
 \eta_u\ge c'M^{1/2}n^{-2/3}\gg n^{-1+\delta},
\)
where $0<\delta<1/3$ is the exponent chosen above.
Since $|z_r-z_u|=\OO(n^{-1})=\oo(\eta_u)$ for
$u\le r\le u+1/n$, the distance bound changes by at most
a fixed factor when $z_u$ is replaced by $z_r$.
Moreover, the Stieltjes integral bounds in
\Cref{lem:stieltjes-bounds} give
\[
 \sup_{u\le r\le u+1/n}
 |m_u^n(z_r)-m_u^n(z_u)|
 \le \frac C{n\eta_u}=\oo(1).
\]
For DMC, combining this with the bound on $\delta_u(z_u)$
and the macroscopic characteristic equation
$\frac{\rd}{\rd r}m_r(z_r)=F_r(z_r)$ gives
\[
 \sup_{u\le r\le u+1/n}
 |m_u^n(z_r)-m_r(z_r)|
 \le |\delta_u(z_u)|+\frac C{n\eta_u}+\frac Cn=\oo(1).
\]
Here we use the uniform bound on $F_r(z_r)$.
Thus the additional DMC condition in
\Cref{lem:s4-characteristic-evolution} holds for
sufficiently large $n$.
Thus the one-step estimates and the Taylor expansions
in \Cref{lem:s4-characteristic-evolution} apply along
the whole segment, with the configuration held at time $u$.

We now accumulate the conditional variances.
For the discrete models, the stopping convention in
\eqref{s4:edge-stopped-martingale} gives
\[
 \langle\Re N^j\rangle_s+\langle\Im N^j\rangle_s
 =
 \sum_{u\in[b_j,s)\cap\mathbb Z_n}
 \mathbf1_{\{u<\tau\}}
 \mathbb E_u|\xi_u(z_{u+1/n})|^2.
\]
By \eqref{s4:moment-input}, \eqref{s4:edge-P},
and \eqref{s4:edge-distance}, this is bounded by
\[
 \frac C{n^2}\int_0^\infty
 \frac{\rd x}
 {(\sqrt{\kappa_s}+x)^3
  (x^2+\eta_s^2/\kappa_s)}
 \le \frac C{n^2\kappa_s\eta_s}.
\]
Here $x=s-u$. The sum is bounded by the displayed
integral up to a fixed factor, since
$\eta_s/\sqrt{\kappa_s}\gtrsim n^{-1/2}\gg n^{-1}$.
Only updates with $u<\tau$ contribute, so the estimate
holds for every deterministic $s\in[b_j,t_j]$,
including $s>\tau$.
For DBM, \eqref{s4:dbm-martingale-bound} gives the same
integral directly. This proves \eqref{s4:edge-bracket}.

Finally, we accumulate the nonmartingale terms.
For the discrete models, \eqref{s4:step-remainder}
and the same bounds for $P_u^n$ and the spectral distance
bound the left-hand side of \eqref{s4:edge-drift-error} by
\[
 \frac Cn\int_0^\infty
 \frac{\rd x}
 {(\sqrt{\kappa_s}+x)^2
  \sqrt{x^2+\eta_s^2/\kappa_s}}
 \le
 \frac{C[1+\log(1+\kappa_s/\eta_s)]}{n\kappa_s}.
\]
The additional $P_u^n(z_u)/n^2$ term is bounded by
the same integral because the spectral distances
are uniformly bounded above.
By \eqref{e: realimage} and \eqref{s4:edge-probe},
\[
 \frac{\kappa_s}{\eta_s}
 \le Cn^{1/2}\sqrt{\kappa_s}\,f(t_j)^{1/4}
 \le C\bigl(n^{1/3}\sqrt{\kappa_s}\bigr)^{3/2},
\]
where we used $f(t_j)\le\kappa_s$.
This proves \eqref{s4:edge-drift-error}.
For DBM, \eqref{s4:dbm-error-decomposition} and
\eqref{s4:dbm-second-order} bound the stated integral
by the same expression.
\end{proof}

\paragraph{Martingale estimates.}
We now use the quadratic-variation bound \eqref{s4:edge-bracket}
to control the stopped martingales $N_s^j$ from
\eqref{s4:edge-stopped-martingale} along the local edge characteristics
$z_s^j$, simultaneously for all $j$ and all $s\in[b_j,t_j]$.

\begin{proposition}\label{p:bddquadrac:DMC}
For each of the three models, for every sufficiently large
fixed $M$ and all $n\ge n_0(M)$, with probability at least
$1-Ce^{-cM^{1/8}}-Cn^{-20}$,
\begin{equation}\label{s4:edge-martingale}
 |N_s^j|
 \le
 \frac{C(n^{1/3}\sqrt{\kappa_s^j})^{1/16}}
      {n\sqrt{\kappa_s^j\eta_s^j}},
 \qquad b_j\le s\le t_j,
\end{equation}
simultaneously for every $j$.
The constants are independent of $M,n$.
\end{proposition}

\begin{proof}
We first consider the discrete models.
By \eqref{s4:moment-input}, \eqref{s4:edge-P}, and
$\operatorname{dist}(z_{u+1/n}^j,I_u)\gtrsim\eta_u^j$,
for each fixed integer $p\ge1$,
\[
 \begin{aligned}
 &(n^2\kappa_u^j\eta_u^j)^p
 \mathbb E_u\left|
 \mathbf1_{\{u<\tau\}}\xi_u(z_{u+1/n}^j)
 \right|^{2p} \le
 C_p\left(\frac{\sqrt{\kappa_u^j}}{n\eta_u^j}\right)^p
 \left[
 1+\left(\frac{\sqrt{\kappa_u^j}}
               {n(\eta_u^j)^2}\right)^{p-1}
 \right]
 \le C_p n^{-p/2}.
 \end{aligned}
\]
Here we used $\xi_u=n^{-1}\Delta\mathcal M_u^n$ and,
by \eqref{e: realimage} and \eqref{s4:edge-probe},
\[
 \frac{n\eta_u^j}{\sqrt{\kappa_u^j}}\gtrsim\sqrt n,
 \qquad
 \frac{n(\eta_u^j)^2}{\sqrt{\kappa_u^j}}\gtrsim1.
\]
The constants may depend on the fixed probe constant.

To apply Freedman's inequality, we first bound the
martingale increments uniformly.
Truncate each stopped increment at
$n^{-1/16}/(n\sqrt{\kappa_u^j\eta_u^j})$ and subtract
its conditional mean, so that the resulting increments
are still martingale differences.
Markov's inequality with $p=64$ bounds the truncation
probability by $Cn^{-24}$ per update.
There are $\OO(n^2)$ path--update pairs, so all truncations
leave the increments unchanged outside an event of
probability $Cn^{-22}$.

Since the original increments have conditional mean zero,
the centering errors are bounded by the conditional
expectations of the discarded tails.
The same moment estimate, summed over the updates, gives
\(
 \frac{Cn^{-23}}{n\sqrt{\kappa_s^j\eta_s^j}}
\)
as an upper bound for the accumulated centering error up
to time $s$. Here we used that $\kappa_u^j$ and $\eta_u^j$
decrease along each local characteristic.
Thus, outside the preceding exceptional event, every
original partial sum differs from the truncated centered
sum by at most this quantity.

Fix $A=2^\ell M$, $\ell\ge0$, and consider a path meeting
\(
 A\le n^{1/3}\sqrt{\kappa_r^j}\le2A.
\)
These times form a deterministic interval.
Let $s$ be its last grid time, omitting intervals with
no grid point.
By \eqref{e: realimage},
\(
 {n\sqrt{\kappa_r^j\eta_r^j}}
 \asymp
 {n\sqrt{\kappa_s^j\eta_s^j}}
\)
throughout this interval.

Truncation and conditional centering do not increase the
conditional second moments.
Consequently, \eqref{s4:edge-bracket} bounds the sum of
the conditional variances of the truncated centered
increments up to time $s$ by
\(
 \frac{C}{n^2\kappa_s^j\eta_s^j}.
\)
Since $\kappa_u^j$ and $\eta_u^j$ decrease in time,
every such increment before $s$ has absolute value at most
\(
 \frac{Cn^{-1/16}}{n\sqrt{\kappa_s^j\eta_s^j}}.
\)
Apply Freedman's maximal inequality to the real and
imaginary parts, from $b_j$ to $s$, with deviation
\(
 \frac{CA^{1/16}}{n\sqrt{\kappa_s^j\eta_s^j}}.
\)
Substituting the preceding variance and increment bounds
gives failure probability
\[
 C\exp\left\{
 -\frac{cA^{1/8}}{1+n^{-1/16}A^{1/16}}
 \right\}
 \le Ce^{-c'A^{1/8}}.
\]
Here $A\le Cn^{1/3}$, since $\kappa_s^j\le\kappa_*$.
The comparison of the deterministic bounds within the
interval, together with the centering-error estimate,
therefore gives \eqref{s4:edge-martingale} at every grid time
in this interval.

For DBM, take $s$ to be the last time in the same
deterministic interval.
The continuous exponential martingale inequality,
with the quadratic-variation bound \eqref{s4:edge-bracket}
and deviation $CA^{1/16}/(n\sqrt{\kappa_s^j\eta_s^j})$,
gives the same failure probability $Ce^{-c'A^{1/8}}$.
No truncation is needed.

It remains to count the paths for each $A$.
If a path meets the displayed interval, then
$n^{1/3}\sqrt{f(t_j)}\le2A$, since $\kappa_r^j$ decreases
to $f(t_j)$ at its end point.
By \eqref{e:deff}, all such end times lie in an interval
of length $\OO(An^{-1/3})$.
The mesh in \eqref{s4:edge-time-mesh} has spacing at least
$c'n^{-1/3}A^{-1/2}$ there, apart from the possible final
capped increment; grid rounding changes this bound only
by a fixed factor.
The number of such paths is therefore at most
\[
 C\left(1+\frac{An^{-1/3}}
                  {n^{-1/3}A^{-1/2}}\right)
 \le CA^{3/2}.
\]
Each path meets each interval at most once.
Taking the union over the paths and over $A=2^\ell M$
gives total failure probability at most
\[
 \sum_{\ell\ge0}
 C(2^\ell M)^{3/2}e^{-c(2^\ell M)^{1/8}}
 +Cn^{-22}
 \le Ce^{-c'M^{1/8}}+Cn^{-20}.
\]
This proves the proposition.
\end{proof}

\paragraph{Estimates along the characteristics.}
We now estimate $\delta_s(z_s^j)$ along each characteristics flow.
We use \Cref{s4:outer-log} to control the entrance error
$\delta_{b_j}(z_{b_j}^j)$ and the transform difference on
$\omega_+(s)$, \Cref{lem:s4-edge-errors} to control the
accumulated error terms, and \Cref{p:bddquadrac:DMC} to
control $N_s^j$.
Combining these bounds with the evolution equation and applying
Gr\"onwall gives the following estimate.

\begin{lemma}\label{lem:s4-edge-improvement}
On $A_n$ and the common event of \Cref{s4:outer-log,p:bddquadrac:DMC},
for every $j$ and $b_j\le s\le t_j\wedge\tau$,
\begin{equation}\label{s4:edge-improvement}
 \begin{aligned}
 |\delta_s(z_s^j)|\le{}&
 \frac{C(\log n)^2}{n\sqrt{\kappa_s^j}}
 +\frac{C(n^{1/3}\sqrt{\kappa_s^j})^{1/16}}
        {n\sqrt{\kappa_s^j\eta_s^j}} +\frac{C\log(2+n^{1/3}\sqrt{\kappa_s^j})}
        {n\kappa_s^j}.
 \end{aligned}
\end{equation}
In particular, for sufficiently large fixed $M$ and then
$n\ge n_0(M)$,
\[
 |\delta_s(z_s^j)|
 \le
 \frac{(n^{1/3}\sqrt{\kappa_s^j})^{1/8}}
      {2n\sqrt{\kappa_s^j\eta_s^j}},
 \qquad b_j\le s\le t_j\wedge\tau.
\]
This improves the bound for $|\delta_s(z_s^j)|$ in the
definition of $\tau$ in \eqref{def:stoptime edge DMC}
by a factor of two.
The constants are independent of $j,M,n$.
\end{lemma}

\begin{proof}
Fix $j$ and suppress it in the path notation.
The characteristic equation gives
\[
 \int_u^s\frac{-\Im\chi_r(z_r)}{\eta_r}\,\rd r
 =\log\frac{\eta_u}{\eta_s}.
\]
We next bound the integral of the nonlinear coefficient in
\eqref{s4:sharp-linear}.
Using the bound for $|\delta_r(z_r)|$ in
\eqref{def:stoptime edge DMC}, together with
\eqref{e: realimage} and \eqref{e:sharp real part DMC}, gives
\[
 \int_{b_j}^s
 \frac{(n^{1/3}\sqrt{\kappa_r})^{1/8}}
      {n\eta_r\sqrt{\kappa_r\eta_r}}\,\rd r
 \le C\bigl(n^{1/3}\sqrt{f(t_j)}\bigr)^{3/4}
       \bigl(n^{1/3}\sqrt{\kappa_s}\bigr)^{-11/8}
 \le CM^{-5/8}.
\]
Consequently,
\begin{equation}\label{s4:propagation}
 \exp\left\{\int_u^s\left(
 -\frac{\Im\chi_r(z_r)}{\eta_r}+C+
 \frac{C(n^{1/3}\sqrt{\kappa_r})^{1/8}}
      {n\eta_r\sqrt{\kappa_r\eta_r}}
 \right)\rd r\right\}
 \le C\frac{\eta_u}{\eta_s}.
\end{equation}
For grid times, the corresponding discrete product satisfies
the same bound.
The Riemann-sum error for the first coefficient is
$\OO((n\sqrt{\kappa_s})^{-1})$, by \eqref{s4:time-velocity}
and \eqref{e: realimage}; the nonlinear integrand is increasing,
so its left Riemann sum is bounded by its integral.

For the discrete models, sum \eqref{sde for DMC 2} along
the path. Using \eqref{s4:sharp-linear},
\eqref{s4:entrance-bound}, \eqref{s4:outer-log-bound},
and \eqref{s4:edge-drift-error} gives
\[
 \begin{aligned}
 |\delta_s(z_s)|\le{}&
 \frac{C(\log n)^2}{n}+|N_s^j|
 +\frac{C\log(2+n^{1/3}\sqrt{\kappa_s})}{n\kappa_s}\\
 &+\frac1n\sum_{u\in[b_j,s)\cap\mathbb Z_n}
 \left(
 -\frac{\Im\chi_u(z_u)}{\eta_u}+C+
 \frac{C(n^{1/3}\sqrt{\kappa_u})^{1/8}}
      {n\eta_u\sqrt{\kappa_u\eta_u}}
 \right)|\delta_u(z_u)|.
 \end{aligned}
\]
The martingale contribution is the stopped partial sum $N_s^j$;
we estimate it by \eqref{s4:edge-martingale} before applying
Gr\"onwall.

The coefficient in the sum is at most $C/\sqrt{\kappa_u}$,
by \eqref{e: realimage} and the calculation leading to
\eqref{s4:edge-P}.
By \eqref{s4:propagation}, the propagation factor is at most
$C\eta_u/\eta_s\asymp C\sqrt{\kappa_u/\kappa_s}$.
Discrete Gr\"onwall, followed by the change of variable
$x=\sqrt{\kappa_u}$, therefore bounds the additional
contribution from the martingale by
\[
 \frac{C(n^{1/3})^{1/16}}{n\sqrt{\kappa_s}}
 \left(\frac{\sqrt{\kappa_s}}{\eta_s}\right)^{1/2}
 \int_{\sqrt{\kappa_s}}^{\sqrt{\kappa_*}}
 x^{-23/16}\,\rd x
 \le
 \frac{C(n^{1/3}\sqrt{\kappa_s})^{1/16}}
      {n\sqrt{\kappa_s\eta_s}}.
\]
The corresponding contribution from the accumulated error
terms is bounded by
\[
 \frac C{n\sqrt{\kappa_s}}
 \int_{\sqrt{\kappa_s}}^{\sqrt{\kappa_*}}
 \frac{\log(2+n^{1/3}x)}{x^2}\,\rd x
 \le
 \frac{C\log(2+n^{1/3}\sqrt{\kappa_s})}{n\kappa_s}.
\]
The entrance error and outer-contour terms contribute at most
$C(\log n)^2/(n\sqrt{\kappa_s})$.
This proves \eqref{s4:edge-improvement} for the discrete models.

For DBM, integrate \eqref{s4:dbm-characteristic-evolution}
and use \eqref{s4:dbm-sharp-linear} together with the preceding
error and martingale bounds.
Integral Gr\"onwall gives the same two integrals.
All coefficients are evaluated strictly before $\tau$, so the
estimates include the discrete update reaching $\tau$.
If $s=b_j=\tau$, use \eqref{s4:entrance-bound} directly.

To compare \eqref{s4:edge-improvement} with the bound for
$|\delta_s(z_s)|$ in \eqref{def:stoptime edge DMC},
divide its three terms by
\(
 \frac{(n^{1/3}\sqrt{\kappa_s})^{1/8}}
      {n\sqrt{\kappa_s\eta_s}}.
\)
The resulting ratios are bounded, respectively, by
\[
 C(\log n)^2n^{-5/24},\qquad
 CM^{-1/16},\qquad
 C\log(2+M)M^{-7/8}.
\]
First choose $M$ sufficiently large and then $n\ge n_0(M)$
to make their sum less than $1/2$.
\end{proof}

We next show that $\tau=t_N$ by comparing these estimates
with the two bounds defining $\tau$.
We then obtain the support bound on $[a_n,t_N]$ and restrict
it to the time window in \Cref{large-eigen: DMC}.

\begin{proof}[Proof of \Cref{prop:s4-edge-full-time}]
Work on $A_n\cap\{\sigma_R>T_0\}$ and the common event of
\Cref{s4:outer-log,p:bddquadrac:DMC}.
Since $A_n\subset\mathcal G_n^0$, \eqref{s4:dbm-confinement}
and the preceding estimates give probability at least
\[
 1-q_n-Ce^{-cM^{1/8}}-Cn^{-20}-Ce^{-c(\log n)^2}.
\]
Here $\sigma_R=\infty$ for the discrete models.
For sufficiently large fixed $M$ and $n\ge n_0(M)$,
\Cref{s4:outer-log,lem:s4-edge-improvement} give
\[
 \sup_{z\in\cC(s)}|\delta_s(z)|
 <\frac{(\log n)^{10}}{2n},
 \qquad a_n\le s\le\tau,
\]
and
\[
 |\delta_s(z_s^j)|
 \le
 \frac{(n^{1/3}\sqrt{\kappa_s^j})^{1/8}}
      {2n\sqrt{\kappa_s^j\eta_s^j}},
 \qquad b_j\le s\le t_j\wedge\tau.
\]
These inequalities also cover paths entering at $\tau$,
by \eqref{s4:entrance-bound}.

If $\tau<t_N$, the definition of $\tau$ and
$\sigma_R>T_0$ imply that either
\(
 \sup_{z\in\cC(\tau)}|\delta_\tau(z)|
 \ge\frac{(\log n)^{10}}n,
\)
or, for some $j$ with $b_j\le\tau\le t_j$,
\(
 |\delta_\tau(z_\tau^j)|
 \ge
 \frac{(n^{1/3}\sqrt{\kappa_\tau^j})^{1/8}}
      {n\sqrt{\kappa_\tau^j\eta_\tau^j}}.
\)
For the discrete models, the preceding upper bounds include
the update reaching $\tau$.
For DBM, continuity on the finitely many active intervals
ensures that one of these inequalities holds at $\tau$.
Both alternatives contradict the preceding upper bounds.
Hence $\tau=t_N$.

The transform estimates therefore hold at every time in
$[a_n,t_N]$.
If the support first reached or crossed $E_s+f(s)$ at some
$a_n<b\le t_N$, the probe and first-crossing argument in the
proof of \Cref{l:beforestop DMC} would apply at $b$.
This contradiction proves \eqref{s4:full-support}.
\end{proof}

\begin{proof}[Proof of \Cref{large-eigen: DMC}]
For $|s-T|\le M'n^{-1/3}$, \eqref{e:deff} gives
\[
 f(s)\le(M+vM')^2n^{-2/3}\le CM^2n^{-2/3},
 \qquad M\ge\max\{1,M'\}.
\]
For sufficiently large $n$, this window is contained in
$[a_n,t_N]$.
Since $x_1(s)\le\sup I_s$, \Cref{prop:s4-edge-full-time}
gives the required upper bound simultaneously throughout
the window.
For each fixed $M$, $q_n=\oo(1)$, so increasing $n_0(M)$
absorbs $q_n+Cn^{-20}+Ce^{-c(\log n)^2}$ into $Ce^{-cM^{1/8}}$.
This proves \eqref{eq:edge-window-bound} with $\delta=1/8$.
\end{proof}

\subsubsection{Bulk estimate}\label{subsec:s4-bulk}

We prove \Cref{bulktight: DMC} for all three models.
We first estimate $\delta_s$ along the characteristic corresponding
to a fixed $w\in R(K)$, using the support bound and exterior estimates
proved above. We then obtain the estimate simultaneously for all
$w\in R(K)$ by a finite mesh and the positivity of the Stieltjes kernel.
We use the time convention from \Cref{subsec:edge-main-results}.

\paragraph{Characteristics and stopping time.}
Fix $J\Subset(s_*,T_0)$ as in \Cref{bulktight: DMC}, choose
$s_*<a<\inf J$, and use $a_n$ from \eqref{s4:working-time}.
Fix a sufficiently large constant $K$, a deterministic $t\in J$,
with $t\in\mathbb Z_n$ for the discrete models, and $w\in R(K)$.
Consider the deterministic characteristic with end point
$z_t=E_t+n^{-2/3}w$, as in \Cref{lem:bulk-characteristics}.
Recall that
\begin{equation}\label{s4:bulk-scales}
 u_s+\ri v_s=\sqrt{z_s-E_s},\qquad
 \eta_s=\Im z_s=2u_sv_s.
\end{equation}
By \eqref{def: bulk region},
\begin{equation}\label{s4:bulk-minimum}
 n^{1/3}v_t=\Im\sqrt w\ge cK,\qquad
 n^{2/3}\eta_t=\Im w\ge K^2,\qquad
 u_t,v_t\le n^{-17/60}.
\end{equation}
Indeed, if $u_t\le2v_t$, the second inequality gives the first;
if $u_t>2v_t$, the additional factor
$\sqrt{1+\max\{\Re w,0\}}$ in \eqref{def: bulk region}
gives $n^{1/3}v_t\ge cK^2$.

Trace the characteristic backwards until time $a_n$ or its first
contact with $\Re z_s-E_s=\kappa_*$.
Denote this time by $b(w)$, taking the smallest grid time not less
than it for the discrete models. We use the resulting path
$(z_s(w))_{b(w)\le s\le t}$ only on this interval.
On $[b(w),t]$, \Cref{lem:bulk-characteristics} gives
$u_s\asymp u_t+t-s$, $v_s\asymp v_t$, and
$-\eta_s'\asymp v_s$.
Since $u_s/v_s$ decreases, the portion $u_s\ge4v_s$ precedes
the portion $u_s\le4v_s$.
If $b(w)>a_n$, then
\[
 \Re z_{b(w)}-E_{b(w)}=\kappa_*+\OO(n^{-1}),
 \qquad \Im z_{b(w)}=\OO(v_t)=\oo(1).
\]
Since the right side of $\cC(b(w))$ has real part
$E_{b(w)}+\kappa_*/4$, the entrance point lies outside this
contour at distance at least $\kappa_*/2$ for sufficiently large $n$.
If $b(w)=a_n$, then \eqref{s4:uv-flow} and
\eqref{s4:bulk-minimum} give, for sufficiently large $n$,
\[
 \Re z_{a_n}(w)-E_{a_n}
 =u_{a_n}^2-v_{a_n}^2\ge c(\inf J-a)^2>0.
\]
Enlarge $\mathcal K_a$ to contain fixed neighborhoods of all these
starting points, uniformly for $t\in J$ and $w\in R(K)$, and of
the edge starting points with $T\in J$.
Choose $\varepsilon$ also smaller than $v^2(\inf J-a)^2/8$.
By \Cref{prop:s4-positive-time-preparation}, choose $A_n$ as in
\eqref{s4:initial-good-event} for this set, with $q_n=\oo(1)$.
These choices are independent of $t$, $w$ and $K$.

Use the control function $f$ and stopping time $\tau$ from the
preceding subsection with $T=t$ and $M=K^{1/2}$.
Let $\mathcal G\subset A_n$ be the common event constructed in
the proofs of \Cref{prop:s4-edge-full-time,s4:outer-log}.
Then $\tau=t_N>t$ on $\mathcal G$, the support bound
\eqref{s4:full-support} holds, and the contour errors are bounded
by $C(\log n)^2/n$. Moreover,
\begin{equation}\label{s4:bulk-common-event}
 \mathbb P(\mathcal G^c)
 \le q_n+Ce^{-cK^{1/16}}+Cn^{-20}+Ce^{-c(\log n)^2}.
\end{equation}
The choices of constants in these results can be made uniformly
for $t\in J$.

When $u_s\ge4v_s$, the preceding geometric estimates give
$u_s\ge cKn^{-1/3}$ and $u_s\ge c(t-s)$.
Taking $K$ sufficiently large and the fixed constant $v$ in
\eqref{e:deff} sufficiently small therefore gives
\[
 \sqrt{f(s)}=K^{1/2}n^{-1/3}+v(t-s)\le u_s/4.
\]
Since $\Re z_s-E_s=u_s^2-v_s^2$, \Cref{l:beforestop DMC}
implies, for $b(w)\le s<t\wedge\tau$,
\begin{equation}\label{s4:bulk-distance}
 \operatorname{dist}(z_s,I_s)\ge
 \begin{cases}
 cu_s^2,&u_s\ge4v_s,\\
 \eta_s,&u_s\le4v_s.
 \end{cases}
\end{equation}
Also, $\eta_s\ge\eta_t\ge K^2n^{-2/3}\gg n^{-1+\delta}$
for the fixed $0<\delta<1/3$ chosen above.

\begin{figure}[htbp]
 \centering
 \includegraphics[width=0.70\linewidth]{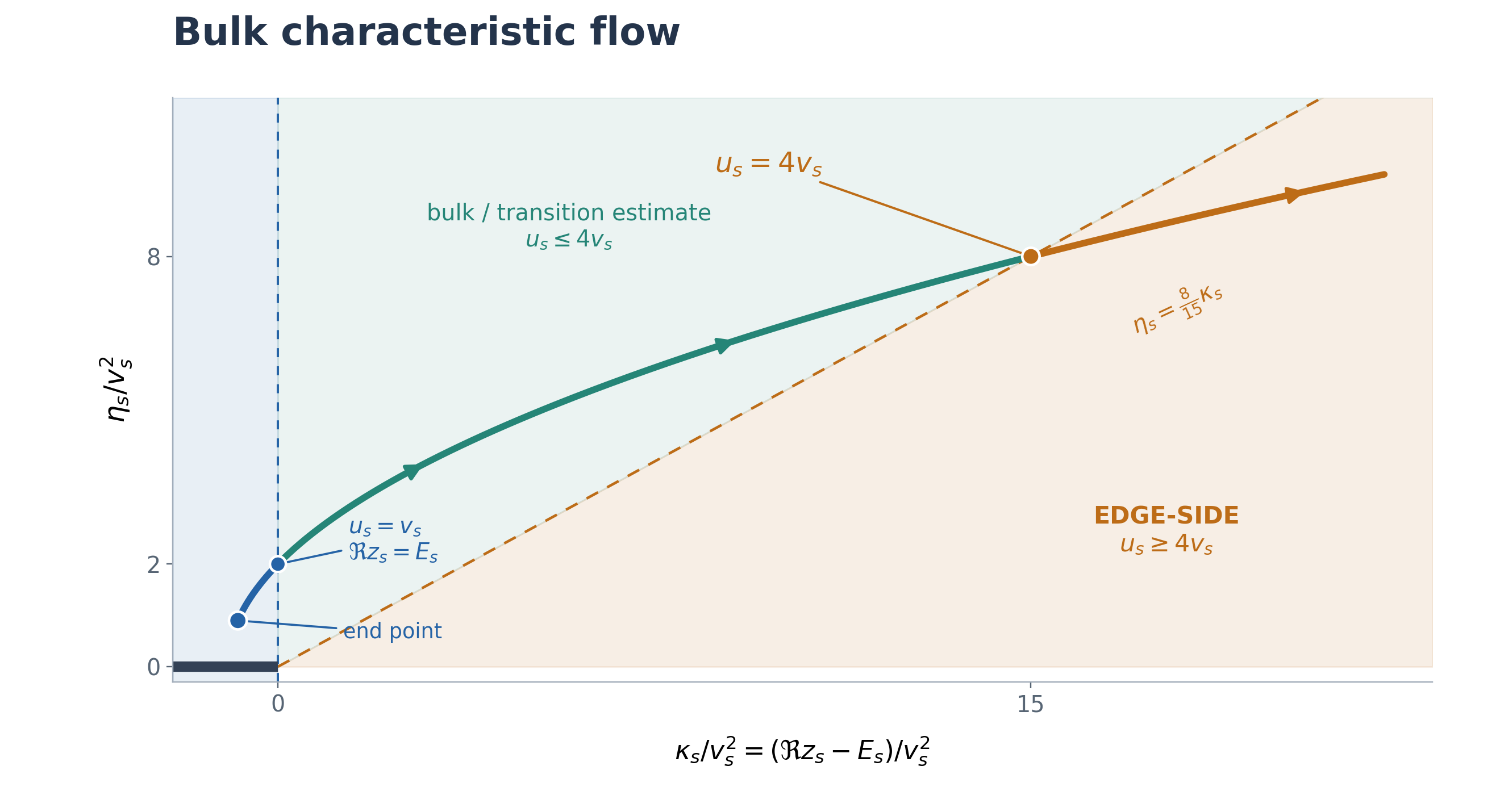}
 \caption{A characteristic traced backwards from the end point
$z_t=E_t+n^{-2/3}w$, with $\Re w<0$.
Here $u_s+\ri v_s=\sqrt{z_s-E_s}$, and the plotted coordinates
$(\kappa_s/v_s^2,\eta_s/v_s^2)$ are normalized at each time $s$.
Arrows indicate decreasing $s$, along which $u_s/v_s$ increases.
The vertical dashed line is $\Re z_s=E_s$, or $u_s=v_s$.
The oblique dashed line is $u_s=4v_s$, separating the two
regions in \eqref{s4:bulk-distance}.
The path continues backwards beyond the displayed portion
to time $b(w)$.}
 \label{fig:bulk-characteristic-flow}
\end{figure}

For this fixed $w$, define
\[
 \tau_w=t\wedge\tau\wedge\inf\left\{
 s\in[b(w),t]:
 |\delta_s(z_s(w))|\ge
 \frac{(n^{1/3}(u_s+v_s))^{1/4}}{nu_s(u_s+v_s)}
 \right\},
\]
where the infimum is over grid times for the discrete models
and $\inf\varnothing=\infty$.
Using the stopping convention in \eqref{s4:edge-stopped-martingale}, set
\begin{equation}\label{s4:bulk-stop}
 N_s(w)=
 \begin{cases}
 \displaystyle\mathbf1_{A_n}\sum_{r\in[b(w),s)\cap\mathbb Z_n}
 \mathbf1_{\{r<\tau_w\}}\xi_r(z_{r+1/n}(w)),
 &\text{DMC and corners},\\[6pt]
 \displaystyle\mathbf1_{A_n}\int_{b(w)}^s
 \mathbf1_{\{r\le\tau_w\}}\,\rd \mathcal M_r^n(z_r(w)),
 &\text{DBM},
 \end{cases}
 \qquad b(w)\le s\le t.
\end{equation}
These are martingales for each deterministic $w$.
We prove their estimates with this stopping and intersect with
$\mathcal G$ only afterwards.

\paragraph{Estimates on a fixed characteristic.}
We estimate the accumulated error terms and the stopped martingale
along $z_s(w)$, then apply Gr\"onwall as in the edge argument.
The following lemma records the resulting bounds and the
pointwise probability estimate needed for the spectral mesh.

\begin{proposition}\label{lem:s4-bulk-characteristic}
For sufficiently large fixed $K$ and $n\ge n_0(K)$, the following
holds for each of the three models and every deterministic $w\in R(K)$.
Outside an event of probability at most
\[
  C\bigl(1+\log(2+|\sqrt w|)\bigr)
          e^{-c|\sqrt w|^{1/8}}+Cn^{-22},
\]
we have on $\mathcal G$ that $\tau_w=t$, and,
for $b(w)\le s\le t$,
\begin{align}
 |\delta_s(z_s)|
 &\le C\left[
 \frac{(n^{1/3}u_s)^{1/16}}{nu_s^2}
 +\frac{(\log n)^2}{nu_s}\right],
 &&u_s\ge4v_s,\label{s4:bulk-exterior-final}\\
 |\delta_s(z_s)|
 &\le\frac C{n\eta_s}
 \left[
 (n^{1/3}v_s)^{1/16}\bigl(1+\log(2+n^{1/3}v_s)\bigr)
 +v_s(\log n)^2\right],
 &&u_s\le4v_s.\label{s4:bulk-interior-final}
\end{align}
In particular,
\begin{equation}\label{s4:bulk-point-prob}
 \mathbb P\left(\mathcal G\cap
 \left\{|\delta_t(E_t+n^{-2/3}w)|>
 \frac{C(\Im\sqrt w)^{1/4}}{n^{1/3}\Im w}\right\}\right)
 \le C\bigl(1+\log(2+|\sqrt w|)\bigr)
          e^{-c|\sqrt w|^{1/8}}+Cn^{-22}.
\end{equation}
The constants are independent of $t\in J$, $w$, $K$ and $n$.
\end{proposition}

\begin{proof}
Fix $w$ and omit it from the path notation.
For $s<\tau_w$, its definition gives
\[
 \frac{|\delta_s(z_s)|}{v_s}\le C
 \begin{cases}
 (n^{1/3}u_s)^{-7/4}/(n^{1/3}v_s),&u_s\ge4v_s,\\
 (n^{1/3}v_s)^{-3/4}/(n^{2/3}\eta_s),&u_s\le4v_s.
 \end{cases}
\]
Both bounds are at most $CK^{-11/4}$, by
\eqref{s4:bulk-minimum} and \eqref{s4:uv-flow}.
In particular, $m_s^n(z_s)-m_s(z_s)=\oo(1)$.
On the outer contour, the definition of $\tau$ and
\Cref{lem:exterior-cauchy} also give an $\oo(1)$ difference.
Taylor expansion of the flux, as in the proof of
\Cref{lem:s4-edge-errors}, therefore gives
\[
 P_s^n(z_s)\le C/u_s\asymp v_s/\eta_s.
\]
For DMC, the Stieltjes integral bounds in
\Cref{lem:stieltjes-bounds} and the macroscopic
characteristic equation give, for $s<\tau_w$,
\[
 \begin{aligned}
 \sup_{s\le r\le s+n^{-1}}
 |m_s^n(z_r)-m_r(z_r)|
 &\le |\delta_s(z_s)|+\frac C{n\eta_s}+\frac Cn\\
 &\le CK^{-11/4}v_s+CK^{-2}n^{-1/3}+\frac Cn
 =\oo(1).
 \end{aligned}
\]
Here we use $v_s\le Cn^{-17/60}$,
$\eta_s\ge K^2n^{-2/3}$, and
$\frac{\rd}{\rd r}m_r(z_r)=F_r(z_r)$
with $F_r(z_r)$ uniformly bounded.
Thus the additional DMC condition in
\Cref{lem:s4-characteristic-evolution} is satisfied.
The one-step comparison in the proof of that lemma and
\eqref{s4:bulk-distance} justify applying the evolution and
moment estimates throughout each one-step segment.

We first bound the accumulated nonmartingale terms in
$\mathcal E_2$ and the additional $CP_s^n(z_s)/n^2$ term
in \eqref{s4:sharp-linear}.
For the discrete models, \eqref{s4:moment-input} and
\eqref{s4:step-remainder} give
\begin{equation}\label{s4:bulk-exterior-errors}
 \langle\Re N(w)\rangle_s+\langle\Im N(w)\rangle_s
 \le C
 \begin{cases}
 (n^2u_s^4)^{-1},&u_s\ge4v_s,\\
 (n^2\eta_s^2)^{-1},&u_s\le4v_s,
 \end{cases}
 \qquad b(w)\le s\le t,
\end{equation}
and
\begin{equation}\label{s4:bulk-interior-errors}
 \sum_{r\in[b(w),s)\cap\mathbb Z_n}
 \left(
 |\mathcal E_2(r)-\xi_r(z_{r+1/n})|
 +\frac{P_r^n(z_r)}{n^2}
 \right)
 \le C
 \begin{cases}
 (nu_s^2)^{-1},&u_s\ge4v_s,\\
 (n\eta_s)^{-1},&u_s\le4v_s,
 \end{cases}
 \qquad s\le t\wedge\tau_w.
\end{equation}
Indeed, on $u_r\ge4v_r$, the two sums are bounded by
\[
 \frac C{n^2}\int_{u_s}^{\OO(1)}x^{-5}\,\rd x,
 \qquad
 \frac Cn\int_{u_s}^{\OO(1)}x^{-3}\,\rd x.
\]
On $u_r\le4v_r$, use $-\eta_r'\asymp v_r$ to bound the
new contributions by
\[
 \frac C{n^2}\int_{\eta_s}^{\infty}y^{-3}\,\rd y,
 \qquad
 \frac Cn\int_{\eta_s}^{\infty}y^{-2}\,\rd y.
\]
At the transition $u_r=4v_r$, $\eta_r\asymp u_r^2$,
so these bounds include the preceding portion of the path.
The discrete sums obey the same bounds because the variables
change by a bounded factor over one step.
Only increments before $\tau_w$ contribute to the bracket,
so \eqref{s4:bulk-exterior-errors} also holds after $\tau_w$.
For DBM, \eqref{s4:dbm-martingale-bound} gives the same bracket bound,
and \eqref{s4:dbm-second-order} gives \eqref{s4:bulk-interior-errors}
with its left-hand side replaced by
\[
 \frac{|1-\theta|}{2\theta n}
 \int_{b(w)}^s|\partial_z^2m_r^n(z_r)|\,\rd r.
\]

We next estimate $N_s(w)$.
For the discrete models, \eqref{s4:moment-input} gives,
for each fixed integer $p\ge1$ and $b(w)\le s<t$,
\begin{equation}\label{s4:bulk-high-moment}
 \mathbb E_s\left|
 \mathbf1_{\{s<\tau_w\}}\xi_s(z_{s+1/n})\right|^{2p}
 \le
 \frac{C_p n^{-p/2}}{[nu_s(u_s+v_s)]^{2p}}.
\end{equation}
To check the two terms, multiply the moment bound by
$[nu_s(u_s+v_s)]^{2p}$.
For $u_s\ge4v_s$, they are bounded by
$(nu_s)^{-p}$ and $(nu_s)^{-p}(nu_s^3)^{1-p}$,
where $nu_s^3\ge cK^3$.
For $u_s\le4v_s$, they are bounded by
$(v_s/(n\eta_s))^p$ and $v_s(n\eta_s)^{1-2p}$.
Using $v_s\le Cn^{-17/60}$ and $n\eta_s\ge K^2n^{1/3}$
proves \eqref{s4:bulk-high-moment}.

Truncate the stopped increments at
$n^{-1/16}/[nu_s(u_s+v_s)]$ and subtract their conditional means.
By the calculation in the proof of \Cref{p:bddquadrac:DMC},
with $p=64$, truncation fails with probability at most $Cn^{-24}$
per update, and the accumulated centering error is negligible.
On the portion $u_s\ge4v_s$, divide the values of $n^{1/3}u_s$
into intervals $[A,2A]$.
The bracket bound and Freedman's maximal inequality, applied
from $b(w)$ to the last time in each interval, give
\begin{equation}\label{s4:bulk-martingale}
 |N_s(w)|\le C
 \begin{cases}
 (n^{1/3}u_s)^{1/16}/(nu_s^2),&u_s\ge4v_s,\\
 (n^{1/3}v_s)^{1/16}/(n\eta_s),&u_s\le4v_s.
 \end{cases}
\end{equation}
For the second bound, divide the decreasing heights $\eta_s$
into intervals in which they vary by at most a factor of two,
and use $v_s\asymp v_t$.
On the first portion, use $C/(n^2u_s^4)$ for the variance
and $Cn^{-1/16}/(nu_s^2)$ for the increment bound; on the
second, use $C/(n^2\eta_s^2)$ and $Cn^{-1/16}/(n\eta_s)$.
Freedman's inequality therefore gives failure probability
$Ce^{-cA^{1/8}}$, with
$A=n^{1/3}u_s$ on the first portion and $A\asymp n^{1/3}v_t$
on the second, exactly as in the proof of
\Cref{p:bddquadrac:DMC}.
All partial sums start at $b(w)$, including on the second portion.

The first collection of intervals sums geometrically.
For the second, \eqref{s4:bulk-minimum} gives
\begin{equation}\label{s4:bulk-logarithm}
 \log\left(2+\frac{v_t^2}{\eta_t}\right)
 \le C\log(2+n^{1/3}v_t),
\end{equation}
which bounds the number of height intervals.
If $u_t\ge4v_t$, the smallest first-portion value of
$n^{1/3}u_s$ is comparable to $|\sqrt w|$.
Otherwise, $n^{1/3}v_t\asymp|\sqrt w|$, and the first
portion ends at $u_s\asymp v_t$.
Consequently, \eqref{s4:bulk-martingale} holds simultaneously for
$b(w)\le s\le t$ outside an event of probability
\[
 C\bigl(1+\log(2+|\sqrt w|)\bigr)e^{-c|\sqrt w|^{1/8}}
 +Cn^{-22}.
\]
For DBM, the continuous exponential martingale inequality
and \eqref{s4:bulk-exterior-errors} give the same conclusion
without truncation.

If $b(w)>a_n$, the entrance point lies outside $\cC(b(w))$
at a fixed positive distance. On $\mathcal G$, the contour
encloses both supports, so \eqref{e:cauchy int for stij} and
the contour bound give
\[
 |\delta_{b(w)}(z_{b(w)})|\le C(\log n)^2/n.
\]
If $b(w)=a_n$, then $z_{a_n}(w)\in\mathcal K_a$, and
\eqref{s4:initial-good-event} gives
$|\delta_{a_n}(z_{a_n}(w))|\le(\log n)/n$ on $\mathcal G$.
Work on $\mathcal G$ outside the local martingale exceptional
event just estimated.
We apply Gr\"onwall to the evolution equation up to $\tau_w$,
using \eqref{s4:sharp-linear} and its DBM counterpart
\eqref{s4:dbm-sharp-linear}.
As in the proof of \Cref{lem:s4-edge-improvement},
\[
 \int_r^s\frac{-\Im\chi_q(z_q)}{\eta_q}\,\rd q
 =\log\frac{\eta_r}{\eta_s}.
\]
The integral of the additional coefficient
$|\delta_q(z_q)|/\eta_q$ is bounded by
$C(n^{1/3}u_s)^{-7/4}/(n^{1/3}v_s)$ on the first portion,
and by $C(n^{1/3}v_s)^{-3/4}/K^2$ on the second,
including the preceding contribution.
These follow by integrating the bound defining $\tau_w$,
using $-u_q'\asymp1$ and $-\eta_q'\asymp v_q$.
The propagation factor is therefore at most $C\eta_r/\eta_s$.
For the discrete models, the corresponding product has the
same bound: the imaginary-velocity derivative calculation in
the proof of \Cref{c:kappa:DMC}, now using
\eqref{s4:im-edge}, bounds the Riemann-sum error by
$C/(n\min_{[b(w),t]}u_s)=\oo(1)$.
Here $u_s=\eta_s/(2v_s)\ge cK^2n^{-23/60}$.
The deterministic bounds used for the nonlinear coefficient vary
by a bounded factor over one step, so they also bound its discrete sums.

Sum \eqref{sde for DMC 2} before taking absolute values of the
martingale contribution, so that it is $N_s(w)$.
Use \eqref{s4:bulk-martingale}, \eqref{s4:bulk-interior-errors}, and
the entrance and contour estimates, and then apply Gr\"onwall.
On $u_s\ge4v_s$, the martingale contribution is bounded by
\[
 \frac C{nu_s}\int_{u_s}^{\OO(1)}
 \frac{(n^{1/3}x)^{1/16}}{x^2}\,\rd x
 \le\frac{C(n^{1/3}u_s)^{1/16}}{nu_s^2}.
\]
On $u_s\le4v_s$, the corresponding bound is
\[
 \frac{C(n^{1/3}v_s)^{1/16}}{n\eta_s}
 \left(1+\int_{\eta_s}^{C v_s^2}\frac{\rd y}{y}\right)
 \le\frac{C(n^{1/3}v_s)^{1/16}
                 (1+\log(2+n^{1/3}v_s))}{n\eta_s}.
\]
The constant upper limit is chosen to exceed the height at the
transition. The accumulated nonmartingale errors obey the same
bounds without the factors $(n^{1/3}u_s)^{1/16}$ and
$(n^{1/3}v_s)^{1/16}$.
The entrance error and contour terms contribute at most
$C(\log n)^2/(nu_s)$, which is comparable to
$Cv_s(\log n)^2/(n\eta_s)$ on the second portion.
This proves \eqref{s4:bulk-exterior-final} and
\eqref{s4:bulk-interior-final} for $s\le\tau_w$.
For DBM, integrate \eqref{s4:dbm-characteristic-evolution}
and use integral Gr\"onwall; the same two integrals occur.

Dividing these estimates by the bound defining $\tau_w$,
the martingale terms have ratios at most
\[
 C(n^{1/3}u_s)^{-3/16},\qquad
 C(n^{1/3}v_s)^{-3/16}\log(2+n^{1/3}v_s),
\]
and the entrance and contour terms have ratios at most
\[
 C(\log n)^2n^{-1/12},\qquad
 C(\log n)^2n^{-71/240},
\]
respectively.
They are less than $1/2$ in total when $K$ is sufficiently
large and then $n\ge n_0(K)$.
These estimates include the first stopping update; for DBM,
they extend to $\tau_w$ by continuity.
Since $\tau>t$ on $\mathcal G$, the definition of $\tau_w$
therefore gives $\tau_w=t$.

At $s=t$ with $u_t\ge4v_t$, multiplying
\eqref{s4:bulk-exterior-final} by $n\eta_t$ gives
\[
 n\eta_t|\delta_t(z_t)|
 \le C\left[
 (n^{1/3}v_t)(n^{1/3}u_t)^{-15/16}
 +v_t(\log n)^2\right]
 \le C(n^{1/3}v_t)^{1/4}.
\]
For $u_t\le4v_t$, the same conclusion follows from
\eqref{s4:bulk-interior-final}.
Since $n^{1/3}v_t=\Im\sqrt w$ and
$n\eta_t=n^{1/3}\Im w$, this proves
\eqref{s4:bulk-point-prob}.
\end{proof}

\paragraph{All spectral parameters.}
We now pass from the pointwise estimate to a common event for
all $w\in R(K)$, using the mesh discretizations.

\begin{proof}[Proof of \Cref{bulktight: DMC}]
For $A_\ell=2^\ell K$, $\ell\ge0$, partition the complex plane
into squares of side $A_\ell^{-1}$.
From each square meeting
$R(K)\cap\{A_\ell\le|\sqrt w|\le2A_\ell\}$,
choose one point of this intersection.
The resulting finite set $\mathcal N_\ell$ satisfies
\begin{equation}\label{s4:bulk-net}
 \#\mathcal N_\ell\le CA_\ell^6,\qquad
 \min_{w_0\in\mathcal N_\ell}|w-w_0|\le CA_\ell^{-1}
 \quad\text{if }w\in R(K),\ A_\ell\le|\sqrt w|\le2A_\ell.
\end{equation}
Since $K\le|\sqrt w|\le n^{1/20}$ on $R(K)$,
the corresponding annuli cover $R(K)$, and the meshes contain
$\OO(n^{3/10})$ points in total.
Apply \eqref{s4:bulk-point-prob} at these deterministic points.
Summing the probabilities gives
\[
 \begin{aligned}
 &\mathbb P\bigl(\mathcal G\cap
 \{\text{the pointwise bound fails at some mesh point}\}\bigr)\\
 &\qquad\le
 \sum_{\ell\ge0}CA_\ell^6
       (1+\log(2+A_\ell))e^{-cA_\ell^{1/8}}+Cn^{-20}
 \le Ce^{-c'K^{1/8}}+Cn^{-20}.
 \end{aligned}
\]

It remains to extend the estimate within each square.
Let $w_0\in\mathcal N_\ell$ and $w\in R(K)$ lie in the same
square and satisfy $A_\ell\le|\sqrt w|\le2A_\ell$.
Then $\Im w\asymp\Im w_0$, $\Im\sqrt w\asymp\Im\sqrt{w_0}$,
and $|w-w_0|\le\frac12\Im w_0$ for large $K$.
At a mesh point satisfying the pointwise bound,
\eqref{s4:kernel-split} and \eqref{s4:im-edge} give
\[
 \int\frac{(\mu_t^n+\mu_t)(\rd r)}
 {|r-E_t-n^{-2/3}w_0|^2}
 \le\frac{Cn^{1/3}\Im\sqrt{w_0}}{\Im w_0}.
\]
The ordinary kernels at $w$ and $w_0$ differ by their spatial
separation divided by the product of their denominators.
By the previous distance comparison, positivity of the measures,
and the bounded derivative of the analytic term in
\eqref{s4:kernel-split}, we obtain
\begin{equation}\label{s4:bulk-net-extension}
 \begin{aligned}
 &|\delta_t(E_t+n^{-2/3}w)-\delta_t(E_t+n^{-2/3}w_0)|\\
 &\qquad\le
 \frac{Cn^{-1/3}\Im\sqrt{w_0}}{\Im w_0}|w-w_0|
 +Cn^{-2/3}|w-w_0|
 \le\frac{Cn^{-1/3}(\Im\sqrt{w_0})^{1/4}}{\Im w_0}.
 \end{aligned}
\end{equation}
For the first term, use
$(\Im\sqrt{w_0})^{3/4}/A_\ell\le C$.
The ratio of the second term to the final bound is at most
\[
 Cn^{-1/3}\frac{\Re\sqrt{w_0}}{A_\ell}
                  (\Im\sqrt{w_0})^{3/4}
 \le Cn^{-71/240}.
\]
Thus the pointwise estimate extends to every $w\in R(K)$.

Adding \eqref{s4:bulk-common-event} once, the total failure
probability is at most
$q_n+Ce^{-cK^{1/16}}+Cn^{-20}+Ce^{-c(\log n)^2}$.
For each fixed $K$, increase $n_0(K)$ to absorb $q_n$ and the last
two terms.
This proves \Cref{bulktight: DMC}, for example with
$\alpha=1/32$ and $\alpha'=1/4$.
All constants are uniform for the deterministic time $t\in J$;
the resulting event may depend on $t$.
\end{proof}

\subsubsection{Uniform bounds on a time window}\label{subsec:s4-window}
We first improve the exterior bound in \Cref{s4:outer-log}
from $C(\log n)^2/n$ to $C/n$, with arbitrarily high probability,
by estimating the martingale simultaneously along all contour
characteristics.

\begin{lemma}\label[lemma]{s4:outer-tight}
Under \Cref{ass:weightfunc,ass:prelim-log-branches,ass:edge-data},
fix $a\in(s_*,T_0)$ and use \eqref{s4:working-time}.
For each of the three models, let $\cC(s)$ and $\omega_+(s)$
be the transported contours from \Cref{subsec:s4-characteristics}.
For every $\varepsilon\in(0,1)$, there exist constants
$C<\infty$ and $n_0$ such that, for all $n\ge n_0$,
\begin{equation}\label{s4:outer-Op}
 \mathbb P\left(
 \sup_{a_n\le s\le T_0}
 \left[
 \sup_{z\in\cC(s)}|\delta_s(z)|
 +\sup_{z\in\omega_+(s)}|\delta_s(z)|
 \right]
 \le \frac Cn
 \right)\ge1-\varepsilon.
\end{equation}
The constants may depend on $\varepsilon,a,T_0$ and the fixed
model, spectral domains and assumption constants, but not on $n$.
The same conclusion holds with the contours replaced by
sufficiently small fixed neighborhoods, chosen independently
of $\varepsilon$ and $n$.
\end{lemma}

\begin{proof}
Fix $T=(a+T_0)/2$ and a sufficiently large constant $M$,
and use the stopping time $\tau$ from \Cref{def:s4-edge-stopping}.
By \Cref{prop:s4-positive-time-preparation}, choose $1\le B<\infty$
such that, for sufficiently large $n$, the events
\[
 A_{n,B}=A_n\cap\left\{
 \sup_{z\in\mathcal K_a}|\delta_{a_n}(z)|\le B/n\right\}
 \in\mathscr F_{a_n}
\]
satisfy $\mathbb P(A_{n,B}^c)\le\varepsilon/4$.
The set $\mathcal K_a$, the events $A_n$ and $B$ can be chosen
independently of $M$.
The localized quantities below are set to zero on $A_{n,B}^c$.
On $A_n$, the support bounds preceding \Cref{s4:outer-log} imply that,
up to $\tau$, the contours $\cC(s)$,
$\widetilde{\cC}_k(s)$ and $\omega_+(s)$, together with
sufficiently small fixed neighborhoods of them, stay a fixed
positive distance from both supports.
The neighborhoods and separation constants can be chosen
independently of $M$, for all sufficiently large $n$.
Use the partition $a_n=u_0<\cdots<u_J=t_N$ from
\Cref{subsec:s4-characteristics}.

For $w\in\cC(u_k)$ and $u_k\le t\le u_{k+1}$, set
\[
 N_t^k(w)=
 \begin{cases}
 \displaystyle\mathbf1_{A_{n,B}}\sum_{s\in[u_k,t)\cap\mathbb Z_n}
 \mathbf1_{\{s<\tau\}}
 \xi_s\bigl(\Phi_{u_k,s+1/n}(w)\bigr),
 &\text{DMC and corners},\\[6pt]
 \displaystyle\mathbf1_{A_{n,B}}\int_{u_k}^{t}\mathbf1_{\{s\le\tau\}}
 \,\rd\mathcal M_s^n\bigl(\Phi_{u_k,s}(w)\bigr),
 &\text{DBM}.
 \end{cases}
\]
Both $N_t^k(w)$
and $\partial_wN_t^k(w)$ are martingales.
For the discrete models, \eqref{s4:moment-input}, Cauchy's
formula on the fixed neighborhoods, and the bounded spatial
derivative of $\Phi_{u_k,s}$ give
\[
 \mathbb E_s\left|
 \partial_w^q\left[
 \mathbf1_{A_{n,B}}\mathbf1_{\{s<\tau\}}
 \xi_s\bigl(\Phi_{u_k,s+1/n}(w)\bigr)
 \right]\right|^2
 \le Cn^{-3},\qquad q=0,1.
\]
For $s<\tau$, we have
$I_{s+1/n}\subset I_s+[0,C/n]$ for DMC and
$I_{s+1/n}\subset I_s$ for corners.
Thus, on the fixed neighborhoods, for sufficiently large $n$,
\(
 \dist(z,I_{s+1/n})
 \ge \dist(z,I_s)-\frac Cn
 \ge \frac12\dist(z,I_s).
\)
For DBM, \eqref{eq:pote_mt} gives quadratic-variation rates
bounded by $C/n^2$ for both $N_t^k(w)$ and
$\partial_wN_t^k(w)$.
The one-dimensional Sobolev inequality on $\cC(u_k)$
and Doob's maximal inequality therefore give
\begin{equation}\label{s4:outer-Sobolev}
 \begin{aligned}
 &\mathbb E\sup_{u_k\le t\le u_{k+1}}
       \sup_{w\in\cC(u_k)}|N_t^k(w)|^2 \le C\oint_{\cC(u_k)}
 \mathbb E\sup_{u_k\le t\le u_{k+1}}
 \left(|N_t^k(w)|^2+|\partial_wN_t^k(w)|^2\right)
 |\rd w|
 \le \frac C{n^2}.
 \end{aligned}
\end{equation}
The constants are uniform in $k$.

On $A_{n,B}$, for $u_k\le\tau$, use \eqref{s4:contour-transfer} and
\eqref{s4:outer-cauchy} in the exterior estimates of
\Cref{lem:s4-characteristic-evolution,lem:s4-dbm-characteristic-evolution}.
As in the derivation of \eqref{s4:outer-gronwall},
Gronwall's inequality gives
\[
 \begin{aligned}
 &\sup_{u_k\le t\le u_{k+1}\wedge\tau}
       \sup_{z\in\widetilde{\cC}_k(t)}|\delta_t(z)| \le C\left(
       \sup_{z\in\cC(u_k)}|\delta_{u_k}(z)|
       +\sup_{u_k\le t\le u_{k+1}}
        \sup_{w\in\cC(u_k)}|N_t^k(w)|
       +\frac1n
       \right).
 \end{aligned}
\]
On $A_{n,B}$, the initial contour is contained in $\mathcal K_a$,
so
\[
 \sup_{z\in\cC(a_n)}|\delta_{a_n}(z)|\le\frac Bn.
\]
Iterating over the $J$ intervals, using
\eqref{s4:contour-transfer} and \eqref{s4:outer-Sobolev},
we obtain
\[
 \begin{aligned}
 &\left\|\mathbf1_{A_{n,B}}\sup_{a_n\le s\le\tau}
       \sup_{z\in\cC(s)}|\delta_s(z)|\right\|_2 \le \frac{CB}n+
 C\sum_{k=0}^{J-1}
 \left\|\sup_{u_k\le t\le u_{k+1}}
       \sup_{w\in\cC(u_k)}|N_t^k(w)|\right\|_2
 \le\frac{C(B+1)}n.
 \end{aligned}
\]

Move the sides of $\cC(s)$ inward by a sufficiently small
fixed amount. The support separation and
\eqref{s4:contour-inward} remain valid.
After shortening the time intervals if necessary, we can repeat
the argument with the same stopping time $\tau$, using the
exterior estimates in
\Cref{lem:s4-characteristic-evolution,lem:s4-dbm-characteristic-evolution},
which do not require smallness of $\delta_s$.
The perturbed initial contours are contained in $\mathcal K_a$.
Cauchy's formula then gives the same truncated $L^2$ bound on fixed
neighborhoods of $\cC(s)$ and $\omega_+(s)$.

By this truncated $L^2$ bound, \Cref{prop:s4-edge-full-time},
$q_n=\oo(1)$ and Markov's inequality, for every $A>0$,
\[
 \limsup_{n\to\infty}
 \mathbb P\left(
 n\sup_{a_n\le s\le T_0}
 \left[
 \sup_{z\in\cC(s)}|\delta_s(z)|
 +\sup_{z\in\omega_+(s)}|\delta_s(z)|
 \right]>A\right)
 \le\frac{\varepsilon}{4}+Ce^{-cM^{1/8}}+\frac{C(B+1)^2}{A^2},
\]
where $C$ is independent of $M$ and $B$.
The same estimate holds on the fixed neighborhoods.
With $B$ fixed, choose $M$ and then $A$ so that the last two terms
are at most $\varepsilon/4$ each. Taking $n$ sufficiently large
proves both assertions.
\end{proof}

\begin{proof}[Proof of \Cref{l: Delta is sml}]
Fix $s_*<a<T$ and use $a_n$ from \eqref{s4:working-time}.
For sufficiently large $n$, the time window is contained in
$[a_n,T_0]$.
The exterior estimate follows from \Cref{s4:outer-tight}.
For the estimate on $\mathcal K$, we first control $\delta_s$
along one characteristic throughout the time window, and then
compare the Stieltjes kernels at the remaining spectral points.

Fix a sufficiently large constant $K$, depending on
$\varepsilon$ and $M'$.
For the discrete models, let
\[
 t=\frac1n\left\lfloor n\bigl(T+M'n^{-1/3}\bigr)\right\rfloor;
\]
for DBM, let $t=T+M'n^{-1/3}$.
Thus $t$ is the last time in the window, with the time convention
of \Cref{subsec:edge-main-results}.
Let $z_s$ be the local characteristic from
\Cref{lem:bulk-characteristics}, traced backwards to $a_n$ or
its first contact with $\Re z_s-E_s=\kappa_*$, with end point
\[
 z_t=E_t+\ri K^2n^{-2/3}.
\]
For sufficiently large $n$, $\ri K^2\in R(K)$ and $t$ belongs
to a fixed compact subinterval of $(a,T_0)$.

At the end point, $u_t=v_t=Kn^{-1/3}/\sqrt2$.
Taking $K$ sufficiently large compared with $M'$, the bounds
for $u_s'$ and $v_s'$ in \Cref{lem:bulk-characteristics}
give, on the local part of the path within the window,
\[
 u_s\le4v_s,\qquad
 |z_s-E_s|\le CK^2n^{-2/3},\qquad
 \eta_s\ge K^2n^{-2/3}.
\]
Since $CK^2n^{-2/3}<\kappa_*$ for sufficiently large $n$,
the backward path cannot reach $\kappa_s=\kappa_*$ within
this window. Hence the entire window lies in the local
part of the characteristic.

Apply \Cref{lem:s4-bulk-characteristic} and
\eqref{s4:bulk-common-event} with $w=\ri K^2$.
By choosing $K$ and then $n$ sufficiently large, with
probability at least $1-\varepsilon/2$,
\(
 \sup_{|s-T|\le M'n^{-1/3}}
 |\delta_s(z_s)|\le Cn^{-1/3}.
\)
On this event, \eqref{s4:im-edge} and
\eqref{s4:kernel-split} give, throughout the window,
\[
 \int\frac{(\mu_s^n+\mu_s)(\rd r)}{|r-z_s|^2}
 =
 \frac{\Im\widehat m_s^n(z_s)+\Im\widehat m_s(z_s)}
      {\eta_s}
 \le Cn^{1/3}.
\]

For $w\in\mathcal K$, put $z=E_s+n^{-2/3}w$.
The preceding bounds and compactness of $\mathcal K$ give
\[
 |z-z_s|\le Cn^{-2/3},\qquad
 |r-z_s|\le C|r-z|,\qquad r\in\mathbb R.
\]
Using the identity
\[
 \frac1{r-z}-\frac1{r-z_s}
 =\frac{z-z_s}{(r-z)(r-z_s)}
\]
and \eqref{s4:kernel-split}, we obtain
\begin{equation}\label{s4:compact-transfer}
 \begin{aligned}
 |\delta_s(z)-\delta_s(z_s)|
 &\le C|z-z_s|
 \left(
 \int\frac{(\mu_s^n+\mu_s)(\rd r)}{|r-z_s|^2}+1
 \right)\le Cn^{-1/3}.
 \end{aligned}
\end{equation}
The final $1$ accounts for the bounded spatial derivative
of $g_{\mathsf b}$ in \eqref{s4:kernel-split}.
All constants may depend on the fixed $K,M'$ and
$\mathcal K$, but not on $n$.

Combining \eqref{s4:compact-transfer} with the estimate
at $z_s$ proves \eqref{eq:edge-compact-window}
simultaneously in $s$ and $w$.
Intersecting this event with the exterior event from
\Cref{s4:outer-tight}, chosen with probability at least
$1-\varepsilon/2$, proves both assertions with probability
at least $1-\varepsilon$.
\end{proof}

\section{Convergence to \texorpdfstring{$\mathrm{ALE}_{\beta}$}{ALE-beta}}
\label{sec:ale-convergence}

In this section, we prove convergence to $\mathrm{ALE}_\beta$
for the three particle systems introduced in \Cref{subsec:models}.
We use the characterization in
\Cref{thm:prelim-airychar} done by \cite{charofairy}.
We first define the edge scaling and state the convergence theorem.
We then derive the evolution of the rescaled Stieltjes transform,
identifying its drift and martingale covariations.
Combining these estimates with the bounds from \Cref{s:eost}, we prove
joint tightness of the rescaled Stieltjes transforms, the associated
martingale processes, and the particle trajectories.
Finally, we verify the hypotheses of \Cref{thm:prelim-airychar}
for every subsequential limit and complete the proof.

\subsection{Main result and scaling}
\label{subsec:s5-scaling}

We center the particles at the deterministic edge and choose the
space and time scales so that the limiting evolution has the
normalization in \eqref{eq:prelim-airy-martingale}.
Recall $E_s$, $\cA(s)$ and $\cF_s$ from
\Cref{subsec:edge-main-results,subsec:prelim-characteristics}.
Fix $T\in(s_*,T_0)$ and put
\begin{equation}\label{e:scaling paras}
 \chi=n^{-2/3}\cA(T)^{-1/3},\qquad
 \zeta=\frac{\theta n^{-1/3}\cA(T)^{-2/3}}
 {\partial_{11}\cF_T(m_T(E_T),E_T)}.
\end{equation}
Both parameters are positive by
\eqref{def: classical measure} and \eqref{eq:edge-flux-nondegeneracy}.
Here $\chi$ is the spatial scale which is different from the characteristic velocity  $\chi_s(z)$.
The normalizations give
\begin{equation}\label{s5:scale-identities}
 n\chi^{3/2}{\cA(T)}^{1/2}=1,\qquad
 \frac{\zeta\,\partial_{11}\cF_T
 (m_T(E_T),E_T)}{\theta n\chi^2}=1.
\end{equation}
For DBM, \eqref{e:F_CP} gives $\zeta=n\chi^2$.

For DBM, write $x_i(s)=\lambda_i(s)$.
For each of the three models, define
\begin{equation}\label{e: scale particle DPS}
 \widetilde x_i^{(n)}(t)
 =\frac{x_i(T+\zeta t)-E_{T+\zeta t}}{\chi}.
\end{equation}
For the discrete models this is first defined on
\begin{equation}\label{jt:grid}
 \mathbb T_n
 =\{t\in\mathbb R:T+\zeta t\in[s_*,T_0]\cap\mathbb Z_n\},
\end{equation}
whose spacing is $(n\zeta)^{-1}$, and then extended by linear
interpolation. DBM retains its continuous paths.
Every fixed compact time interval is contained in the rescaled
time domain for all sufficiently large $n$.
For corners, the same holds for the existence of every fixed
particle rank, since $T<T_0<1$.

\begin{theorem}\label{thm: main thm}
For each of the three models, assume
\Cref{ass:weightfunc,ass:prelim-log-branches,ass:edge-data}.
With the scaling in \eqref{e:scaling paras} and \eqref{e: scale particle DPS},
for every fixed $\sfT>0$ and $k\ge1$,
\begin{equation}\label{e:scaling result}
 (\widetilde x_1^{(n)},\ldots,\widetilde x_k^{(n)})
 \Longrightarrow
 (\mathcal A_1^\beta,\ldots,\mathcal A_k^\beta)
 \quad\text{in }C([-\sfT,\sfT],\mathbb R^k),
\end{equation}
where $(\mathcal A_i^\beta)_{i\ge1}$ is the
$\mathrm{ALE}_\beta$ from \Cref{subsec:prelim-overview}.
\end{theorem}

\begin{remark}
The arguments in \Cref{sec:dynamical_loop,s:eost,sec:ale-convergence}
also apply to the related continuous-time chains described in
\Cref{sec:models}, under the corresponding analytic and edge
regularity assumptions. We omit the details.
\end{remark}

The square-root behavior assumption in \Cref{ass:edge-data} gives the same local geometry at the edge
in all three models. In suitable coordinates, the leading relation
is $y^2=x$, the \textit{Airy spectral curve}. In Eynard--Orantin topological
recursion, this curve is the local model of a regular branch point
\cite[Section~10.5]{EynardOrantin2007}. This common geometry gives
a reason to expect Airy statistics at the edge across different
particle systems under the square root assumption.

In units of $1/n$ for both space and time, the scales $\chi$ and
$\zeta$ correspond to particle fluctuations of order $n^{1/3}$ over
a window of order $n^{2/3}$ around an observation time of order $n$.
These are the KPZ $1{:}2{:}3$ exponents \cite{Corwin2012KPZ}, with
the parameter along the line ensemble playing the role of the KPZ
spatial coordinate.
The square-root behavior also enters the characterization of the
Airy$_\beta$ line ensemble \cite{charofairy}.

To apply \Cref{thm:prelim-airychar}, we associate a
particle-generated Nevanlinna function with the rescaled particles.
Recall $\delta_s=m_s^n-m_s$ and the ordinary transforms
$\widehat m_s^n,\widehat m_s$ from
\eqref{s4:ordinary-transforms}.
Set
\begin{align}
 \Delta_t^{(n)}(w)
 &=n\chi\,\delta_{T+\zeta t}
       (E_{T+\zeta t}+\chi w),
 \label{e: def of del}\\
 Y_t^{(n)}(w)
 &=\sum_i\frac1{\widetilde x_i^{(n)}(t)-w}
       -n\chi\widehat m_{T+\zeta t}(E_{T+\zeta t}).
 \label{jt:Y}
\end{align}
The sum is over the current particles; for corners, the normalization
remains the top-row size $n$, as in \eqref{eq:corner-level-measures}.
Since $\widehat m_s(E_s)$ is real, subtracting the
constant in \eqref{jt:Y} preserves the Nevanlinna property
and leaves the poles unchanged, so $Y_t^{(n)}$ is a
particle-generated Nevanlinna function with poles
$\widetilde x_i^{(n)}(t)$.
For the discrete models, we linearly interpolate $\Delta^{(n)}$ and
$Y^{(n)}$ between grid times; the atomic representation in \eqref{jt:Y} is used only
on the grid.
For DMC, $\Delta^{(n)}$ is defined on the corresponding spectral
domain, which contains every fixed rescaled compact subset of
$\mathbb H$ for all sufficiently large $n$.

The estimates of Section~4 control $\Delta^{(n)}$, whereas
$Y^{(n)}$ records the particle locations.
The following comparison allows us to use those estimates in
the characterization of the limit.

\begin{lemma}\label{lem: expan particle stij}
Under either the assumptions of \Cref{thm: main thm}, for every fixed
$\sfT>0$ and $\mathcal K\Subset\mathbb H$,
\begin{equation}\label{jt:bridge}
 \sup_{|t|\le \sfT}\sup_{w\in\mathcal K}
 \left|Y_t^{(n)}(w)-\Delta_t^{(n)}(w)-\sqrt w\right|
 \xrightarrow[n\to\infty]{\mathbb P}0.
\end{equation}
\end{lemma}

\begin{proof}
We first work at the times of the model and write
$s=T+\zeta t$ and $z=E_s+\chi w$.
By \Cref{lem:imastij,lem:edge-interior-expansion},
\eqref{s4:kernel-split}, and \eqref{s5:scale-identities},
uniformly for $|t|\le \sfT$ and $w\in\mathcal K$,
\begin{equation}\label{jt:det-bridge}
 n\chi\bigl(\widehat m_s(z)-\widehat m_s(E_s)\bigr)
 =\sqrt{\frac{\cA(s)}{\cA(T)}}\sqrt w
       +\OO_{\mathcal K}(n^{-1/3}).
\end{equation}
By the Lipschitz continuity of $\sqrt{\cA}$ in
\Cref{lem:imastij} and $|s-T|\le\zeta \sfT$,
the first term differs from $\sqrt w$ by
$\OO_{\mathcal K,\sfT}(n^{-1/3})$. For corners and DBM, subtracting
\eqref{e: def of del} from \eqref{jt:Y} gives the left-hand side
of \eqref{jt:det-bridge}, proving the assertion.

For DMC, we also compare the generalized and ordinary kernels,
and replace the smoothed empirical measure by its particles.
On the support events used in the proof of \Cref{s4:outer-tight},
the contours $\cC(s)$ enclose both measures.
The residue formula and \eqref{s4:kernel-split} give
\begin{equation}\label{jt:test-contour}
 \int g_{\mathsf b}(r,z)(\mu_s^n-\mu_s)(\rd r)
 =-\frac1{2\pi\ri}\oint_{\cC(s)}
       g_{\mathsf b}(u,z)\delta_s(u)\,\rd u.
\end{equation}
The kernel is bounded on these contours. Therefore
\Cref{l: Delta is sml} implies
\begin{equation}\label{jt:field-compare}
 n\chi\left|\delta_s(z)
       -\bigl(\widehat m_s^n(z)-\widehat m_s(z)\bigr)\right|
 \le Cn\chi\sup_{u\in\cC(s)}|\delta_s(u)|
 =\OO_{\mathbb P}(\chi),
\end{equation}
uniformly over the same time and spectral sets.
The support events have probabilities arbitrarily close to one,
so the bound holds in probability without this restriction.

Since $\Im z\asymp\chi\gg n^{-1}$, the denominators
$x_i(s)-z$ and $r-z$ are comparable for
$r\in[x_i(s),x_i(s)+\theta/n]$.
Integrating their resolvent difference gives
\begin{equation}\label{jt:desmoothing}
 \left|\frac1n\sum_i\frac1{x_i(s)-z}-\widehat m_s^n(z)\right|
 \le\frac Cn\frac{\Im\widehat m_s^n(z)}{\Im z}.
\end{equation}
By \Cref{l: Delta is sml}, \eqref{s4:im-edge}, and
$|\Im g_{\mathsf b}(r,z)|\le C\Im z$,
$\Im\widehat m_s^n(z)=\OO_{\mathbb P}(n^{-1/3})$ uniformly.
Thus the right-hand side of \eqref{jt:desmoothing}, multiplied
by $n\chi$, is $\OO_{\mathbb P}(n^{-1/3})$.
Combining \eqref{jt:det-bridge}, \eqref{jt:field-compare},
and \eqref{jt:desmoothing}
proves \eqref{jt:bridge} on the grid.
For the discrete models, apply these bounds on a slightly enlarged
time window. Linear interpolation preserves the bound for
$Y^{(n)}-\Delta^{(n)}-\sqrt{\,\cdot\,}$, completing the proof.
\end{proof}

\subsection{Evolution of the rescaled Stieltjes transform}
\label{subsec:s5-evolution}

In this part, we derive the drift and martingale covariance in the edge scaling.
We also record the increment bounds used to prove tightness used in \Cref{sec:jt-revised}.

Fix $\sfT>0$ and consider the time interval $[-\sfT,\sfT]$ for the rescaled particle systems.
For DBM, set $a=-\sfT$. For the discrete models, let $a$ be
the largest point of $\mathbb T_n$ not exceeding $-\sfT$.
In the discrete models, we work from $a$ through the first grid point
at or after $\sfT$ and do the linear interpolation to $[-\sfT,\sfT]$.

For the discrete models, let $t,t'\in\mathbb T_n$ be successive
times, and write $s=T+\zeta t$, so that
$t'-t=(n\zeta)^{-1}$ and $T+\zeta t'=s+1/n$.
Define the martingale increment by
\begin{equation}\label{jt:ordinary-noise}
 \begin{aligned}
 \Delta M_t^{(n)}(w)=n\chi\Bigl\{&
 (\widehat m_{s+1/n}^n-\widehat m_s^n)(E_{s+1/n}+\chi w) -\mathbb E_s
 (\widehat m_{s+1/n}^n-\widehat m_s^n)(E_{s+1/n}+\chi w)
 \Bigr\}.
 \end{aligned}
\end{equation}
For the discrete models, set $M_a^{(n)}=0$ and define
$M_t^{(n)}$ at later grid times by summing these increments
over $u\in\mathbb T_n$ with $a\le u<t$.
For DBM, the martingale part of $\Delta_t^{(n)}(w)$ obtained from
\eqref{eq:pote_mt} is
\[
 M_t^{(n)}(w)
 =-\frac{\chi}{\sqrt{\theta n}}\sum_{i=1}^n
 \int_{T+\zeta a}^{T+\zeta t}
 \frac{\rd B_i(s)}{(\lambda_i(s)-E_s-\chi w)^2}.
\]
Thus $M_a^{(n)}=0$ for all three models.
We linearly interpolate $\Delta^{(n)}$, $Y^{(n)}$, and
$M^{(n)}$ between successive grid times when discussing
continuous-path convergence.

Choose compact neighborhoods
$\mathcal K\Subset\operatorname{Int}\mathcal K_1
\Subset\operatorname{Int}\mathcal K_2\Subset\mathbb H$.
Recall the time $\sigma_{R}$ introduced
before \eqref{s4:dbm-confinement}, with $R$ fixed
and sufficiently large.
Also, set $\sigma_{R}=\infty$ for DMC and corners.

\begin{definition}\label{def:s5-local-stopping}
For $R>0$, let $\tau_R^{(n)}$ be the first time $t$,
starting at $a$, at which at least one of the following
conditions holds, where $s=T+\zeta t$:
\begin{equation}\label{s5:local-stop}
 \begin{gathered}
 \sup_{w\in\mathcal K_2}|\Delta_t^{(n)}(w)|>R,
 \qquad
 n\sup_{z\in\omega_+(s)}|\delta_s(z)|>R,
 \qquad s\ge\sigma_{R},\\
 I_s\not\subset
 \begin{cases}
  [\ell_0,E_s+\kappa_*/8],&\text{DMC and corners},\\
  [L_s-\kappa_*/8,E_s+\kappa_*/8],&\text{DBM}.
 \end{cases}
 \end{gathered}
\end{equation}
The time $t$ ranges over $\mathbb T_n$ for DMC and corners and over
continuous times for DBM.
If a condition holds at $a$, then $\tau_R^{(n)}=a$;
if no condition occurs, set $\tau_R^{(n)}=\infty$.
\end{definition}

By the contour construction in
\Cref{subsec:s4-characteristics}, for $t<\tau_R^{(n)}$
and $s=T+\zeta t$, both $\cC(s)$ and $\omega_+(s)$
enclose $\operatorname{supp}\mu_s^n\cup\operatorname{supp}\mu_s$
and remain a fixed positive distance from this union.
Since $I_{s+1/n}\subset I_s+[0,1/n]$ for DMC and
$I_{s+1/n}\subset I_s$ for corners, the same contour
bounds hold at time $s+1/n$.
We use the stopping convention in
\eqref{s4:edge-stopped-martingale}; in particular, the
discrete martingale sums include the increment ending at
$\tau_R^{(n)}$.

The following lemma controls the stopping probability
and the difference between $Y^{(n)}$ and
$\Delta^{(n)}+\sqrt{\,\cdot\,}$ before the stopping time.

\begin{lemma}\label{lem:s5-stopped-comparison}
One has
\(
 \lim_{R\to\infty}\limsup_{n\to\infty}
 \mathbb P(\tau_R^{(n)}\le \sfT)=0.
\)
Moreover, for every fixed $R>0$ and all sufficiently
large $n$, at every time of the model with $|t|\le \sfT$
and $t<\tau_R^{(n)}$, we have
\begin{equation}\label{s5:stopped-bridge}
 \sup_{w\in\mathcal K_1}
 \left|Y_t^{(n)}(w)-\Delta_t^{(n)}(w)-\sqrt w\right|
 \le C_{R,\mathcal K,\sfT}n^{-1/3}.
\end{equation}
For each fixed integer $q\ge1$, we also have
\[
 \sup_{w\in\mathcal K}
 \left|
 \partial_w^q\bigl(
 Y_t^{(n)}(w)-\Delta_t^{(n)}(w)-\sqrt w
 \bigr)
 \right|
 \le C_{R,\mathcal K,\sfT,q}n^{-1/3}.
\]
The constants may depend on the fixed neighborhoods,
but not on $n$ or the current configuration.
\end{lemma}

\begin{proof}
By \eqref{e:deff},
$f(T+\zeta t)=\OO_{M,\sfT}(n^{-2/3})$.
Thus \Cref{prop:s4-edge-full-time,lem:b-bounds},
including the left-edge bound for DBM, show that the
probability of failure of the support condition in
\eqref{s5:local-stop} tends to zero, by taking first
$n\to\infty$ and then $M\to\infty$.
The first assertion now follows from
\Cref{l: Delta is sml,s4:outer-tight} and
\eqref{s4:dbm-confinement}.

For the comparison, we use the estimates in the proof
of \Cref{lem: expan particle stij}.
Write $s=T+\zeta t$ and $z=E_s+\chi w$.
By \eqref{jt:det-bridge} and the Lipschitz continuity
of $\sqrt{\cA}$ in \Cref{lem:imastij},
\[
 n\chi\bigl(\widehat m_s(z)-\widehat m_s(E_s)\bigr)
 =\sqrt w+\OO_{\mathcal K,\sfT}(n^{-1/3})
\]
uniformly for $|t|\le \sfT$ and $w\in\mathcal K_1$.
For corners and DBM, this proves
\eqref{s5:stopped-bridge} directly.

For DMC, applying \eqref{jt:test-contour} on
$\omega_+(s)$ and using \eqref{s5:local-stop} gives
\(
 n\chi\left|
 \delta_s(z)-\bigl(\widehat m_s^n(z)-\widehat m_s(z)\bigr)
 \right|
 \le CR\chi.
\)
The local bound in \eqref{s5:local-stop}, together with
\eqref{s4:im-edge} and \eqref{s4:kernel-split}, also gives
\(
 \Im\widehat m_s^n(z)
 \le C_{R,\mathcal K,\sfT}n^{-1/3}.
\)
Since $\Im z\asymp\chi$, \eqref{jt:desmoothing} implies
\[
 n\chi\left|
 \frac1n\sum_i\frac1{x_i(s)-z}-\widehat m_s^n(z)
 \right|
 \le C_{R,\mathcal K,\sfT}n^{-1/3}.
\]
Combining these estimates and using
$\chi=\OO(n^{-2/3})$ proves \eqref{s5:stopped-bridge}.
Cauchy's formula on $\mathcal K_1$ gives the stated
bounds for the $w$-derivatives on $\mathcal K$.
\end{proof}

We next derive the evolution of the rescaled transforms
and estimate their martingale increments.
The following proposition provides the bounds needed
for tightness and identifies the drift and covariance
that appear in the limiting equation.
\begin{proposition}\label{p:dismark scale sde}
Under the assumptions of \Cref{thm: main thm}, before
$\tau_R^{(n)}$, the discrete models satisfy
\begin{equation}\label{s5:scaled-increments}
 \sup_{w\in\mathcal K_1}
 \left|\mathbb E_s\bigl[
 \Delta_{t'}^{(n)}(w)-\Delta_t^{(n)}(w)\bigr]\right|
 \le\frac{C_R}{n\zeta},
\end{equation}
and, for each fixed integer $p\ge1$,
\begin{equation}\label{s5:scaled-moments}
 \sup_{w\in\mathcal K_1}\left\{
 \mathbb E_s|\Delta M_t^{(n)}(w)|^{2p}
 +\mathbb E_s\left|
 \chi\Delta\mathcal M_s^n(E_{s+1/n}+\chi w)\right|^{2p}
 \right\}
 \le C_{R,p}(n\zeta)^{-p}.
\end{equation}
For DBM, the drift of $\Delta^{(n)}$ and the sum of the
quadratic-variation rates of the real and imaginary parts of $M^{(n)}$
are bounded by $C_R$ on $\mathcal K_1$.

For each of the three models, with $a$ chosen as above,
\begin{equation}\label{s5:scaled-evolution}
 \begin{aligned}
 Y_t^{(n)}(w)-Y_a^{(n)}(w)
 ={}&M_t^{(n)}(w)-M_a^{(n)}(w)\\
 &+\int_a^t\left[
 \frac{2-\beta}{2\beta}\partial_w^2Y_u^{(n)}(w)
 +\frac12\partial_w\bigl(Y_u^{(n)}(w)\bigr)^2
 -\frac12\right]\rd u+\oo_{\mathbb P}(1),
 \end{aligned}
\end{equation}
uniformly for $a\le t\le \sfT$ and $w\in\mathcal K$.
For the discrete models, replace the integral by its left-endpoint
Riemann sum with step $(n\zeta)^{-1}$; the identity is first
asserted at grid times. The corresponding stopped identity also holds.
Uniformly before $\tau_R^{(n)}$,
\begin{equation}\label{jt:bracket-density}
 \begin{aligned}
 (n\zeta)\mathbb E_s
 [\Delta M_t^{(n)}(w)\Delta M_t^{(n)}(v)]
 &=\frac2\beta\partial_w\partial_v
       \frac{Y_t^{(n)}(w)-Y_t^{(n)}(v)}{w-v}+\oo(1),
 \end{aligned}
\end{equation}
for the discrete models. For DBM, the right-hand side of
\eqref{jt:bracket-density} without the error term equals
$\frac{\rd}{\rd t}\langle M^{(n)}(w),M^{(n)}(v)\rangle_t$.
The identity includes case $v=\bar w$. At $v=w$, its right-hand side is
$(3\beta)^{-1}\partial_w^3Y_t^{(n)}(w)$.
\end{proposition}

\begin{proof}
We write $s=T+\zeta t$ and $z=E_s+\chi w$.
All estimates below hold before \eqref{s5:local-stop}, on slightly
larger spectral compacts when derivatives are taken.
By the edge expansion and Cauchy's formula,
\begin{equation}\label{jt:localized-jets}
 \begin{gathered}
 |m_s^n(z)-m_s(E_s)|\le C_Rn^{-1/3},\\
 |\partial_z^qm_s^n(z)|+|\partial_z^qm_s(z)|
 \le C_{R,q}n^{(2q-1)/3},\quad q\ge1,\qquad
 P_s^n(z)\le C_Rn^{1/3}.
 \end{gathered}
\end{equation}
The last bound follows from \eqref{eq:main-P} and the velocity
formula \eqref{eq:velocity-from-flux}.
Also $\operatorname{dist}(z,I_s)\ge\Im z\ge c\chi$.
Since $|E_{s+1/n}-E_s|=\OO(n^{-1})=\oo(\chi)$, these bounds hold
at $E_{s+1/n}+\chi w$ and along the intervening segment.
Thus the one-step estimates of Section~3 apply, with the spectral
exponent chosen smaller than $1/3$.

\paragraph{The martingale increments.}
Inserting \eqref{jt:localized-jets} into
\eqref{s4:moment-input} bounds each of its two terms, after
multiplication by $\chi^{2p}$, by $C_{R,p}n^{-2p/3}$.
This proves \eqref{s5:scaled-moments} for the martingale increments
$\chi\Delta\mathcal M_s^n(E_{s+1/n}+\chi w)$.
For DMC, comparison with $\Delta M_t^{(n)}(w)$ gives
\begin{equation}\label{jt:noise-compare-increment}
 \chi\Delta\mathcal M_s^n(E_{s+1/n}+\chi w)
       -\Delta M_t^{(n)}(w)
 =-\frac{\chi}{2\pi\ri}\oint_{\omega_+(s)}
 g_{\mathsf b}(u,E_{s+1/n}+\chi w)
             \Delta\mathcal M_s^n(u)\,\rd u.
\end{equation}

By \Cref{def:s5-local-stopping},
$\omega_+(s)$ encloses $I_s+[0,1/n]$ before
$\tau_R^{(n)}$ for all sufficiently large $n$.
Applying \eqref{jt:test-contour} at time $s+1/n$
and subtracting the conditional expectation at time $s$
gives \eqref{jt:noise-compare-increment}.
At fixed spectral distance, \eqref{s4:moment-input} bounds the
conditional $2p$-th moment of this difference by
$C_p\chi^{2p}n^{-p}$. This also proves
\eqref{s5:scaled-moments} for $\Delta M_t^{(n)}(w)$.
Doob's inequality and Cauchy's formula further give
\begin{equation}\label{jt:noise-compare}
 \mathbb E\sup_t\sup_{w\in\mathcal K}
 \left|\sum_{u<t}\mathbf1_{\{u<\tau_R^{(n)}\}}
 \left[\chi\Delta\mathcal M_{T+\zeta u}^n
 (E_{T+\zeta u+1/n}+\chi w)-\Delta M_u^{(n)}(w)\right]
 \right|^2
 \le C_{R,\mathcal K,\sfT}\zeta\chi^2=\oo(1).
\end{equation}
Here and below such sums run over the grid points of the fixed
time window. For corners and DBM the two kernels agree.

\paragraph{The leading drift.}
The identity $a_s^n=-\chi_s^n\partial_zm_s^n+F_s^n$ and
$E_s'=\chi_s(E_s)$ give the exact cancellation
\begin{equation}\label{jt:drift-cancellation}
 \begin{aligned}
 a_s^n-\partial_sm_s+E_s'\partial_z\delta_s
 ={}&(\chi_s(E_s)-\chi_s^n)\partial_z\delta_s
 +(\chi_s-\chi_s^n)\partial_zm_s+F_s^n-F_s.
 \end{aligned}
\end{equation}
Each velocity difference is $\OO_R(n^{-1/3})$, while the
transform derivatives are $\OO_R(n^{1/3})$.
The outer-contour difference is bounded by
$C\sup_{\omega_+(s)}|\delta_s|$ after integration by parts
in \eqref{def of F_s}.
These estimates give a bounded leading drift after multiplication
by $n\chi\zeta$.

To identify it, expand $\cF_s$ in its first variable around
$m_s(E_s)$. The linear term cancels the last transport term on
the left of \eqref{jt:drift-cancellation}. The remaining local
quadratic term is
\[
 \frac{\partial_{11}\cF_s(m_s(E_s),E_s)}{2\theta}
 \partial_z\left[
 (m_s^n(z)-m_s(E_s))^2-(m_s(z)-m_s(E_s))^2\right].
\]
The cubic terms are $\OO_R(n^{-1/3})$ by
\eqref{jt:localized-jets}; replacing the second flux argument by
$E_s$ has the same bound. The finite-$n$ change of the corners
flux in \eqref{e:F_CP} contributes $\OO_R(n^{-2/3})$.
Using \eqref{jt:det-bridge} and \eqref{s5:scale-identities}, we
therefore obtain
\begin{equation}\label{s5:quadratic-drift}
 n\chi\zeta
 [a_s^n-\partial_sm_s+E_s'\partial_z\delta_s](z)
 =\frac12\partial_w
 \left[(\Delta_t^{(n)}(w)+\sqrt w)^2-w\right]+\OO_R(n^{-1/3}).
\end{equation}
The error is uniform before $\tau_R^{(n)}$ on the fixed
spectral compacts.
Here \Cref{lem:imastij} gives
$|\sqrt{\cA(s)}-\sqrt{\cA(T)}|\le C_Tn^{-1/3}$.

\paragraph{The second-order drift.}
For the discrete models we must retain the second-order mean from
Section~3, as well as the change of the evaluation point.
Recall that $E_s'=\chi_s(E_s)$.
Taylor expansion of \eqref{eq:main-drift} at $z$ gives
\begin{equation}\label{s5:moving-edge-mean}
 \begin{aligned}
 &\mathbb E_s[m_{s+1/n}^n(E_{s+1/n}+\chi w)-m_s^n(z)]\\
 &\quad=\frac1n[a_s^n+E_s'\partial_zm_s^n](z)
 +\frac1{n^2}
 \left[c_s^n+E_s'\partial_za_s^n+\frac{(E_s')^2}{2}\partial_z^2m_s^n\right](z)
 +\OO_R(n^{-4/3}).
 \end{aligned}
\end{equation}
Indeed, the third spatial derivative is $\OO_R(n^{5/3})$,
and \eqref{eq:main-second} gives the same remainder bound.
The error $\OO(n^{-2})$ in the edge displacement, from
\eqref{s4:edge-time}, contributes only $\OO_R(n^{-5/3})$.

Deforming the contours in \eqref{t:loopstudy2} and
\eqref{eq:corners-second-order-mean} gives
\begin{equation}\label{s5:second-mean-coefficients}
 c_s^n(z)=
 \begin{cases}
 \displaystyle
 \frac{(1-\theta)\chi_s^n(z)+\theta(\chi_s^n(z))^2}{2}
       \partial_z^2m_s^n(z)+\OO_R(n^{2/3}),&\text{DMC},\\[5pt]
 \displaystyle
 \frac{(\chi_s^n(z))^2}{2\theta}
       \partial_z^2m_s^n(z)+\OO_R(n^{2/3}),&\text{corners}.
 \end{cases}
\end{equation}
For DMC, the integral term
in \eqref{err:secondder} equals $-2\partial_zm_s^n$ plus a
holomorphic term with bounded derivatives. Using also \eqref{e:hz} and
$\partial_z[h_-(\mathsf b(z))]=-\theta a_s^n(z)$, its residues give
\[
 c_s^n=-\frac12\partial_z^2\chi_s^n
 -\theta\partial_z(\chi_s^na_s^n)
 -\theta\partial_z(\chi_s^n\partial_zm_s^n)
 -\frac1{2\theta}\partial_z^2\log\cB_s^n+\OO_R(n^{1/3}).
\]
For corners, the residue is
$-\theta^{-1}\partial_z
(\chi_s^na_s^n-\frac12\partial_z\chi_s^n)$,
up to a bounded outer-contour term.
Inserting the explicit fluxes from \eqref{e:F_CP} gives
\eqref{s5:second-mean-coefficients}; all products of two first
spatial derivatives are $\OO_R(n^{2/3})$.

Moreover,
$\partial_za_s^n=-\chi_s^n\partial_z^2m_s^n+\OO_R(n^{2/3})$
and $\chi_s^n(z)-v=\OO_R(n^{-1/3})$.
The two models consequently have the same coefficient after
the moving-edge correction:
\begin{equation}\label{s5:common-second-drift}
 c_s^n+v\partial_za_s^n+\frac{v^2}{2}\partial_z^2m_s^n
 =\frac{1-\theta}{2\theta^2}
 \partial_{11}\cF_s(m_s(E_s),E_s)\partial_z^2m_s^n
 +\OO_R(n^{2/3}).
\end{equation}
For DMC this uses
$\partial_{11}\cF_s=\theta^2\chi_s(1-\chi_s)$ at the edge;
for corners it uses $\partial_{11}\cF_s=\theta\chi_s^2$.
After rescaling, \eqref{s5:scale-identities} gives
\[
 \frac{\zeta}{n\chi^2}\,
 \frac{1-\theta}{2\theta^2}
 \partial_{11}\cF_T(m_T(E_T),E_T)
 =\frac{1-\theta}{2\theta}=\frac{2-\beta}{2\beta}.
\]

For the deterministic transform, its equation along the moving
edge is
\[
 (\partial_s+E_s'\partial_z)m_s(z)
 =F_s(z)+(E_s'-\chi_s(z))\partial_zm_s(z).
\]
Differentiating this identity once more along the moving edge,
using \eqref{s4:edge-time} and \eqref{jt:localized-jets}, bounds
the second derivative by $Cn^{1/3}$ almost everywhere.
Its time-evaluation error in \eqref{s5:moving-edge-mean} is
therefore negligible. After multiplication by $n\chi$, the Taylor and mean
error terms above are $\OO_R(n^{-1})$ per update.
There are $\OO(n\zeta)=\OO(n^{2/3})$ updates on a fixed
microscopic window, so their sum is $\OO_R(n^{-1/3})$.
The error in \eqref{s5:quadratic-drift} has the same
accumulated bound after multiplication by the time step
$(n\zeta)^{-1}$.
The same estimates, without identifying the coefficient, prove
\eqref{s5:scaled-increments}.

For DBM, \eqref{eq:pote_mt} already contains the second-order
drift $(1-\theta)(2\theta n)^{-1}\partial_z^2m_s^n$.
Since $\zeta=n\chi^2$, it gives the same limiting coefficient.
Its leading drift is covered by \eqref{jt:drift-cancellation} and
\eqref{s5:quadratic-drift}, including the potential terms through
$F_s^n-F_s$.

Summing the resulting evolution for $\Delta^{(n)}$, using
\eqref{jt:noise-compare}, and then applying \eqref{jt:bridge}
gives \eqref{s5:scaled-evolution} before each fixed stop.
The comparison error in \eqref{jt:bridge} enters as the difference
of its two end point values; we do not differentiate it in time.
Cauchy's formula on larger compacts handles its spatial
derivatives. Increasing $R$ removes the stop and proves the
asserted identity in probability.

\paragraph{The martingale covariations.}
For the discrete models, use \eqref{eq:main-cov} at
$E_{s+1/n}+\chi w$ and $E_{s+1/n}+\chi v$.
The velocity expansion on the stopped neighborhood is
\begin{equation}\label{jt:velocity-expansion}
 \chi_s^n(z)=\chi_s(E_s)
 -\frac{\partial_{11}\cF_T(m_T(E_T),E_T)}{\theta}
       (m_s^n(z)-m_s(E_s))+\OO_R(n^{-2/3}).
\end{equation}
For an ordinary kernel, contour deformation gives
\[
 \frac1{2\pi\ri}\oint_{\omega_-}
 \frac{f(r)\,\rd r}{(r-z)^2(r-u)^2}
 =\frac1{2\pi\ri}\oint_{\omega_+}
 \frac{f(r)\,\rd r}{(r-z)^2(r-u)^2}
 -\partial_z\partial_u\frac{f(z)-f(u)}{z-u}.
\]
The constant velocity contributes zero to the divided difference.
The negative sign of the linear term in
\eqref{jt:velocity-expansion} cancels the negative residue sign.
The resulting coefficient after scaling is
\[
 \frac{\zeta\,\partial_{11}\cF_T
 (m_T(E_T),E_T)}{n\chi^2\theta^2}
 =\frac1\theta=\frac2\beta.
\]
The analytic parts of the generalized kernels and the separated
outer-contour integral tend to zero after scaling.
The remainder in \eqref{jt:velocity-expansion} and the
one-step covariance error contribute $\OO_R(n^{-1/3})$; for the
latter, explicitly,
\[
 (n\zeta)\chi^2\frac{P_s^n(z)}{n^2d^3}
 \le C_Rn^{-1/3},\qquad d\ge c\chi.
\]
By \eqref{s5:stopped-bridge} and Cauchy's formula, replacing
$\Delta_t^{(n)}+\sqrt{\,\cdot\,}$ by $Y_t^{(n)}$ in the
divided-difference derivatives contributes $\OO_R(n^{-1/3})$.
Together with \eqref{jt:noise-compare-increment},
\eqref{s5:scaled-moments}, and conditional Cauchy--Schwarz,
this proves \eqref{jt:bracket-density} for the increments $\Delta M_t^{(n)}$.
Here moving the evaluation point back from $E_{s+1/n}+\chi w$
to $E_s+\chi w$ changes the rescaled parameter by $\OO(n^{-1/3})$;
the same Cauchy bounds control this error.
For DBM it follows directly from \eqref{eq:dbm-bracket} and
$\zeta=n\chi^2$.
Finally,
\[
 \left.\partial_w\partial_v
 \frac{f(w)-f(v)}{w-v}\right|_{v=w}=\frac16f'''(w),
\]
which gives the diagonal coefficient and completes the proof.
\end{proof}

\subsection{Joint tightness of the rescaled transforms and particle trajectories}
\label{sec:jt-revised}

We combine the estimates of the previous subsection with
\Cref{l: Delta is sml} to prove tightness in time and in the spatial
variable. We then bound each fixed ranked particle from below and
recover the particle trajectories from the ordered poles of $Y_t$
as $t$ varies.

\begin{theorem}[Joint tightness]\label{prop:joint-tightness-revised}
Under the assumptions of \Cref{thm: main thm}, for every fixed
$\sfT>0$ and $k\ge1$, the processes
\[
 \bigl(Y^{(n)},M^{(n)},
       \widetilde x_1^{(n)},\ldots,\widetilde x_k^{(n)}\bigr)
\]
are jointly tight in
$C([-\sfT,\sfT],\operatorname{Hol}(\mathbb H))^2
\times C([-\sfT,\sfT],\mathbb R^k)$, where the topology is
locally uniform convergence. 

Every joint subsequential limit satisfies
$Y_t(w)=\Delta_t(w)+\sqrt w$. At every time $t$, $Y_t$ is
particle-generated with infinitely many poles, and its poles,
ordered decreasingly and counted with multiplicity, are the
limiting particle positions $x_i(t)$.
Here the convergence of $\Delta^{(n)}$ is on each fixed spectral
compact, whose scaled image belongs to the spectral domain for
all sufficiently large $n$.
\end{theorem}

\begin{proof}
Fix $\sfT>0$ and choose $a$ as in the previous subsection.
We apply the estimates on a slightly larger fixed time interval,
which contains all grid intervals meeting $[-\sfT,\sfT]$ for the discrete models.
By \Cref{large-eigen: DMC,l: Delta is sml,s4:outer-tight}, for every
$\mathcal K\Subset\mathbb H$,
\begin{equation}\label{jt:window-input}
 \begin{gathered}
 \sup_{|t|\le \sfT}\sup_{w\in\mathcal K}
       |\Delta_t^{(n)}(w)|=\OO_{\mathbb P}(1),\qquad
 \left(\sup_{|t|\le \sfT}\widetilde x_1^{(n)}(t)\right)_+
       =\OO_{\mathbb P}(1),\\
 n\sup_{|t|\le \sfT}\sup_{z\in\omega_+(T+\zeta t)}
       |\delta_{T+\zeta t}(z)|=\OO_{\mathbb P}(1).
 \end{gathered}
\end{equation}
For discrete models the displayed bounds are first applied at grid
times, with the interpolation convention of the preceding subsection.

\paragraph{Tightness of the rescaled transforms.}
We follow the proof of \cite[Lemma~7.12]{charofairy}, using the
stopped increment estimates \eqref{s5:scaled-increments} and
\eqref{s5:scaled-moments} in place of the estimates there.
Choose slightly larger spectral compacts and use the common stopping
time $\tau_R^{(n)}$ from those estimates.
For $N$ consecutive discrete increments, the Burkholder--Davis--Gundy
inequality and \eqref{s5:scaled-moments} bound the $2p$-th moment
of their stopped sum by $C_{R,p}(N/(n\zeta))^p$.
The bounded drift in \eqref{s5:scaled-increments} gives the same
time-increment bound for $\Delta^{(n)}$.
Cauchy's formula on the larger compacts therefore yields
\begin{equation}\label{jt:time-moments}
 \mathbb E\sup_{w\in\mathcal K}
 \left|X_{v\wedge\tau_R^{(n)}}^{(n)}(w)
       -X_{u\wedge\tau_R^{(n)}}^{(n)}(w)\right|^{2p}
 \le C_{R,p,\mathcal K,\sfT}|v-u|^p,
 \qquad X=\Delta\ \text{or }M.
\end{equation}
For $\Delta^{(n)}$ we use the bounds for the rescaled martingale
increments $\chi\Delta\mathcal M_s^n$ in \eqref{s5:scaled-moments};
for $M^{(n)}$ we use those for $\Delta M_t^{(n)}$.
For discrete models, linear interpolation gives the same estimate
between grid times; for DBM, it follows from the continuous-time
Burkholder--Davis--Gundy inequality.

Kolmogorov's criterion gives tightness in time, and Cauchy's formula
and Montel's theorem give compactness in the spectral variable.
The initial values of $\Delta^{(n)}$ are tight by \eqref{jt:window-input},
and $M^{(n)}$ starts from zero on the enlarged time window.
Removing the stopping times by that estimate and taking a compact
exhaustion of $\mathbb H$ proves tightness of $\Delta^{(n)}$ and $M^{(n)}$;
\eqref{jt:bridge} gives tightness of $Y^{(n)}$ and the relation
$Y=\Delta+\sqrt{\,\cdot\,}$.
The fourth-moment estimate also gives
\[
 \mathbb E\sum_{t\in[-\sfT,\sfT)\cap\mathbb T_n}\sup_{w\in\mathcal K}
   |\mathbf1_{\{t<\tau_R^{(n)}\}}\Delta M_t^{(n)}(w)|^4
   \le \frac{C_{R,\mathcal K,\sfT}}{n\zeta}
   \longrightarrow0.
\]
Thus the step versions and linear interpolations have the same
limits, and the stopped increments have the uniform integrability
needed below.

\paragraph{Lower bounds for the particles.}
We show that each fixed ranked particle stays bounded below
throughout $[-\sfT,\sfT]$ with arbitrarily high probability.
Fix a small $\varepsilon_0>0$, and stop at the first failure of
\begin{equation}\label{jt:high-stop}
 |\Delta_t^{(n)}(\ri bH)|\le\varepsilon_0\sqrt{bH}
 \quad(b=1,2),\qquad
 n\sup_{z\in\omega_+(T+\zeta t)}
       |\delta_{T+\zeta t}(z)|\le R.
\end{equation}
As in the definition of $\tau_R^{(n)}$ in
\eqref{s5:local-stop}, we also stop if the empirical
support leaves the fixed interval enclosed by the contours.
The same time convention applies.
Here $H,R$ are fixed before taking $n\to\infty$.
Applying \eqref{jt:det-bridge} and \eqref{jt:desmoothing}
at $w=\ri bH$, and using \eqref{jt:test-contour} on
$\omega_+(T+\zeta t)$, gives, before this stopping time,
\[
 \left|Y_t^{(n)}(\ri bH)-\Delta_t^{(n)}(\ri bH)-\sqrt{\ri bH}\right|
 \le C_{H,R,\sfT}n^{-1/3},
 \qquad b=1,2.
\]
Here we use \eqref{jt:high-stop} and the Lipschitz
continuity of $\sqrt{\cA}$ from \Cref{lem:imastij}.
Consequently, $|Y_t^{(n)}(\ri bH)|\le C\sqrt H$
for all sufficiently large $n$.
By \cite[Lemma~2.2]{charofairy}, followed by Cauchy's formula,
\begin{equation}\label{jt:high-jets}
 |Y_t^{(n)}(w)|\le C\sqrt H,\qquad
 |\partial_w^qY_t^{(n)}(w)|\le C_qH^{1/2-q},
 \qquad |w-\ri bH|\le bH/4.
\end{equation}
The ordinary-kernel comparison transfers these bounds to
$\Delta^{(n)}$ on slightly smaller disks.

We apply the calculation giving \eqref{s5:scaled-increments} at
these spectral points. At $s=T+\zeta t$ and
$z=E_s+\chi \ri bH$, the bounds are
\[
 \operatorname{dist}(z,I_s)\ge c\chi H,\qquad
 P_s^n(z)\le Cn^{1/3}H^{-1/2}.
\]
The two factors in the leading drift difference are
$\OO(n^{-1/3}\sqrt H)$ and $\OO(n^{1/3}H^{-1/2})$.
The second-order terms contribute $\OO(H^{-3/2})$ to the drift
rate, and the outer-contour difference contributes $\OO(R/n)$.
Consequently, for $t'=t+(n\zeta)^{-1}$ and $b=1,2$,
\begin{equation}\label{jt:high-moments}
 \begin{aligned}
 \left|\mathbb E_s\bigl[
   \Delta_{t'}^{(n)}(\ri bH)-\Delta_t^{(n)}(\ri bH)\bigr]\right|
 &\le \frac C{n\zeta},\\
 \mathbb E_s\left|
   \Delta_{t'}^{(n)}(\ri bH)
     -\mathbb E_s\Delta_{t'}^{(n)}(\ri bH)\right|^2
 &\le \frac C{n\zeta H^{5/2}}.
 \end{aligned}
\end{equation}
The variance estimate is \eqref{s4:moment-input} with the displayed
values of $P_s^n$ and the spectral distance.
For DBM, \eqref{eq:pote_mt} and \eqref{eq:dbm-bracket} give the corresponding
drift bound $C$ and quadratic-variation rate $CH^{-5/2}$.
The constants are independent of $H,R$, with the lower bound on
$n$ allowed to depend on them.

At time $a$, \Cref{bulktight: DMC}
gives $|\Delta_a^{(n)}(\ri bH)|\le\varepsilon_0\sqrt{bH}/4$
outside an event whose limiting probability is at most
$Ce^{-cH^{\alpha/2}}$.
For $H$ sufficiently large compared with $(1+\sfT)^2$, the accumulated
drift in \eqref{jt:high-moments} is less than
$\varepsilon_0\sqrt H/4$.
Doob's maximal inequality bounds the probability of a larger
martingale increment by $C_{\sfT,\varepsilon_0}H^{-7/2}$.
The increment reaching the stopping time is included.
Using \eqref{jt:window-input} to let $R\to\infty$ after
$n\to\infty$, and then letting $H\to\infty$, proves
\begin{equation}\label{jt:high-propagation}
 \lim_{H\to\infty}\limsup_{n\to\infty}
 \mathbb P\left(
  \sup_{|t|\le \sfT}\max_{b=1,2}
  \frac{|Y_t^{(n)}(\ri bH)-\sqrt{\ri bH}|}{\sqrt H}
    >3\varepsilon_0\right)=0.
\end{equation}

At each grid time, the real centering of $Y^{(n)}$ gives
\begin{equation}\label{jt:positive-count}
 2\Im Y_t^{(n)}(\ri H)-\Im Y_t^{(n)}(2\ri H)
 =\sum_i\frac{6H^3}
 {((\widetilde x_i^{(n)}(t))^2+H^2)
  ((\widetilde x_i^{(n)}(t))^2+4H^2)}.
\end{equation}
For sufficiently small $\varepsilon_0$, \eqref{jt:high-propagation}
makes the left side at least $c\sqrt H$, simultaneously in time.
For $|x|>LH$, the summand is at most
$6L^{-2}H/(x^2+H^2)$. Since
$\Im Y_t^{(n)}(\ri H)\le C\sqrt H$, choosing $L$ sufficiently large
makes the total contribution of these particles at most
$c\sqrt H/2$.
Each remaining particle contributes at most $3/(2H)$, and hence
\begin{equation}\label{jt:count-bound}
 \#\{i:|\widetilde x_i^{(n)}(t)|\le LH\}
 \ge c'H^{3/2}
\end{equation}
at every grid time on this event. For every fixed $k$, this bounds
$\widetilde x_k^{(n)}$ below uniformly in time with arbitrarily
high probability. Linear interpolation preserves the bound.
For DBM, the identity and the argument hold at every time.

\paragraph{Convergence of the particle trajectories.}
We use the pole convergence argument in the proof of
\cite[Proposition~7.6]{charofairy}, with the particle bounds just
proved. To apply it simultaneously in time, include the real
random variables
\[
 \left(\sup_{|t|\le \sfT}\widetilde x_1^{(n)}(t)\right)_+,
 \qquad
 \left(\sup_{|t|\le \sfT}-\widetilde x_i^{(n)}(t)\right)_+,
 \quad i\ge1,
\]
in the joint state with $\Delta^{(n)}$, $Y^{(n)}$, and $M^{(n)}$.
Missing ranks may be assigned zero in these auxiliary coordinates;
each fixed rank is present for all sufficiently large $n$.
These variables are tight by \eqref{jt:window-input} and
\eqref{jt:count-bound}. After taking a subsequence and a
Skorokhod representation, they converge to finite values.
Thus all particles have a common upper bound on the coupled
sequence, and each fixed rank has a lower bound, both uniform
throughout the time window.

By \cite[Lemmas~2.7--2.8]{charofairy}, locally uniform convergence
of particle-generated Nevanlinna functions gives convergence of
their ordered poles, provided these poles do not escape to infinity.
For any $t$ and any sequence of grid times $t_n\to t$, convergence of
$Y^{(n)}$ uniformly on compact time and spectral sets gives
$Y_{t_n}^{(n)}\to Y_t$ locally uniformly.
The particle bounds exclude escape of every fixed rank, so $Y_t$
is particle-generated with infinitely many poles and
$\widetilde x_i^{(n)}(t_n)\to x_i(t)$ for every fixed $i$,
where $x_i(t)$ is its $i$-th pole.
Applying the same pole convergence to $Y_{t_j}\to Y_t$ as
$t_j\to t$ shows that each $x_i$ is continuous.

If convergence of a particle trajectory were not uniform at grid
times, compactness would give a sequence $t_n\to t$ contradicting
the preceding convergence and continuity of $x_i$.
Linear interpolation therefore gives
\[
 \sup_{|t|\le \sfT}|\widetilde x_i^{(n)}(t)-x_i(t)|
 \longrightarrow0.
\]
For DBM the argument uses arbitrary times directly.
This proves the asserted joint tightness and identifies the
limiting trajectories. These supreme bound are used only for the coupling and are omitted
from the filtration below.
\end{proof}

\subsection{Identification of the limit}
\label{subsec:s5-identification}

We verify the Airy-like and martingale conditions in
\Cref{thm:prelim-airychar} for every subsequential limit from
\Cref{prop:joint-tightness-revised}.
Passing to a diagonal subsequence, we obtain $(\Delta,Y)$ and
compatible limiting increments of $M^{(n)}$ on all compact time
intervals. On each $[-\sfT,\sfT]$, write $M$ for the resulting process,
normalized by $M_{-\sfT}=0$.

\paragraph{The Airy-like condition.}
Following the proof of \cite[Proposition~7.7]{charofairy}, we pass
the bounds in \Cref{large-eigen: DMC,bulktight: DMC} to the limit
using \eqref{jt:bridge}. For every $\varepsilon>0$, there is a
deterministic $C<\infty$, independent of $t$, such that, at each
fixed $t\in\mathbb R$, with probability at least $1-\varepsilon$,
all poles of $Y_t$ are at most $C$ and
\[
 |Y_t(w)-\sqrt w|
 \le \frac{C(\Im\sqrt w)^{\alpha'}}{\Im w},
 \qquad \Im w\ge C\sqrt{1+\max\{\Re w,0\}}.
\]
Here $\alpha'<1$ is from \Cref{bulktight: DMC}.
The uniformity in $t$ follows from the uniformity of the finite-$n$
estimates for $s$ near $T$, since
$T+\zeta t\to T$ for each fixed $t$.
Taking $t_j=-j$ and letting $C_{*,j}$ be the least positive integer
for which these bounds hold gives a tight family $(C_{*,j})_{j\ge1}$.
Thus the Airy-like condition holds with $\mathfrak d=1-\alpha'$.

\begin{proposition}\label{s4:joint-martingale}
Every joint subsequential limit in
\Cref{prop:joint-tightness-revised} satisfies
\eqref{eq:prelim-airy-martingale} and
\eqref{eq:prelim-airy-bracket} relative to a common filtration.
For every $\sfT>0$ and $w\in\mathbb H$, the process
$(M_t(w))_{-\sfT\le t\le \sfT}$, normalized by $M_{-\sfT}(w)=0$,
is a continuous square-integrable martingale.
\end{proposition}

\begin{proof}
We follow the proof of \cite[Proposition~7.6]{charofairy}, using
\eqref{s5:scaled-evolution}, \eqref{jt:bracket-density} and
\eqref{jt:time-moments}.

Fix an integer $\sfT\ge1$.
Set $\tau_R^{(n)}=\sfT+1$ if there is no exit by time $\sfT$.
Passing to a further subsequence, we may assume that these times
converge jointly with the processes, for all integer $R,\sfT$ and
an increasing sequence of compact subsets covering $\mathbb H$.
Write $\tau_R$ for their limits on $[-\sfT,\sfT]$.
We use the filtration generated by $\Delta_s,Y_s$, increments
$M_v-M_u$ with $u\le v\le t$, and the events
$\{\tau_R\le s\}$ with $s\le t$, for all these limiting times.

The fourth-moment bound in \eqref{jt:time-moments} gives uniform
integrability. We can therefore pass the stopped martingale
identities to the limit and identify their quadratic covariations
using \eqref{jt:bracket-density}.
For the discrete models, we first work at grid times.
The estimate in the proof of \Cref{prop:joint-tightness-revised}
shows that the largest martingale jump tends to zero, so linear
interpolation gives the same limit.

Cauchy's formula gives convergence of the spatial derivatives.
Passing to the limit in \eqref{s5:scaled-evolution}, with sums
replaced by time integrals for the discrete models, gives
\eqref{eq:prelim-airy-martingale} up to $\tau_R$.
By \Cref{lem:s5-stopped-comparison}, almost surely
$\tau_R=\sfT+1$ for all sufficiently large $R$.
Thus the evolution equation holds on $[-\sfT,\sfT]$, and $M$ is a
continuous local martingale with
\begin{align}
 \rd\langle M(w),M(v)\rangle_t
 &=\frac2\beta\partial_w\partial_v
       \frac{Y_t(w)-Y_t(v)}{w-v}\,\rd t,
       &&w\ne v,\label{jt:limit-bracket}\\
 \rd\langle M(w),M(w)\rangle_t
 &=\frac1{3\beta}\partial_w^3Y_t(w)\,\rd t.
       \label{jt:diag-bracket}
\end{align}

It remains to prove square-integrability.
By \Cref{bulktight: DMC} and \eqref{jt:bridge},
\[
 \mathbb P\bigl(|Y_t(\ri H)|>C\sqrt H\bigr)
 \le Ce^{-cH^{\alpha/2}},
\]
uniformly for $|t|\le \sfT$.
On the complementary event, Schwarz--Pick gives
$\sup_{w\in\mathcal K}|Y_t(w)|\le C_{\mathcal K}H^{3/2}$.
Hence $\sup_{|t|\le \sfT}\mathbb E\Im Y_t(w)<\infty$.
The bracket identity and the Nevanlinna representation yield
\[
 \mathbb E\langle M(w),M(\bar w)\rangle_\sfT
 \le \frac{2}{\beta(\Im w)^3}
       \int_{-\sfT}^{\sfT}\mathbb E\Im Y_t(w)\,\rd t<\infty.
\]
Therefore $M(w)$ is a square-integrable martingale.
The filtration is the same for all time intervals.
\end{proof}

\begin{proof}[Proof of \Cref{thm: main thm}]
By \Cref{prop:joint-tightness-revised}, every subsequence has a
further jointly convergent subsequence, whose limiting particle
positions are the ordered poles of $Y$.
The Airy-like bounds and \Cref{s4:joint-martingale} verify the two
conditions of \Cref{thm:prelim-airychar}, so these particle trajectories
have law $\mathrm{ALE}_\beta$.
Every subsequential limit has this law, proving the claimed convergence.
\end{proof}

\appendix
\section{One-step dynamical estimates}
\label{s:dyn-loop-beta-corners}

\subsection{The corners loop equation}

We prove the fixed-distance corners estimates used in
\Cref{thm:loopstudy} and then record the DMC covariance estimate.
Throughout this appendix, we assume
\Cref{ass:weightfunc,ass:prelim-log-branches} and use the fixed-distance
spectral sets and contour conventions of
\Cref{subsec:prelim-fixed-inputs}.
All estimates are uniform in time, in the admissible current
configuration and on these spectral sets. Constants may depend on
$\theta$, $T$ and the fixed domains, but not on $n$.

Fix $s=s_j=j/n$ and suppress the time index, writing
$I=I_s$, $\Lambda=\Lambda_s$ and $\mathscr D=\mathscr D_s$.
All expectations are conditional on the current configuration $\bmx$.
Recall \eqref{eq:corners-forward-increment} and
\Cref{subsec:prelim-fixed-inputs}; we write
\begin{equation}\label{e:Dn-definition}
 D(z)=D_j^n(z),\qquad
 \Delta\cM(z)=\Delta\cM_j^n(z)=D(z)-\bE D(z).
\end{equation}
For corners, the current and next rows have $k+1=n-j$ and $k=n-j-1$
particles, respectively. Thus $k\asymp n$ for $s\le T<1$, and
\eqref{e:B_function} gives
\begin{equation}\label{e:finite-step-slope}
 \cB(z)=\cB_s^n(z)=\frac{k}{n}W'(z)-\theta m_s^n(z).
\end{equation}
In particular, \Cref{ass:prelim-log-branches} gives
\begin{equation}\label{eq:corners-annulus}
 c\le |\cB(z)|\le C\quad(z\in\mathscr D),\qquad
 \operatorname{wind}_{\omega_-}\cB=-1.
\end{equation}
We use $R=R_{\bmy}$ from \eqref{eq:corners-ratio} and
$\cA$ from the zero-insertion case of \eqref{e:phi-tilde-def-Corner}:
\begin{equation}\label{e:cA-n-definition}
 \cA(z)=\bE R(z).
\end{equation}

\begin{theorem}[\cite{Gorin:2022qky}, Corollary~3.6]
\label{lem:dyn-loop-one-step}
For the transition \eqref{eq:beta-corners}, the function
\begin{equation}\label{e:cC-n-definition}
 \cC(z)=\bE\!\left[
 R(z)\left(\frac{k}{n}W'(z)
 +\frac{\theta-1}{n}\sum_{a=1}^{k+1}\frac1{z-x_a}
 +\frac1n\sum_{i=1}^k\frac1{z-y_i}\right)\right]
\end{equation}
extends holomorphically to $\Lambda$.
\end{theorem}

Here the potential in the cited result is $kW$, and we divide its
observable by $n$. Since $R'/R=D$ by \eqref{eq:corners-ratio},
\begin{equation}\label{eq:corners-scalar-loop}
 \cC=\cB\cA-\frac1n\cA'.
\end{equation}

\subsection{Properties of \texorpdfstring{$\cA(z)$}{A}}

We apply the Cauchy factorization argument of
\cite[Propositions~4.11 and~4.13]{Gorin:2022qky} to
\eqref{eq:corners-scalar-loop}.
Fix $e\in I$. The function $(z-e)\cB(z)$ has a single-valued
holomorphic logarithm on $\mathscr D$. For $z$ outside $\omega_-$, put
\begin{align}
 h_-(z)&=\frac1{2\pi\ri}\oint_{\omega_-}
 \frac{\log((w-e)\cB(w))}{w-z}\,\rd w,\notag\\
 h_-'(z)+\frac1{z-e}
 &=\frac1{2\pi\ri}\oint_{\omega_-}
 \frac{\cB'(w)}{\cB(w)(w-z)}\,\rd w.
 \label{eq:corners-factor-derivative}
\end{align}
Let $h_+$ be the interior Cauchy part, so that
$\log((z-e)\cB(z))=h_+(z)-h_-(z)$ on the working annulus.
On the integration contours, we use their holomorphic continuations
within this annulus.

\begin{proposition}
\label{prop:first-order-cA}\label{prop:second-order-cA}
For all sufficiently large $n$,
\begin{equation}\label{e:first-order-log-cA}
 \cA(z)=(z-e)e^{h_-(z)}\bigl(1+\OO(n^{-1})\bigr).
\end{equation}
Moreover, $\cC$ has no zeros inside $\omega_-$, and
\begin{equation}\label{eq:corners-factor-second}
 \frac{\cA(z)}{z-e}e^{-h_-(z)}
 =1-\frac1{2\pi\ri n}\oint_{\omega_-}
 \frac{h_-'(w)+\dfrac1{w-e}}{\cB(w)(w-z)}\,\rd w
 +\OO(n^{-2}).
\end{equation}
These statements, and the exact identity
\eqref{eq:corners-factor-exact} below, also hold for complex measures
of mass one and uniformly bounded total variation, provided
\eqref{eq:corners-annulus} holds and $\cC$ in
\eqref{eq:corners-scalar-loop} extends holomorphically inside the contour.
\end{proposition}

\begin{proof}
Interlacing gives
\[
 R(z)^{-1}=\sum_{a=1}^{k+1}\frac{a_a}{z-x_a},\qquad
 a_a=\frac{\prod_i(x_a-y_i)}{\prod_{b\ne a}(x_a-x_b)}>0,
 \qquad \sum_a a_a=1.
\]
Consequently, $R$, $R^{-1}$ and $D$ are uniformly bounded on the
fixed spectral sets. In particular, $\cA$ and its fixed-order
derivatives are bounded there, also for complex measures of bounded
total variation. Mass one gives $\cA(z)=z+\OO(1)$ at infinity.

Multiplying \eqref{eq:corners-scalar-loop} by $e^{-h_+}$ and taking
the exterior Cauchy projection, as in
\cite[(121)--(126)]{Gorin:2022qky}, gives
\begin{equation}\label{eq:corners-factor-exact}
 \frac{\cA(z)}{z-e}e^{-h_-(z)}
 =1-\frac1{2\pi\ri n}\oint_{\omega_-}
 \frac{\cA'(w)e^{-h_+(w)}}{w-z}\,\rd w.
\end{equation}
The integral is bounded, proving \eqref{e:first-order-log-cA}.
Apply this estimate on a slightly larger fixed neighborhood and
substitute its derivative into \eqref{eq:corners-factor-exact};
this gives \eqref{eq:corners-factor-second}.

Finally, \eqref{e:first-order-log-cA} gives
$\operatorname{wind}_{\omega_-}\cA=1$, while
$\cC=\cB\cA(1+\OO(n^{-1}))$ on the contour.
Thus $\operatorname{wind}_{\omega_-}\cC=0$, and the argument principle
shows that $\cC$ has no zeros inside. This proof uses only upper
bounds on the complex measure and therefore also proves the stated extension.
\end{proof}

\subsection{One-step mean estimate}
We derive the second-order mean expansion for the corners transition.
An inverse insertion expresses $\bE D(z)$ as a logarithmic derivative
of the deformed $\cA$, to which we apply
\Cref{prop:second-order-cA}.

For the inverse insertion in \eqref{e:phi-tilde-def-Corner},
\begin{equation}\label{eq:corners-inverse-tilt}
 \cA^{(v,1)}(u)
 =\frac{\bE[R(u)R(v)^{-1}]}{\bE R(v)^{-1}},
\end{equation}
and its coefficient is
\begin{equation}\label{eq:corners-tilted-observables}
 \cB^{(v,1)}(w)=\cB(w)+\frac1{n(w-v)}.
\end{equation}
The residue representation in the proof of
\Cref{prop:first-order-cA} gives $|\bE R(v)^{-1}|\ge c$ on the fixed
spectral sets. Hence this normalized complex measure has bounded
total variation.

\begin{proposition}[One-step mean estimate]
\label{prop:one-step-mean}\label{lem:corners-inverse-tilt}
With the insertion parameter $v$ held fixed during differentiation,
\begin{equation}\label{eq:corners-tilt-diagonal}
 \cA^{(v,1)}(v)=\frac1{\bE R(v)^{-1}},\qquad
 \left.\partial_u\log\cA^{(v,1)}(u)\right|_{u=v}=\bE D(v).
\end{equation}
Moreover,
\begin{equation}\label{eq:corners-explicit-mean}
\begin{aligned}
 \bE D(z)
 ={}&\frac1{2\pi\ri}\oint_{\omega_-}
 \frac{\cB'(w)}{\cB(w)(w-z)}\,\rd w\\
 &+\frac1{2\pi\ri n}\oint_{\omega_-'}\frac1{(w-z)^2}
 \left[-\frac1{\cB(w)}\frac1{2\pi\ri}\oint_{\omega_-}
 \frac{\cB'(u)}{\cB(u)(u-w)}\,\rd u
 -\frac{\cB'(w)}{2\cB(w)^2}\right]\rd w+\OO(n^{-2}).
\end{aligned}
\end{equation}
Here $\omega_-$ lies strictly inside $\omega_-'$; both contours
enclose $I$ and exclude $z$, with fixed positive separation.
\end{proposition}

\begin{proof}
The identities \eqref{eq:corners-tilt-diagonal} follow by evaluating
\eqref{eq:corners-inverse-tilt} and its $u$-derivative at $u=v$.
For real $v>\sup I$, the tilted transition has potential
$W+k^{-1}\log(v-\cdot)$. Apply \Cref{lem:dyn-loop-one-step} to
this potential and continue the unnormalized loop identity
holomorphically in $v$. The nonzero normalization and
\eqref{eq:corners-tilted-observables} then allow us to apply
\Cref{prop:second-order-cA} for every admissible complex $v$.

Following the inverse-insertion argument of
\cite[Proposition~4.16]{Gorin:2022qky}, take the logarithmic
$u$-derivative of \eqref{eq:corners-factor-second} for the tilted
transition and set $u=v=z$. By \eqref{eq:corners-factor-derivative},
\begin{equation}\label{eq:corners-general-mean}
\begin{aligned}
 \bE D(z)
 ={}&\frac1{2\pi\ri}\oint_{\omega_-}
 \frac{(\cB^{(z,1)})'(w)}{\cB^{(z,1)}(w)(w-z)}\,\rd w\\
 &-\frac1{(2\pi\ri)^2n}\oint_{\omega_-'}
 \frac{\rd w}{\cB^{(z,1)}(w)(w-z)^2}
 \oint_{\omega_-}
 \frac{(\cB^{(z,1)})'(u)}{\cB^{(z,1)}(u)(u-w)}\,\rd u
 +\OO(n^{-2}).
\end{aligned}
\end{equation}
Expand \eqref{eq:corners-tilted-observables} in this expression.
The insertion contributes
\[
 \frac1{2\pi\ri n}\oint_{\omega_-'}
 \frac{\rd w}{\cB(w)(w-z)^3}
 =-\frac1{4\pi\ri n}\oint_{\omega_-'}
 \frac{\cB'(w)\,\rd w}{\cB(w)^2(w-z)^2}.
\]
This gives \eqref{eq:corners-explicit-mean}.
All expansions hold on fixed neighborhoods, so Cauchy's formula
preserves the $\OO(n^{-2})$ error after differentiation.
\end{proof}

\subsection{One-step covariance estimate}
We obtain the conditional covariance and Laplace estimates from
the mean expansion under exponential deformations.
The Laplace estimate then gives the Gaussian asymptotics of
all joint moments.

\begin{proposition}[Covariance and Laplace estimates]
\label{prop:corners-covariance}
For admissible $z,u$,
\begin{equation}\label{e:one-step-covariance}
 n\bE[\Delta\cM(z)\Delta\cM(u)]
 =-\frac1{2\pi\ri}\oint_{\omega_-}
 \frac{\rd w}{\cB(w)(w-z)^2(w-u)^2}+\OO(n^{-1}).
\end{equation}
For a fixed number of conjugate-paired parameters $(a_j,z_j)$ with
bounded coefficients and $\sum_j a_j/(z_j-w)\in\mathbb R$ for $w\in I$,
\begin{equation}\label{eq:corners-laplace}
 \log\bE\exp\!\left(-\sqrt n\sum_j a_j\Delta\cM(z_j)\right)
 =-\frac12\sum_{i,j}a_i a_j\frac1{2\pi\ri}\oint_{\omega_-}
 \frac{\rd w}{\cB(w)(w-z_i)^2(w-z_j)^2}
 +\OO(n^{-1/2}).
\end{equation}
These estimates are uniform on fixed admissible compact spectral sets;
coincident and conjugate spectral parameters are allowed.
\end{proposition}

\begin{proof}
For covariance, use the potential $W_t(w)=W(w)+t/[k(u-w)]$.
Its normalized expectation is
$\bE^{(t)}f=\bE[f e^{-tD(u)}]/\bE e^{-tD(u)}$, and
\[
 \cB_t(w)=\cB(w)+\frac{t}{n(w-u)^2},\qquad
 \left.\partial_t\bE^{(t)}D(z)\right|_{t=0}
 =-\operatorname{Cov}(D(z),D(u)).
\]
Since $D(u)$ is bounded, the exponential and inverse-insertion
normalizations stay nonzero on a fixed disk $|t|\le t_0$.
Thus \eqref{eq:corners-explicit-mean} applies holomorphically and
uniformly in $t$, giving
\[
 \bE^{(t)}D(z)-\bE D(z)
 =\frac{t}{2\pi\ri n}\oint_{\omega_-}
 \frac{\rd w}{\cB(w)(w-z)^2(w-u)^2}+\OO(n^{-2}).
\]
Cauchy's formula differentiates the error at $t=0$ and proves
\eqref{e:one-step-covariance}.

For \eqref{eq:corners-laplace}, use the real exponential deformation
of \cite[Proposition~4.18, (160)--(168)]{Gorin:2022qky}.
Put $h(w)=\sum_j a_j/(z_j-w)$ and $W_t=W+t\sqrt n\,h/k$,
$0\le t\le1$. The deformed transition is a probability measure,
and $\cB_t=\cB+tn^{-1/2}h'$ preserves the annular bounds and winding.
Expanding \eqref{eq:corners-explicit-mean} gives
\begin{equation}\label{eq:corners-gaussian-response}
 \bE^{(t)}D(z)-\bE D(z)
 =\frac{t}{2\pi\ri\sqrt n}\oint_{\omega_-}
 \frac{h'(w)}{\cB(w)(w-z)^2}\,\rd w+\OO(n^{-1}).
\end{equation}
Integrate the logarithmic Laplace derivative
$-\sqrt n\sum_j a_j(\bE^{(t)}D(z_j)-\bE D(z_j))$
over $t\in[0,1]$ to obtain \eqref{eq:corners-laplace}.
\end{proof}

\begin{corollary}
\label{prop:one-step-covariance}
For every fixed $r\ge1$,
\begin{equation}\label{eq:corners-wick-moments}
 \bE\prod_{j=1}^r\sqrt n\,\Delta\cM(z_j)
 =\sum_{\pi\in\mathcal P_2(r)}\prod_{\{i,j\}\in\pi}
 n\bE[\Delta\cM(z_i)\Delta\cM(z_j)]+\oo(1),
\end{equation}
uniformly on fixed admissible compact spectral sets.
Here $\mathcal P_2(r)$ is the set of pairings; the sum is empty
for odd $r$. Mixed conjugate moments are included by replacing
the corresponding arguments by $\bar z_j$.
Along every covariance-convergent subsequence, the joint moments
converge to those of a centered Gaussian vector, possibly degenerate.
\end{corollary}

\begin{proof}
Apply \cite[Lemma~4.19]{Gorin:2022qky} to the real and imaginary
parts, using \eqref{eq:corners-laplace}.
The covariance matrices are uniformly bounded by
\eqref{e:one-step-covariance}. Extracting a convergent subsequence
of these finite matrices gives Gaussian moment convergence and
\eqref{eq:corners-wick-moments}, also when the parameters vary in
the fixed spectral sets.
\end{proof}

\begin{theorem}
\label{t:loop_beta_ctscorner}
We have
\begin{equation}\label{eq:appendix-corners-expansion}
 D(z)=\frac1{2\pi\ri}\oint_{\omega_-}
 \frac{\cB'(w)}{\cB(w)(w-z)}\,\rd w
 +\Delta\cM(z)+\OO(n^{-1}).
\end{equation}
Its second-order mean, covariance and Laplace estimates are
\eqref{eq:corners-explicit-mean}, \eqref{e:one-step-covariance}
and \eqref{eq:corners-laplace}, respectively.
\end{theorem}

\begin{proof}
Combine \eqref{e:Dn-definition} with
\Cref{prop:one-step-mean,prop:corners-covariance}.
\end{proof}

\subsection{Proof of \texorpdfstring{\Cref{thm:loopstudy}}{Proposition \getrefnumber{thm:loopstudy}}}
\label{subsec:proof-prop-2.5}
We recall the DMC covariance estimate in \Cref{thm:loopstudy}
by differentiating the second-order mean expansion under an exponential
tilt and then finish its proof.

\begin{lemma}\label{lem:dmc-fixed-covariance}
For DMC, uniformly for $z,u$ in the compact subset of spectral domain in \Cref{subsec:prelim-fixed-inputs}, which have positive distance from the cut,
\[
 n\bE[\Delta\cM_j^n(z)\Delta\cM_j^n(u)\mid\mathscr F_j]
 =\frac1{2\pi\ri\theta}\oint_{\omega_-}
 \chi_{s_j}^n(w)K_{\sfb}(w,z)K_{\sfb}(w,u)\,\rd w
 +\OO(n^{-1}).
\]
\end{lemma}

\begin{proof}
We use the second-order mean expansion of
\cite[Proposition~4.16]{Gorin:2022qky} under a small exponential tilt.
Condition on $\bmx$ and suppress the time index.
By the smoothed measure in \Cref{subsec:prelim-measures} and
\eqref{eq:kernel-common}, the change in $nm(u)$ caused by one
particle moving from $w$ to $w+1/n$ is
\[
 H_u(w)=-\frac n\theta\int_0^{\theta/n}\int_0^{1/n}
 K_{\sfb}(w+r+t,u)\,\rd t\,\rd r
 =-\frac1nK_{\sfb}(w,u)+\OO(n^{-2}).
\]
Thus $D(u)=\sum_i e_iH_u(x_i)$ and $|D(u)|\le C$.
Tilting by $e^{tD(u)}$ replaces $\varphi^+$ by
$\varphi^+e^{tH_u}$ and leaves $\varphi^-$ unchanged.
The normalization is nonzero and the total variation bounded on
a fixed disk $|t|\le t_0$. The bounds and logarithm branch on the
working annulus persist since $\cB_t-\cB=\OO(n^{-1})$.

The complex-weight argument in
\cite[Propositions~4.11 and~4.13, Claim~4.17]{Gorin:2022qky}
makes the second-order mean expansion uniform and holomorphic in $t$.
Its leading term is the DMC drift in \eqref{e:dmg}, and
$\left.\partial_t\log\cB_t\right|_{t=0}=\chi H_u$.
The order-$n^{-1}$ coefficient, written in our notation in
\eqref{err:secondder}, changes by $\OO(n^{-1})$: every occurrence
of $\chi$, $\cB'/\cB$ and the Cauchy part of $\log\cB$ has this
bound under the tilt, while the empirical measure is fixed.
Cauchy's formula differentiates the uniform $\OO(n^{-2})$ error.
Consequently,
\[
 \operatorname{Cov}(D(z),D(u))
 =-\frac1{2\pi\ri\theta}\oint_{\omega_-}
 \chi(w)H_u(w)K_{\sfb}(w,z)\,\rd w+\OO(n^{-2}).
\]
Substitute the expansion of $H_u$ and multiply by $n$.
\end{proof}

\begin{proof}[Proof of \Cref{thm:loopstudy}]
For DMC, the fixed-distance assumptions allow us to apply
\cite[Theorem~4.5, (89)]{Gorin:2022qky}, which gives
\eqref{e:dmg}. The covariance estimate \eqref{e:covT}, including its
$\OO(n^{-1})$ remainder, is \Cref{lem:dmc-fixed-covariance}.

For corners, \eqref{e:dmg} and \eqref{e:covT} follow from
\Cref{t:loop_beta_ctscorner,prop:corners-covariance}.
Appendix~A uses the same normalization
$D=D_j^n$ and $\cB=\cB_{s_j}^n$; see
\eqref{e:Dn-definition}--\eqref{e:finite-step-slope}.
Substituting $\cF_{s_j}^n=-\theta\log\cB_{s_j}^n$
and $\chi_{s_j}^n=-\theta/\cB_{s_j}^n$ gives the displayed formulas,
with the flux differentiated through its logarithmic derivative.

For DBM, It\^o's formula gives \eqref{eq:pote_mt};
see \cite[(3.19)]{MR4009708}.
The martingale coefficient gives the first equality in
\eqref{eq:dbm-bracket}, and Cauchy's formula gives the second.
\end{proof}

\section{Sufficient conditions for the analytic assumptions}
\label{sec:appendix-b-verification}

We verify \Cref{ass:weightfunc,ass:prelim-log-branches} and the
contour conditions preceding \Cref{thm:loopstudy} on $[0,T]$.
The deterministic evolution and the separate edge regularity
condition in \Cref{ass:edge-data} are discussed in
\Cref{app:deterministic-proofs}.
Retain the notation of \Cref{subsec:prelim-measures} and
\Cref{def: function F}, and assume $I_0\subset[\ell_0,r_0]$,
where $\ell_0<r_0$ are fixed. Take $T<1$ for corners.
Discrete coefficient statements are at $s\in[0,T]\cap\bZ_n$;
measure bounds also hold under the interpolation in
\Cref{subsec:prelim-measures}. All conclusions hold for sufficiently
large $n$, with constants depending only on $T$, the fixed domains
and the stated strict margins. The coefficient bounds concern fixed
exterior domains; the mesoscopic assumptions in \Cref{sec:main}
remain separate. 

\subsection{Support and mass bounds}
\label{subsec:b-support}

The following bounds allow us to choose common working domains.

\begin{lemma}\label[lemma]{lem:b-bounds}
For DMC, $\mu_s^n$ has mass one and density at most $\theta^{-1}$,
and $I_s\subset[\ell_0,r_0+s]$ at every grid time $s\le T$.
At intermediate times its support lies in
$[\ell_0,r_0+s+n^{-1}]$.
For corners, for every $s\in[0,T]$,
\[
 \mu_s^n(\bR)=1-s,\qquad
 \supp\mu_s^n\subset[\ell_0,r_0],\qquad
 0\le\mu_0^n(({-\infty},x])-\mu_s^n(({-\infty},x])\le s.
\]
These bounds hold for every realization.
\end{lemma}

\begin{proof}
For DMC, the smoothed particle intervals have disjoint interiors,
and each particle moves by $0$ or $n^{-1}$ per step.
For corners, interlacing bounds each successive CDF difference
by $n^{-1}$. Summation and interpolation give the assertions.
\end{proof}

\subsection{Discrete Markov chains}
\label{subsec:b-dmc}

We give sufficient conditions for
\Cref{ass:weightfunc,ass:prelim-log-branches} and the contour conditions
preceding \Cref{thm:loopstudy}, and verify them for the discrete-time
models listed in the introduction. We first treat general analytic
weights, then use the explicit coefficients of the tiling and descending
transitions. Throughout, the bounds concern fixed exterior domains and
are uniform for $s\in[0,T]\cap\bZ_n$ and sufficiently large $n$.

\paragraph{A sufficient condition for general weights.}
Assume \Cref{ass:weightfunc} on a fixed complex neighborhood $\Lambda$
of $[\ell_0,r_0+T]$, with both weights uniformly positive on its real
part. Recall $g_{\sfb}$ from \eqref{s4:kernel-split}.

\begin{proposition}\label[proposition]{prop:b-dmc}
Suppose that, for some $\eta>0$ and every real $x$ in this neighborhood,
\begin{equation}\tag{D}\label{cond:b-D}
  \partial_x\log\frac{\varphi^+(x)}{\varphi^-(x)}
  -\theta\max_{r\in[\ell_0,r_0+T]}
       \partial_x g_{\sfb}(r,x)
  \ge-\frac4\theta+\eta.
\end{equation}
After shrinking $\Lambda$ uniformly, $\cB_s^n$ has a single-valued
holomorphic logarithm on $\Lambda\setminus I_s$. On every fixed exterior
compact set $\mathscr D\Subset\Lambda\setminus[\ell_0,r_0+T]$,
\begin{equation}\label{eq:b-dmc-bounds}
  c\le |\cB_s^n(z)|\le C,
  \qquad z\in\mathscr D,\quad 0\le s\le T,
  \qquad \operatorname{wind}_{\cin}\cB_s^n=0,
\end{equation}
where $\cin\subset\mathscr D$ is any positively oriented enclosing
contour. These conclusions hold for every reachable configuration.
\end{proposition}

\begin{proof}
By \eqref{s4:kernel-split}, \eqref{cond:b-D} bounds the real derivative
of the analytic part of $\log(\varphi^+/\varphi^-)-\theta m_s^n$ from
below by $-4/\theta+\eta$. The mass and density bounds in
\Cref{lem:b-bounds} give
\[
  \theta\int\frac{y}{(r-x)^2+y^2}\,\mu_s^n(\rd r)
  \le 2\arctan\frac{\theta}{2y}
  =\pi-\frac{4y}{\theta}+\OO(y^3),\qquad y>0.
\]
Taylor expansion of the analytic part and Schwarz symmetry therefore
show that, on a sufficiently thin common neighborhood,
\[
  \left|\Im\left(\log\frac{\varphi^+}{\varphi^-}
                         -\theta m_s^n\right)\right|<\pi.
\]
Thus $\cB_s^n/\varphi^-$ avoids $(-\infty,0]$ and is positive on the
real complement of $I_s$. Its principal logarithm, together with a
holomorphic logarithm of $\varphi^-$, gives the required branch.
The strict margin and the fixed support distance give
\eqref{eq:b-dmc-bounds}.
\end{proof}

In particular, \eqref{cond:b-D} is preserved by sufficiently small
holomorphic perturbations with a common strict margin.

\paragraph{Explicit transition weights.}
For the tiling and descending models below, we verify the same
analytic conditions directly from the transition coefficients.
At time $s$, the model prescribes deterministic spatial bounds
\[
  a_s\le x_n(s),\qquad x_1(s)\le b_s-\frac{\theta}{n}.
\]
We call $[a_s,b_s]$ the \textbf{allowed support interval}.
Under the smoothing convention of \Cref{subsec:prelim-measures},
each particle $x_i(s)$ is represented by the interval
$[x_i(s),x_i(s)+\theta/n]$. Thus
\[
  I_s=[x_n(s),x_1(s)+\theta/n]\subset[a_s,b_s].
\]
The endpoints $a_s,b_s$ are part of the model data, whereas $I_s$
depends on the current particle configuration. Their values are
given in each example below.

Recall that $I_0\subset[\ell_0,r_0]$. Since each particle moves by
$0$ or $1/n$ at each step, \Cref{lem:b-bounds} gives
$I_s\subset[\ell_0,r_0+s]$ at every grid time.
We assume that this deterministic support bound stays a fixed
positive distance from both endpoints of the allowed interval:
\begin{equation}\label{eq:b-initial-wall-clearance}
  \inf_{0\le s\le T}
  \min\{\ell_0-a_s,\ b_s-r_0-s\}>0.
\end{equation}
Equivalently, there is a constant $c>0$, independent of $s$ and $n$,
such that
\[
  [\ell_0,r_0+s]\subset[a_s+c,b_s-c],\qquad 0\le s\le T.
\]
We refer to \eqref{eq:b-initial-wall-clearance} as the
\textbf{uniform boundary-separation condition}.
Together with the parameter restrictions below, it gives uniform
positive lower bounds for both transition weights on all reachable
supports and allows us to choose working domains away from their
boundary zeros. We use it as a sufficient condition for these
analytic estimates. All parameters below are fixed, and all strict
inequalities have a uniform positive margin over $[0,T]$.

We use three forms of the coefficients. Common nonvanishing holomorphic
factors in the two weights can be removed, since they do not change
the normalized transition.

\textbf{Affine coefficients:}
\[
  \sfb(z)=z,\qquad
  \varphi_s^{-,n}(z)=z-a_s,\qquad
  \varphi_s^{+,n}(z)=b_s-\theta/n-z.
\]

\textbf{Exponential coefficients:} for $\gamma>0$ and
$D_s=D_0e^{\gamma s}$, with $D_0\in\bR$, take
\begin{equation}\label{eq:b-explicit-exponential-data}
\begin{aligned}
  \sfb_s(z)&=e^{\gamma z}+D_s e^{-\gamma z},\\
  \varphi_s^{-,n}(z)
    &=(1-e^{\gamma(a_s-z)})
      (1-D_s e^{-\gamma(b_s+z)}),\\
  \varphi_s^{+,n}(z)
    &=-e^{\gamma(a_s-b_s+\theta/n)}
      (1-e^{\gamma(b_s-z-\theta/n)})
      (1-D_s e^{-\gamma(a_s+z+\theta/n)}).
\end{aligned}
\end{equation}
We require one of
\[
  D_s\le0,\qquad
  0<D_s<e^{2\gamma a_s},\qquad
  D_s>e^{2\gamma b_s}
\]
to hold throughout $[0,T]$. In the last case, multiply both weights
by $-1$ to make them positive on the reachable supports.

\textbf{Trigonometric coefficients:} take $\gamma>0$ and
\begin{equation}\label{eq:b-explicit-trigonometric-data}
\begin{aligned}
  \sfb_s(z)&=2\cos\bigl(\delta+\gamma(z-s/2)\bigr),\\
  \varphi_s^{-,n}(z)
    &=\sin\frac{\gamma(z-a_s)}2\,
      \sin\left(\delta+\frac{\gamma(z+b_s-s)}2\right),\\
  \varphi_s^{+,n}(z)
    &=\sin\frac{\gamma(b_s-\theta/n-z)}2\,
      \sin\left(\delta+\frac{\gamma(z+a_s+\theta/n-s)}2\right).
\end{aligned}
\end{equation}
We require
$[\delta+\gamma(a_s-s/2),\delta+\gamma(b_s-s/2)]$
to lie in a fixed compact subinterval of $(0,\pi)$.

\begin{corollary}\label[corollary]{cor:b-affine-weights}
Under \eqref{eq:b-initial-wall-clearance} and the stated parameter
restrictions, each of the three coefficient families satisfies
\Cref{ass:weightfunc,ass:prelim-log-branches} on common fixed exterior
domains throughout $[0,T]$. The contour conditions preceding
\Cref{thm:loopstudy} can be chosen uniformly in time.
\end{corollary}

\begin{proof}
The explicit coefficients and boundary separation give the
uniform weight bounds and a common neighborhood on which
$\sfb_s$ is biholomorphic. To verify the logarithm condition,
we use the zero--pole argument of
\cite[Proposition~5.12]{Gorin:2022qky}.

Fix a time and merge adjacent occupied intervals into $k$
components. For affine coefficients,
\[
  \cB_s^n(z)=z-a_s+(b_s-\theta/n-z)
       \prod_i\frac{z-x_i}{z-x_i-\theta/n}.
\]
After cancellation, the numerator has degree $k$, with
leading coefficient $b_s-a_s-\theta-\theta/n>c$.
The endpoint signs give one zero in each occupied component,
matching its pole. Hence there are no other zeros.
For exponential coefficients with $D_s=0$, the same argument
applies after the change of variables $u=e^{\gamma z}$.

For $D_s\ne0$, the rational expression in $u=e^{\gamma z}$
has poles in the occupied intervals and their images under
$u\mapsto D_s/u$. The endpoint signs account for all zeros
except two, which, after removing the finite-$n$ shifts in
the coefficients, are at $u=\pm\sqrt{D_s}$.
Their derivatives are nonzero multiples of
\[
  \int_{a_s}^{b_s}
  \frac{\gamma e^{\gamma r}}{(u-e^{\gamma r})^2}
  \left(1-\theta\frac{\rd\mu_s^n}{\rd r}(r)\right)\rd r.
\]
This is uniformly nonzero: the kernel has a fixed positive
sign for $D_s>0$, and a fixed nonzero imaginary sign for
$D_s<0$. The boundary-separation condition gives a uniform
positive vacancy mass.

For the trigonometric coefficients, the substitution
$u=\tan((\delta+\gamma(z-s/2))/2)$ gives the analogous
zero--pole count, with the two additional zeros at
$u=0,\infty$. Their simplicity follows from the positive
vacancy integrals against $u^{-2}$ and $1$, respectively.
The strict parameter margins and Rouch\'e's theorem preserve
these zero locations under the finite-$n$ shifts.

Thus each occupied component contains equally many zeros
and poles, while the remaining zeros lie outside the working
domain. Consequently, $\cB_s^n$ is nonzero and has a
single-valued holomorphic logarithm on $\Lambda\setminus I_s$.
The fixed-distance upper bound follows directly from the
kernel estimates. Compactness and Hurwitz's theorem give the
uniform lower bound, and the contours can be chosen with
fixed positive separations.
\end{proof}

\paragraph{Time-independent examples.}
\textbf{Nonintersecting Bernoulli random walks} have $\theta=1$,
$\sfb(z)=z$, $\varphi^-=1$ and $\varphi^+=c>0$
\cite[Section~3.1]{Gorin:2022qky}. The same coefficients give the
\textbf{ascending Jack process} for arbitrary $\theta>0$, by the Jack
limit following \cite[Claim~3.10]{Gorin:2022qky}. Both satisfy
\eqref{cond:b-D} directly.

For the \textbf{ascending Macdonald process} with the principal and dual
specializations of \cite[Claim~3.10, (61)--(63)]{Gorin:2022qky}, take
$q_n=e^{-\gamma/n}$ and $t_n=q_n^\theta$, with $\gamma>0$. Then
\[
 \sfb(z)=e^{-\gamma z},\qquad \varphi^-=1,\qquad \varphi^+=c>0.
\]
On a sufficiently thin injective strip,
\[
 \partial_xg_{\sfb}(r,x)
 =\frac{\gamma^2}{4\sinh^2(\gamma(r-x)/2)}
       -\frac1{(r-x)^2}\le0,
\]
with diagonal value $-\gamma^2/12$. Thus \eqref{cond:b-D} holds.
The same proposition covers spatially varying weights satisfying
\eqref{cond:b-D}.

\paragraph{Lozenge tilings.}
Let $L$ be the terminal time in the normalization $s=j/n$. Here
$\theta=1$, $a_s=s-L$ and $b_s=1$. We impose
\eqref{eq:b-initial-wall-clearance} on the prescribed interval $[0,T]$.
The transition coefficients are given by
\cite[Theorem~3.14]{Gorin:2022qky}; the positivity regimes are described
in \cite[Section~1.3]{Gorin:2022qky} and
\cite[Section~2.2]{borodin2010q}.

For \textbf{uniform lozenge tilings}, the coefficients are affine:
\[
 \sfb(z)=z,\qquad
 \varphi_s^{-,n}(z)=z+L-s,\qquad
 \varphi_s^{+,n}(z)=1-n^{-1}-z.
\]
Here $\varphi_s^{-,n}$ vanishes at $a_s=s-L$, whereas
$\varphi_s^{+,n}$ vanishes at $b_s-1/n=1-1/n$. The model restricts
the particle positions to
\[
  s-L\le x_n(s),\qquad x_1(s)\le1-\frac1n,
\]
so the allowed support interval is $[s-L,1]$.
In this case, \eqref{eq:b-initial-wall-clearance} reads
\[
  T<\min\{L+\ell_0,\,1-r_0\}.
\]
Indeed, the two distances in that condition are
$\ell_0+L-s$ and $1-r_0-s$. For $0\le s\le T$ and
$x\in[\ell_0,r_0+s]$, the same inequalities give
\[
  \varphi_s^{-,n}(x)\ge \ell_0+L-T,\qquad
  \varphi_s^{+,n}(x)\ge 1-r_0-T-\frac1n.
\]
Both weights therefore have uniform positive lower bounds for
sufficiently large $n$. The nonvanishing, logarithm and contour
conditions follow from \Cref{cor:b-affine-weights}.

For \textbf{$(q,\kappa)$-distributed lozenge tilings}, take
$q_n=Q^{1/n}$, $Q=e^{-\gamma}$, and $D_s=\kappa^2Q^{-s}$.
After a common sign normalization, the coefficients of
\cite[Theorem~3.14]{Gorin:2022qky} are exactly
\eqref{eq:b-explicit-exponential-data}. In the real regime,
$0<\kappa<Q^L$ gives $0<D_s<e^{2\gamma a_s}$, while
$\kappa>Q^{-1}$ gives $D_s>e^{2\gamma b_s}$. In the imaginary regime,
$\kappa=\ri e^\delta$, with $\delta\in\bR$, gives $D_s<0$.
Thus \Cref{cor:b-affine-weights} applies in both regimes. In the
setting of \cite[Theorem~5.1 and Assumption~5.11]{Gorin:2022qky},
one may instead invoke \cite[Proposition~5.12]{Gorin:2022qky}
for the fixed-domain bounds and logarithm; \eqref{eq:b-initial-wall-clearance} supplies
the uniform weight lower bounds required here.

For \textbf{trigonometric weights}, take
$Q=e^{\ri\gamma}$, $\kappa=e^{\ri\delta}$ and
\[
  \gamma>0,\qquad \gamma L<\delta<\pi-\gamma.
\]
Removing the common phase from the weights and rescaling the
interaction map by a nonzero constant give
\eqref{eq:b-explicit-trigonometric-data}; see
\cite[Section~2.2]{borodin2010q} for the positive trigonometric regime.
The displayed inequalities keep
$\delta+\gamma(x-s/2)$ in a fixed compact subinterval of $(0,\pi)$
for $a_s\le x\le b_s$. The trigonometric case of
\Cref{cor:b-affine-weights} therefore verifies the assumptions.

\paragraph{Descending transitions.}
Let $K$ be the initial rank, $M$ a bound on the number of columns,
and $\lambda'$ the conjugate partition. Use the particle--hole
coordinates of \cite[Section~2.2]{Gorin:2022qky},
\[
 x_i=\frac{-\lambda'_{M+1-i}+\theta^{-1}(M-i)}{M},
 \qquad 1\le i\le M.
\]
At time $s=j/M$, the particle number is $M$ and the interaction
parameter is $\theta^{-1}$. The spatial bounds in these coordinates are
\[
  -K/M+s\le x_M(s),\qquad
  x_1(s)\le\theta^{-1}-\frac{1}{\theta M}.
\]
Each smoothing interval now has length $1/(\theta M)$, so the allowed
support interval is $[-K/M+s,\theta^{-1}]$. Thus
$a_s=-K/M+s$ and $b_s=\theta^{-1}$.
The preceding criteria are used with $n=M$ and interaction parameter
$\theta^{-1}$. We assume uniform bounds on $K/M$ and impose
\eqref{eq:b-initial-wall-clearance} with these endpoints.

For \textbf{descending Jack}, the Jack limit of
\cite[Claim~3.12]{Gorin:2022qky} gives
\[
 \sfb(z)=z,\qquad
 \varphi_s^{-,M}(z)=z+K/M-s,\qquad
 \varphi_s^{+,M}(z)=\theta^{-1}(1-M^{-1})-z.
\]
These are the affine coefficients, so \Cref{cor:b-affine-weights}
applies.

For \textbf{descending Macdonald}, use the principal specializations
of \cite[Claim~3.12, (64)--(66)]{Gorin:2022qky}, with
$q_M=e^{-\gamma/(\theta M)}$ and $t_M=e^{-\gamma/M}$.
Particle--hole duality and the adjacent $q$-Gamma ratio give
\[
\begin{aligned}
 \sfb(z)&=e^{\gamma z},\\
 \varphi_s^{-,M}(z)&=1-e^{\gamma(a_s-z)},\\
 \varphi_s^{+,M}(z)&=e^{\gamma(a_s-z)}
                  -e^{\gamma(a_s-\theta^{-1}+1/(\theta M))}.
\end{aligned}
\]
These are the $D_s=0$ case of
\eqref{eq:b-explicit-exponential-data}, with the stated substitution
of parameters. Hence the same corollary applies.

The \textbf{Koornwinder quasi-branching transitions} are given in
\cite[Claim~3.17]{Gorin:2022qky}. For this principal specialization,
the same particle--hole coordinates and scaling give the coefficients
\eqref{eq:b-explicit-exponential-data}, with
\[
 \sfb_s(z)=e^{\gamma z}+D_s e^{-\gamma z},\qquad
 D_s=D_0e^{\gamma s}.
\]
Indeed, if $v$ is the parameter in \cite[(75)]{Gorin:2022qky}, then
$D_s=q_M^{v+2\theta}$; the shift comes from the hole-coordinate
translation. The quasi-branching parameter change multiplies $D_s$
by $t_M^{-1}$ at each step. The weights follow from the adjacent
$q$-Gamma ratio, as for descending Macdonald, with the two additional
factors in \eqref{eq:b-explicit-exponential-data}.
Under
\[
 D_s\le0,\qquad
 0<D_s<e^{2\gamma a_s},\qquad\text{or}\qquad
 D_s>e^{2\gamma/\theta},
\]
with one alternative holding throughout $[0,T]$ and the stated strict
margins, \Cref{cor:b-affine-weights} gives the required conditions.

\begin{remark}\label[remark]{rem:b-time-dependent-map}
For the interaction map $\sfb_s$ is time-dependent examples, the preceding verifications apply
to $\sfb_s$ and $\varphi_s^{\pm,n}$ at each time, with common bounds
throughout $[0,T]$. When $\sfb_s$ is $C^2$ in time, the transform
evolution has the additional term
\[
 \int\partial_s\left[
   \frac{\sfb_s'(z)}{\sfb_s(r)-\sfb_s(z)}
 \right]\mu_s^n(\rd r).
\]
In discrete time, this is obtained by expanding the change of the
kernel over one step. Applying \Cref{s:eost,sec:ale-convergence}
to these examples requires this term and the corresponding uniform
time-regularity estimates. Time-dependent weights with a fixed
interaction map do not contribute a kernel derivative.
\end{remark}

\subsection{Weighted \texorpdfstring{$\beta$}{beta}-corners}
\label{subsec:b-corners}

We verify the analytic assumptions for the weighted Dixon--Anderson
transitions of \cite[Section~3.5 and Corollary~3.6]{Gorin:2022qky}.
We retain the normalization \eqref{eq:corner-level-measures} and
coefficient \eqref{e:B_function}. The first condition is uniform over
initial configurations in a fixed interval. The second uses the
limiting initial measure to obtain the same analytic bounds.

Fix $0<T<1$ and assume $I_0\subset[\ell_0,r_0]$.
By \Cref{lem:b-bounds}, interlacing gives
\[
  I_s\subset[\ell_0,r_0],\qquad
  \mu_s^n(\bR)=1-s\ge1-T,
  \qquad s\in[0,T]\cap\bZ_n.
\]
Thus one interval contains every subsequent row, and the mass stays
bounded away from zero. Assume that $W$ is holomorphic on a complex
neighborhood of this interval and satisfies Schwarz symmetry.
Recall that
\[
  \cB_s^n(z)=\left(1-s-\frac1n\right)W'(z)-\theta m_s^n(z).
\]
All bounds in this subsection concern fixed exterior domains.

\paragraph{Sufficient conditions.}
\begin{proposition}\label[proposition]{prop:b-corners}
Suppose that either of the following conditions holds.
\begin{enumerate}[label=\textup{(\roman*)},leftmargin=*,itemsep=4pt]
  \item For some $\eta>0$,
  \begin{equation}\label{eq:b-corners-derivative}
    \frac{\rd}{\rd x}\left[(x-\ell_0)(x-r_0)W'(x)\right]
    \ge-\theta+\eta,
    \qquad \ell_0\le x\le r_0.
  \end{equation}

  \item The initial measures satisfy $\mu_0^n\Rightarrow\mu_0$ in
  probability for a deterministic probability measure $\mu_0$.
  There is a fixed annulus $\mathscr D$ surrounding $[\ell_0,r_0]$
  at a positive distance, with its closure and enclosed region in
  the analytic domain of $W$, such that, for
  $\cB_0=W'-\theta m_0$,
  \begin{equation}\tag{C$_0$}\label{cond:b-C0}
    \operatorname{wind}_{\cin}\cB_0=-1,
    \qquad
    \inf_{z\in\mathscr D}\left(
      |\cB_0(z)|-\theta T\int_{\ell_0}^{r_0}
                       \frac{\rd r}{|r-z|^2}
    \right)>0,
  \end{equation}
  where $\cin\subset\mathscr D$ is a positively oriented contour
  enclosing $[\ell_0,r_0]$.
\end{enumerate}
Under \textup{(i)}, one can choose such a fixed annulus
$\mathscr D$; under \textup{(ii)}, use the specified annulus.
There are $c,C>0$ and nested positively oriented enclosing contours
$\cin,\Gamma\subset\mathscr D$ such that
\begin{equation}\label{eq:b-corners-bounds}
  c\le|\cB_s^n(z)|\le C,\qquad z\in\mathscr D,
  \qquad
  \operatorname{wind}_{\cin}\cB_s^n
  =\operatorname{wind}_{\Gamma}\cB_s^n=-1,
\end{equation}
uniformly for $s\in[0,T]\cap\bZ_n$ and sufficiently large $n$.
For every $e\in I_s$, $(z-e)\cB_s^n(z)$ has a single-valued
holomorphic logarithm on $\mathscr D$.
These conclusions hold for every realization under \textup{(i)}.
Under \textup{(ii)}, they hold on events
$\mathcal G_n\in\mathscr F_0$ with $\bP(\mathcal G_n)\to1$,
simultaneously for all subsequent realizations.
Consequently, \Cref{ass:weightfunc,ass:prelim-log-branches} and the
contour conditions preceding \Cref{thm:loopstudy} hold on these
fixed domains, with the same event convention.
\end{proposition}

\begin{proof}
For \textup{(i)}, decrease $\eta$ so that $0<\eta<\theta$.
For $x\in[\ell_0,r_0]$ and $\Im z>0$,
\[
  \frac{\Im\bigl[(z-\ell_0)(z-r_0)/(z-x)\bigr]}{\Im z}
  =1+\frac{(x-\ell_0)(r_0-x)}{|z-x|^2}\ge1.
\]
Using the mass $1-s$, Schwarz symmetry and
\eqref{eq:b-corners-derivative}, we obtain on a sufficiently thin
common neighborhood
\[
  \frac{\Im[(z-\ell_0)(z-r_0)\cB_s^n(z)]}{\Im z}
  \ge\left(1-s-\frac1n\right)
       \left(-\theta+\frac\eta2\right)+\theta(1-s)
  \ge c_T>0.
\]
At the endpoints, the same condition gives
\[
  W'(\ell_0)\le\frac{\theta-\eta}{r_0-\ell_0},
  \qquad
  W'(r_0)\ge\frac{-\theta+\eta}{r_0-\ell_0}.
\]
Together with the Stieltjes integral, these inequalities give
$\cB_s^n<-c_T$ at a fixed point just left of $\ell_0$ and
$\cB_s^n>c_T$ at a fixed point just right of $r_0$.
The imaginary-part estimate and these real bounds give a fixed
annulus on which $\cB_s^n$ is uniformly nonzero.
Replacing $W$ by $uW$, $0\le u\le1$, preserves these bounds.
The winding therefore equals that of $-\theta m_s^n$, namely $-1$.

For \textup{(ii)}, interlacing and the mass bounds give
\[
  \sup_x\left|
    \frac{\mu_s^n(({-\infty},x])}{1-s}
    -\mu_0^n(({-\infty},x])
  \right|\le s.
\]
Since the two normalized measures have equal mass, integration by
parts yields, for $z\in\mathscr D$,
\[
  \left|\frac{\cB_s^n(z)}{1-s}-\cB_0(z)\right|
  \le\theta T\int_{\ell_0}^{r_0}\frac{\rd r}{|r-z|^2}
     +\frac{C}{n(1-T)}
     +\theta\sup_{\mathscr D}|m_0^n-m_0|.
\]
The final two terms tend to zero in probability, since the initial
supports are uniformly bounded and $\mathscr D$ stays a fixed
distance from them. On initial events of probability tending to
one, \eqref{cond:b-C0} therefore gives a uniform positive lower
bound and preserves the winding of $\cB_0$, for every subsequent
row and all $s\le T$.

In both cases, the upper bound follows from the fixed support
distance. Multiplication by $z-e$ changes the winding to zero,
giving the asserted logarithm. The separated contours are chosen
as in \Cref{cor:b-contours}.
\end{proof}

\paragraph{Initial measure condition.}
The following endpoint test provides a nonvanishing annulus for
\eqref{cond:b-C0}. Suppose that $\mu_0$ has support
$[\ell_0,r_0]$ and a continuous density $\rho_0$ that is strictly
positive on $(\ell_0,r_0)$. Assume that $m_0$ extends continuously from each
half-plane to this interval, with finite real endpoint values.
If
\begin{equation}\label{eq:b-corners-endpoint-signs}
  W'(\ell_0)<\theta m_0(\ell_0),\qquad
  W'(r_0)>\theta m_0(r_0),
\end{equation}
then a fixed annulus $\mathscr D$ can be chosen so that
\eqref{cond:b-C0} holds for all sufficiently small fixed $T>0$.

Indeed, for $0\le u\le1$ and $\ell_0<x<r_0$,
\[
  \Im\bigl[uW'(x)-\theta m_0(x+\ri0)\bigr]
  =-\theta\pi\rho_0(x)<0.
\]
Also $m_0(\ell_0)>0>m_0(r_0)$, so
\eqref{eq:b-corners-endpoint-signs} keeps the left endpoint value
negative and the right endpoint value positive throughout this
homotopy. Continuity and compactness give a common thin
nonvanishing annulus. The winding is that of $-\theta m_0$, hence
$-1$. On a smaller fixed annulus, $\inf|\cB_0|>0$ and the integral
in \eqref{cond:b-C0} is uniformly bounded, which proves the claim.
For a prescribed time horizon, the inequality in
\eqref{cond:b-C0} determines whether this annulus can be used.

\paragraph{Examples.}
The following examples concern the transition potential $W$.
Under the derivative condition, the proposition applies for every
fixed $T<1$. For the initial-measure condition, $T$ must satisfy
\eqref{cond:b-C0}.

\textbf{Unweighted $\beta$-corners.}
For $W=0$, \eqref{eq:b-corners-derivative} holds with
$\eta=\theta/2$. Thus the proposition applies to every initial
configuration in the fixed interval and every fixed $T<1$.

\textbf{Gaussian $\beta$-corners.}
The top row has Gaussian $\beta$-ensemble distribution, while
the conditional transitions have $W=0$; see
\cite[Definitions~1.2--1.3]{GorinXuZhang2024Airy}.
By \cite[Theorem~7.9(1)]{charofairy}, the normalized top row
has compact support with probability tending to one.
Interlacing gives the same bound for all subsequent rows, so \Cref{prop:b-corners} applies on this event. 

\textbf{Linear potentials.}
For $W(x)=\lambda x$, $\lambda\in\bR$,
\eqref{eq:b-corners-derivative} holds whenever
$|\lambda|(r_0-\ell_0)<\theta$.

\textbf{Quadratic potentials.}
On $[\ell_0,r_0]=[e-R_0,e+R_0]$, the potential
$W(x)=\lambda(x-e)^2/2$ satisfies
\eqref{eq:b-corners-derivative} if
\[
  -\frac{\theta}{2R_0^2}<\lambda<
  \frac{\theta}{R_0^2}.
\]
The initial-measure condition gives a larger range.
For $I_0\subset[-1,1]$ and $\mu_0^n\Rightarrow\mu_0$
in probability, where
\[
  \mu_0(\rd x)=\frac2\pi\sqrt{1-x^2}\,
  \mathbf1_{[-1,1]}(x)\,\rd x,
\]
we have $m_0(-1)=2$ and $m_0(1)=-2$.
Thus \eqref{eq:b-corners-endpoint-signs} holds for
$W(x)=\lambda x^2/2$ whenever $\lambda>-2\theta$,
giving \eqref{cond:b-C0} for sufficiently small fixed $T>0$.

\textbf{Laguerre potentials.}
Let $0<\ell_0<r_0$ and
$W(x)=\alpha\log x-\lambda x$, with $\alpha,\lambda\ge0$,
on an analytic neighborhood avoiding zero.
The derivative in \eqref{eq:b-corners-derivative} is concave,
so checking its endpoint values shows that the condition holds if
\[
  (r_0-\ell_0)\max\left\{
    \frac{\alpha}{\ell_0}-\lambda,\,
    \lambda-\frac{\alpha}{r_0}
  \right\}<\theta.
\]
For initial measures satisfying the endpoint test, the alternative
conditions are
\[
  \frac{\alpha}{\ell_0}-\lambda<\theta m_0(\ell_0),
  \qquad
  \frac{\alpha}{r_0}-\lambda>\theta m_0(r_0),
\]
which imply \eqref{cond:b-C0} for sufficiently small fixed $T>0$.

\textbf{Jacobi potentials.}
Let $0<\ell_0<r_0<1$ and
$W(x)=\alpha\log x+\nu\log(1-x)$, with $\alpha,\nu\ge0$,
on an analytic neighborhood avoiding $0$ and $1$.
The same endpoint calculation gives
\[
  (r_0-\ell_0)\max\left\{
    \frac{\alpha}{\ell_0}-\frac{\nu}{1-\ell_0},\,
    \frac{\nu}{1-r_0}-\frac{\alpha}{r_0}
  \right\}<\theta.
\]
Alternatively, under the initial-measure regularity of the endpoint
test, it suffices that
\[
  \frac{\alpha}{\ell_0}-\frac{\nu}{1-\ell_0}
       <\theta m_0(\ell_0),\qquad
  \frac{\alpha}{r_0}-\frac{\nu}{1-r_0}
       >\theta m_0(r_0)
\]
for sufficiently small fixed $T>0$.

\textbf{Polynomial potentials and analytic perturbations.}
For the semicircle initial measure on $[-1,1]$ above,
$W(x)=\lambda x^4/4$ satisfies the endpoint test for every
$\lambda\ge0$. More generally, any real polynomial $W$ with
\[
  W'(-1)<2\theta,\qquad W'(1)>-2\theta
\]
satisfies \eqref{cond:b-C0} for sufficiently small fixed $T>0$.
Both criteria in \Cref{prop:b-corners} are preserved by sufficiently
small Schwarz-symmetric holomorphic perturbations of $W$.

The derivative criterion also gives the potential condition in
\Cref{thm:app-edge-formation}.

\subsection{Dyson Brownian motion and related continuous-time chains}
\label{subsec:b-dbm}

We verify the fixed-domain analytic assumptions for Dyson Brownian
motion and the hypotheses of the dynamical loop equations for the
related Poisson chains. For DBM, the logarithm is explicit, so a common
compact support bound suffices. For the jump processes, we use the
generator identity of \cite[Section~3.2]{Gorin:2022qky} directly.
The quantum-walk examples at the end are identified with these
Poisson or Bernoulli transitions.

\paragraph{A sufficient condition.}
We retain the potential assumptions and the convention for solutions
of \Cref{subsec:prelim-dbm-model}. Fix $T>0$ and a deterministic
compact interval $[-R,R]$ in the analytic domain of $V$.
We work on the event that $I_s\subset[-R,R]$ for $0\le s\le T$,
or before the corresponding exit time $\sigma_R$ from
\Cref{subsec:s4-upper}.

\begin{proposition}\label[proposition]{prop:b-dbm}
Assume \Cref{ass:weightfunc} in the preceding DBM setting.
On the support event above, $\cB_s^n$ has the holomorphic logarithm
\[
  \log\cB_s^n(z)
  =\frac\theta2\left(m_s^n(z)+\frac{V'(z)}2\right)^2
\]
on $\Lambda\setminus I_s$. For every fixed compact exterior domain
$\mathscr D\Subset\Lambda\setminus[-R,R]$, there are constants
$c,C>0$, independent of $s$ and $n$, such that
\[
  c\le |\cB_s^n(z)|\le C,
  \qquad z\in\mathscr D,\quad 0\le s\le T.
\]
Thus \Cref{ass:prelim-log-branches} holds on fixed exterior annuli.
The same conclusions hold up to $T\wedge\sigma_R$, and throughout
$[0,T]$ with probability at least $1-Ce^{-cn}$ whenever
\eqref{s4:dbm-confinement} applies.
\end{proposition}

\begin{proof}
By \eqref{e:B_function}, the displayed expression is a holomorphic
logarithm of $\cB_s^n$. Since the empirical measure has mass one,
\[
  |m_s^n(z)|\le\frac1{\dist(z,[-R,R])},\qquad z\in\mathscr D.
\]
This bound and boundedness of $V'$ on $\mathscr D$ give
$|\log\cB_s^n|\le C$. Exponentiating its real part gives the
claimed upper and lower bounds.
\end{proof}

The DBM loop equation is given by
\cite[Corollary~3.3]{Gorin:2022qky}; its form in the normalization
of this paper is recorded in \Cref{thm:loopstudy}.

\paragraph{Examples.}
\textbf{Dyson Brownian motion without pontential.}
For $V=0$, the logarithm reduces to $\theta(m_s^n)^2/2$.
The proposition therefore applies on every common compact support
event, including for the process started from zero.

\textbf{Quadratic potentials and Gaussian Dyson Brownian motion.}
For $V(x)=\lambda x^2/2$, with $\lambda>0$,
\[
  \log\cB_s^n(z)
  =\frac\theta2\left(m_s^n(z)+\frac{\lambda z}{2}\right)^2.
\]
The same verification applies to both deterministic and random
initial data on the stated support event. In particular, it includes
the stationary Gaussian, or Hermite, process.

\textbf{General analytic potentials.}
Every potential satisfying the assumptions of
\Cref{subsec:prelim-dbm-model,ass:weightfunc} is covered on the
same support event. For example,
\[
  V(x)=\frac g4x^4+\frac a2x^2,
  \qquad g>0,\quad a\in\bR,
\]
is entire and satisfies $V''(x)=3gx^2+a\ge a$ and
$xV'(x)=gx^4+ax^2\ge-|a|x^2$.
Thus this example satisfies the global potential conditions,
including when $a<0$. The fixed-domain verification
uses only analyticity and the support bound.

\paragraph{Exterior contours.}
\begin{corollary}\label[corollary]{cor:b-contours}
Under the preceding sufficient conditions, the contours
$\cin,\Gamma$ and spectral domains required in
\Cref{thm:loopstudy} can be chosen with uniform positive
separations for $s\le T$, on the events or up to the
stopping times specified for each model.
\end{corollary}

\paragraph{Related continuous-time chains.}
We verify the analytic hypotheses for the Poisson processes of
\cite[Section~3.2]{Gorin:2022qky}. In our spatial normalization,
a particle $x_i$ jumps to $x_i+1/n$ at rate
\[
  n\varphi(x_i)
  \prod_{j\ne i}
  \left(1+\frac{\theta}{n(x_i-x_j)}\right).
\]
Assume that $\varphi$ is positive on the real particle domain
and holomorphic on a common complex neighborhood of the particle
positions and their $1/n$-translates. On a common compact support
event, these conditions hold uniformly in time on fixed working
domains. The generator loop equation then follows from the
residue argument of \cite[Corollary~3.2]{Gorin:2022qky}.

\textbf{$\theta$-nonintersecting Poisson random walks.}
For $\varphi=1$, the analytic condition holds automatically.
The dynamical loop equation is given by
\cite[Corollary~3.2]{Gorin:2022qky}.
For compact initial supports, the particle bound of
\cite[Proposition~3.6]{huang2017beta}, together with monotonicity
of the paths, gives a common compact support event throughout
$[0,T]$ with probability at least $1-Ce^{-cn}$.

\textbf{Spatially dependent Poisson rates.}
For positive analytic $\varphi$, the same loop equation applies
on the preceding domains. If $\varphi\le C$ on $\bR$, comparison
with the homogeneous process at speed $C$ gives the same support
bound, by the coupling argument of
\cite[Proposition~3.6]{huang2017beta}.
For example, $\varphi(x)=1+\varepsilon\cos x$,
$|\varepsilon|<1$, satisfies these conditions.

\textbf{Quantum random walks.}
At $\theta=1$, the Pieri rule identifies the
defining-representation walk of
\cite[Section~5.2]{bufetov2019fourier}
with the embedded homogeneous Poisson chain.
Its transition probability is
\[
  \frac1n\prod_{j\ne i}
  \left(1+\frac1{n(x_i-x_j)}\right).
\]
Thus its Poissonization has the loop equation verified above.
Similarly, the binomial mixture of exterior-power transitions
with weights $\binom nk c^k/(1+c)^n$ gives the Bernoulli
coefficients $\varphi^-=1$, $\varphi^+=c$.
These satisfy \eqref{cond:b-D}, as verified in
\Cref{subsec:b-dmc}.

\section{Evolution and formation of regular edges}
\label{app:deterministic-proofs}

We prove Proposition~\ref{prop:prelim-deterministic-limit} and give sufficient
conditions on the initial measure for a regular square-root edge
to form and persist. We also compute waiting times in the scalar
cases and establish the time regularity and uniformity needed to
verify the deterministic edge assumptions in
\Cref{s:eost,sec:ale-convergence}.
Throughout, we retain the notation of
\Cref{subsec:prelim-measures,subsec:prelim-characteristics}.

\subsection{The deterministic evolution}
\label{subsec:app-deterministic-limit}

We first treat the discrete models using the one-step estimates of
\Cref{thm:loopstudy}. For DBM, a comparison between different values
of $\beta$ extends the law of large numbers to every $\beta>0$.

\begin{lemma}[Deterministic limit for the discrete models]
\label[lemma]{lem:discrete-deterministic-limit}
Under the assumptions of Proposition~\ref{prop:prelim-deterministic-limit},
its conclusions hold for DMC and weighted $\beta$-corners.
The transforms and their spatial derivatives of every fixed order
converge locally uniformly on $\Lambda\setminus I$, uniformly in
time and in probability. For DMC, $\mu_s\le\theta^{-1}\mathrm dr$.
\end{lemma}

\begin{proof}
We follow the compactness and martingale argument of
\cite[Section~5.5, (243)--(245)]{Gorin:2022qky}.
Fix a time horizon $T$ covered by
Proposition~\ref{prop:prelim-deterministic-limit}, and linearly interpolate
the empirical measures between consecutive grid times.
Bounded particle speed for DMC and interlacing for corners give
$d_{\mathrm{BL}}(\mu_s^n,\mu_t^n)\le(1+|I|)|s-t|$, where
$d_{\mathrm{BL}}$ is the bounded-Lipschitz metric.
The interpolation changes the original paths by $\OO(n^{-1})$
in this metric. Together with the common compact support from
\Cref{lem:b-bounds}, this gives compactness of the empirical paths.

At each fixed spectral point outside the common contour,
\eqref{e:covT}, with $u=\bar z$, bounds the conditional second
moment of $\Delta\cM_j^n(z)$ by $C_z/n$.
Doob's inequality therefore gives
\begin{equation}\label{eq:discrete-limit-martingale}
 \bE\max_{\ell\le\lfloor nT\rfloor}
 \left|\frac1n\sum_{j<\ell}\Delta\cM_j^n(z)\right|^2
 \le \frac{C_{z,T}}{n^2}.
\end{equation}
Summing \eqref{e:dmg}, the accumulated drift error is
$\OO(n^{-1})$. Hence every subsequential limit satisfies
\eqref{e:burger}, with the initial measure prescribed in
Proposition~\ref{prop:prelim-deterministic-limit}.
The coefficient bounds in \Cref{ass:prelim-log-branches}
justify passage to the logarithmic derivatives on the fixed
contours. The support, mass, and DMC density bounds of
\Cref{lem:b-bounds} pass to weak limits.
The logarithm conditions pass by local uniform convergence
on the prescribed domains and, for DMC, Hurwitz's theorem.

For uniqueness, we use the analytic-scale argument in
\cite[Chapter~III, Section~2, proof of Theorem~2.1]{Nishida1978}.
Let $f_s=m_s-\widetilde m_s$ be the difference of two candidate
transforms that agree at time $a$.
Choose parallel smooth contours $\Gamma_r$, $0\le r\le\rho$,
increasing outwards in a fixed collar in $\mathscr D$ enclosing
$\omega_-$, and write $\|f\|_r=\sup_{\Gamma_r}|f|$.
For a short time after $a$, the coefficient lower bounds in
\Cref{ass:prelim-log-branches} and continuity imply that
$\cB_s/\widetilde\cB_s$ stays near $1$ on this collar. Thus
\(
 \sup_{\Gamma_r}
 \left|\operatorname{Log}\frac{\cB_s}{\widetilde\cB_s}\right|
 \le C\|f_s\|_r.
\)
For corners, using the quotient avoids taking a logarithm of
either coefficient separately.
By \Cref{def: function F} and \eqref{e:F_CP} with $n^{-1}=0$,
subtracting \eqref{e:burger} and integrating by parts yields
\begin{equation}\label{eq:discrete-limit-uniqueness}
 \|f_t\|_R
 \le \frac{C}{R-r}\int_a^t\|f_s\|_r\,\mathrm ds,
 \qquad 0\le r<R\le\rho.
\end{equation}
Here the factor $(R-r)^{-1}$ follows from the Cauchy kernel
estimate between the two contours.
Iterating over $k$ equally spaced contours gives
\[
 \|f_t\|_\rho
 \le
 \frac{[Ck(t-a)/\rho]^k}{k!}
 \sup_{a\le s\le t}\|f_s\|_0.
\]
The right-hand side tends to zero for $t-a<\rho/(eC)$.
The transforms therefore agree on $\Gamma_\rho$ and, by
analytic continuation, throughout their common domain.
The uniform bounds allow this argument to be restarted,
giving uniqueness on $[0,T]$.
Equality of the transforms determines the measures, since
$m_s(z)/\sfb'(z)$ is the Cauchy transform of
$\sfb_*\mu_s$ evaluated at $\sfb(z)$.

Every subsequential limit is consequently the same deterministic
path, which proves convergence in probability.
Uniform kernel bounds and Cauchy's formula give the asserted
local uniform convergence of the transforms and their derivatives.
\end{proof}

When the working bounds hold on the initial-data events
$\mathcal G_n$ of \Cref{prop:b-corners}, apply the preceding
argument conditionally on $\mathcal G_n$.
Its constants are uniform on these events, and
$\mathbb P(\mathcal G_n^c)\to0$ gives the same conclusion.

\begin{lemma}[Comparison of DBM parameters]
\label[lemma]{lem:dbm-beta-comparison}
Under the DBM hypotheses of
Proposition~\ref{prop:prelim-deterministic-limit}, let $\beta>0$ and
$T<\infty$. Couple the solutions $\lambda^\beta$ and $\lambda^1$
of \eqref{eq:DBM-potential} with the same initial configuration
and Brownian motions. Then
\[
 \bE\sup_{0\le s\le T}
 \frac1n\sum_{i=1}^n
 |\lambda_i^\beta(s)-\lambda_i^1(s)|^2
 \le \frac{C_{\beta,K,T}}n,
\]
where $V''\ge-2K$.
The constant is independent of the initial configuration,
and the bound also holds for a common random initial configuration.
\end{lemma}

\begin{proof}
We first justify the calculation for $0<\beta<1$, when collisions
may occur.
For globally Lipschitz drift, the ordered strong solution is
given by \cite[Theorem~3.1]{CepaLepingle1997}.
Its inverse gaps are locally time-integrable, collision times
have Lebesgue measure zero, and there is no boundary local-time
term; see \cite[Lemmas~3.3, 3.5 and~3.8]{CepaLepingle1997}.
For the locally Lipschitz drift $-V'/2$, the same statements
follow by localization. Indeed, It\^o's formula gives
\[
 \mathrm d\left(\frac1n\sum_i\lambda_i^2\right)
 =
 \mathrm dL_s+
 \left[
   \frac2{\beta n}+\frac{n-1}{n}
   -\frac1n\sum_i\lambda_iV'(\lambda_i)
 \right]\mathrm ds,
\]
where $L$ is a local martingale.
The lower bound on $V''$ implies
$xV'(x)\ge-C(1+x^2)$, so the stopped moment estimate rules
out explosion.

Write
$\Delta_i=\lambda_i^\beta-\lambda_i^1$,
$R=n^{-1}\sum_i\Delta_i^2$, and
$c_\beta=(\sqrt{2/\beta}-\sqrt2)^2$.
Pairing the interaction terms shows that their contribution
to $\mathrm dR$ is
\[
 -\frac2{n^2}\sum_{i<j}
 \frac{
   [(\lambda_i^\beta-\lambda_j^\beta)
    -(\lambda_i^1-\lambda_j^1)]^2
 }{
   (\lambda_i^\beta-\lambda_j^\beta)
   (\lambda_i^1-\lambda_j^1)
 }
 \,\mathrm ds
 \le0.
\]
The inverse-gap integrability justifies this identity at almost
every time, including when $\beta<1$.
Using $V''\ge-2K$, It\^o's formula therefore yields
\[
 \mathrm dR
 \le \mathrm dM+2KR\,\mathrm ds+\frac{c_\beta}{n}\,\mathrm ds,
 \qquad
 \mathrm d\langle M\rangle_s
 =\frac{4c_\beta}{n^2}R_s\,\mathrm ds.
\]
Since $R_0=0$, localization followed by the
Burkholder--Davis--Gundy, Young, and Gr\"onwall inequalities gives
$\bE\sup_{s\le T}R_s\le C_{\beta,K,T}/n$.
Fatou's lemma removes the stopping.
Conditioning on the common initial configuration proves
the assertion for random initial data.
\end{proof}

\begin{proof}[Proof of Proposition~\ref{prop:prelim-deterministic-limit}]
For DMC and weighted $\beta$-corners, the result follows from
\Cref{lem:discrete-deterministic-limit}.
For DBM with $\beta=1$, the law of large numbers follows from
\cite[Theorem~2.4]{MR4009708}, since compact support implies
the required logarithmic moment bound.
For random initial data, condition on the initial configurations
and apply the deterministic-initial-data result along subsequences
on which the initial convergence and compact containment hold
almost surely.

For general $\beta>0$, use the coupling of
\Cref{lem:dbm-beta-comparison}.
Writing $\mu_s^{n,\beta}$ for the corresponding empirical measure,
Cauchy--Schwarz gives
\[
 \sup_{0\le s\le T}
 d_{\mathrm{BL}}(\mu_s^{n,\beta},\mu_s^{n,1})^2
 \le
 \sup_{0\le s\le T}
 \frac1n\sum_{i=1}^n
 |\lambda_i^\beta(s)-\lambda_i^1(s)|^2.
\]
The right-hand side converges to zero in probability by
\Cref{lem:dbm-beta-comparison}.
Thus the same law of large numbers, with the same deterministic
McKean--Vlasov limit, holds for every fixed $\beta>0$.
\end{proof}

\begin{remark}
For \textit{uniform lozenge tilings} and the \textit{descending Jack}
and \textit{descending Macdonald transitions} in particle--hole coordinates,
the weights $\varphi_s^{\pm,n}(z)$ in \eqref{eq:disc-mc} depend on $s$.
The interaction map is fixed, with $\sfb(z)=z$ for uniform tilings and
descending Jack, and $\sfb(z)=e^{\gamma z}$ for descending Macdonald.
For \textit{$(q,\kappa)$-distributed lozenge tilings} in the real and
imaginary regimes and for \textit{Koornwinder quasi-branching transitions},
the interaction map is $\sfb_s(z)=e^{\gamma z}+D_s e^{-\gamma z}$
with $D_s=D_0e^{\gamma s}$.
For \textit{trigonometric lozenge tilings}, it is
$\sfb_s(z)=2\cos(\delta+\gamma(z-s/2))$.

The convergence statement of Proposition~\ref{prop:prelim-deterministic-limit}
holds for all these models and for the \textit{$\theta$-nonintersecting
Poisson walks}, including those with spatially dependent rates, in
\Cref{subsec:b-dbm}. We retain the initial convergence, common compact
support, fixed domains and contours of that proposition and the analytic
bounds in \Cref{subsec:b-dmc,subsec:b-dbm}. For the discrete chains, we
assume local uniform convergence of the coefficients and their first two
time derivatives, with common analytic bounds. When $\sfb_s$ varies with
time, the limiting transform equation contains the additional term in
\Cref{rem:b-time-dependent-map}.
\end{remark}

\subsection{Sufficient conditions for the formation of a square-root edge}
\label{subsec:app-edge-formation}

We give sufficient conditions on various models and the initial measure
for the deterministic right edge to satisfy
\Cref{ass:edge-data}\textup{(i)--(ii)} after a waiting time.
These are the deterministic edge conditions used in
\Cref{s:eost,sec:ale-convergence}. The initial measure may have an
irregular edge.
We first prove a general criterion for edge formation and persistence.
We then check it for the initial measures and models in the introduction
and \Cref{sec:appendix-b-verification}. Finally, we compute exact waiting
times when the characteristic equations have explicit solutions.

\subsubsection{Criteria for square-root edge formation}

We use the deterministic flows of \Cref{subsec:app-deterministic-limit}
on $[0,T]$, with the support and analytic bounds stated there.
Write $E_s=\sup\supp\mu_s$ and let $\widehat m_s$ be the
ordinary Stieltjes transform in \eqref{s4:ordinary-transforms}.
For DMC and Poisson, we assume $\mu_0\le\theta^{-1}\rd x$.
The particle spacing gives the same bound for $\mu_s$.
For DMC, sufficient conditions for the fixed-domain estimates are given
in \Cref{prop:b-dmc,cor:b-affine-weights}. For corners, take $T<1$.
For DBM, short-time persistence of an initially regular square-root
edge was proved in \cite[Proposition~3.4]{adhikari2020dyson}.

\begin{theorem}[Formation of a regular edge]
\label{thm:app-edge-formation}
Under the preceding assumptions, the following conditions are sufficient
for the right edge to be regular. For corners, assume
\eqref{eq:b-corners-derivative} in parts~\textup{(i)} and~\textup{(ii)}.

\begin{enumerate}[label=\textup{(\roman*)},leftmargin=*,itemsep=4pt]
\item \textbf{Regularity at positive times.}
Suppose that
\begin{equation}\label{eq:app-initial-edge-condition}
 \begin{cases}
 \displaystyle\lim_{u\downarrow E_0}\widehat m_0'(u)=\infty,
     &\textnormal{DBM and Poisson},\\[3pt]
 \displaystyle\lim_{u\downarrow E_0}
       e^{\theta\widehat m_0(u)}\widehat m_0'(u)=\infty,
     &\textnormal{DMC},\\[3pt]
 \displaystyle\lim_{u\downarrow E_0}
       \frac{\widehat m_0'(u)}{\widehat m_0(u)^2}=\infty,
     &\textnormal{corners}.
 \end{cases}
\end{equation}
Then the conclusions below hold on $J=[\varepsilon,T]$ for every
$0<\varepsilon<T$.

\item \textbf{Regularity after a waiting time.}
Let $r_0$ be the upper support bound at time zero. There are
$h_0,C,\delta>0$ such that, if some $u_0\in(E_0,r_0+\delta)$ satisfies
\begin{equation}\label{eq:app-finite-compression}
 -\chi_0'(u_0)>C/h_0,\qquad s_*=\frac{C}{-\chi_0'(u_0)},
\end{equation}
then the same conclusions hold on $J=[\varepsilon,T]$ for every
$s_*<\varepsilon<T$. The constants $h_0,C,\delta$ depend only on the
coefficient, domain and support bounds, and for corners also on $1-T$
and the strict margin in \eqref{eq:b-corners-derivative}.

\item \textbf{Short-time regularity for corners.}
In place of \eqref{eq:b-corners-derivative}, assume \eqref{cond:b-C0} and
\begin{equation}\label{eq:app-corners-local-condition}
 m_0(E_0+)>-\infty,\qquad m_0'(E_0+)=\infty,\qquad \cB_0(E_0)>0.
\end{equation}
There is a $T_1\in(0,T]$ such that the conclusions hold on
$J=[\varepsilon,T_1]$ for every $0<\varepsilon<T_1$. If $m_0(E_0+\zeta)$
already has a convergent expansion in $\sqrt\zeta$ with positive
linear coefficient, they hold on $J=[0,T_1]$.
\end{enumerate}

On each such interval $J$, there are $r,c,C>0$ such that, uniformly
for $s\in J$,
\begin{equation}\label{eq:app-generated-edge}
 \rho_s(E_s-x)=\frac{\sqrt{\cA(s)}}{\pi}\sqrt{x}\,(1+\OO(x)),
 \qquad 0<x<r,\qquad c\le\cA(s)\le C.
\end{equation}
With the square root positive for $\zeta>0$, $m_s(E_s+\zeta)$ has a
convergent expansion in $\sqrt\zeta$ for $0<|\zeta|<r$ off the
negative real axis, with real coefficients and leading terms
\begin{equation}\label{eq:app-generated-transform}
 m_s(E_s+\zeta)=m_s(E_s)+\sqrt{\cA(s)}\sqrt\zeta+\OO(\zeta).
\end{equation}
The functions $E_s$, $m_s(E_s)$ and $\sqrt{\cA(s)}$ are $C^1$ on $J$.
We have $E_s'=\chi_s(E_s)$ and \eqref{s4:edge-time}, with the
additional term $H_s(E_s)$ from \eqref{eq:app-time-dependent-forcing}
when $\sfb_s$ depends on time.

For each $s\in J$, the flux $\cF_s$ is holomorphic and uniformly bounded
in a ball of fixed radius about $(m_s(E_s),E_s)$, and
\eqref{eq:edge-flux-nondegeneracy} holds uniformly on $J$.
\end{theorem}

Thus \Cref{ass:edge-data}\textup{(ii)} holds on $J$ for a time-independent
flux, and on sufficiently short closed subintervals of $J$ for a
time-dependent flux. For DMC and corners it holds on all of $J$ if the
coefficients are holomorphic and uniformly bounded on a fixed neighborhood of
$\{E_t:t\in J\}$ for every $s\in J$, and
\begin{equation}\label{eq:app-common-flux-domain}
 \begin{cases}
 \displaystyle\inf_{s,t\in J}
 \bigl(\varphi_s^-(E_t)+\varphi_s^+(E_t)e^{-\theta m_t(E_t)}\bigr)>0,
 &\textnormal{DMC},\\[3pt]
 \displaystyle\inf_{s,t\in J}
 \bigl((1-s)W'(E_t)-\theta m_t(E_t)\bigr)>0,
 &\textnormal{corners}.
 \end{cases}
\end{equation}

\begin{proof}[Proof of \Cref{thm:app-edge-formation}]
We first estimate the derivatives of the characteristic velocity. We
then show that the characteristic map has a simple critical point
whose image is the right edge.

\paragraph{The characteristic equations.}
Fix $a<T$ and start an exterior characteristic at $z_a=u>E_a$ with
$p_a=m_a(u)$. Write $v_s=\dot z_s$. Denote the contour integral in
\eqref{def of F_s}, including the factor $(2\pi\ri\theta)^{-1}$, by
$S_s$. On the fixed working domains, $S_s$ is continuous in time and
holomorphic in space with bounded derivatives. If the interaction map
depends on time, \Cref{rem:b-time-dependent-map} and
\eqref{s4:kernel-split} give
\begin{equation}\label{eq:app-time-dependent-forcing}
 \partial_s m_s+\chi_s\partial_z m_s=F_s+H_s,\qquad
 H_s(z)=\int\partial_s g_{\sfb_s}(r,z)\,\mu_s(\rd r).
\end{equation}
Here $F_s$ uses the coefficients at time $s$. The function $H_s$ is
holomorphic near the support and vanishes for a fixed map.
For the Poisson chains of \Cref{subsec:b-dbm}, we use $\sfb(z)=z$ and
the smoothing convention of \Cref{subsec:prelim-measures}. Then
\begin{equation}\label{eq:app-poisson-flux}
 \cF_s(p,z)=\varphi(z)e^{-\theta p},\qquad
 \chi_s(z)=\varphi(z)e^{-\theta m_s(z)}.
\end{equation}
The characteristic equations are
\[
 v_s=-\theta^{-1}\partial_1\cF_s(p_s,z_s),\qquad
 \dot p_s=\theta^{-1}\partial_2\cF_s(p_s,z_s)-S_s(z_s)+H_s(z_s).
\]
As long as $z_s>E_s$, the transport equation gives $p_s=m_s(z_s)$.
Differentiating $v_s$ along the characteristic gives
\begin{equation}\label{eq:app-acceleration}
 \dot v_s=
 \begin{cases}
 \displaystyle \tfrac14 V'(z_s)V''(z_s)
     -\int g(z_s,x)\,\mu_s(\rd x),&\textnormal{DBM},\\[3pt]
 \displaystyle v_s(1-v_s)
     \left[\partial_s\log(\varphi_s^+/\varphi_s^-)
     -(\log\varphi_s^-)'+\theta(S_s-H_s)\right](z_s),
     &\textnormal{DMC},\\[3pt]
 -v_s^2\bigl[W'(z_s)/\theta-S_s(z_s)\bigr],&\textnormal{corners},\\[3pt]
 \theta S_s(z_s)v_s,&\textnormal{Poisson}.
 \end{cases}
\end{equation}
For DBM this follows from \eqref{equ:MVequ}, with $g$ given by
\eqref{eq:def_g_func}. For DMC, differentiate
$\log(v_s/(1-v_s))=\log(\varphi_s^+/\varphi_s^-)(z_s)-\theta p_s$.
The spatial derivative of the weight ratio cancels with its term in
$\partial_2\cF_s$. For corners, differentiate
$1/v_s=p_s-(1-s)W'(z_s)/\theta$. For Poisson, differentiating
$\log v_s=\log\varphi(z_s)-\theta p_s$ cancels the terms containing
$\varphi'$.

\paragraph{Derivatives of the velocity.}
We verify that \eqref{eq:app-initial-edge-condition} implies
$-\chi_0'(u)\to\infty$ as $u\downarrow E_0$. We also show that at any
starting time $a<T$, for exterior real $u$ in the working neighborhood,
\begin{equation}\label{eq:app-curvature}
 \chi_a''(u)\ge c(-\chi_a'(u))^{3/2}-C|\chi_a'(u)|-C,
\end{equation}
whenever $-\chi_a'(u)$ is sufficiently large, uniformly over the admissible measures.

\textbf{Dyson Brownian motion.}
By \eqref{eq:velocity-from-flux},
\[
 -\chi_a'=m_a'+V''/2,\qquad
 \chi_a''=-m_a''-V'''/2.
\]
Since $-m_a''\ge2(m_a')^{3/2}$ by H\"older's inequality and the
derivatives of $V$ are bounded,
\eqref{eq:app-curvature} follows. At $a=0$, the DBM condition in
\eqref{eq:app-initial-edge-condition} also gives $-\chi_0'\to\infty$.

\smallskip
\noindent\textbf{Discrete Markov chains.}
By \eqref{s4:kernel-split}, $m_a-\widehat m_a$ has bounded analytic
derivatives. Put
\[
 L=\bigl(\log(\varphi_a^+/\varphi_a^-)\bigr)'-\theta m_a'.
\]
Then \eqref{eq:velocity-from-flux} gives
$\chi_a'=\chi_a(1-\chi_a)L$. Set $t=(u-x)^{-1}$ and
$f(t)=\theta\rho_a(u-1/t)$ for $t>0$. Then $0\le f\le1$,
$\int f=\theta\widehat m_a'$ and $2\int tf=-\theta\widehat m_a''$.
Among such $f$ of given mass, $\int tf$ is minimized by
$\mathbf1_{[0,\int f]}$. Hence $-\widehat m_a''\ge\theta(\widehat m_a')^2$, and
$L'\ge L^2-C|L|-C$. Differentiating $\chi_a'$ gives
\begin{equation}\label{eq:app-dmc-curvature}
 \chi_a''\ge\frac{2(\chi_a')^2}{\chi_a}-C|\chi_a'|-C
       \ge2(\chi_a')^2-C|\chi_a'|-C.
\end{equation}
This proves \eqref{eq:app-curvature}, including when $\chi_a$ approaches $1$.
Moreover, $\chi_a(1-\chi_a)\asymp e^{\theta\widehat m_a}$
near $E_a$, since the weights are bounded above and away from zero.
The DMC condition in
\eqref{eq:app-initial-edge-condition} therefore implies
$-\chi_0'\to\infty$.

\smallskip
\noindent\textbf{Weighted $\beta$-corners.}
Enlarge the support interval slightly to $[\ell,r]$. By the strict
margin, \eqref{eq:b-corners-derivative} still holds with these
endpoints and a smaller margin. Integrating from $u$ to $r$ gives
$W'(u)/\theta+(u-\ell)^{-1}\ge c>0$. Since $\mu_a$ has mass
$1-a\ge1-T$, it follows that $\cB_a/\theta\asymp-m_a\ge c$.
Cauchy--Schwarz for the positive measure
$((u-x)^{-1}-(u-\ell)^{-1})\mu_a(\rd x)$ gives
\[
 -\frac{\cB_a m_a''}{2\theta}-(m_a')^2
 \ge c(-m_a'')-C|m_a|m_a'.
\]
Since $\chi_a=-\theta/\cB_a$ and
$-m_a''\ge2(m_a')^{3/2}/\sqrt{1-a}$, differentiation yields
\[
 \chi_a''\ge c\frac{-m_a''}{|m_a|^3}
             -C\frac{m_a'}{|m_a|^2}-C
 \ge c(-\chi_a')^{3/2}-C|\chi_a'|-C,
\]
where we used $-\chi_a'=\chi_a^2(m_a'-(1-a)W''/\theta)$.
This proves \eqref{eq:app-curvature}. The same identity at $a=0$ shows that
the corners condition in \eqref{eq:app-initial-edge-condition}
implies $-\chi_0'\to\infty$.

\smallskip
\noindent\textbf{Poisson chains.}
The density bound and the Markov--Krein correspondence
\cite[Theorem~2.1]{huang2017beta}, after rescaling space by
$\theta^{-1}$, give
\begin{equation}\label{eq:app-poisson-markov-krein}
 e^{-\theta m_a(z)}=1+\int\frac{\nu(\rd r)}{z-r},
 \qquad \nu\ge0,\qquad \nu(\mathbb R)=\theta,
\end{equation}
with $\nu$ supported in the convex hull of $\supp\mu_a$.
Cauchy--Schwarz and H\"older's inequality give
\[
 (e^{-\theta m_a}-1)^2\le-\theta\partial_u e^{-\theta m_a},
 \qquad
 \partial_u^2e^{-\theta m_a}
 \ge\frac2{\sqrt\theta}
       (-\partial_u e^{-\theta m_a})^{3/2}.
\]
In $\chi_a=\varphi e^{-\theta m_a}$, the first inequality bounds
$|\varphi' e^{-\theta m_a}|$ by
$C(1+\sqrt{-\partial_u e^{-\theta m_a}})$. Since $\varphi\ge c$
and $e^{-\theta m_a}\ge1$, this gives the divergence at $a=0$ and
$\chi_a^2\le C(1-\chi_a')$ whenever $-\chi_a'$ is large.
Differentiating $\varphi e^{-\theta m_a}$ twice and using the second
inequality gives \eqref{eq:app-curvature}.

\paragraph{The critical point of the characteristic map.}
Fix $a<T$ and take $0<h\le T-a$. Let $\Phi_h(x,v)$ be the position at time $a+h$ of the
solution of \eqref{eq:app-acceleration} with $z_a=x$ and $v_a=v$.
Thus $z_{a+h}(u)=\Phi_h(u,\chi_a(u))$. For DMC and corners the initial
velocities are bounded, and for DBM the right-hand side of
\eqref{eq:app-acceleration} depends only on the position and time.
Differentiating the equation in $x,v$ gives, as long as the positions
stay in the working domain,
\begin{equation}\label{eq:app-flow-derivatives}
 \Phi_{h,x}=1+\OO(h^2),\qquad \Phi_{h,v}=h+\OO(h^2),\qquad
 \Phi_{h,xx},\ \Phi_{h,xv},\ \Phi_{h,vv}=\OO(h^2).
\end{equation}
For Poisson, take $0<v\le Ch^{-1/2}$ and put $Z(t)=z_{a+ht}$.
For $0\le t\le1$,
\[
 Z''(t)=h\theta S_{a+ht}(Z(t))Z'(t),\qquad Z(0)=x,\quad Z'(0)=hv.
\]
In particular $Z'(t)=\OO(\sqrt h)$. Differentiating in $x$ and in
$w=hv$, and using $\partial_v=h\partial_w$, gives
\begin{equation}\label{eq:app-poisson-flow-derivatives}
 \begin{aligned}
 \Phi_{h,x}&=1+\OO(h^{3/2}),&
 \Phi_{h,v}&=h+\OO(h^2),\\
 \Phi_{h,xx}&=\OO(h^{3/2}),&
 \Phi_{h,xv}&=\OO(h^2),\qquad \Phi_{h,vv}=\OO(h^3).
 \end{aligned}
\end{equation}
Choose $u_+$ a fixed distance to the right of the support bound at
time $a$, within the working domain, so that $z_{a+t}(u_+)>E_{a+t}$
for $0\le t\le \min\{h_0,T-a\}$, with $h_0>0$ uniformly small.
For DMC, \Cref{lem:b-bounds} gives $E_{a+t}\le r_0+a+t$,
and \eqref{eq:b-initial-wall-clearance} keeps the tiling weights positive
in a fixed neighborhood of $r_0+a$ for small $t$.
For the other models we can take $u_+$ beyond the common support bound.
After decreasing $h_0$, we have $\partial_u z_{a+h}(u_+)>0$.

Suppose first that $-\chi_a'(u)\to\infty$ as $u\downarrow E_a$.
For DMC and corners, \eqref{eq:app-flow-derivatives} gives
$\partial_u z_{a+h}(u_-)<0$ for $u_-$ sufficiently close to $E_a$.
For DBM take $m_a'(u_-)=2/h$. Since $m_a'$ decreases and
$|m_a|\le\sqrt{m_a'}$, the initial velocities on $[u_-,u_+]$
are $\OO(h^{-1/2})$. The positions change by $\OO(\sqrt h)$ and
\eqref{eq:app-flow-derivatives} again gives the negative sign.
For Poisson choose $-(e^{-\theta m_a})'(u_-)=C_1/h$ with $C_1$
sufficiently large. By \eqref{eq:app-poisson-markov-krein} this
derivative in absolute value decreases on $[u_-,u_+]$, and
$e^{-\theta m_a}=\OO(h^{-1/2})$ there. Hence
$\chi_a'(u_-)<-2/h$ and \eqref{eq:app-poisson-flow-derivatives}
gives $\partial_u z_{a+h}(u_-)<0$.

For part~\textup{(ii)}, take $a=0$, choose $\delta$ so that
$r_0+\delta<u_+$, and take
$C/(-\chi_0'(u_0))<h\le h_0$. For DMC and corners,
$u_-=u_0$ gives $\partial_u z_h(u_-)<0$.
For DBM, $m_0'(u_0)>2/h$ after increasing $C$.
For Poisson, the bound $(e^{-\theta m_0}-1)^2\le
-\theta(e^{-\theta m_0})'$ and the bounded derivatives of $\varphi$
give $-\chi_0'\le C_2(1-(e^{-\theta m_0})')$.
Thus for both models the specified $u_-$ can be chosen in $[u_0,u_+)$.

On the interval $[u_-,u_+]$ just constructed, the bounds
\eqref{eq:app-flow-derivatives} or \eqref{eq:app-poisson-flow-derivatives}
apply. At every zero of $\partial_u z_{a+h}$ on this interval,
\[
 \chi_a'=-\frac{\Phi_{h,x}}{\Phi_{h,v}}=-h^{-1}(1+\OO(h)),\qquad
 \partial_u^2 z_{a+h}=\Phi_{h,v}\chi_a''+\OO(1)
                 \ge ch^{-1/2}-C>0
\]
by the chain rule and \eqref{eq:app-curvature}, after decreasing $h_0$.
The derivative is negative at $u_-$ and positive at $u_+$, so there is
exactly one $u_h\in(u_-,u_+)$ with
\begin{equation}\label{eq:app-simple-fold}
 \partial_u z_{a+h}(u_h)=0,\qquad \partial_u^2z_{a+h}(u_h)>0.
\end{equation}

\paragraph{The right edge and its coefficient.}
The flow $(x,v)\mapsto(\Phi_h,\partial_h\Phi_h)$ has positive
Jacobian determinant. At $u_h$,
\[
 \partial_u v_{a+h}(u_h)=
 -\frac{\det D_{(x,v)}(\Phi_h,\partial_h\Phi_h)}{\Phi_{h,v}}<0.
\]
The velocity stays in $(0,1)$ for DMC, negative for corners and positive
for Poisson by \eqref{eq:app-acceleration}. The flux formulas
\eqref{e:F_CP} and \eqref{eq:app-poisson-flux} therefore give
$\partial_{11}\cF_s>0$ along these real characteristics, as for DBM.
Differentiating $v_s=-\theta^{-1}\partial_1\cF_s(p_s,z_s)$ at $u_h$
now gives $\partial_u p_{a+h}(u_h)>0$.

Near $u_+$ we have $p_{a+h}(u)=m_{a+h}(z_{a+h}(u))$.
On $(u_h,u_+]$, $\partial_u z_{a+h}>0$. The local inverses of
$u\mapsto z_{a+h}(u)$ therefore continue $m_{a+h}$ from the upper
half-plane across $(z_{a+h}(u_h),z_{a+h}(u_+))$, where it is real.
This interval carries no mass. At $u_h$, \eqref{eq:app-simple-fold}
gives a convergent expansion of $u-u_h$ in powers of
$\sqrt{z-z_{a+h}(u_h)}$, with positive leading coefficient.
Substituting into $p_{a+h}(u)$ gives positive density immediately to
the left of $z_{a+h}(u_h)$. Since $z_{a+h}(u_+)>E_{a+h}$, we obtain
$E_{a+h}=z_{a+h}(u_h)$ and
\begin{equation}\label{eq:app-fold-coefficient}
 \cA(a+h)=
 \frac{2(\partial_u p_{a+h}(u_h))^2}{\partial_u^2 z_{a+h}(u_h)}>0.
\end{equation}
The coefficients of the expansion are real. At $z=E_{a+h}-x+\ri0$
the even powers are real, and the next imaginary term after $\sqrt x$
is $\OO(x^{3/2})$. By \eqref{s4:kernel-split} this proves
\eqref{eq:app-generated-edge} and \eqref{eq:app-generated-transform}.

A regular square-root edge has finite edge transform and
$\widehat m_s'(E_s+)=\infty$, so \eqref{eq:app-initial-edge-condition}
holds at each such $s$. We can thus restart the argument at any time
in the first interval of regularity and continue in steps of $h_0/2$
up to $T$. By \eqref{eq:app-simple-fold} and the implicit function
theorem, $u_h$ is $C^1$ in $h>0$ wherever this critical point was
constructed. The characteristic equations and \eqref{eq:app-fold-coefficient}
then give $C^1$ dependence of $E_s$, $m_s(E_s)$ and $\sqrt{\cA(s)}$.
A finite cover of $[\varepsilon,T]$ gives a common expansion radius
and uniform coefficient bounds. At $u=u_{s-a}$, differentiating $\partial_u z_s(u)=0$ in $s$ gives
$\frac{\rd}{\rd s}u_{s-a}=-\partial_u v_s/\partial_u^2 z_s$.
Since $\partial_u v_s=-\theta^{-1}\partial_{11}\cF_s\,\partial_u p_s$
at this point, differentiating $E_s=z_s(u_{s-a})$ and
$m_s(E_s)=p_s(u_{s-a})$ gives $E_s'=\chi_s(E_s)$ and
\eqref{s4:edge-time}, with $H_s(E_s)$ added when $\sfb_s$ varies.

For part~\textup{(iii)}, it suffices to prove the velocity estimates in
a right neighborhood of $E_0$ at time zero. Finiteness of $m_0(E_0+)$
gives
\[
 \lim_{u\downarrow E_0}\frac{m_0'(u)^2}{-m_0''(u)}=0.
\]
Indeed, split the integral for $m_0'$ at $E_0-\delta$. The first part
is bounded for fixed $\delta$, while $-m_0''\to\infty$ since
$m_0'\to\infty$. By Cauchy--Schwarz, the square of the
second is at most $(-m_0''/2)\int_{E_0-\delta}^{E_0}
(E_0-x)^{-1}\mu_0(\rd x)$, and the last integral tends to zero with
$\delta$. Since $\cB_0(E_0)>0$, differentiating $\chi_0=-\theta/\cB_0$
now gives $-\chi_0'\asymp m_0'\to\infty$ and
$\chi_0''\ge c(-m_0'')\ge c(-\chi_0')^{3/2}$ near $E_0$.
Choose $u_+>E_0$ in this neighborhood. By interlacing,
$E_s\le E_0<z_s(u_+)$ for $0\le s\le T_1$ if $T_1$ is small.
The same sign change and inversion give \eqref{eq:app-simple-fold} and
\eqref{eq:app-fold-coefficient} for $0<s\le T_1$, with uniform bounds
on each $[\varepsilon,T_1]$.

If the initial edge is regular, use $q=\sqrt{u-E_0}$ as parameter.
At $s=0$, the position has first and second $q$-derivatives $0$ and
$2$ at $q=0$, and the transform has first derivative
$\sqrt{\cA(0)}>0$. The implicit function theorem continues the critical
point in this parameter from $s=0$. Its initial velocity in $q$ is
$-\partial_q v_0/2>0$, so it lies on the exterior branch for small
positive $s$. The expansion and the time regularity proved above
therefore extend to $[0,T_1]$, after decreasing $T_1$ if needed.

Finally, at the edge $\partial_{11}\cF_s$ equals $\theta$ for DBM,
$\theta^2v_s(1-v_s)$ for DMC, $\theta v_s^2$ for corners, and
$\theta^2v_s$ for Poisson, with $v_s=\chi_s(E_s)$. These quantities
are positive and continuous on $J$. The explicit fluxes in
\eqref{e:F_CP} and \eqref{eq:app-poisson-flux} thus give the claimed
bound and a common radius about the respective edge points.
For a time-independent flux these neighborhoods cover the whole edge
curve. For a time-dependent flux, continuity gives a common
neighborhood on each sufficiently short time interval.
Under \eqref{eq:app-common-flux-domain}, the argument of the logarithm
in \eqref{e:F_CP} has positive real part on one neighborhood of the
whole curve for all $s\in J$. This gives
\Cref{ass:edge-data}\textup{(ii)} on $J$.
\end{proof}
\subsubsection{Application to the models}
\label{rem:app-initial-tests}

We retain the parameter, support and analyticity conditions of
\Cref{sec:appendix-b-verification}. In each family,
\eqref{eq:app-initial-edge-condition} gives $s_*=0$.
Under the full-interval hypotheses of \Cref{thm:app-edge-formation},
\eqref{eq:app-finite-compression} gives a sufficient waiting time
whenever it holds at some $u_0>E_0$ in the fixed neighborhood.

\paragraph{Discrete Markov chains.}
For DMC, suppose that $\widehat m_0'(E_0+)=\infty$. If
$\widehat m_0(E_0+)>-\infty$, the DMC condition in
\eqref{eq:app-initial-edge-condition} follows directly. It also holds
if $\mu_0\le(1-\kappa)\theta^{-1}\rd x$ near $E_0$ for some
$0<\kappa<1$. Indeed, the density cap gives
\[
 \bigl(e^{\theta\widehat m_0}\bigr)''
 =\theta e^{\theta\widehat m_0}
       \bigl(\widehat m_0''+\theta(\widehat m_0')^2\bigr)\le0.
\]
If $\widehat m_0(E_0+)=-\infty$, the strict local cap gives
$e^{\theta\widehat m_0(E_0+d)}\ge c d^{1-\kappa}$ near $d=0$.
Concavity then gives $(e^{\theta\widehat m_0})'(E_0+)=\infty$.

Thus $s_*=0$ for a density comparable to $(E_0-x)^\gamma$ with
$0<\gamma\le1$, or with a positive limit strictly below $\theta^{-1}$
at $E_0$. This applies to \textbf{DMC with general analytic weights}
under \eqref{cond:b-D}, including the perturbations of \Cref{subsec:b-dmc},
and to the \textbf{nonintersecting Bernoulli walks}, \textbf{ascending
Jack process}, and \textbf{ascending Macdonald process} of that subsection.
For \textbf{uniform}, \textbf{$(q,\kappa)$-distributed} and
\textbf{trigonometric lozenge tilings}, use the wall and parameter
conditions of \Cref{cor:b-affine-weights}, including the real and
imaginary regimes of the $(q,\kappa)$ weights. Time-dependent weights
enter through \eqref{eq:app-acceleration}, and a time-dependent
interaction map gives the term $H_s$ in \eqref{eq:app-time-dependent-forcing}.

For \textbf{descending Jack}, \textbf{descending Macdonald}, and
\textbf{Koornwinder quasi-branching transitions}, the test is applied
at the right edge of the normalized hole measure. Its mass is one
and its density cap is $\theta$. The condition is therefore
$e^{\widehat m_0/\theta}\widehat m_0'\to\infty$ as $u\downarrow E_0$,
with density cap $\theta$ in the two sufficient tests above.

\paragraph{Weighted $\beta$-corners.}
For the \textbf{weighted Dixon--Anderson transitions}, the corners
condition in \eqref{eq:app-initial-edge-condition} is equivalent to
\[
 \mu_0\{E_0\}=0,\qquad \int(E_0-x)^{-2}\,\mu_0(\rd x)=\infty.
\]
Indeed, the integral is $m_0'(E_0+)$. If it diverges and the end point
carries no atom, split the integral for $m_0(u)$ at $E_0-\delta$. The part
away from $E_0$ stays bounded. By Cauchy--Schwarz, the part near $E_0$,
divided by $\sqrt{m_0'(u)}$, is bounded in absolute value by
$\sqrt{\mu_0([E_0-\delta,E_0])}$. Letting $u\downarrow E_0$ and then
$\delta\downarrow0$ gives $m_0^2/m_0'\to0$. Conversely, a finite
$m_0'(E_0+)$ gives a finite ratio $m_0'/m_0^2$ at the edge,
and an atom of mass $p>0$ gives $m_0'/m_0^2\to1/p$.

For an initial law with no atom at $E_0$, a density comparable to
$(E_0-x)^\gamma$ with $-1<\gamma\le1$ gives $s_*=0$ under \eqref{eq:b-corners-derivative}. This includes
\textbf{unweighted $\beta$-corners} for $T<1$ and the semicircle
initial measure for \textbf{Gaussian $\beta$-corners}. For
\textbf{linear}, \textbf{quadratic}, \textbf{Laguerre}, and
\textbf{Jacobi potentials}, the same conclusion holds in the derivative
ranges of \Cref{subsec:b-corners}. A rightmost atom gives a different
waiting time.

\begin{proposition}[A rightmost atom for weighted corners]
\label[proposition]{prop:app-corners-atom}
Under the working hypotheses of this subsection for corners, including
\eqref{eq:b-corners-derivative}, suppose
$\mu_0\{E_0\}=p\in(0,1]$. Then
\[
 \mu_s\{E_0\}=p-s,\qquad 0\le s\le\min\{p,T\}.
\]
If $p<T$, the conclusions of \Cref{thm:app-edge-formation} hold
on every $[\varepsilon,T]$ with $p<\varepsilon<T$. In this case,
$s_*=p$ is the exact waiting time for a regular right edge.
If $p=1$, then
$\mu_s=(1-s)\delta_{E_0}$ for $0\le s\le T$.
\end{proposition}

\begin{proof}
The expansion $\chi_0(u)=-(u-E_0)/p+\oo(u-E_0)$ and
\eqref{eq:app-acceleration} give
\[
 z_s(u)-E_0=(1-s/p)(u-E_0)+\oo(u-E_0),\qquad
 \partial_s z_s(u)=-(u-E_0)/p+\oo(u-E_0).
\]
For $s<p$, interlacing and
$m_s(z_s)=(\partial_s z_s)^{-1}+(1-s)W'(z_s)/\theta$
therefore give an atom of mass $p-s$ at $E_0$.
If $p\le T$, the right-tail bounds in \Cref{lem:b-bounds}, applied
between two times, imply that this atom vanishes at $s=p$.
Restarting at $a<p$ close to $p$ gives
$\partial_u z_{a+h}(E_0+)=1-h/(p-a)<0$ for $p-a<h\le h_0$.
Since $\partial_u z_{a+h}(u)>0$ at a fixed $u>E_0$ for small $h$,
\eqref{eq:app-curvature} gives the critical point in
\eqref{eq:app-simple-fold}. The edge identification and continuation
in the proof of \Cref{thm:app-edge-formation} yield regularity for
the allowed times $s>p$. The Dirac case follows from interlacing.
\end{proof}

Now assume the endpoint conditions of \Cref{subsec:b-corners}
and $m_0'(E_0+)=\infty$. The finite edge transform and the strict
endpoint signs give \eqref{eq:app-corners-local-condition}. For example, this
gives regularity on $(0,T_1]$ for the Laguerre and Jacobi potentials.
For the semicircle law on $[-1,1]$, the initial square-root edge persists
on $[0,T_1]$ for the quadratic examples with $\lambda>-2\theta$ and
the quartic examples with $\lambda\ge0$, since
$\cB_0(1)=\lambda+2\theta>0$. The same conclusion holds for the
\textbf{polynomial potentials and analytic perturbations} preserving
the strict endpoint conditions. Perturbations preserving
\eqref{eq:b-corners-derivative} retain its full time interval.

\paragraph{Dyson Brownian motion.}
For DBM, \eqref{eq:app-initial-edge-condition} reads
$\int(E_0-x)^{-2}\,\mu_0(\rd x)=\infty$. A rightmost atom, including
a Dirac initial law, or an edge density comparable to $(E_0-x)^\gamma$
with $-1<\gamma\le1$ gives $s_*=0$. This applies to \textbf{additive
DBM}, \textbf{DBM with a quadratic potential}, and the \textbf{general
analytic potentials} of \Cref{subsec:b-dbm}, including its quartic
examples with $g>0$ and $a\in\bR$. The semicircle law of the stationary
\textbf{Gaussian (Hermite) process} satisfies the density test.

\begin{remark}[A restriction for general potentials]
\label[remark]{rem:app-stationary-obstruction}
For general analytic $V$, compact support alone does not imply
eventual square-root behavior.
For example,
\[
 V(x)=\frac{16}{15}x^6-4x^4+6x^2,\qquad
 \rho_0(x)=\frac{16}{5\pi}(1-x^2)^{5/2}\mathbf1_{\{|x|\le1\}}
\]
satisfy $V''\ge-6$, and
$m_0(z)=-V'(z)/2+(16/5)(z^2-1)^{5/2}$, with the branch behaving as
$z^5$ at infinity. Its boundary real part is $-V'/2$, so the
McKean--Vlasov velocity vanishes on the support. This stationary
measure retains a $5/2$ right edge for all times.
\end{remark}

\paragraph{Poisson chains and quantum walks.}
For \textbf{$\theta$-nonintersecting Poisson walks} and the
\textbf{spatially dependent Poisson rates} of \Cref{subsec:b-dbm},
use the DBM integral test together with $\mu_0\le\theta^{-1}\rd x$.
A density comparable to $(E_0-x)^\gamma$ with $0\le\gamma\le1$ gives
$s_*=0$, including the packed interval of density $\theta^{-1}$.
For the \textbf{defining-representation quantum walk} at $\theta=1$,
the clock $s=k/n^2$ gives the homogeneous Poisson evolution, where
$k$ is the number of transitions. The \textbf{binomial mixture of
exterior-power transitions} uses $s=k/n$ and the Bernoulli test at
$\theta=1$, as identified in \Cref{subsec:b-dbm}.

\subsubsection{Explicitly solvable models}
For additive DBM, unweighted corners, homogeneous Poisson walks and
the constant-ratio DMC below, we can compute the waiting time even
when the divergence condition \eqref{eq:app-initial-edge-condition}
fails. We retain the common support and analytic bounds.

\begin{proposition}[Explicit waiting times]
\label[proposition]{prop:app-scalar-waiting}
Consider DBM with $V=0$, corners with $W=0$, Poisson with $\varphi=1$,
or DMC with constant
$\varphi^+/\varphi^-$ and either $\sfb(z)=z$ or
$\sfb(z)=e^{\kappa z}$, $\kappa\in\mathbb R\setminus\{0\}$.
Let $\mu_0$ be a compactly supported probability measure, with
$\mu_0\le\theta^{-1}\rd x$ for DMC and Poisson, and $\mu_0$
non-Dirac for corners.
Define
\begin{equation}\label{eq:app-scalar-threshold}
 s_*=\frac1{-\chi_0'(E_0+)}.
\end{equation}
Then $-\chi_0'$ is positive and strictly decreasing on $(E_0,\infty)$,
and $s_*<\infty$, with $s_*<1$ for corners.
For every $s>s_*$ in the allowed time range, there is a unique
$u_s>E_0$ such that $1+s\chi_0'(u_s)=0$, and
\begin{equation}\label{eq:app-scalar-edge}
 E_s=u_s+s\chi_0(u_s),\qquad
 \cA(s)=\frac{2m_0'(u_s)^2}{s\chi_0''(u_s)}.
\end{equation}
The regularity conclusions of \Cref{thm:app-edge-formation} hold
uniformly on every $[\varepsilon,T]$ with $s_*<\varepsilon<T$.
There is no regular square-root right edge for $0<s<s_*$ in the
allowed time range. No assertion is made at $s=s_*$.
\end{proposition}

\begin{proof}
For these models, $F_s=0$ in \eqref{def of F_s}, after removing
the common analytic factor from the DMC weights. Thus
\Cref{def:prelim-characteristic-flow} gives
\begin{equation}\label{eq:app-straight-characteristics}
 z_s(u)=u+s\chi_0(u),\qquad m_s(z_s(u))=m_0(u).
\end{equation}
Only the initial velocity and its derivatives need to be computed.

\textbf{Dyson Brownian motion without potential.}
For $V=0$, \eqref{eq:velocity-from-flux} gives $\chi_0=-m_0$, so
$-\chi_0'=m_0'$ strictly decreases to zero and
\[
 s_*=\frac1{m_0'(E_0+)}.
\]
This is the critical-point calculation for free convolution used
in the proof of \cite[Lemma~2.14]{edgerigiditydysonbrownian}.
In particular, a Dirac initial measure has $s_*=0$.

\smallskip
\noindent\textbf{Nonintersecting Poisson walks.}
For the homogeneous process of \Cref{subsec:b-dbm},
\eqref{eq:app-poisson-flux} gives $\chi_0=e^{-\theta m_0}$.
By \eqref{eq:app-poisson-markov-krein}, $-\chi_0'$ strictly
decreases to zero, and
\[
 s_*=\left(
 \theta\lim_{u\downarrow E_0}e^{-\theta m_0(u)}m_0'(u)
 \right)^{-1}.
\]
For the packed initial measure on $[a,a+\theta]$, one has
$e^{-\theta m_0(u)}=1+\theta/(u-E_0)$, so $s_*=0$ and
\[
 E_s=E_0+s+2\sqrt{\theta s},\qquad
 \cA(s)=\frac1{\sqrt{\theta s}(\theta+\sqrt{\theta s})^2}.
\]

\smallskip
\noindent\textbf{Unweighted $\beta$-corners.}
For $W=0$, \eqref{eq:velocity-from-flux} gives $\chi_0=1/m_0$, and hence
\[
 s_*=\left(
 \lim_{u\downarrow E_0}\frac{m_0'(u)}{m_0(u)^2}
 \right)^{-1}.
\]
Cauchy--Schwarz shows that $m_0'/m_0^2$ strictly decreases to $1$
at infinity for a non-Dirac measure, so $s_*<1$.
If $\mu_0$ has an atom at $E_0$, then $s_*=\mu_0\{E_0\}$,
in agreement with \Cref{prop:app-corners-atom}.

\smallskip
\noindent\textbf{Nonintersecting Bernoulli walks and the ascending Jack process.}
Here $\sfb(z)=z$ and the weight ratio is constant, as in
\Cref{subsec:b-dmc}. Write
$p=\varphi^+/(\varphi^-+\varphi^+)\in(0,1)$. By
\eqref{eq:velocity-from-flux},
\[
 \chi_0(z)=
 \frac{p e^{-\theta m_0(z)}}{1-p+p e^{-\theta m_0(z)}}.
\]
The density bound and the Nevanlinna representation
\cite[Definition~1.3 and (1.6)]{charofairy} give
\begin{equation}\label{eq:app-dmc-velocity-measure}
 \chi_0(z)=p+\int\frac{\nu(\rd x)}{z-x},\qquad
 \nu\ge0,\qquad \nu(\mathbb R)=\theta p(1-p),
\end{equation}
with $\nu$ supported in the convex hull of $\supp\mu_0$.
Thus $-\chi_0'(u)=\int(u-x)^{-2}\nu(\rd x)$ strictly decreases
to zero. For the packed initial measure
$\mu_0=\theta^{-1}\mathbf1_{[a,a+\theta]}\rd x$,
direct substitution gives
\[
 s_*=\frac{\theta p}{1-p}.
\]

\smallskip
\noindent\textbf{Ascending Macdonald process and exponential coordinates.}
Take $\sfb(z)=e^{\kappa z}$ with constant weight ratio and $p$ as above;
this includes the ascending Macdonald example in \Cref{subsec:b-dmc}.
Applying the Nevanlinna representation used in
\eqref{eq:app-dmc-velocity-measure} in the coordinate $y=e^{\kappa u}$
gives a positive measure $\nu$ on the convex hull of the transformed
support such that
\[
 \begin{aligned}
 -\chi_0'(u)&=|\kappa|y\int\frac{\nu(\rd t)}{(y-t)^2},\\
 \chi_0''(u)&=\kappa^2y
     \int\frac{y+t}{|y-t|^3}\,\nu(\rd t)>0.
 \end{aligned}
\]
Also $\chi_0'(u)\to0$ as $u\to\infty$.
For the packed initial measure on $[a,a+\theta]$,
\[
 e^{-\theta m_0(u)}
 =e^{\kappa\theta}
   \frac{e^{\kappa u}-e^{\kappa a}}
        {e^{\kappa u}-e^{\kappa(a+\theta)}},\qquad
 s_*=\frac{p(e^{\kappa\theta}-1)}{\kappa(1-p)}.
\]

In each of the scalar models, $\chi_0''>0$, so $1+s\chi_0'(u_s)=0$ gives
a simple critical point for $s>s_*$. Formula
\eqref{eq:app-fold-coefficient} and the edge identification in the proof of
\Cref{thm:app-edge-formation} give \eqref{eq:app-scalar-edge}
and its uniformity on $[\varepsilon,T]$.
For $0<s<s_*$, $1+s\chi_0'(u)>0$ for $u>E_0$.
If $m_0'(E_0+)<\infty$, differentiating \eqref{eq:app-straight-characteristics} gives
\[
 m_s'(E_s+)=
 \frac{m_0'(E_0+)}{1+s\chi_0'(E_0+)}<\infty.
\]
The remaining cases retain a corners atom by
\Cref{prop:app-corners-atom}, or an infinite DMC edge value along
\eqref{eq:app-straight-characteristics}, which also excludes a
regular square-root right edge.
\end{proof}

For packed initial data, the identifications in \Cref{subsec:b-dbm} give
$s_*=0$ for the defining-representation quantum walk and
$s_*=p/(1-p)=c$ for the binomial mixture of exterior-power transitions,
on the respective time scales specified above.

\bibliography{ref}
\bibliographystyle{alpha}

\end{document}